# SUPER FUZZY MATRICES AND SUPER FUZZY MODELS FOR SOCIAL SCIENTISTS

**W. B. Vasantha Kandasamy**
e-mail: vasanthakandasamy@gmail.com
web: http://mat.iitm.ac.in/~wbv
www.vasantha.net

**Florentin Smarandache**
e-mail: smarand@unm.edu

**K. Amal**
e-mail: kamalsj@gmail.com

**2008**

# SUPER FUZZY MATRICES AND SUPER FUZZY MODELS FOR SOCIAL SCIENTISTS

**W. B. Vasantha Kandasamy**
**Florentin Smarandache**
**K. Amal**

**2008**

# CONTENTS

# PREFACE

The concept of supermatrix for social scientists was first introduced by Paul Horst. The main purpose of his book was to introduce this concept to social scientists, students, teachers and research workers who lacked mathematical training. He wanted them to be equipped in a branch of mathematics that was increasingly valuable for the analysis of scientific data.

This book introduces the concept of fuzzy super matrices and operations on them. The author has provided only those operations on fuzzy supermatrices that are essential for developing super fuzzy multi expert models. We do not indulge in labourious use of suffixes or superfixes and difficult notations; instead we illustrate the working by simple examples. This book will be highly useful to social scientists who wish to work with multi expert models.

An important feature of this book is its simple approach. Illustrations are given to make the method of approach to the problems easily understandable. Super fuzzy models using Fuzzy Cognitive Maps, Fuzzy Relational maps, Bidirectional Associative Memories and Fuzzy Associative Memories are defined here. Every model is a multi expert model. This book will certainly be a boon not only to social scientists but also to engineers, students, doctors and researchers.

The authors introduce thirteen multi expert models using the notion of fuzzy supermatrices. These models are also described by illustrative examples.

This book has three chapters. In the first chapter we recall some basic concepts about supermatrices and fuzzy matrices. Chapter two introduces the notion of fuzzy supermatrices and their properties. Chapter three introduces many super fuzzy multi expert models.

The authors deeply acknowledge the unflinching support of Dr.K.Kandasamy, Meena and Kama.

W.B.VASANTHA KANDASAMY
FLORENTIN SMARANDACHE
AMAL. K

**Chapter One**

# BASIC CONCEPTS

In this chapter we just recall the definition of supermatrix and some of its basic properties which comprises the section 1. In section 2 fuzzy matrices are introduced.

## 1.1 Supermatrices

The general rectangular or square array of numbers such as

$$A = \begin{bmatrix} 2 & 3 & 1 & 4 \\ -5 & 0 & 7 & -8 \end{bmatrix}, \; B = \begin{bmatrix} 1 & 2 & 3 \\ -4 & 5 & 6 \\ 7 & -8 & 11 \end{bmatrix},$$

$$C = [3, 1, 0, -1, -2] \text{ and } D = \begin{bmatrix} -7/2 \\ 0 \\ \sqrt{2} \\ 5 \\ -41 \end{bmatrix}$$

are known as matrices.

We shall call them as simple matrices [92]. By a simple matrix we mean a matrix each of whose elements are just an ordinary number or a letter that stands for a number. In other words, the elements of a simple matrix are scalars or scalar quantities.

A supermatrix on the other hand is one whose elements are themselves matrices with elements that can be either scalars or other matrices. In general the kind of supermatrices we shall deal with in this book, the matrix elements which have any scalar for their elements. Suppose we have the four matrices;

$$a_{11} = \begin{bmatrix} 2 & -4 \\ 0 & 1 \end{bmatrix}, \; a_{12} = \begin{bmatrix} 0 & 40 \\ 21 & -12 \end{bmatrix}$$

$$a_{21} = \begin{bmatrix} 3 & -1 \\ 5 & 7 \\ -2 & 9 \end{bmatrix} \text{ and } a_{22} = \begin{bmatrix} 4 & 12 \\ -17 & 6 \\ 3 & 11 \end{bmatrix}.$$

One can observe the change in notation $a_{ij}$ denotes a matrix and not a scalar of a matrix ($1 \leq i, j \leq 2$).

Let

$$a = \begin{bmatrix} a_{11} & a_{12} \\ a_{21} & a_{22} \end{bmatrix};$$

we can write out the matrix a in terms of the original matrix elements i.e.,

$$a = \left[\begin{array}{cc|cc} 2 & -4 & 0 & 40 \\ 0 & 1 & 21 & -12 \\ \hline 3 & -1 & 4 & 12 \\ 5 & 7 & -17 & 6 \\ -2 & 9 & 3 & 11 \end{array}\right].$$

Here the elements are divided vertically and horizontally by thin lines. If the lines were not used the matrix a would be read as a simple matrix.

Thus far we have referred to the elements in a supermatrix as matrices as elements. It is perhaps more usual to call the elements of a supermatrix as submatrices. We speak of the submatrices within a supermatrix. Now we proceed on to define the order of a supermatrix.

The order of a supermatrix is defined in the same way as that of a simple matrix. The height of a supermatrix is the number of rows of submatrices in it. The width of a supermatrix is the number of columns of submatrices in it.

All submatrices with in a given row must have the same number of rows. Likewise all submatrices with in a given column must have the same number of columns.

A diagrammatic representation is given by the following figure.

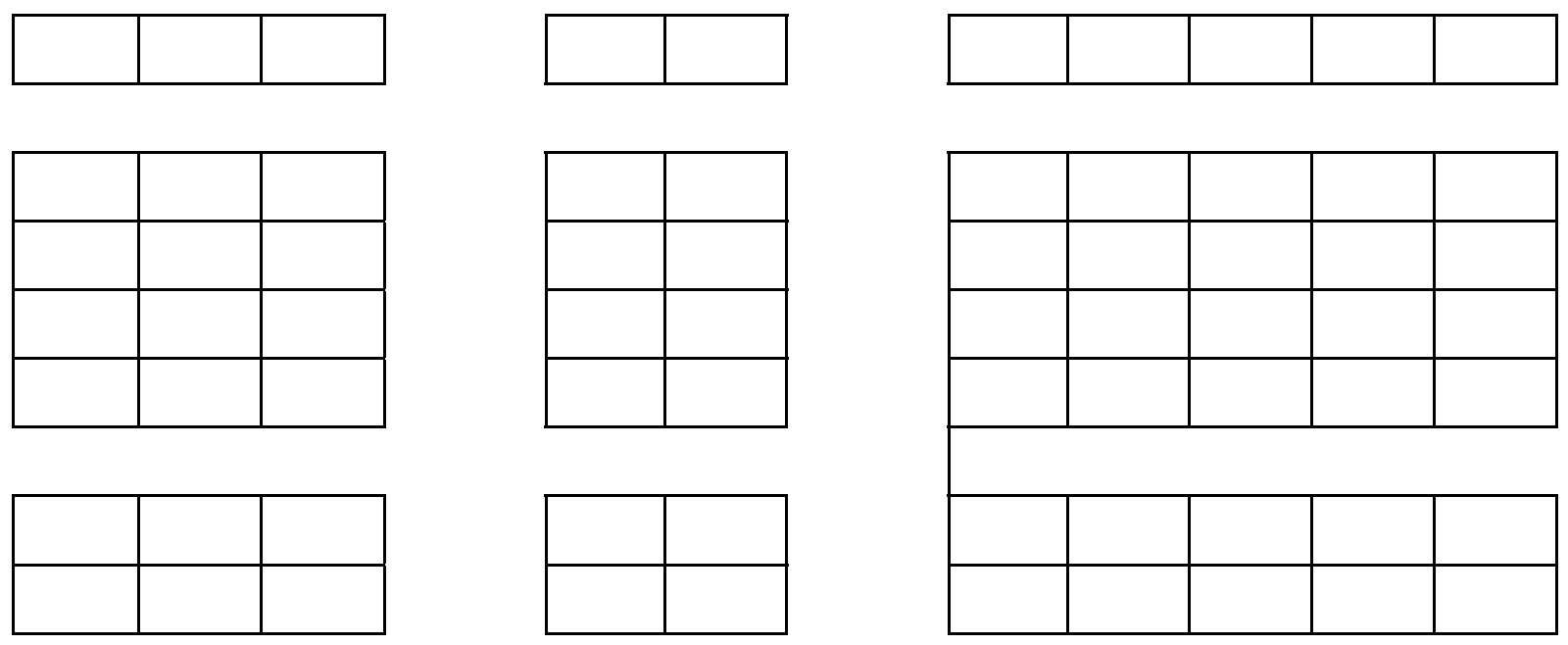

In the first row of rectangles we have one row of a square for each rectangle; in the second row of rectangles we have four rows of squares for each rectangle and in the third row of rectangles we have two rows of squares for each rectangle. Similarly for the first column of rectangles three columns of squares for each rectangle. For the second column of rectangles we have two column of squares for each rectangle, and for the third column of rectangles we have five columns of squares for each rectangle.

Thus we have for this supermatrix 3 rows and 3 columns.

One thing should now be clear from the definition of a supermatrix. The super order of a supermatrix tells us nothing about the simple order of the matrix from which it was obtained by partitioning. Furthermore, the order of supermatrix tells us nothing about the orders of the submatrices within that supermatrix.

Now we illustrate the number of rows and columns of a supermatrix.

***Example 1.1.1:*** Let

$$a = \left[\begin{array}{c|cccc} 3 & 3 & 0 & 1 & 4 \\ -1 & 2 & 1 & -1 & 6 \\ \hline 0 & 3 & 4 & 5 & 6 \\ 1 & 7 & 8 & -9 & 0 \\ 2 & 1 & 2 & 3 & -4 \end{array}\right].$$

a is a supermatrix with two rows and two columns.

Now we proceed on to define the notion of partitioned matrices. It is always possible to construct a supermatrix from any simple matrix that is not a scalar quantity.

The supermatrix can be constructed from a simple matrix this process of constructing supermatrix is called the partitioning.

A simple matrix can be partitioned by dividing or separating the matrix between certain specified rows, or the procedure may be reversed. The division may be made first between rows and then between columns.

We illustrate this by a simple example.

***Example 1.1.2:*** Let

$$A = \begin{bmatrix} 3 & 0 & 1 & 1 & 2 & 0 \\ 1 & 0 & 0 & 3 & 5 & 2 \\ 5 & -1 & 6 & 7 & 8 & 4 \\ 0 & 9 & 1 & 2 & 0 & -1 \\ 2 & 5 & 2 & 3 & 4 & 6 \\ 1 & 6 & 1 & 2 & 3 & 9 \end{bmatrix}$$

is a $6 \times 6$ simple matrix with real numbers as elements.

$$A_1 = \left[\begin{array}{cc|cccc} 3 & 0 & 1 & 1 & 2 & 0 \\ 1 & 0 & 0 & 3 & 5 & 2 \\ 5 & -1 & 6 & 7 & 8 & 4 \\ 0 & 9 & 1 & 2 & 0 & -1 \\ 2 & 5 & 2 & 3 & 4 & 6 \\ 1 & 6 & 1 & 2 & 3 & 9 \end{array}\right].$$

Now let us draw a thin line between the 2nd and 3rd columns.

This gives us the matrix $A_1$. Actually $A_1$ may be regarded as a supermatrix with two matrix elements forming one row and two columns.

Now consider

$$A_2 = \left[\begin{array}{cccccc} 3 & 0 & 1 & 1 & 2 & 0 \\ 1 & 0 & 0 & 3 & 5 & 2 \\ 5 & -1 & 6 & 7 & 8 & 4 \\ 0 & 9 & 1 & 2 & 0 & -1 \\ \hline 2 & 5 & 2 & 3 & 4 & 6 \\ 1 & 6 & 1 & 2 & 3 & 9 \end{array}\right]$$

Draw a thin line between the rows 4 and 5 which gives us the new matrix $A_2$. $A_2$ is a supermatrix with two rows and one column.

Now consider the matrix

$$A_3 = \left[\begin{array}{cc|cccc} 3 & 0 & 1 & 1 & 2 & 0 \\ 1 & 0 & 0 & 3 & 5 & 2 \\ 5 & -1 & 6 & 7 & 8 & 4 \\ 0 & 9 & 1 & 2 & 0 & -1 \\ \hline 2 & 5 & 2 & 3 & 4 & 6 \\ 1 & 6 & 1 & 2 & 3 & 9 \end{array}\right],$$

$A_3$ is now a second order supermatrix with two rows and two columns. We can simply write $A_3$ as

$$\begin{bmatrix} a_{11} & a_{12} \\ a_{21} & a_{22} \end{bmatrix}$$

where

$$a_{11} = \begin{bmatrix} 3 & 0 \\ 1 & 0 \\ 5 & -1 \\ 0 & 9 \end{bmatrix},$$

$$a_{12} = \begin{bmatrix} 1 & 1 & 2 & 0 \\ 0 & 3 & 5 & 2 \\ 6 & 7 & 8 & 4 \\ 1 & 2 & 0 & -1 \end{bmatrix},$$

$$a_{21} = \begin{bmatrix} 2 & 5 \\ 1 & 6 \end{bmatrix} \text{ and } a_{22} = \begin{bmatrix} 2 & 3 & 4 & 6 \\ 1 & 2 & 3 & 9 \end{bmatrix}.$$

The elements now are the submatrices defined as $a_{11}$, $a_{12}$, $a_{21}$ and $a_{22}$ and therefore $A_3$ is in terms of letters.

According to the methods we have illustrated a simple matrix can be partitioned to obtain a supermatrix in any way that happens to suit our purposes.

The natural order of a supermatrix is usually determined by the natural order of the corresponding simple matrix. Further more we are not usually concerned with natural order of the submatrices within a supermatrix.

Now we proceed on to recall the notion of symmetric partition, for more information about these concepts please refer [92]. By a symmetric partitioning of a matrix we mean that the rows and columns are partitioned in exactly the same way. If the matrix is partitioned between the first and second column and between the third and fourth column, then to be symmetrically partitioning, it must also be partitioned between the first and second rows and third and fourth rows. According to this rule of symmetric partitioning only square simple matrix can be symmetrically partitioned. We give an example of a symmetrically partitioned matrix $a_s$,

***Example 1.1.3:*** Let

$$a_s = \left[\begin{array}{c|cc|c} 2 & 3 & 4 & 1 \\ \hline 5 & 6 & 9 & 2 \\ 0 & 6 & 1 & 9 \\ \hline 5 & 1 & 1 & 5 \end{array}\right].$$

Here we see that the matrix has been partitioned between columns one and two and three and four. It has also been partitioned between rows one and two and rows three and four.

Now we just recall from [92] the method of symmetric partitioning of a symmetric simple matrix.

***Example 1.1.4:*** Let us take a fourth order symmetric matrix and partition it between the second and third rows and also between the second and third columns.

$$a = \left[\begin{array}{cc|cc} 4 & 3 & 2 & 7 \\ 3 & 6 & 1 & 4 \\ \hline 2 & 1 & 5 & 2 \\ 7 & 4 & 2 & 7 \end{array}\right].$$

We can represent this matrix as a supermatrix with letter elements.

$$a_{11} = \begin{bmatrix} 4 & 3 \\ 3 & 6 \end{bmatrix}, \; a_{12} = \begin{bmatrix} 2 & 7 \\ 1 & 4 \end{bmatrix}$$

$$a_{21} = \begin{bmatrix} 2 & 1 \\ 7 & 4 \end{bmatrix} \text{ and } a_{22} = \begin{bmatrix} 5 & 2 \\ 2 & 7 \end{bmatrix},$$

so that

$$a = \begin{bmatrix} a_{11} & a_{12} \\ a_{21} & a_{22} \end{bmatrix}.$$

The diagonal elements of the supermatrix a are $a_{11}$ and $a_{22}$. We also observe the matrices $a_{11}$ and $a_{22}$ are also symmetric matrices.

The non diagonal elements of this supermatrix a are the matrices $a_{12}$ and $a_{21}$. Clearly $a_{21}$ is the transpose of $a_{12}$.

The simple rule about the matrix element of a symmetrically partitioned symmetric simple matrix are (1) The diagonal submatrices of the supermatrix are all symmetric matrices. (2) The matrix elements below the diagonal are the transposes of the corresponding elements above the diagonal.

The forth order supermatrix obtained from a symmetric partitioning of a symmetric simple matrix a is as follows.

$$a = \begin{bmatrix} a_{11} & a_{12} & a_{13} & a_{14} \\ a'_{12} & a_{22} & a_{23} & a_{24} \\ a'_{13} & a'_{23} & a_{33} & a_{34} \\ a'_{14} & a'_{24} & a'_{34} & a_{44} \end{bmatrix}.$$

How to express that a symmetric matrix has been symmetrically partitioned (i) $a_{11}$ and $a^t_{11}$ are equal. (ii) $a^t_{ij}$ $(i \neq j)$; $a^t_{ij} = a_{ji}$ and $a^t_{ji} = a_{ij}$. Thus the general expression for a symmetrically partitioned symmetric matrix;

$$a = \begin{bmatrix} a_{11} & a_{12} & \dots & a_{1n} \\ a'_{12} & a_{22} & \dots & a_{2n} \\ \vdots & \vdots & & \vdots \\ a'_{1n} & a'_{2n} & \dots & a_{nn} \end{bmatrix}.$$

If we want to indicate a symmetrically partitioned simple diagonal matrix we would write

$$D = \begin{bmatrix} D_1 & 0 & \dots & 0 \\ 0' & D_2 & \dots & 0 \\ & & & \\ 0' & 0' & \dots & D_n \end{bmatrix}$$

0' only represents the order is reversed or transformed. We denote $a_{ij}^t$ = a'$_{ij}$ just the ' means the transpose.

D will be referred to as the super diagonal matrix. The identity matrix

$$I = \begin{bmatrix} I_s & 0 & 0 \\ 0 & I_t & 0 \\ 0 & 0 & I_r \end{bmatrix}$$

s, t and r denote the number of rows and columns of the first second and third identity matrices respectively (zeros denote matrices with zero as all entries).

***Example 1.1.5:*** We just illustrate a general super diagonal matrix d;

$$d = \left[\begin{array}{ccc|cc} 3 & 1 & 2 & 0 & 0 \\ 5 & 6 & 0 & 0 & 0 \\ \hline 0 & 0 & 0 & 2 & 5 \\ 0 & 0 & 0 & -1 & 3 \\ 0 & 0 & 0 & 9 & 10 \end{array}\right]$$

$$\text{i.e., } d = \begin{bmatrix} m_1 & 0 \\ 0 & m_2 \end{bmatrix}.$$

An example of a super diagonal matrix with vector elements is given, which can be useful in experimental designs.

***Example 1.1.6:*** Let

$$\left[\begin{array}{c|c|c|c} 1 & 0 & 0 & 0 \\ 1 & 0 & 0 & 0 \\ 1 & 0 & 0 & 0 \\ \hline 0 & 1 & 0 & 0 \\ 0 & 1 & 0 & 0 \\ \hline 0 & 0 & 1 & 0 \\ 0 & 0 & 1 & 0 \\ 0 & 0 & 1 & 0 \\ \hline 0 & 0 & 1 & 0 \\ \hline 0 & 0 & 0 & 1 \\ 0 & 0 & 0 & 1 \\ 0 & 0 & 0 & 1 \\ 0 & 0 & 0 & 1 \end{array}\right].$$

Here the diagonal elements are only column unit vectors. In case of supermatrix [92] has defined the notion of partial triangular matrix as a supermatrix.

***Example 1.1.7:*** Let

$$u = \left[\begin{array}{ccc|cc} 2 & 1 & 1 & 3 & 2 \\ 0 & 5 & 2 & 1 & 1 \\ 0 & 0 & 1 & 0 & 2 \end{array}\right]$$

u is a partial upper triangular supermatrix.

***Example 1.1.8:*** Let

$$L = \left[\begin{array}{ccccc} 5 & 0 & 0 & 0 & 0 \\ 7 & 2 & 0 & 0 & 0 \\ 1 & 2 & 3 & 0 & 0 \\ 4 & 5 & 6 & 7 & 0 \\ 1 & 2 & 5 & 2 & 6 \\ \hline 1 & 2 & 3 & 4 & 5 \\ 0 & 1 & 0 & 1 & 0 \end{array}\right];$$

L is partial upper triangular matrix partitioned as a supermatrix.

Thus $T = \left[\begin{array}{c} T \\ \hline a' \end{array}\right]$ where T is the lower triangular submatrix, with

$$T = \begin{bmatrix} 5 & 0 & 0 & 0 & 0 \\ 7 & 2 & 0 & 0 & 0 \\ 1 & 2 & 3 & 0 & 0 \\ 4 & 5 & 6 & 7 & 0 \\ 1 & 2 & 5 & 2 & 6 \end{bmatrix} \text{ and } a' = \begin{bmatrix} 1 & 2 & 3 & 4 & 5 \\ 0 & 1 & 0 & 1 & 0 \end{bmatrix}.$$

We proceed on to define the notion of supervectors i.e., Type I column supervector. A simple vector is a vector each of whose elements is a scalar. It is nice to see the number of different types of supervectors given by [92].

***Example 1.1.9:*** Let

$$v = \left[\begin{array}{c} 1 \\ 3 \\ 4 \\ \hline 5 \\ 7 \end{array}\right].$$

This is a type I i.e., type one column supervector.

$$v = \begin{bmatrix} v_1 \\ v_2 \\ \vdots \\ v_n \end{bmatrix}$$

where each $v_i$ is a column subvectors of the column vector v.

Type I row supervector is given by the following example.

***Example 1.1.10:*** $v^1$ = [2 3 1 | 5 7 8 4] is a type I row supervector. i.e., v' = [$v'_1$, $v'_2$, …, $v'_n$] where each $v'_i$ is a row subvector; $1 \leq i \leq n$.

Next we recall the definition of type II supervectors.

Type II column supervectors.

**DEFINITION 1.1.1:** *Let*

$$a = \begin{bmatrix} a_{11} & a_{12} & \dots & a_{1m} \\ a_{21} & a_{22} & \dots & a_{2m} \\ \dots & \dots & \dots & \dots \\ a_{n1} & a_{n2} & \dots & a_{nm} \end{bmatrix}$$

$$\begin{aligned} a_1^1 &= [a_{11} \dots a_{1m}] \\ a_2^1 &= [a_{21} \dots a_{2m}] \\ &\dots \\ a_n^1 &= [a_{n1} \dots a_{nm}] \end{aligned}$$

*i.e.,*
$$a = \begin{bmatrix} a_1^1 \\ a_2^1 \\ \vdots \\ a_n^1 \end{bmatrix}_m$$

*is defined to be the type II column supervector.*
*Similarly if*

$$a^1 = \begin{bmatrix} a_{11} \\ a_{21} \\ \vdots \\ a_{n1} \end{bmatrix}, \ a^2 = \begin{bmatrix} a_{12} \\ a_{22} \\ \vdots \\ a_{n2} \end{bmatrix}, \dots, \ a^m = \begin{bmatrix} a_{1m} \\ a_{2m} \\ \vdots \\ a_{nm} \end{bmatrix}.$$

*Hence now* $a = [a^1\ a^2 \dots a^m]_n$ *is defined to be the type II row supervector.*

*Clearly*

$$a = \begin{bmatrix} a_1^1 \\ a_2^1 \\ \vdots \\ a_n^1 \end{bmatrix}_m = [a^1\ a^2 \ldots a^m]_n$$

*the equality of supermatrices.*

***Example 1.1.11:*** Let

$$A = \begin{bmatrix} 3 & 6 & 0 & 4 & 5 \\ 2 & 1 & 6 & 3 & 0 \\ 1 & 1 & 1 & 2 & 1 \\ 0 & 1 & 0 & 1 & 0 \\ 2 & 0 & 1 & 2 & 1 \end{bmatrix}$$

be a simple matrix. Let a and b the supermatrix made from A.

$$a = \left[\begin{array}{ccc|cc} 3 & 6 & 0 & 4 & 5 \\ 2 & 1 & 6 & 3 & 0 \\ 1 & 1 & 1 & 2 & 1 \\ \hline 0 & 1 & 0 & 1 & 0 \\ 2 & 0 & 1 & 2 & 1 \end{array}\right]$$

where

$$a_{11} = \begin{bmatrix} 3 & 6 & 0 \\ 2 & 1 & 6 \\ 1 & 1 & 1 \end{bmatrix},\ a_{12} = \begin{bmatrix} 4 & 5 \\ 3 & 0 \\ 2 & 1 \end{bmatrix},$$

$$a_{21} = \begin{bmatrix} 0 & 1 & 0 \\ 2 & 0 & 1 \end{bmatrix} \text{ and } a_{22} = \begin{bmatrix} 1 & 0 \\ 2 & 1 \end{bmatrix}.$$

i.e.,
$$a = \begin{bmatrix} a_{11} & a_{12} \\ a_{21} & a_{22} \end{bmatrix}.$$

$$b = \left[\begin{array}{cccc|c} 3 & 6 & 0 & 4 & 5 \\ 2 & 1 & 6 & 3 & 0 \\ 1 & 1 & 1 & 2 & 1 \\ 0 & 1 & 0 & 1 & 0 \\ \hline 2 & 0 & 1 & 2 & 1 \end{array}\right] = \begin{bmatrix} b_{11} & b_{12} \\ b_{21} & b_{22} \end{bmatrix}$$

where

$$b_{11} = \begin{bmatrix} 3 & 6 & 0 & 4 \\ 2 & 1 & 6 & 3 \\ 1 & 1 & 1 & 2 \\ 0 & 1 & 0 & 1 \end{bmatrix}, \; b_{12} = \begin{bmatrix} 5 \\ 0 \\ 1 \\ 0 \end{bmatrix},$$

$b_{21}$ = [2 0 1 2 ] and $b_{22}$ = [1].

$$a = \left[\begin{array}{ccc|cc} 3 & 6 & 0 & 4 & 5 \\ 2 & 1 & 6 & 3 & 0 \\ 1 & 1 & 1 & 2 & 1 \\ \hline 0 & 1 & 0 & 1 & 0 \\ 2 & 0 & 1 & 2 & 1 \end{array}\right]$$

and

$$b = \left[\begin{array}{cccc|c} 3 & 6 & 0 & 4 & 5 \\ 2 & 1 & 6 & 3 & 0 \\ 1 & 1 & 1 & 2 & 1 \\ 0 & 1 & 0 & 1 & 0 \\ \hline 2 & 0 & 1 & 2 & 1 \end{array}\right].$$

We see that the corresponding scalar elements for matrix a and matrix b are identical. Thus two supermatrices are equal if and only if their corresponding simple forms are equal.

Now we give examples of type III supervector for more refer [92].

***Example 1.1.12:***

$$a = \left[\begin{array}{ccc|cc} 3 & 2 & 1 & 7 & 8 \\ 0 & 2 & 1 & 6 & 9 \\ 0 & 0 & 5 & 1 & 2 \end{array}\right] = [T' \mid a']$$

and

$$b = \left[\begin{array}{ccc} 2 & 0 & 0 \\ 9 & 4 & 0 \\ 8 & 3 & 6 \\ \hline 5 & 2 & 9 \\ 4 & 7 & 3 \end{array}\right] = \left[\begin{array}{c} T \\ \hline b' \end{array}\right]$$

are type III supervectors.

One interesting and common example of a type III supervector is a prediction data matrix having both predictor and criterion attributes.

The next interesting notion about supermatrix is its transpose. First we illustrate this by an example before we give the general case.

***Example 1.1.13:*** Let

$$a = \left[\begin{array}{ccc|cc} 2 & 1 & 3 & 5 & 6 \\ 0 & 2 & 0 & 1 & 1 \\ 1 & 1 & 1 & 0 & 2 \\ \hline 2 & 2 & 0 & 1 & 1 \\ 5 & 6 & 1 & 0 & 1 \\ \hline 2 & 0 & 0 & 0 & 4 \\ 1 & 0 & 1 & 1 & 5 \end{array}\right]$$

$$= \begin{bmatrix} a_{11} & a_{12} \\ a_{21} & a_{22} \\ a_{31} & a_{32} \end{bmatrix}$$

where

$$a_{11} = \begin{bmatrix} 2 & 1 & 3 \\ 0 & 2 & 0 \\ 1 & 1 & 1 \end{bmatrix}, \; a_{12} = \begin{bmatrix} 5 & 6 \\ 1 & 1 \\ 0 & 2 \end{bmatrix},$$

$$a_{21} = \begin{bmatrix} 2 & 2 & 0 \\ 5 & 6 & 1 \end{bmatrix}, \; a_{22} = \begin{bmatrix} 1 & 1 \\ 0 & 1 \end{bmatrix},$$

$$a_{31} = \begin{bmatrix} 2 & 0 & 0 \\ 1 & 0 & 1 \end{bmatrix} \text{ and } a_{32} = \begin{bmatrix} 0 & 4 \\ 1 & 5 \end{bmatrix}.$$

The transpose of a

$$a^t = a' = \left[\begin{array}{ccc|cc|cc} 2 & 0 & 1 & 2 & 5 & 2 & 1 \\ 1 & 2 & 1 & 2 & 6 & 0 & 0 \\ 3 & 0 & 1 & 0 & 1 & 0 & 1 \\ \hline 5 & 1 & 0 & 1 & 0 & 0 & 1 \\ 6 & 1 & 2 & 1 & 1 & 4 & 5 \end{array}\right].$$

Let us consider the transposes of $a_{11}$, $a_{12}$, $a_{21}$, $a_{22}$, $a_{31}$ and $a_{32}$.

$$a'_{11} = a^t_{11} = \begin{bmatrix} 2 & 0 & 1 \\ 1 & 2 & 1 \\ 3 & 0 & 1 \end{bmatrix}$$

$$a'_{12} = a^t_{12} = \begin{bmatrix} 5 & 1 & 0 \\ 6 & 1 & 2 \end{bmatrix}$$

$$a'_{21} = a^t_{21} = \begin{bmatrix} 2 & 5 \\ 2 & 6 \\ 0 & 1 \end{bmatrix}$$

$$a'_{31} = a^t_{31} = \begin{bmatrix} 2 & 1 \\ 0 & 0 \\ 0 & 1 \end{bmatrix}$$

$$a'_{22} = a^t_{22} = \begin{bmatrix} 1 & 0 \\ 1 & 1 \end{bmatrix}$$

$$a'_{32} = a^t_{32} = \begin{bmatrix} 0 & 1 \\ 4 & 5 \end{bmatrix}.$$

$$a' = \begin{bmatrix} a'_{11} & a'_{21} & a'_{31} \\ a'_{12} & a'_{22} & a'_{32} \end{bmatrix}.$$

Now we describe the general case. Let

$$a = \begin{bmatrix} a_{11} & a_{12} & \cdots & a_{1m} \\ a_{21} & a_{22} & \cdots & a_{2m} \\ \vdots & \vdots & & \vdots \\ a_{n1} & a_{n2} & \cdots & a_{nm} \end{bmatrix}$$

be a n × m supermatrix. The transpose of the supermatrix a denoted by

$$a' = \begin{bmatrix} a'_{11} & a'_{21} & \cdots & a'_{n1} \\ a'_{12} & a'_{22} & \cdots & a'_{n2} \\ \vdots & \vdots & & \vdots \\ a'_{1m} & a'_{2m} & \cdots & a'_{nm} \end{bmatrix}.$$

a' is a m by n supermatrix obtained by taking the transpose of each element i.e., the submatrices of a.

Now we will find the transpose of a symmetrically partitioned symmetric simple matrix. Let a be the symmetrically partitioned symmetric simple matrix.

Let a be a m × m symmetric supermatrix i.e.,

$$a = \begin{bmatrix} a_{11} & a_{21} & \cdots & a_{m1} \\ a_{12} & a_{22} & \cdots & a_{m2} \\ \vdots & \vdots & & \vdots \\ a_{1m} & a_{2m} & \cdots & a_{mm} \end{bmatrix}$$

the transpose of the supermatrix is given by a'

$$a' = \begin{bmatrix} a'_{11} & (a'_{12})' & \cdots & (a'_{1m})' \\ a'_{12} & a'_{22} & \cdots & (a'_{2m})' \\ \vdots & \vdots & & \vdots \\ a'_{1m} & a'_{2m} & \cdots & a'_{mm} \end{bmatrix}$$

The diagonal matrix $a_{11}$ are symmetric matrices so are unaltered by transposition. Hence

$$a'_{11} = a_{11},\ a'_{22} = a_{22},\ \ldots,\ a'_{mm} = a_{mm}.$$

Recall also the transpose of a transpose is the original matrix. Therefore

$$(a'_{12})' = a_{12},\ (a'_{13})' = a_{13},\ \ldots,\ (a'_{ij})' = a_{ij}.$$

Thus the transpose of supermatrix constructed by symmetrically partitioned symmetric simple matrix a of a' is given by

$$a' = \begin{bmatrix} a_{11} & a_{12} & \cdots & a_{1m} \\ a'_{21} & a_{22} & \cdots & a_{2m} \\ \vdots & \vdots & & \vdots \\ a'_{1m} & a'_{2m} & \cdots & a_{mm} \end{bmatrix}.$$

Thus a = a'.
Similarly transpose of a symmetrically partitioned diagonal matrix is simply the original diagonal supermatrix itself;

i.e., if

$$D = \begin{bmatrix} d_1 & & & \\ & d_2 & & \\ & & \ddots & \\ & & & d_n \end{bmatrix}$$

$$D' = \begin{bmatrix} d_1' & & & \\ & d_2' & & \\ & & \ddots & \\ & & & d_n' \end{bmatrix}$$

$d'_1 = d_1$, $d'_2 = d_2$ etc. Thus D = D'.

Now we see the transpose of a type I supervector.

***Example 1.1.14:*** Let

$$V = \left[\begin{array}{c} 3 \\ 1 \\ 2 \\ \hline 4 \\ 5 \\ 7 \\ \hline 5 \\ 1 \end{array}\right]$$

The transpose of V denoted by V' or $V^t$ is

$$V' = [3\ 1\ 2\ |\ 4\ 5\ 7\ |\ 5\ 1].$$

If

$$V = \begin{bmatrix} v_1 \\ v_2 \\ v_3 \end{bmatrix}$$

where

$$v_1 = \begin{bmatrix} 3 \\ 1 \\ 2 \end{bmatrix}, v_2 = \begin{bmatrix} 4 \\ 5 \\ 7 \end{bmatrix} \text{ and } v_3 = \begin{bmatrix} 5 \\ 1 \end{bmatrix}$$

$$V' = [v'_1\ v'_2\ v'_3].$$

Thus if

$$V = \begin{bmatrix} v_1 \\ v_2 \\ \vdots \\ v_n \end{bmatrix}$$

then

$$V' = [v'_1\ v'_2 \ldots v'_n].$$

***Example 1.1.15:*** Let

$$t = \left[\begin{array}{cccc|cc} 3 & 0 & 1 & 1 & 5 & 2 \\ 4 & 2 & 0 & 1 & 3 & 5 \\ 1 & 0 & 1 & 0 & 1 & 6 \end{array}\right]$$

= [T | a ]. The transpose of t

$$\text{i.e., } t' = \left[\begin{array}{ccc} 3 & 4 & 1 \\ 0 & 2 & 0 \\ 1 & 0 & 1 \\ 1 & 1 & 0 \\ \hline 5 & 3 & 1 \\ 2 & 5 & 6 \end{array}\right] = \left[\begin{array}{c} T' \\ \hline a' \end{array}\right].$$

The addition of supermatrices may not be always be defined.

***Example 1.1.16:*** For instance let

$$a = \begin{bmatrix} a_{11} & a_{12} \\ a_{21} & a_{22} \end{bmatrix}$$

and

$$b = \begin{bmatrix} b_{11} & b_{12} \\ b_{21} & b_{22} \end{bmatrix}$$

where

$$a_{11} = \begin{bmatrix} 3 & 0 \\ 1 & 2 \end{bmatrix}, \quad a_{12} = \begin{bmatrix} 1 \\ 7 \end{bmatrix}$$

$$a_{21} = [4 \;\; 3], \qquad a_{22} = [6].$$

$$b_{11} = [2], \qquad b_{12} = [1 \;\; 3]$$

$$b_{21} = \begin{bmatrix} 5 \\ 2 \end{bmatrix} \text{ and } b_{22} = \begin{bmatrix} 4 & 1 \\ 0 & 2 \end{bmatrix}.$$

It is clear both a and b are second order square supermatrices but here we cannot add together the corresponding matrix elements of a and b because the submatrices do not have the same order.

Now we proceed onto recall the definition of minor product of two supervectors.

Suppose

$$v_a = \begin{bmatrix} v_{a_1} \\ v_{a_2} \\ \vdots \\ v_{a_n} \end{bmatrix} \text{ and } v_b = \begin{bmatrix} v_{b_1} \\ v_{b_2} \\ \vdots \\ v_{b_n} \end{bmatrix}.$$

The minor product of these two supervectors $v_a$ and $v_b$ is given by

$$= v'_a v_b = \left[v'_{a_1} v'_{a_2} \cdots v'_{a_n}\right] \begin{bmatrix} v_{b_1} \\ v_{b_2} \\ \vdots \\ v_{b_n} \end{bmatrix}$$

$$= v'_{a_1} v_{b_1} + v'_{a_2} v_{b_2} + \cdots + v'_{a_n} v_{b_n} .$$

We illustrate this by the following example.

***Example 1.1.17:*** Let $V_a$ and $V_b$ be two type I supervectors where

$$V_a = \begin{bmatrix} v_{a_1} \\ v_{a_2} \\ v_{a_3} \end{bmatrix}$$

with

$$v_{a_1} = \begin{bmatrix} 0 \\ 1 \\ 2 \end{bmatrix},\ v_{a_2} = \begin{bmatrix} 4 \\ 0 \\ 1 \\ -1 \end{bmatrix} \text{ and } v_{a_3} = \begin{bmatrix} 1 \\ 2 \end{bmatrix}.$$

Let

$$V_b = \begin{bmatrix} v_{b_1} \\ v_{b_2} \\ v_{b_3} \end{bmatrix}$$

where

$$v_{b_1} = \begin{bmatrix} 1 \\ -1 \\ 0 \end{bmatrix},\ v_{b_2} = \begin{bmatrix} -4 \\ 1 \\ 2 \\ 0 \end{bmatrix} \text{ and } v_{b_3} = \begin{bmatrix} -1 \\ 1 \end{bmatrix}.$$

$$V'_a V_b = \left[ v'_{a_1}\, v'_{a_2}\, v'_{a_3} \right] \begin{bmatrix} v_{b_1} \\ v_{b_2} \\ v_{b_3} \end{bmatrix}$$

$$= \quad v'_{a_1} v_{b_1} + v'_{a_2} v_{b_2} + \cdots + v'_{a_n} v_{b_n}$$

$$= \quad \begin{bmatrix} 0 & 1 & 2 \end{bmatrix} \begin{bmatrix} 1 \\ -1 \\ 0 \end{bmatrix} + \begin{bmatrix} 4 & 0 & 1 & -1 \end{bmatrix} \begin{bmatrix} -4 \\ 1 \\ 2 \\ 0 \end{bmatrix} + \begin{bmatrix} 1 & 2 \end{bmatrix} \begin{bmatrix} -1 \\ 1 \end{bmatrix}$$

$= \quad -1 + (-16+2) + (-1+2)$
$= \quad -1 - 16 + 2 - 1 + 2$
$= \quad -14.$

It is easily proved $V'_a V_b = V'_b V_a$.

Now we proceed on to recall the definition of major product of type I supervectors.

Suppose

$$V_a = \begin{bmatrix} v_{a_1} \\ v_{a_2} \\ \vdots \\ v_{a_n} \end{bmatrix} \text{ and } V_b = \begin{bmatrix} v_{b_1} \\ v_{b_2} \\ \vdots \\ v_{b_m} \end{bmatrix}$$

be any two supervectors of type I. The major product is defined as

$$V_a V'_b = \begin{bmatrix} v_{a_1} \\ v_{a_2} \\ \vdots \\ v_{a_n} \end{bmatrix} \cdot \left[ v'_{b_1} v'_{b_2} \cdots v'_{b_m} \right]$$

$$= \begin{bmatrix} v_{a_1} v'_{b_1} & v_{a_1} v'_{b_2} & \cdots & v_{a_1} v'_{b_m} \\ v_{a_2} v'_{b_1} & v_{a_2} v'_{b_2} & \cdots & v_{a_2} v'_{b_m} \\ \vdots & & & \\ v_{a_n} v'_{b_1} & v_{a_n} v'_{b_2} & \cdots & v_{a_n} v'_{b_m} \end{bmatrix}.$$

Now we illustrate this by the following example.

***Example 1.1.18:*** Let

$$V_a = \begin{bmatrix} v_{a_1} \\ v_{a_2} \\ v_{a_3} \end{bmatrix} \text{ and } V_b = \begin{bmatrix} v_{b_1} \\ v_{b_2} \\ v_{b_3} \\ v_{b_4} \end{bmatrix}$$

where

$$v_{a_1} = [2],\ v_{a_2} = \begin{bmatrix} 1 \\ -1 \end{bmatrix} \text{ and } v_{a_3} = \begin{bmatrix} 1 \\ 2 \\ 0 \end{bmatrix}$$

and

$$v_{b_1} = \begin{bmatrix} 3 \\ 1 \\ 2 \end{bmatrix},\ v_{b_2} = \begin{bmatrix} 1 \\ 2 \end{bmatrix},\ v_{b_3} = \begin{bmatrix} 3 \\ 4 \\ -1 \\ 0 \end{bmatrix} \text{ and } v_{b_4} = [5].$$

$$V_a V'_b = \begin{bmatrix} 2 \\ \hline 1 \\ -1 \\ \hline 1 \\ 2 \\ 0 \end{bmatrix} \left[ 3 \quad 1 \quad 2 \mid 1 \quad 2 \mid 3 \quad 4 \quad -1 \quad 0 \mid 5 \right]$$

$$= \left[\begin{array}{c|c|c|c} [2][3 \quad 1 \quad 2] & [2][1 \quad 2] & [2][3 \quad 4 \quad -1 \quad 0] & [2][5] \\ \hline \begin{bmatrix} 1 \\ -1 \end{bmatrix}[3 \quad 1 \quad 2] & \begin{bmatrix} 1 \\ -1 \end{bmatrix}[1 \quad 2] & \begin{bmatrix} 1 \\ -1 \end{bmatrix}[3 \quad 4 \quad -1 \quad 0] & \begin{bmatrix} 1 \\ -1 \end{bmatrix}[5] \\ \hline \begin{bmatrix} 1 \\ 2 \\ 0 \end{bmatrix}[3 \quad 1 \quad 2] & \begin{bmatrix} 1 \\ 2 \\ 0 \end{bmatrix}[1 \quad 2] & \begin{bmatrix} 1 \\ 2 \\ 0 \end{bmatrix}[3 \quad 4 \quad -1 \quad 0] & \begin{bmatrix} 1 \\ 2 \\ 0 \end{bmatrix}[5] \end{array}\right]$$

$$= \left[\begin{array}{ccc|cc|cccc|c} 6 & 2 & 4 & 2 & 4 & 6 & 4 & -2 & 0 & 10 \\ \hline 3 & 1 & 2 & 1 & 2 & 3 & 4 & -1 & 0 & 5 \\ -3 & -1 & -2 & -1 & -2 & -3 & -4 & 1 & 0 & -5 \\ \hline 3 & 1 & 2 & 1 & 2 & 3 & 4 & -1 & 0 & 5 \\ 6 & 2 & 4 & 2 & 4 & 6 & 8 & -2 & 0 & 10 \\ 0 & 0 & 0 & 0 & 0 & 0 & 0 & 0 & 0 & 0 \end{array}\right].$$

We leave it for the reader to verify that $(V_a V'_b)' = V_b V'_a$.

***Example 1.1.19:*** We just recall if

$$v = \begin{bmatrix} 3 \\ 4 \\ 7 \end{bmatrix}$$

is a column vector and v' the transpose of v is a row vector then we have

$$v'v = [3 \quad 4 \quad 7] \begin{bmatrix} 3 \\ 4 \\ 7 \end{bmatrix}$$
$$= 3^2 + 4^2 + 7^2 = 74.$$

Thus if

$$V'_x = [x_1 \; x_2 \; \ldots \; x_n]$$

$$V'_x V_x = [x_1\ x_2\ \ldots\ x_n] \begin{bmatrix} x_1 \\ x_2 \\ \vdots \\ x_n \end{bmatrix}$$

$$= x_1^2 + x_2^2 + \ldots + x_n^2 .$$

Also

$$[1\ 1\ \ldots\ 1\ ] \begin{bmatrix} x_1 \\ x_2 \\ \vdots \\ x_n \end{bmatrix} = [x_1 + x_2 + \ldots + x_n]$$

and

$$[x_1\ x_2\ \ldots\ x_n] \begin{bmatrix} 1 \\ 1 \\ \vdots \\ 1 \end{bmatrix} = [x_1 + x_2 + \ldots + x_n];$$

$$\text{i.e., } 1'v_x = v'_x 1 = \sum x_i$$

where

$$v_x = \begin{bmatrix} x_1 \\ x_2 \\ \vdots \\ x_n \end{bmatrix}$$

and

$$\sum x_i = x_1 + x_2 + \ldots + x_n.$$

We have the following types of products defined.

***Example 1.1.20:*** We have

$$[0\ 1\ 0\ 0]\begin{bmatrix}0\\1\\0\\0\end{bmatrix} = 1,$$

$$[0\ 1\ 0\ 0]\begin{bmatrix}1\\0\\0\\0\end{bmatrix} = 0,$$

$$[0\ 1\ 0\ 0]\begin{bmatrix}1\\1\\1\\1\end{bmatrix} = 1$$

and

$$\begin{bmatrix}0\\1\\0\\0\end{bmatrix}[1\ 0\ 0] = \begin{bmatrix}0&0&0\\1&0&0\\0&0&0\\0&0&0\end{bmatrix}.$$

Recall

$$a = \begin{bmatrix}a_{11} & a_{12} & \cdots & a_{1m}\\ a_{21} & a_{22} & \cdots & a_{2m}\\ a_{n1} & a_{n2} & \cdots & a_{nm}\end{bmatrix}$$

we have

$$a = \begin{bmatrix}a_1^1\\ a_2^1\\ \vdots\\ a_n^1\end{bmatrix}_m \qquad (1)$$

and

$$a = [a^1\ a^2\ \ldots\ a^m]_n\ . \qquad (2)$$

Now transpose of

$$a = \begin{bmatrix} a_1^1 \\ a_2^1 \\ \vdots \\ a_n^1 \end{bmatrix}_m$$

is given by the equation

$$a' = \begin{bmatrix} (a_1^1)' \, (a_2^1)' \cdots (a_n^1)' \end{bmatrix}_m$$

$$a' = \begin{bmatrix} (a^1)' \\ (a^2)' \\ \vdots \\ (a^m)' \end{bmatrix}_n .$$

The matrix

$$b = \begin{bmatrix} b_{11} & b_{12} & \cdots & b_{1s} \\ b_{21} & b_{22} & \cdots & b_{2s} \\ \vdots & \cdots & & \vdots \\ b_{t1} & b_{t2} & \cdots & b_{ts} \end{bmatrix}$$

row supervector of b is

$$b = [b_1 \; b_2 \ldots b_s]_t = [b^1 \; b^2 \ldots b^s]_t .$$

Column supervector of b is

$$b = \begin{bmatrix} b_1^1 \\ b_2^1 \\ \vdots \\ b_t^1 \end{bmatrix}_s .$$

Transpose of b;

$$b' = \begin{bmatrix} b_1^1 \\ b_2^1 \\ \vdots \\ b_s^1 \end{bmatrix}_t$$

$$b' = [b_1\ b_2\ \ldots\ b_t]_s.$$

The product of two matrices as a minor product of type II supervector.

$$ab = [a^1\ a^2\ \ldots\ a^m]_n \begin{bmatrix} b_1^1 \\ b_2^1 \\ \vdots \\ b_t^1 \end{bmatrix}_s$$

$$= \left[ a_1 b_1^1 + a_2 b_2^1 + \ldots + a_m b_t^1 \right]_{ns}.$$

How ever to make this point clear we give an example.

***Example 1.1.21:*** Let

$$a = \begin{array}{c} \begin{matrix} a_1^1 & a_2^1 \end{matrix} \\ \begin{bmatrix} 2 & 1 \\ 3 & 5 \\ 6 & 1 \end{bmatrix} \end{array}$$

and

$$b = \begin{bmatrix} 1 & 2 \\ 3 & 1 \end{bmatrix} \begin{matrix} b^1 \\ b^2 \end{matrix}.$$

$$ab = \begin{bmatrix} 2 \\ 3 \\ 6 \end{bmatrix} \begin{bmatrix} 1 & 2 \end{bmatrix} + \begin{bmatrix} 1 \\ 5 \\ 1 \end{bmatrix} \begin{bmatrix} 3 & 1 \end{bmatrix}$$

$$= \begin{bmatrix} 2 & 4 \\ 3 & 6 \\ 6 & 12 \end{bmatrix} + \begin{bmatrix} 3 & 1 \\ 15 & 5 \\ 3 & 1 \end{bmatrix}$$

$$= \begin{bmatrix} 5 & 5 \\ 18 & 11 \\ 9 & 13 \end{bmatrix}.$$

It is easily verified that if the major product of the type II supervector is computed between a and b, then the major product coincides with the minor product. From the above example.

$$ab = \begin{bmatrix} [2 \quad 1] \begin{bmatrix} 1 \\ 3 \end{bmatrix} & [2 \quad 1] \begin{bmatrix} 2 \\ 1 \end{bmatrix} \\ [3 \quad 5] \begin{bmatrix} 1 \\ 3 \end{bmatrix} & [3 \quad 5] \begin{bmatrix} 2 \\ 1 \end{bmatrix} \\ [6 \quad 1] \begin{bmatrix} 1 \\ 3 \end{bmatrix} & [6 \quad 1] \begin{bmatrix} 2 \\ 1 \end{bmatrix} \end{bmatrix}$$

$$= \begin{bmatrix} 2\times1+1\times3 & 2\times2+1\times1 \\ 3\times1+5\times3 & 3\times2+5\times1 \\ 6\times1+1\times3 & 6\times2+1\times1 \end{bmatrix}$$

$$= \begin{bmatrix} 5 & 5 \\ 18 & 11 \\ 9 & 13 \end{bmatrix}.$$

We can find the minor and major product of supervectors by reversing the order of the factors. Since the theory of multiplication of supermatrices involves lots of notations we have resolved to explain these concepts by working out these concepts with numerical illustrations, which we feel is easy for the grasp of the reader. Now we give the numerical illustration of the minor product of Type III vectors.

***Example 1.1.22:*** Let

$$X = \left[\begin{array}{cc|c|ccc} 2 & 3 & 4 & 2 & 2 & 2 \\ -1 & 1 & 1 & 1 & 0 & 1 \\ 0 & 0 & 2 & -4 & 0 & 0 \end{array}\right]$$

and

$$Y = \left[\begin{array}{cc} 2 & 0 \\ 1 & 1 \\ \hline 2 & 1 \\ \hline 5 & 3 \\ 1 & -1 \\ 0 & 2 \end{array}\right]$$

be two type III supervectors. To find the product XY.

$$X\,Y = \left[\begin{array}{cc|c|ccc} 2 & 3 & 4 & 2 & 2 & 2 \\ -1 & 1 & 1 & 1 & 0 & 1 \\ 0 & 0 & 2 & -4 & 0 & 0 \end{array}\right] \left[\begin{array}{cc} 2 & 0 \\ 1 & 1 \\ \hline 2 & 1 \\ \hline 5 & 3 \\ 1 & -1 \\ 0 & 2 \end{array}\right]$$

$$= \begin{bmatrix} 2 & 3 \\ -1 & 1 \\ 0 & 0 \end{bmatrix} \begin{bmatrix} 2 & 0 \\ 1 & 1 \end{bmatrix} + \begin{bmatrix} 4 \\ 1 \\ 2 \end{bmatrix} \begin{bmatrix} 2 & 1 \end{bmatrix} + \begin{bmatrix} 2 & 2 & 2 \\ 1 & 0 & 1 \\ -4 & 0 & 0 \end{bmatrix} \begin{bmatrix} 5 & 3 \\ 1 & -1 \\ 0 & 2 \end{bmatrix}$$

$$= \begin{bmatrix} 7 & 3 \\ -1 & 1 \\ 0 & 0 \end{bmatrix} + \begin{bmatrix} 8 & 4 \\ 2 & 1 \\ 4 & 2 \end{bmatrix} + \begin{bmatrix} 12 & 8 \\ 5 & 5 \\ -20 & -12 \end{bmatrix}$$

$$= \begin{bmatrix} 27 & 15 \\ 6 & 7 \\ -16 & -10 \end{bmatrix}.$$

$$Y^t X^t = \left[\begin{array}{cc|c|ccc} 2 & 1 & 2 & 5 & 1 & 0 \\ 0 & 1 & 1 & 3 & -1 & 2 \end{array}\right] \left[\begin{array}{ccc} 2 & -1 & 0 \\ 3 & 1 & 0 \\ \hline 4 & 1 & 2 \\ \hline 2 & 1 & -4 \\ 2 & 0 & 0 \\ 2 & 1 & 0 \end{array}\right]$$

$$= \begin{bmatrix} 2 & 1 \\ 0 & 1 \end{bmatrix} \begin{bmatrix} 2 & -1 & 0 \\ 3 & 1 & 0 \end{bmatrix} + \begin{bmatrix} 2 \\ 1 \end{bmatrix} \begin{bmatrix} 4 & 1 & 2 \end{bmatrix} + \begin{bmatrix} 5 & 1 & 0 \\ 3 & -1 & 2 \end{bmatrix} \begin{bmatrix} 2 & 1 & -4 \\ 2 & 0 & 0 \\ 2 & 1 & 0 \end{bmatrix}$$

$$= \begin{bmatrix} 7 & -1 & 0 \\ 3 & 1 & 0 \end{bmatrix} + \begin{bmatrix} 8 & 2 & 4 \\ 4 & 1 & 2 \end{bmatrix} + \begin{bmatrix} 12 & 5 & -20 \\ 8 & 5 & -12 \end{bmatrix}$$

$$= \begin{bmatrix} 27 & 6 & -16 \\ 15 & 7 & -10 \end{bmatrix}.$$

From this example it is very clear.

$$(XY)^t = Y^t X^t.$$

Now we illustrate the minor product moment of type III row supervector by an example.

***Example 1.1.23:*** Let

$$X = \left[\begin{array}{cc|c|cccc} 2 & 3 & 4 & 3 & 4 & 5 & 0 \\ 1 & 4 & 1 & 1 & 1 & -1 & 6 \\ 2 & 1 & 2 & 0 & 2 & 1 & 1 \end{array}\right].$$

Consider

$$XX' = \left[\begin{array}{cc|c|cccc} 2 & 3 & 4 & 3 & 4 & 5 & 0 \\ 1 & 4 & 1 & 1 & 1 & -1 & 6 \\ 2 & 1 & 2 & 0 & 2 & 1 & 1 \end{array}\right] \left[\begin{array}{ccc} 2 & 1 & 2 \\ 3 & 4 & 1 \\ \hline 4 & 1 & 2 \\ \hline 3 & 1 & 0 \\ 4 & 1 & 2 \\ 5 & -1 & 1 \\ 0 & 6 & 1 \end{array}\right]$$

$$= \begin{bmatrix} 2 & 3 \\ 1 & 4 \\ 2 & 1 \end{bmatrix} \begin{bmatrix} 2 & 1 & 2 \\ 3 & 4 & 1 \end{bmatrix} + \begin{bmatrix} 4 \\ 1 \\ 2 \end{bmatrix} \begin{bmatrix} 4 & 1 & 2 \end{bmatrix} +$$

$$\begin{bmatrix} 3 & 4 & 5 & 0 \\ 1 & 1 & -1 & 6 \\ 0 & 2 & 1 & 1 \end{bmatrix} \begin{bmatrix} 3 & 1 & 0 \\ 4 & 1 & 2 \\ 5 & -1 & 1 \\ 0 & 6 & 1 \end{bmatrix}$$

$$= \begin{bmatrix} 13 & 14 & 7 \\ 14 & 17 & 6 \\ 7 & 6 & 5 \end{bmatrix} + \begin{bmatrix} 16 & 4 & 8 \\ 4 & 1 & 2 \\ 8 & 2 & 4 \end{bmatrix} + \begin{bmatrix} 50 & 2 & 13 \\ 2 & 39 & 7 \\ 13 & 7 & 6 \end{bmatrix}$$

$$= \begin{bmatrix} 79 & 20 & 28 \\ 20 & 57 & 15 \\ 28 & 15 & 15 \end{bmatrix}.$$

Minor product of Type III column supervector is illustrated by the following example.

***Example 1.1.24:*** Let

$$Y^t = \left[\begin{array}{ccc|c|ccccc} 2 & 3 & 1 & 0 & 1 & 2 & 1 & 5 & 1 \\ 0 & 1 & 5 & 2 & 0 & 3 & 0 & 1 & 0 \end{array}\right]$$

where Y is the column supervector

$$Y^tY = \left[\begin{array}{ccc|c|ccccc} 2 & 3 & 1 & 0 & 1 & 2 & 1 & 5 & 1 \\ 0 & 1 & 5 & 2 & 0 & 3 & 0 & 1 & 0 \end{array}\right] \left[\begin{array}{cc} 2 & 0 \\ 3 & 1 \\ 1 & 5 \\ \hline 0 & 2 \\ \hline 1 & 0 \\ 2 & 3 \\ 1 & 0 \\ 5 & 1 \\ 1 & 0 \end{array}\right]$$

$$= \begin{bmatrix} 2 & 3 & 1 \\ 0 & 1 & 5 \end{bmatrix} \begin{bmatrix} 2 & 0 \\ 3 & 1 \\ 1 & 5 \end{bmatrix} + \begin{bmatrix} 0 \\ 2 \end{bmatrix} \begin{bmatrix} 0 & 2 \end{bmatrix} + \begin{bmatrix} 1 & 2 & 1 & 5 & 1 \\ 0 & 3 & 0 & 1 & 0 \end{bmatrix} \begin{bmatrix} 1 & 0 \\ 2 & 3 \\ 1 & 0 \\ 5 & 1 \\ 1 & 0 \end{bmatrix}$$

$$= \begin{bmatrix} 14 & 8 \\ 8 & 26 \end{bmatrix} + \begin{bmatrix} 0 & 0 \\ 0 & 4 \end{bmatrix} + \begin{bmatrix} 32 & 11 \\ 11 & 10 \end{bmatrix} = \begin{bmatrix} 46 & 19 \\ 19 & 40 \end{bmatrix}.$$

Next we proceed on to illustrate the major product of Type III vectors.

***Example 1.1.25:*** Let

$$X = \left[\begin{array}{ccc} 3 & 1 & 6 \\ 2 & 0 & -1 \\ \hline 1 & 2 & 3 \\ \hline 6 & 3 & 0 \\ 4 & 2 & 1 \\ 5 & 1 & -1 \end{array}\right]$$

and

$$Y = \left[\begin{array}{c|ccc} 3 & 5 & 2 & 0 \\ 1 & 1 & 2 & 2 \\ 0 & 3 & 1 & -2 \end{array}\right].$$

$$XY = \left[\begin{array}{ccc} 3 & 1 & 6 \\ 2 & 0 & -1 \\ \hline 1 & 2 & 3 \\ \hline 6 & 3 & 0 \\ 4 & 2 & 1 \\ 5 & 1 & -1 \end{array}\right] \left[\begin{array}{c|ccc} 3 & 5 & 2 & 0 \\ 1 & 1 & 2 & 2 \\ 0 & 3 & 1 & -2 \end{array}\right]$$

$$= \left[\begin{array}{c|c} \begin{bmatrix} 3 & 1 & 6 \\ 2 & 0 & -1 \end{bmatrix} \begin{bmatrix} 3 \\ 1 \\ 0 \end{bmatrix} & \begin{bmatrix} 3 & 1 & 6 \\ 2 & 0 & -1 \end{bmatrix} \begin{bmatrix} 5 & 2 & 0 \\ 1 & 2 & 2 \\ 3 & 1 & -2 \end{bmatrix} \\ \hline \begin{bmatrix} 1 & 2 & 3 \end{bmatrix} \begin{bmatrix} 3 \\ 1 \\ 0 \end{bmatrix} & \begin{bmatrix} 1 & 2 & 3 \end{bmatrix} \begin{bmatrix} 5 & 2 & 0 \\ 1 & 2 & 2 \\ 3 & 1 & -2 \end{bmatrix} \\ \hline \begin{bmatrix} 6 & 3 & 0 \\ 4 & 2 & 1 \\ 5 & 1 & -1 \end{bmatrix} \begin{bmatrix} 3 \\ 1 \\ 0 \end{bmatrix} & \begin{bmatrix} 6 & 3 & 0 \\ 4 & 2 & 1 \\ 5 & 1 & -1 \end{bmatrix} \begin{bmatrix} 5 & 2 & 0 \\ 1 & 2 & 2 \\ 3 & 1 & -2 \end{bmatrix} \end{array}\right]$$

$$= \left[\begin{array}{c|ccc} 10 & 34 & 14 & -10 \\ 6 & 7 & 3 & 2 \\ \hline 5 & 16 & 9 & -2 \\ \hline 21 & 33 & 18 & 6 \\ 14 & 25 & 13 & 2 \\ 16 & 23 & 11 & 4 \end{array}\right].$$

Now minor product of type IV vector is illustrated by the following example.

***Example 1.1.26:*** Let

$$X = \left[\begin{array}{c|ccc|cc} 1 & 3 & 1 & 2 & 5 & 1 \\ 2 & 1 & 1 & 1 & 2 & 0 \\ \hline 1 & 5 & 1 & 1 & 1 & 2 \\ 4 & 1 & 0 & 2 & 2 & 1 \\ 3 & 2 & 1 & 0 & 1 & 1 \\ -1 & 0 & 1 & 1 & 0 & 1 \\ \hline 4 & 2 & 1 & 3 & 1 & 1 \end{array}\right]$$

and

$$Y = \left[\begin{array}{cccc|ccc|cc} 1 & 1 & 0 & 1 & 3 & 1 & 2 & 1 & 2 \\ \hline 2 & 0 & 1 & 0 & 1 & 2 & 0 & 0 & 1 \\ 1 & 1 & 0 & 2 & 3 & 0 & 1 & 1 & 4 \\ 1 & 0 & 1 & 1 & 2 & 1 & 1 & 2 & 0 \\ \hline 1 & 2 & 0 & 1 & 1 & 0 & 1 & 1 & 2 \\ 0 & 1 & 1 & 0 & 1 & 1 & 0 & 2 & 1 \end{array}\right].$$

$$XY = \left[\begin{array}{c|ccc|cc} 1 & 3 & 1 & 2 & 5 & 1 \\ 2 & 1 & 1 & 1 & 2 & 0 \\ \hline 1 & 5 & 1 & 1 & 1 & 2 \\ 4 & 1 & 0 & 2 & 2 & 1 \\ 3 & 2 & 1 & 0 & 1 & 1 \\ -1 & 0 & 1 & 1 & 0 & 1 \\ \hline 4 & 2 & 1 & 3 & 1 & 1 \end{array}\right] \times$$

$$\left[\begin{array}{cccc|ccc|cc} 1 & 1 & 0 & 1 & 3 & 1 & 2 & 1 & 2 \\ \hline 2 & 0 & 1 & 0 & 1 & 2 & 0 & 0 & 1 \\ 1 & 1 & 0 & 2 & 3 & 0 & 1 & 1 & 4 \\ 1 & 0 & 1 & 1 & 2 & 1 & 1 & 2 & 0 \\ \hline 1 & 2 & 0 & 1 & 1 & 0 & 1 & 1 & 2 \\ 0 & 1 & 1 & 0 & 1 & 1 & 0 & 2 & 1 \end{array}\right]$$

$$= \left[\begin{array}{c} 1 \\ 2 \\ \hline 1 \\ 4 \\ 3 \\ -1 \\ \hline 4 \end{array}\right] [1\ 1\ 0\ 1 \mid 3\ 1\ 2 \mid 1\ 2] +$$

$$\left[\begin{array}{ccc} 3 & 1 & 2 \\ 1 & 1 & 1 \\ \hline 5 & 1 & 1 \\ 1 & 0 & 2 \\ 2 & 1 & 0 \\ 0 & 1 & 1 \\ \hline 2 & 1 & 3 \end{array}\right] \left[\begin{array}{cccc|ccc|cc} 2 & 0 & 1 & 0 & 1 & 2 & 0 & 0 & 1 \\ 1 & 1 & 0 & 2 & 3 & 0 & 1 & 1 & 4 \\ 1 & 0 & 1 & 1 & 2 & 1 & 1 & 2 & 0 \end{array}\right]$$

$$+ \left[\begin{array}{cc} 5 & 1 \\ 2 & 0 \\ \hline 1 & 2 \\ 2 & 1 \\ 1 & 1 \\ 0 & 1 \\ \hline 1 & 1 \end{array}\right] \left[\begin{array}{cccc|ccc|cc} 1 & 2 & 0 & 1 & 1 & 0 & 1 & 1 & 2 \\ 0 & 1 & 1 & 0 & 1 & 1 & 0 & 2 & 1 \end{array}\right]$$

$$= \left[\begin{array}{c|c|c} \begin{bmatrix} 1 \\ 2 \end{bmatrix}\begin{bmatrix} 1 & 1 & 0 & 1 \end{bmatrix} & \begin{bmatrix} 1 \\ 2 \end{bmatrix}\begin{bmatrix} 3 & 1 & 2 \end{bmatrix} & \begin{bmatrix} 1 \\ 2 \end{bmatrix}\begin{bmatrix} 1 & 2 \end{bmatrix} \\ \hline \begin{bmatrix} 1 \\ 4 \\ 3 \\ -1 \end{bmatrix}\begin{bmatrix} 1 & 1 & 0 & 1 \end{bmatrix} & \begin{bmatrix} 1 \\ 4 \\ 3 \\ -1 \end{bmatrix}\begin{bmatrix} 3 & 1 & 2 \end{bmatrix} & \begin{bmatrix} 1 \\ 4 \\ 3 \\ -1 \end{bmatrix}\begin{bmatrix} 1 & 2 \end{bmatrix} \\ \hline [4]\begin{bmatrix} 1 & 1 & 0 & 1 \end{bmatrix} & [4]\begin{bmatrix} 3 & 1 & 2 \end{bmatrix} & [4]\begin{bmatrix} 1 & 2 \end{bmatrix} \end{array}\right] +$$

$$
\left[\begin{array}{c|c|c}
\begin{bmatrix}3&1&2\\1&1&1\end{bmatrix}\begin{bmatrix}2&0&1&0\\1&1&0&2\\1&0&1&1\end{bmatrix} & \begin{bmatrix}3&1&2\\1&1&1\end{bmatrix}\begin{bmatrix}1&2&0\\3&0&1\\2&1&1\end{bmatrix} & \begin{bmatrix}3&1&2\\1&1&1\end{bmatrix}\begin{bmatrix}0&1\\1&4\\2&0\end{bmatrix} \\ \hline
\begin{bmatrix}5&1&1\\1&0&2\\2&1&0\\0&1&1\end{bmatrix}\begin{bmatrix}2&0&1&0\\1&1&0&2\\1&0&1&1\end{bmatrix} & \begin{bmatrix}5&1&1\\1&0&2\\2&1&0\\0&1&1\end{bmatrix}\begin{bmatrix}1&2&0\\3&0&1\\2&1&1\end{bmatrix} & \begin{bmatrix}5&1&1\\1&0&2\\2&1&0\\0&1&1\end{bmatrix}\begin{bmatrix}0&1\\1&4\\2&0\end{bmatrix} \\ \hline
\begin{bmatrix}2&1&3\end{bmatrix}\begin{bmatrix}2&0&1&0\\1&1&0&2\\1&0&1&1\end{bmatrix} & \begin{bmatrix}2&1&3\end{bmatrix}\begin{bmatrix}1&2&0\\3&0&1\\2&1&1\end{bmatrix} & \begin{bmatrix}2&1&3\end{bmatrix}\begin{bmatrix}0&1\\1&4\\2&0\end{bmatrix}
\end{array}\right]
$$

$$
+\left[\begin{array}{c|c|c}
\begin{bmatrix}5&1\\2&0\end{bmatrix}\begin{bmatrix}1&2&0&1\\0&1&1&0\end{bmatrix} & \begin{bmatrix}5&1\\2&0\end{bmatrix}\begin{bmatrix}1&0&1\\1&1&0\end{bmatrix} & \begin{bmatrix}5&1\\2&0\end{bmatrix}\begin{bmatrix}1&2\\2&1\end{bmatrix} \\ \hline
\begin{bmatrix}1&2\\2&1\\1&1\\0&1\end{bmatrix}\begin{bmatrix}1&2&0&1\\0&1&1&0\end{bmatrix} & \begin{bmatrix}1&2\\2&1\\1&1\\0&1\end{bmatrix}\begin{bmatrix}1&0&1\\1&1&0\end{bmatrix} & \begin{bmatrix}1&2\\2&1\\1&1\\0&1\end{bmatrix}\begin{bmatrix}1&2\\2&1\end{bmatrix} \\ \hline
\begin{bmatrix}1&1\end{bmatrix}\begin{bmatrix}1&2&0&1\\0&1&1&0\end{bmatrix} & \begin{bmatrix}1&1\end{bmatrix}\begin{bmatrix}1&0&1\\1&1&0\end{bmatrix} & \begin{bmatrix}1&1\end{bmatrix}\begin{bmatrix}1&2\\2&1\end{bmatrix}
\end{array}\right]
$$

$$
=\left[\begin{array}{c|c|c}
\begin{bmatrix}1&1&0&1\\2&2&0&2\end{bmatrix} & \begin{bmatrix}3&1&2\\6&2&4\end{bmatrix} & \begin{bmatrix}1&2\\2&4\end{bmatrix} \\ \hline
\begin{bmatrix}1&1&0&1\\4&4&0&4\\3&3&0&3\\-1&-1&0&-1\end{bmatrix} & \begin{bmatrix}3&1&2\\12&4&8\\9&3&6\\-3&-1&-2\end{bmatrix} & \begin{bmatrix}1&2\\4&8\\3&6\\-1&-2\end{bmatrix} \\ \hline
\begin{bmatrix}4&4&0&4\end{bmatrix} & \begin{bmatrix}12&4&8\end{bmatrix} & \begin{bmatrix}4&8\end{bmatrix}
\end{array}\right]+
$$

$$\left[\begin{array}{cccc|ccc|cc} 9 & 1 & 5 & 4 & 10 & 8 & 3 & 5 & 7 \\ 4 & 1 & 2 & 3 & 6 & 3 & 2 & 3 & 5 \\ \hline 12 & 1 & 6 & 3 & 10 & 11 & 2 & 3 & 9 \\ 4 & 0 & 3 & 2 & 5 & 4 & 2 & 4 & 1 \\ 5 & 1 & 2 & 2 & 5 & 4 & 1 & 1 & 6 \\ 2 & 1 & 1 & 3 & 5 & 1 & 2 & 3 & 4 \\ \hline 8 & 1 & 5 & 5 & 11 & 7 & 4 & 7 & 6 \end{array}\right]$$

$$+\left[\begin{array}{cccc|ccc|cc} 5 & 11 & 1 & 5 & 6 & 1 & 5 & 7 & 11 \\ 2 & 4 & 0 & 2 & 2 & 0 & 2 & 2 & 4 \\ \hline 1 & 4 & 2 & 1 & 3 & 2 & 1 & 5 & 4 \\ 2 & 5 & 1 & 2 & 3 & 1 & 2 & 4 & 5 \\ 1 & 3 & 1 & 1 & 2 & 1 & 1 & 3 & 3 \\ 0 & 1 & 1 & 0 & 1 & 1 & 0 & 2 & 1 \\ \hline 1 & 3 & 1 & 1 & 2 & 1 & 1 & 3 & 3 \end{array}\right]$$

$$=\left[\begin{array}{cccc|ccc|cc} 15 & 13 & 6 & 10 & 19 & 10 & 10 & 13 & 20 \\ 8 & 7 & 2 & 7 & 14 & 5 & 8 & 7 & 13 \\ \hline 14 & 6 & 8 & 5 & 16 & 14 & 5 & 9 & 15 \\ 10 & 9 & 4 & 8 & 20 & 9 & 12 & 12 & 14 \\ 9 & 7 & 3 & 6 & 16 & 8 & 8 & 7 & 15 \\ 1 & 1 & 2 & 2 & 3 & 1 & 0 & 4 & 3 \\ \hline 13 & 8 & 6 & 10 & 25 & 12 & 13 & 14 & 17 \end{array}\right].$$

We now illustrate minor product moment of type IV row vector

***Example 1.1.27:*** Let

$$X=\left[\begin{array}{ccc|cc|c} 1 & 1 & 1 & 1 & 0 & 1 \\ 2 & 1 & 2 & 2 & 1 & 2 \\ \hline 0 & 1 & 1 & 3 & 1 & 1 \\ 1 & 0 & 1 & 1 & 3 & 2 \\ 5 & 1 & 0 & 2 & 1 & 3 \\ 1 & 1 & 0 & 1 & 2 & 4 \\ \hline 2 & 1 & 1 & 5 & 0 & 2 \end{array}\right].$$

$$
XX^t = \left[\begin{array}{ccc|cc|c}
1 & 1 & 1 & 1 & 0 & 1 \\
2 & 1 & 2 & 2 & 1 & 2 \\ \hline
0 & 1 & 1 & 3 & 1 & 1 \\
1 & 0 & 1 & 1 & 3 & 2 \\
5 & 1 & 0 & 2 & 1 & 3 \\
1 & 1 & 0 & 1 & 2 & 4 \\ \hline
2 & 1 & 1 & 5 & 0 & 2
\end{array}\right]
\left[\begin{array}{cc|cccc|c}
1 & 2 & 0 & 1 & 5 & 1 & 2 \\
1 & 1 & 1 & 0 & 1 & 1 & 1 \\
1 & 2 & 1 & 1 & 0 & 0 & 1 \\ \hline
1 & 2 & 3 & 1 & 2 & 1 & 5 \\
0 & 1 & 1 & 3 & 1 & 2 & 0 \\ \hline
1 & 2 & 1 & 2 & 3 & 4 & 2
\end{array}\right]
$$

$$
= \left[\begin{array}{ccc}
1 & 1 & 1 \\
2 & 1 & 2 \\ \hline
0 & 1 & 1 \\
1 & 0 & 1 \\
5 & 1 & 0 \\
1 & 1 & 0 \\ \hline
2 & 1 & 1
\end{array}\right]
\left[\begin{array}{cc|cccc|c}
1 & 2 & 0 & 1 & 5 & 1 & 2 \\
1 & 1 & 1 & 0 & 1 & 1 & 1 \\
1 & 2 & 1 & 1 & 0 & 0 & 1
\end{array}\right] +
$$

$$
\left[\begin{array}{cc}
1 & 0 \\
2 & 1 \\ \hline
3 & 1 \\
1 & 3 \\
2 & 1 \\
1 & 2 \\ \hline
5 & 0
\end{array}\right]
\left[\begin{array}{cc|cccc|c}
1 & 2 & 3 & 1 & 2 & 1 & 5 \\
0 & 1 & 1 & 3 & 1 & 2 & 0
\end{array}\right] +
\left[\begin{array}{c}
1 \\
2 \\ \hline
1 \\
2 \\
3 \\
4 \\ \hline
2
\end{array}\right]
[1\ 2 \mid 1\ 2\ 3\ 4 \mid 2] =
$$

$$
\left[\begin{array}{c|c|c}
\begin{bmatrix}1&1&1\\2&1&2\end{bmatrix}\begin{bmatrix}1&2\\1&1\\1&2\end{bmatrix} & \begin{bmatrix}1&1&1\\2&1&2\end{bmatrix}\begin{bmatrix}0&1&5&1\\1&0&1&1\\1&1&0&0\end{bmatrix} & \begin{bmatrix}1&1&1\\2&1&2\end{bmatrix}\begin{bmatrix}2\\1\\1\end{bmatrix} \\
\hline
\begin{bmatrix}0&1&1\\1&0&1\\5&1&0\\1&1&0\end{bmatrix}\begin{bmatrix}1&2\\1&1\\1&2\end{bmatrix} & \begin{bmatrix}0&1&1\\1&0&1\\5&1&0\\1&1&0\end{bmatrix}\begin{bmatrix}0&1&5&1\\1&0&1&1\\1&1&0&0\end{bmatrix} & \begin{bmatrix}0&1&1\\1&0&1\\5&1&0\\1&1&0\end{bmatrix}\begin{bmatrix}2\\1\\1\end{bmatrix} \\
\hline
\begin{bmatrix}2&1&1\end{bmatrix}\begin{bmatrix}1&2\\1&1\\1&2\end{bmatrix} & \begin{bmatrix}2&1&1\end{bmatrix}\begin{bmatrix}0&1&5&1\\1&0&1&1\\1&1&0&0\end{bmatrix} & \begin{bmatrix}2&1&1\end{bmatrix}\begin{bmatrix}2\\1\\1\end{bmatrix}
\end{array}\right]
$$

$$
+\left[\begin{array}{c|c|c}
\begin{bmatrix}1&0\\2&1\end{bmatrix}\begin{bmatrix}1&2\\0&1\end{bmatrix} & \begin{bmatrix}1&0\\2&1\end{bmatrix}\begin{bmatrix}3&1&2&1\\1&3&1&2\end{bmatrix} & \begin{bmatrix}1&0\\2&1\end{bmatrix}\begin{bmatrix}5\\0\end{bmatrix} \\
\hline
\begin{bmatrix}3&1\\1&3\\2&1\\1&2\end{bmatrix}\begin{bmatrix}1&2\\0&1\end{bmatrix} & \begin{bmatrix}3&1\\1&3\\2&1\\1&2\end{bmatrix}\begin{bmatrix}3&1&2&1\\1&3&1&2\end{bmatrix} & \begin{bmatrix}3&1\\1&3\\2&1\\1&2\end{bmatrix}\begin{bmatrix}5\\0\end{bmatrix} \\
\hline
\begin{bmatrix}5&0\end{bmatrix}\begin{bmatrix}1&2\\0&1\end{bmatrix} & \begin{bmatrix}5&0\end{bmatrix}\begin{bmatrix}3&1&2&1\\1&3&1&2\end{bmatrix} & \begin{bmatrix}5&0\end{bmatrix}\begin{bmatrix}5\\0\end{bmatrix}
\end{array}\right]
$$

$$
+\left[\begin{array}{c|c|c}
\begin{bmatrix}1\\2\end{bmatrix}\begin{bmatrix}1&2\end{bmatrix} & \begin{bmatrix}1\\2\end{bmatrix}\begin{bmatrix}1&2&3&4\end{bmatrix} & \begin{bmatrix}1\\2\end{bmatrix}\begin{bmatrix}2\end{bmatrix} \\
\hline
\begin{bmatrix}1\\2\\3\\4\end{bmatrix}\begin{bmatrix}1&2\end{bmatrix} & \begin{bmatrix}1\\2\\3\\4\end{bmatrix}\begin{bmatrix}1&2&3&4\end{bmatrix} & \begin{bmatrix}1\\2\\3\\4\end{bmatrix}\begin{bmatrix}2\end{bmatrix} \\
\hline
\begin{bmatrix}2\end{bmatrix}\begin{bmatrix}1&2\end{bmatrix} & \begin{bmatrix}2\end{bmatrix}\begin{bmatrix}1&2&3&4\end{bmatrix} & \begin{bmatrix}2\end{bmatrix}\begin{bmatrix}2\end{bmatrix}
\end{array}\right]
$$

$$= \left[\begin{array}{cc|cccc|c} 3 & 5 & 2 & 2 & 6 & 2 & 4 \\ 5 & 9 & 3 & 4 & 11 & 3 & 7 \\ \hline 2 & 3 & 2 & 1 & 1 & 1 & 2 \\ 2 & 4 & 1 & 2 & 5 & 1 & 3 \\ 6 & 11 & 1 & 5 & 26 & 6 & 11 \\ 2 & 3 & 1 & 1 & 6 & 2 & 3 \\ \hline 4 & 7 & 2 & 3 & 11 & 3 & 6 \end{array}\right] +$$

$$\left[\begin{array}{cc|cccc|c} 1 & 2 & 3 & 1 & 2 & 1 & 5 \\ 2 & 5 & 7 & 5 & 5 & 4 & 10 \\ \hline 3 & 7 & 10 & 6 & 7 & 5 & 15 \\ 1 & 5 & 6 & 10 & 5 & 7 & 5 \\ 2 & 5 & 7 & 5 & 5 & 4 & 10 \\ 1 & 4 & 5 & 7 & 4 & 5 & 5 \\ \hline 5 & 10 & 15 & 5 & 10 & 5 & 25 \end{array}\right] + \left[\begin{array}{cc|cccc|c} 1 & 2 & 1 & 2 & 3 & 4 & 2 \\ 2 & 4 & 2 & 4 & 6 & 8 & 4 \\ \hline 1 & 2 & 1 & 2 & 3 & 4 & 2 \\ 2 & 4 & 2 & 4 & 6 & 8 & 4 \\ 3 & 6 & 3 & 6 & 9 & 12 & 6 \\ 4 & 8 & 4 & 8 & 12 & 16 & 8 \\ \hline 2 & 4 & 2 & 4 & 6 & 8 & 4 \end{array}\right]$$

$$= \left[\begin{array}{cc|cccc|c} 5 & 9 & 6 & 5 & 11 & 7 & 11 \\ 9 & 18 & 12 & 13 & 22 & 15 & 21 \\ \hline 6 & 12 & 13 & 9 & 11 & 10 & 19 \\ 5 & 13 & 9 & 16 & 16 & 16 & 12 \\ 11 & 22 & 11 & 16 & 40 & 22 & 27 \\ 7 & 15 & 10 & 16 & 22 & 23 & 16 \\ \hline 11 & 21 & 19 & 12 & 27 & 16 & 35 \end{array}\right].$$

The minor product moment of type IV column vector is illustrated for the same X just given in case of row product.

***Example 1.1.28:*** Let

$$X = \left[\begin{array}{ccc|cc|c} 1 & 1 & 1 & 1 & 0 & 1 \\ 2 & 1 & 2 & 2 & 1 & 2 \\ \hline 0 & 1 & 1 & 3 & 1 & 1 \\ 1 & 0 & 1 & 1 & 3 & 2 \\ 5 & 1 & 0 & 2 & 1 & 3 \\ 1 & 1 & 0 & 1 & 2 & 4 \\ \hline 2 & 1 & 1 & 5 & 0 & 2 \end{array}\right].$$

$$X^t X = \left[\begin{array}{cc|cccc|c} 1 & 2 & 0 & 1 & 5 & 1 & 2 \\ 1 & 1 & 1 & 0 & 1 & 1 & 1 \\ 1 & 2 & 1 & 1 & 0 & 0 & 1 \\ \hline 1 & 2 & 3 & 1 & 2 & 1 & 5 \\ 0 & 1 & 1 & 3 & 1 & 2 & 0 \\ \hline 1 & 2 & 1 & 2 & 3 & 4 & 2 \end{array}\right] \left[\begin{array}{ccc|cc|c} 1 & 1 & 1 & 1 & 0 & 1 \\ 2 & 1 & 2 & 2 & 1 & 2 \\ \hline 0 & 1 & 1 & 3 & 1 & 1 \\ 1 & 0 & 1 & 1 & 3 & 2 \\ 5 & 1 & 0 & 2 & 1 & 3 \\ 1 & 1 & 0 & 1 & 2 & 4 \\ \hline 2 & 1 & 1 & 5 & 0 & 2 \end{array}\right]$$

$$= \left[\begin{array}{cc} 1 & 2 \\ 1 & 1 \\ 1 & 2 \\ \hline 1 & 2 \\ 0 & 1 \\ \hline 1 & 2 \end{array}\right] \left[\begin{array}{ccc|cc|c} 1 & 1 & 1 & 1 & 0 & 1 \\ 2 & 1 & 2 & 2 & 1 & 2 \end{array}\right]$$

$$+ \left[\begin{array}{cccc} 0 & 1 & 5 & 1 \\ 1 & 0 & 1 & 1 \\ 1 & 1 & 0 & 0 \\ \hline 3 & 1 & 2 & 1 \\ 1 & 3 & 1 & 2 \\ \hline 1 & 2 & 3 & 4 \end{array}\right] \left[\begin{array}{ccc|cc|c} 0 & 1 & 1 & 3 & 1 & 1 \\ 1 & 0 & 1 & 1 & 3 & 2 \\ 5 & 1 & 0 & 2 & 1 & 3 \\ 1 & 1 & 0 & 1 & 2 & 4 \end{array}\right] + \left[\begin{array}{c} 2 \\ 1 \\ 1 \\ \hline 5 \\ 0 \\ \hline 2 \end{array}\right] [2\ 1\ 1 \mid 5\ 0 \mid 2]$$

$$
= \left[\begin{array}{c|c|c}
\begin{bmatrix}1 & 2\\ 1 & 1\\ 1 & 2\end{bmatrix}\begin{bmatrix}1 & 1 & 1\\ 2 & 1 & 2\end{bmatrix} & \begin{bmatrix}1 & 2\\ 1 & 1\\ 1 & 2\end{bmatrix}\begin{bmatrix}1 & 0\\ 2 & 1\end{bmatrix} & \begin{bmatrix}1 & 2\\ 1 & 1\\ 1 & 2\end{bmatrix}\begin{bmatrix}1\\ 2\end{bmatrix} \\ \hline
\begin{bmatrix}1 & 2\\ 0 & 1\end{bmatrix}\begin{bmatrix}1 & 1 & 1\\ 2 & 1 & 2\end{bmatrix} & \begin{bmatrix}1 & 2\\ 0 & 1\end{bmatrix}\begin{bmatrix}1 & 0\\ 2 & 1\end{bmatrix} & \begin{bmatrix}1 & 2\\ 0 & 1\end{bmatrix}\begin{bmatrix}1\\ 2\end{bmatrix} \\ \hline
\begin{bmatrix}1 & 2\end{bmatrix}\begin{bmatrix}1 & 1 & 1\\ 2 & 1 & 2\end{bmatrix} & \begin{bmatrix}1 & 2\end{bmatrix}\begin{bmatrix}1 & 0\\ 2 & 1\end{bmatrix} & \begin{bmatrix}1 & 2\end{bmatrix}\begin{bmatrix}1\\ 2\end{bmatrix}
\end{array}\right] +
$$

$$
\left[\begin{array}{c|c|c}
\begin{bmatrix}0 & 1 & 5 & 1\\ 1 & 0 & 1 & 1\\ 1 & 1 & 0 & 0\end{bmatrix}\begin{bmatrix}0 & 1 & 1\\ 1 & 0 & 1\\ 5 & 1 & 0\\ 1 & 1 & 0\end{bmatrix} & \begin{bmatrix}0 & 1 & 5 & 1\\ 1 & 0 & 1 & 1\\ 1 & 1 & 0 & 0\end{bmatrix}\begin{bmatrix}3 & 1\\ 1 & 3\\ 2 & 1\\ 1 & 2\end{bmatrix} & \begin{bmatrix}0 & 1 & 5 & 1\\ 1 & 0 & 1 & 1\\ 1 & 1 & 0 & 0\end{bmatrix}\begin{bmatrix}1\\ 2\\ 3\\ 4\end{bmatrix} \\ \hline
\begin{bmatrix}3 & 1 & 2 & 1\\ 1 & 3 & 1 & 2\end{bmatrix}\begin{bmatrix}0 & 1 & 1\\ 1 & 0 & 1\\ 5 & 1 & 0\\ 1 & 1 & 0\end{bmatrix} & \begin{bmatrix}3 & 1 & 2 & 1\\ 1 & 3 & 1 & 2\end{bmatrix}\begin{bmatrix}3 & 1\\ 1 & 3\\ 2 & 1\\ 1 & 2\end{bmatrix} & \begin{bmatrix}3 & 1 & 2 & 1\\ 1 & 3 & 1 & 2\end{bmatrix}\begin{bmatrix}1\\ 2\\ 3\\ 4\end{bmatrix} \\ \hline
\begin{bmatrix}1 & 2 & 3 & 4\end{bmatrix}\begin{bmatrix}0 & 1 & 1\\ 1 & 0 & 1\\ 5 & 1 & 0\\ 1 & 1 & 0\end{bmatrix} & \begin{bmatrix}1 & 2 & 3 & 4\end{bmatrix}\begin{bmatrix}3 & 1\\ 1 & 3\\ 2 & 1\\ 1 & 2\end{bmatrix} & \begin{bmatrix}1 & 2 & 3 & 4\end{bmatrix}\begin{bmatrix}1\\ 2\\ 3\\ 4\end{bmatrix}
\end{array}\right]
$$

$$
+ \left[\begin{array}{c|c|c}
\begin{bmatrix}2\\ 1\\ 1\end{bmatrix}\begin{bmatrix}2 & 1 & 1\end{bmatrix} & \begin{bmatrix}2\\ 1\\ 1\end{bmatrix}\begin{bmatrix}5 & 0\end{bmatrix} & \begin{bmatrix}2\\ 1\\ 1\end{bmatrix}\begin{bmatrix}2\end{bmatrix} \\ \hline
\begin{bmatrix}5\\ 0\end{bmatrix}\begin{bmatrix}2 & 1 & 1\end{bmatrix} & \begin{bmatrix}5\\ 0\end{bmatrix}\begin{bmatrix}5 & 0\end{bmatrix} & \begin{bmatrix}5\\ 0\end{bmatrix}\begin{bmatrix}2\end{bmatrix} \\ \hline
\begin{bmatrix}2\end{bmatrix}\begin{bmatrix}2 & 1 & 1\end{bmatrix} & \begin{bmatrix}2\end{bmatrix}\begin{bmatrix}5 & 0\end{bmatrix} & \begin{bmatrix}2\end{bmatrix}\begin{bmatrix}2\end{bmatrix}
\end{array}\right] =
$$

$$\left[\begin{array}{ccc|cc|c} 5 & 3 & 5 & 5 & 2 & 5 \\ 3 & 2 & 3 & 3 & 1 & 3 \\ 5 & 3 & 5 & 5 & 2 & 5 \\ \hline 5 & 3 & 5 & 5 & 2 & 5 \\ 2 & 1 & 2 & 2 & 1 & 2 \\ \hline 5 & 3 & 5 & 5 & 2 & 5 \end{array}\right] + \left[\begin{array}{ccc|cc|c} 27 & 6 & 1 & 12 & 10 & 21 \\ 6 & 3 & 1 & 6 & 4 & 8 \\ 1 & 1 & 2 & 4 & 4 & 3 \\ \hline 12 & 6 & 4 & 15 & 10 & 15 \\ 10 & 4 & 4 & 10 & 15 & 18 \\ \hline 21 & 8 & 3 & 15 & 18 & 30 \end{array}\right] +$$

$$\left[\begin{array}{ccc|cc|c} 4 & 2 & 2 & 10 & 0 & 4 \\ 2 & 1 & 1 & 5 & 0 & 2 \\ 2 & 1 & 1 & 5 & 0 & 2 \\ \hline 10 & 5 & 5 & 25 & 0 & 10 \\ 0 & 0 & 0 & 0 & 0 & 0 \\ \hline 4 & 2 & 2 & 10 & 0 & 4 \end{array}\right]$$

$$= \left[\begin{array}{ccc|cc|c} 36 & 11 & 8 & 27 & 12 & 30 \\ 11 & 6 & 5 & 14 & 5 & 13 \\ 8 & 5 & 8 & 14 & 6 & 10 \\ \hline 27 & 14 & 14 & 45 & 12 & 30 \\ 12 & 5 & 6 & 12 & 16 & 20 \\ \hline 30 & 13 & 10 & 30 & 20 & 39 \end{array}\right].$$

Now we proceed on to illustrate the major product of type IV vectors

***Example 1.1.29:*** Let

$$X = \left[\begin{array}{c} \left[\begin{array}{c|cc|ccc} 1 & 2 & 1 & 1 & 2 & 3 \\ 3 & 1 & 2 & 3 & 1 & 1 \end{array}\right] \\ \left[\begin{array}{c|cc|ccc} 1 & 1 & 3 & 1 & 1 & 1 \\ 2 & 3 & 1 & 2 & 0 & 1 \\ 3 & 4 & 2 & 0 & 1 & 0 \\ 4 & 2 & 4 & 1 & 0 & 0 \end{array}\right] \\ \left[\begin{array}{c|cc|ccc} 5 & 0 & 1 & 1 & 1 & 1 \end{array}\right] \end{array}\right]$$

and

$$Y = \left[\begin{bmatrix} 1 & 1 & 2 & 1 \\ \hline 1 & 0 & 2 & 4 \\ 0 & 1 & 0 & 3 \\ \hline 1 & 1 & 0 & 0 \\ 1 & 0 & 1 & 1 \\ 0 & 1 & 0 & 1 \end{bmatrix} \begin{bmatrix} 2 & 1 \\ \hline 3 & 1 \\ 1 & 0 \\ \hline 2 & 1 \\ 1 & 2 \\ 1 & 1 \end{bmatrix} \begin{bmatrix} 3 & 1 & 0 \\ \hline 4 & 1 & 1 \\ 1 & 2 & 1 \\ \hline 1 & 2 & 1 \\ 2 & 1 & 2 \\ 1 & -1 & 0 \end{bmatrix}\right].$$

Now we find the major product of XY. The product of the first row of X with first column of Y gives

$$\left[\begin{array}{c|cc|ccc} 1 & 2 & 1 & 1 & 2 & 3 \\ 3 & 1 & 2 & 3 & 1 & 1 \end{array}\right] \begin{bmatrix} 1 & 1 & 2 & 1 \\ \hline 1 & 0 & 2 & 4 \\ 0 & 1 & 0 & 3 \\ \hline 1 & 1 & 0 & 0 \\ 1 & 0 & 1 & 1 \\ 0 & 1 & 0 & 1 \end{bmatrix}$$

$$= \begin{bmatrix} 1 \\ 3 \end{bmatrix} \begin{bmatrix} 1 & 1 & 2 & 1 \end{bmatrix} + \begin{bmatrix} 2 & 1 \\ 1 & 2 \end{bmatrix} \begin{bmatrix} 1 & 0 & 2 & 4 \\ 0 & 1 & 0 & 3 \end{bmatrix} +$$

$$\begin{bmatrix} 1 & 2 & 3 \\ 3 & 1 & 1 \end{bmatrix} \begin{bmatrix} 1 & 1 & 0 & 0 \\ 1 & 0 & 1 & 1 \\ 0 & 1 & 0 & 1 \end{bmatrix}$$

$$= \begin{bmatrix} 1 & 1 & 2 & 1 \\ 3 & 3 & 6 & 3 \end{bmatrix} + \begin{bmatrix} 2 & 1 & 4 & 11 \\ 1 & 2 & 2 & 10 \end{bmatrix} + \begin{bmatrix} 3 & 4 & 2 & 5 \\ 4 & 4 & 1 & 2 \end{bmatrix}$$

$$= \begin{bmatrix} 6 & 6 & 8 & 17 \\ 8 & 9 & 9 & 15 \end{bmatrix}.$$

Now

$$\left[\begin{array}{c|cc|ccc} 1 & 2 & 1 & 1 & 2 & 3 \\ 3 & 1 & 2 & 3 & 1 & 1 \end{array}\right]\left[\begin{array}{cc} 2 & 1 \\ \hline 3 & 1 \\ 1 & 0 \\ \hline 2 & 1 \\ 1 & 2 \\ 1 & 1 \end{array}\right]$$

$$= \begin{bmatrix} 1 \\ 3 \end{bmatrix}\begin{bmatrix} 2 & 1 \end{bmatrix} + \begin{bmatrix} 2 & 1 \\ 1 & 2 \end{bmatrix}\begin{bmatrix} 3 & 1 \\ 1 & 0 \end{bmatrix} + \begin{bmatrix} 1 & 2 & 3 \\ 3 & 1 & 1 \end{bmatrix}\begin{bmatrix} 2 & 1 \\ 1 & 2 \\ 1 & 1 \end{bmatrix}$$

$$= \begin{bmatrix} 2 & 1 \\ 6 & 3 \end{bmatrix} + \begin{bmatrix} 7 & 2 \\ 5 & 1 \end{bmatrix} + \begin{bmatrix} 7 & 8 \\ 8 & 6 \end{bmatrix}$$

$$= \begin{bmatrix} 16 & 11 \\ 19 & 10 \end{bmatrix}.$$

Consider the product of first row with the 3[rd] column.

$$\left[\begin{array}{c|cc|ccc} 1 & 2 & 1 & 1 & 2 & 3 \\ 3 & 1 & 2 & 3 & 1 & 1 \end{array}\right]\left[\begin{array}{ccc} 3 & 1 & 0 \\ \hline 4 & 1 & 1 \\ 1 & 2 & 1 \\ \hline 1 & 2 & 1 \\ 2 & 1 & 2 \\ 1 & -1 & 0 \end{array}\right]$$

$$= \begin{bmatrix} 1 \\ 3 \end{bmatrix}\begin{bmatrix} 3 & 1 & 0 \end{bmatrix} + \begin{bmatrix} 2 & 1 \\ 1 & 2 \end{bmatrix}\begin{bmatrix} 4 & 1 & 1 \\ 1 & 2 & 1 \end{bmatrix} + \begin{bmatrix} 1 & 2 & 3 \\ 3 & 1 & 1 \end{bmatrix}\begin{bmatrix} 1 & 2 & 1 \\ 2 & 1 & 2 \\ 1 & -1 & 0 \end{bmatrix}$$

$$= \begin{bmatrix} 3 & 1 & 0 \\ 9 & 3 & 0 \end{bmatrix} + \begin{bmatrix} 9 & 4 & 3 \\ 6 & 5 & 3 \end{bmatrix} + \begin{bmatrix} 8 & 1 & 5 \\ 6 & 6 & 5 \end{bmatrix}$$

$$= \begin{bmatrix} 20 & 6 & 8 \\ 21 & 14 & 8 \end{bmatrix}.$$

The product of 2$^{nd}$ row of X with first column of Y gives

$$\left[\begin{array}{c|cc|ccc} 1 & 1 & 3 & 1 & 1 & 1 \\ 2 & 3 & 1 & 2 & 0 & 1 \\ 3 & 4 & 2 & 0 & 1 & 0 \\ 4 & 2 & 4 & 1 & 0 & 0 \end{array}\right] \left[\begin{array}{cccc} 1 & 1 & 2 & 1 \\ \hline 1 & 0 & 2 & 4 \\ 0 & 1 & 0 & 3 \\ \hline 1 & 1 & 0 & 0 \\ 1 & 0 & 1 & 1 \\ 0 & 1 & 0 & 1 \end{array}\right] =$$

$$\begin{bmatrix} 1 \\ 2 \\ 3 \\ 4 \end{bmatrix} \begin{bmatrix} 1 & 1 & 2 & 1 \end{bmatrix} + \begin{bmatrix} 1 & 3 \\ 3 & 1 \\ 4 & 2 \\ 2 & 4 \end{bmatrix} \begin{bmatrix} 1 & 0 & 2 & 4 \\ 0 & 1 & 0 & 3 \end{bmatrix} + \begin{bmatrix} 1 & 1 & 1 \\ 2 & 0 & 1 \\ 0 & 1 & 0 \\ 1 & 0 & 0 \end{bmatrix} \begin{bmatrix} 1 & 1 & 0 & 0 \\ 1 & 0 & 1 & 1 \\ 0 & 1 & 0 & 1 \end{bmatrix}$$

$$= \begin{bmatrix} 1 & 1 & 2 & 1 \\ 2 & 2 & 4 & 2 \\ 3 & 3 & 6 & 3 \\ 4 & 4 & 8 & 4 \end{bmatrix} + \begin{bmatrix} 1 & 3 & 2 & 13 \\ 3 & 1 & 6 & 15 \\ 4 & 2 & 8 & 22 \\ 2 & 4 & 4 & 20 \end{bmatrix} + \begin{bmatrix} 2 & 2 & 1 & 2 \\ 2 & 3 & 0 & 1 \\ 1 & 0 & 1 & 1 \\ 1 & 1 & 0 & 0 \end{bmatrix}$$

$$= \begin{bmatrix} 4 & 6 & 5 & 16 \\ 7 & 6 & 10 & 18 \\ 8 & 5 & 15 & 26 \\ 7 & 9 & 12 & 24 \end{bmatrix}.$$

The product of 3$^{rd}$ row of X with the 3$^{rd}$ column of Y.

$$\left[5 \mid 0 \quad 1 \mid 1 \quad 1 \quad 1\right]\begin{bmatrix} 3 & 1 & 0 \\ \hline 4 & 1 & 1 \\ 1 & 2 & 1 \\ \hline 1 & 2 & 1 \\ 2 & 1 & 2 \\ 1 & -1 & 0 \end{bmatrix}$$

$$= \quad [5][3 \quad 1 \quad 0] + [0 \quad 1]\begin{bmatrix} 4 & 1 & 1 \\ 1 & 2 & 1 \end{bmatrix} + [1 \quad 1 \quad 1]\begin{bmatrix} 1 & 2 & 1 \\ 2 & 1 & 2 \\ 1 & -1 & 0 \end{bmatrix}$$

= [15 5 0] + [1 2 1] + [4 2 3]
= [20 9 4].

The product of second row of X with second column of Y.

$$\left[\begin{array}{c|cc|ccc} 1 & 1 & 3 & 1 & 1 & 1 \\ 2 & 3 & 1 & 2 & 0 & 1 \\ 3 & 4 & 2 & 0 & 1 & 0 \\ 4 & 2 & 4 & 1 & 0 & 0 \end{array}\right]\begin{bmatrix} 2 & 1 \\ \hline 3 & 1 \\ 1 & 0 \\ \hline 2 & 1 \\ 1 & 2 \\ 1 & 1 \end{bmatrix}$$

$$= \begin{bmatrix} 1 \\ 2 \\ 3 \\ 4 \end{bmatrix}[2 \quad 1] + \begin{bmatrix} 1 & 3 \\ 3 & 1 \\ 4 & 2 \\ 2 & 4 \end{bmatrix}\begin{bmatrix} 3 & 1 \\ 1 & 0 \end{bmatrix} + \begin{bmatrix} 1 & 1 & 1 \\ 2 & 0 & 1 \\ 0 & 1 & 0 \\ 1 & 0 & 0 \end{bmatrix}\begin{bmatrix} 2 & 1 \\ 1 & 2 \\ 1 & 1 \end{bmatrix}$$

$$= \begin{bmatrix} 2 & 1 \\ 4 & 2 \\ 6 & 3 \\ 8 & 4 \end{bmatrix} + \begin{bmatrix} 6 & 1 \\ 10 & 3 \\ 14 & 4 \\ 10 & 2 \end{bmatrix} + \begin{bmatrix} 4 & 4 \\ 5 & 3 \\ 1 & 2 \\ 2 & 1 \end{bmatrix}$$

$$= \begin{bmatrix} 12 & 6 \\ 19 & 8 \\ 21 & 9 \\ 20 & 7 \end{bmatrix}.$$

The product of the 2nd row with the last column of Y.

$$\left[\begin{array}{c|cc|ccc} 1 & 1 & 3 & 1 & 1 & 1 \\ 2 & 3 & 1 & 2 & 0 & 1 \\ 3 & 4 & 2 & 0 & 1 & 0 \\ 4 & 2 & 4 & 1 & 0 & 0 \end{array}\right] \left[\begin{array}{ccc} 3 & 1 & 0 \\ \hline 4 & 1 & 1 \\ 1 & 2 & 1 \\ \hline 1 & 2 & 1 \\ 2 & 1 & 2 \\ 1 & -1 & 0 \end{array}\right]$$

$$= \begin{bmatrix} 1 \\ 2 \\ 3 \\ 4 \end{bmatrix} \begin{bmatrix} 3 & 1 & 0 \end{bmatrix} + \begin{bmatrix} 1 & 3 \\ 3 & 1 \\ 4 & 2 \\ 2 & 4 \end{bmatrix} \begin{bmatrix} 4 & 1 & 1 \\ 1 & 2 & 1 \end{bmatrix} + \begin{bmatrix} 1 & 1 & 1 \\ 2 & 0 & 1 \\ 0 & 1 & 0 \\ 1 & 0 & 0 \end{bmatrix} \begin{bmatrix} 1 & 2 & 1 \\ 2 & 1 & 2 \\ 1 & -1 & 0 \end{bmatrix}$$

$$= \begin{bmatrix} 3 & 1 & 0 \\ 6 & 2 & 0 \\ 9 & 3 & 0 \\ 12 & 4 & 0 \end{bmatrix} + \begin{bmatrix} 7 & 7 & 4 \\ 13 & 5 & 4 \\ 18 & 8 & 6 \\ 12 & 10 & 6 \end{bmatrix} + \begin{bmatrix} 4 & 2 & 3 \\ 3 & 3 & 2 \\ 2 & 1 & 2 \\ 1 & 2 & 1 \end{bmatrix}$$

$$= \begin{bmatrix} 14 & 10 & 7 \\ 22 & 10 & 6 \\ 29 & 12 & 8 \\ 25 & 16 & 7 \end{bmatrix}.$$

The product of 3rd row of X with 1st column of Y

$$\left[5 \mid 0 \quad 1 \mid 1 \quad 1 \quad 1\right]\left[\begin{array}{cccc} 1 & 1 & 2 & 1 \\ \hline 1 & 0 & 2 & 4 \\ 0 & 1 & 0 & 3 \\ \hline 1 & 1 & 0 & 0 \\ 1 & 0 & 1 & 1 \\ 0 & 1 & 0 & 1 \end{array}\right] =$$

$$[5][1 \quad 1 \quad 2 \quad 1] + [0 \quad 1]\begin{bmatrix} 1 & 0 & 2 & 4 \\ 0 & 1 & 0 & 3 \end{bmatrix} + [1 \quad 1 \quad 1]\begin{bmatrix} 1 & 1 & 0 & 0 \\ 1 & 0 & 1 & 1 \\ 0 & 1 & 0 & 1 \end{bmatrix}$$

= [5 5 10 5] + [0 1 0 3] + [2 2 1 2]

= [7 8 11 10].

The product of 3rd row of X with 2nd column of Y.

$$\left[5 \mid 0 \quad 1 \mid 1 \quad 1 \quad 1\right]\left[\begin{array}{cc} 2 & 1 \\ \hline 3 & 1 \\ 1 & 0 \\ \hline 2 & 1 \\ 1 & 2 \\ 1 & 1 \end{array}\right] =$$

$$[5]\ [2 \ 1] + [0 \ 1]\begin{bmatrix} 3 & 1 \\ 1 & 0 \end{bmatrix} + [1 \quad 1 \quad 1]\begin{bmatrix} 2 & 1 \\ 1 & 2 \\ 1 & 1 \end{bmatrix}$$

= [10 5] + [1 0] + [4 4]

= [15 9].

$$XY = \left[\begin{array}{cccc|cc|ccc} 6 & 6 & 8 & 17 & 16 & 11 & 20 & 6 & 8 \\ 8 & 8 & 9 & 15 & 19 & 10 & 21 & 14 & 8 \\ \hline 4 & 6 & 5 & 16 & 12 & 6 & 14 & 10 & 7 \\ 7 & 6 & 10 & 18 & 19 & 8 & 22 & 10 & 6 \\ 8 & 5 & 15 & 26 & 21 & 9 & 29 & 12 & 8 \\ 7 & 9 & 12 & 24 & 20 & 7 & 25 & 16 & 7 \\ \hline 7 & 8 & 11 & 10 & 15 & 9 & 20 & 9 & 4 \end{array}\right]_{7\times 9}$$

On similar lines we can find the transpose of major product of Type IV vectors.

Now we proceed on to just show the major product moment of a type IV vector.

***Example 1.1.30:*** Suppose

$$X = \left[\begin{array}{c|cc|ccc} 1 & 2 & 1 & 3 & 2 & 1 \\ 2 & 3 & 1 & 2 & 1 & 2 \\ \hline 1 & 4 & 2 & 3 & 2 & 2 \\ 4 & 1 & 3 & 2 & 1 & 1 \\ 2 & 3 & 2 & 3 & 2 & 3 \\ 3 & 4 & 1 & 1 & 4 & 2 \\ \hline 2 & 1 & 2 & 2 & 1 & 3 \end{array}\right]$$

and

$$X^t = \left[\begin{array}{cc|cccc|c} 1 & 2 & 1 & 4 & 2 & 3 & 2 \\ \hline 2 & 3 & 4 & 1 & 3 & 4 & 1 \\ 1 & 1 & 2 & 3 & 2 & 1 & 2 \\ \hline 3 & 2 & 3 & 2 & 3 & 1 & 2 \\ 2 & 1 & 2 & 1 & 2 & 4 & 1 \\ 1 & 2 & 2 & 1 & 3 & 2 & 3 \end{array}\right].$$

$$X^t X = \left[\begin{array}{cc|cccc|c} 1 & 2 & 1 & 4 & 2 & 3 & 2 \\ \hline 2 & 3 & 4 & 1 & 3 & 4 & 1 \\ 1 & 1 & 2 & 3 & 2 & 1 & 2 \\ \hline 3 & 2 & 3 & 2 & 3 & 1 & 2 \\ 2 & 1 & 2 & 1 & 2 & 4 & 1 \\ 1 & 2 & 2 & 1 & 3 & 2 & 3 \end{array}\right] \times$$

$$\left[\begin{array}{c|cc|ccc} 1 & 2 & 1 & 3 & 2 & 1 \\ 2 & 3 & 1 & 2 & 1 & 2 \\ \hline 1 & 4 & 2 & 3 & 2 & 2 \\ 4 & 1 & 3 & 2 & 1 & 1 \\ 2 & 3 & 2 & 3 & 2 & 3 \\ 3 & 4 & 1 & 1 & 4 & 2 \\ \hline 2 & 1 & 2 & 2 & 1 & 3 \end{array}\right].$$

Product of 1st row of $X^t$ with 1st column of X

$$\left[1 \quad 2 \mid 1 \quad 4 \quad 2 \quad 3 \mid 2\right] \left[\begin{array}{c} 1 \\ 2 \\ \hline 1 \\ 4 \\ 2 \\ 3 \\ \hline 2 \end{array}\right]$$

$$= \quad [1 \quad 2]\begin{bmatrix} 1 \\ 2 \end{bmatrix} + [1 \quad 4 \quad 2 \quad 3]\begin{bmatrix} 1 \\ 4 \\ 2 \\ 3 \end{bmatrix} + [2][2]$$

$= \quad 5 + 30 + 4$

$= \quad 39.$

Product of 1st row of $X^t$ with 2nd column of X.

$$\left[1 \quad 2 \mid 1 \quad 4 \quad 2 \quad 3 \mid 2\right]\left[\begin{array}{cc} 2 & 1 \\ 3 & 1 \\ \hline 4 & 2 \\ 1 & 3 \\ 3 & 2 \\ 4 & 1 \\ \hline 1 & 2 \end{array}\right]$$

$$= \begin{bmatrix}1 & 2\end{bmatrix}\begin{bmatrix}2 & 1 \\ 3 & 1\end{bmatrix} + \begin{bmatrix}1 & 4 & 2 & 3\end{bmatrix}\begin{bmatrix}4 & 2 \\ 1 & 3 \\ 3 & 2 \\ 4 & 1\end{bmatrix} + \begin{bmatrix}2\end{bmatrix}\begin{bmatrix}1 & 2\end{bmatrix}$$

$= [8 \;\; 3] + [26 \;\; 21] + [2 \;\; 4]$

$= [36 \;\; 28]$.

The product of 1st row of $X^t$ with 3rd column of X.

$$\left[1 \quad 2 \mid 1 \quad 4 \quad 2 \quad 3 \mid 2\right]\left[\begin{array}{ccc} 3 & 2 & 1 \\ 2 & 1 & 2 \\ \hline 3 & 2 & 2 \\ 2 & 1 & 1 \\ 3 & 2 & 3 \\ 1 & 4 & 2 \\ \hline 2 & 1 & 3 \end{array}\right]$$

$$= \begin{bmatrix}1 & 2\end{bmatrix}\begin{bmatrix}3 & 2 & 1 \\ 2 & 1 & 2\end{bmatrix} + \begin{bmatrix}1 & 4 & 2 & 3\end{bmatrix}\begin{bmatrix}3 & 2 & 2 \\ 2 & 1 & 1 \\ 3 & 2 & 3 \\ 1 & 4 & 2\end{bmatrix} + \begin{bmatrix}2\end{bmatrix}\begin{bmatrix}2 & 1 & 3\end{bmatrix}$$

$$= \quad [7\ 4\ 5] + [20\ \ 22\ \ 18] + [4\ \ 2\ \ 6]$$

$$= \quad [31\ \ 28\ \ 29].$$

The product of 2nd row of $X^t$ with 1st column of $X$.

$$\left[\begin{array}{cc|cccc|c} 2 & 3 & 4 & 1 & 3 & 4 & 1 \\ 1 & 1 & 2 & 3 & 2 & 1 & 2 \end{array}\right] \left[\begin{array}{c} 1 \\ 2 \\ \hline 1 \\ 4 \\ 2 \\ 3 \\ \hline 2 \end{array}\right]$$

$$= \begin{bmatrix} 2 & 3 \\ 1 & 1 \end{bmatrix} \begin{bmatrix} 1 \\ 2 \end{bmatrix} + \begin{bmatrix} 4 & 1 & 3 & 4 \\ 2 & 3 & 2 & 1 \end{bmatrix} \begin{bmatrix} 1 \\ 4 \\ 2 \\ 3 \end{bmatrix} + \begin{bmatrix} 1 \\ 2 \end{bmatrix} [2]$$

$$= \begin{bmatrix} 8 \\ 3 \end{bmatrix} + \begin{bmatrix} 26 \\ 21 \end{bmatrix} + \begin{bmatrix} 2 \\ 4 \end{bmatrix} = \begin{bmatrix} 36 \\ 28 \end{bmatrix}.$$

The product of 2nd row of $X^t$ with 2nd column of $X$.

$$\left[\begin{array}{cc|cccc|c} 2 & 3 & 4 & 1 & 3 & 4 & 1 \\ 1 & 1 & 2 & 3 & 2 & 1 & 2 \end{array}\right] \left[\begin{array}{cc} 2 & 1 \\ 3 & 1 \\ \hline 4 & 2 \\ 1 & 3 \\ 3 & 2 \\ 4 & 1 \\ \hline 1 & 2 \end{array}\right]$$

$$= \begin{bmatrix} 2 & 3 \\ 1 & 1 \end{bmatrix} \begin{bmatrix} 2 & 1 \\ 3 & 1 \end{bmatrix} + \begin{bmatrix} 4 & 1 & 3 & 4 \\ 2 & 3 & 2 & 1 \end{bmatrix} \begin{bmatrix} 4 & 2 \\ 1 & 3 \\ 3 & 2 \\ 4 & 1 \end{bmatrix} + \begin{bmatrix} 1 \\ 2 \end{bmatrix} \begin{bmatrix} 1 & 2 \end{bmatrix}$$

$$= \begin{bmatrix} 13 & 5 \\ 5 & 2 \end{bmatrix} + \begin{bmatrix} 42 & 21 \\ 21 & 18 \end{bmatrix} + \begin{bmatrix} 1 & 2 \\ 2 & 4 \end{bmatrix}$$

$$= \begin{bmatrix} 56 & 28 \\ 28 & 24 \end{bmatrix}.$$

The product of 2nd row of $X^t$ with 3rd column of X.

$$\left[\begin{array}{cc|cccc|c} 2 & 3 & 4 & 1 & 3 & 4 & 1 \\ 1 & 1 & 2 & 3 & 2 & 1 & 2 \end{array}\right] \left[\begin{array}{ccc} 3 & 2 & 1 \\ 2 & 1 & 2 \\ \hline 3 & 2 & 2 \\ 2 & 1 & 1 \\ 3 & 2 & 3 \\ 1 & 4 & 2 \\ \hline 2 & 1 & 3 \end{array}\right]$$

$$= \begin{bmatrix} 2 & 3 \\ 1 & 1 \end{bmatrix} \begin{bmatrix} 3 & 2 & 1 \\ 2 & 1 & 2 \end{bmatrix} + \begin{bmatrix} 4 & 1 & 3 & 4 \\ 2 & 3 & 2 & 1 \end{bmatrix} \begin{bmatrix} 3 & 2 & 2 \\ 2 & 1 & 1 \\ 3 & 2 & 3 \\ 1 & 4 & 2 \end{bmatrix} + \begin{bmatrix} 1 \\ 2 \end{bmatrix} \begin{bmatrix} 2 & 1 & 3 \end{bmatrix}$$

$$= \begin{bmatrix} 12 & 7 & 8 \\ 5 & 3 & 3 \end{bmatrix} + \begin{bmatrix} 27 & 31 & 26 \\ 19 & 15 & 15 \end{bmatrix} + \begin{bmatrix} 2 & 1 & 3 \\ 4 & 2 & 6 \end{bmatrix}$$

$$= \begin{bmatrix} 41 & 39 & 37 \\ 28 & 20 & 24 \end{bmatrix}.$$

The product of $3^{rd}$ row of $X^t$ with $1^{st}$ column of X.

$$\left[\begin{array}{cc|cccc|c} 3 & 2 & 3 & 2 & 3 & 1 & 2 \\ 2 & 1 & 2 & 1 & 2 & 4 & 1 \\ 1 & 2 & 2 & 1 & 3 & 2 & 3 \end{array}\right]\left[\begin{array}{c} 1 \\ 2 \\ \hline 1 \\ 4 \\ 2 \\ 3 \\ \hline 2 \end{array}\right]$$

$$= \begin{bmatrix} 7 \\ 4 \\ 5 \end{bmatrix} + \begin{bmatrix} 20 \\ 22 \\ 18 \end{bmatrix} + \begin{bmatrix} 4 \\ 2 \\ 6 \end{bmatrix} = \begin{bmatrix} 31 \\ 28 \\ 29 \end{bmatrix}.$$

The product of $3^{rd}$ row of $X^t$ with $2^{nd}$ column of X.

$$\left[\begin{array}{cc|cccc|c} 3 & 2 & 3 & 2 & 3 & 1 & 2 \\ 2 & 1 & 2 & 1 & 2 & 4 & 1 \\ 1 & 2 & 2 & 1 & 3 & 2 & 3 \end{array}\right]\left[\begin{array}{cc} 2 & 1 \\ 3 & 1 \\ \hline 4 & 2 \\ 1 & 3 \\ 3 & 2 \\ 4 & 1 \\ \hline 1 & 2 \end{array}\right]$$

$$= \begin{bmatrix} 12 & 5 \\ 7 & 3 \\ 8 & 3 \end{bmatrix} + \begin{bmatrix} 27 & 19 \\ 31 & 15 \\ 26 & 15 \end{bmatrix} + \begin{bmatrix} 2 & 4 \\ 1 & 2 \\ 3 & 6 \end{bmatrix}$$

$$= \begin{bmatrix} 41 & 28 \\ 39 & 20 \\ 37 & 24 \end{bmatrix}.$$

The product 3rd row of $X^t$ with 3rd column of X.

$$\left[\begin{array}{cc|cccc|c} 3 & 2 & 3 & 2 & 3 & 1 & 2 \\ 2 & 1 & 2 & 1 & 2 & 4 & 1 \\ 1 & 2 & 2 & 1 & 3 & 2 & 3 \end{array}\right] \left[\begin{array}{ccc} 3 & 2 & 1 \\ 2 & 1 & 2 \\ \hline 3 & 2 & 2 \\ 2 & 1 & 1 \\ 3 & 2 & 3 \\ 1 & 4 & 2 \\ \hline 2 & 1 & 3 \end{array}\right]$$

$$= \begin{bmatrix} 13 & 8 & 7 \\ 8 & 5 & 4 \\ 7 & 4 & 5 \end{bmatrix} + \begin{bmatrix} 23 & 18 & 19 \\ 18 & 25 & 19 \\ 19 & 19 & 18 \end{bmatrix} + \begin{bmatrix} 4 & 2 & 6 \\ 2 & 1 & 3 \\ 6 & 3 & 9 \end{bmatrix}$$

$$= \begin{bmatrix} 40 & 28 & 32 \\ 28 & 31 & 26 \\ 32 & 26 & 32 \end{bmatrix}.$$

$$X^t X = \left[\begin{array}{c|cc|ccc} 39 & 36 & 28 & 31 & 28 & 29 \\ \hline 36 & 56 & 28 & 41 & 39 & 37 \\ 28 & 28 & 24 & 28 & 20 & 24 \\ \hline 31 & 41 & 28 & 40 & 28 & 32 \\ 28 & 39 & 20 & 28 & 31 & 26 \\ 29 & 37 & 24 & 32 & 26 & 32 \end{array}\right].$$

On similar lines interested reader can find the major product moment of type IV column vector.

For more about supermatrix and their properties please refer [92]. We have not given the general method for the notations seems to be little difficult as given by this book [92]. They have also defined new form of products of supermatrices etc. Most of the properties pertain to applications to models mostly needed by social scientists.

## 1.2 Introduction to Fuzzy Matrices

Just we recall the definition of fuzzy matrix. A fuzzy matrix is a matrix which has its elements from [0, 1]. We as in case of matrix have rectangular fuzzy matrix, fuzzy square matrix, fuzzy row matrix and fuzzy column matrix. We just illustrate them by the following examples.

***Examples 1.2.1:*** Let A = [1 0.4 0.6 0 1 0.7 0.1]; A is a row fuzzy matrix. In fact A is a $1 \times 7$ row fuzzy matrix.

$A_1$ = [0 1] is a $1 \times 2$ row fuzzy matrix.

$A_2$ = [1 0 1 0 0.6 0.2 110] is a $1 \times 9$ row fuzzy matrix.

***Example 1.2.2:*** Let

$$B = \begin{bmatrix} 1 \\ 0 \\ 0.14 \\ 0 \\ 1 \\ 1 \end{bmatrix},$$

B is $6 \times 1$ column fuzzy matrix.

$$B_1 = \begin{bmatrix} 1 \\ 0 \end{bmatrix}$$

is a $2 \times 1$ column fuzzy matrix.

$$B_2 = \begin{bmatrix} 0.06 \\ 0.2 \\ 0.14 \\ 0 \\ 1 \\ 0.03 \\ 0.12 \\ 0.31 \\ 0 \end{bmatrix}$$

is a $9 \times 1$ column fuzzy matrix.

Now having seen a fuzzy column matrix and fuzzy row matrix we proceed on to define the fuzzy square matrix.

***Example 1.2.3:*** Consider the fuzzy matrix

$$A = \begin{bmatrix} 0 & 0.1 & 1 & 0.6 & 0.7 \\ 0.4 & 1 & 0 & 1 & 0 \\ 1 & 0.5 & 0.6 & 0.2 & 1 \\ 0 & 1 & 0 & 0.1 & 0 \\ 0.2 & 0.6 & 1 & 1 & 0.2 \end{bmatrix}$$

is a square fuzzy matrix. Infact A is a 5 × 5 square fuzzy matrix.

***Example 1.2.4:*** Consider the fuzzy matrix

$$B = \begin{bmatrix} 1 & 0 & 0 \\ 0 & 1 & 1 \\ 1 & 1 & 1 \end{bmatrix},$$

B is a square 3 × 3 fuzzy matrix.

***Example 1.2.5:*** Consider the fuzzy matrix

$$C = \begin{bmatrix} 1 & 0.3 \\ 0.1 & 0.04 \end{bmatrix},$$

C is a 2 × 2 fuzzy square matrix.

Now we proceed into illustrate fuzzy rectangular matrices.

***Example 1.2.6:*** Let

$$A = \begin{bmatrix} 0 & 0.3 & 1 & 0.2 \\ 0.1 & 1 & 0 & 0.6 \end{bmatrix}$$

A is a 2 × 4 fuzzy rectangular matrix.

Let

$$T = \begin{bmatrix} 0 & 1 \\ 0.3 & 0.1 \\ 1 & 1 \\ 0.8 & 0.71 \\ 0.5 & 0.11 \\ 0 & 0.1 \end{bmatrix};$$

T is a 6 × 2 fuzzy rectangular matrix.

$$V = \begin{bmatrix} 0.3 & 1 & 0 & 1 \\ 1 & 0.31 & 0.6 & 0.7 \\ 0 & 0.11 & 0.2 & 0 \\ 0.2 & 0.5 & 0.1 & 1 \\ 0.1 & 0.14 & 1 & 0.8 \end{bmatrix}$$

is a 5 × 4 fuzzy rectangular matrix.

Now we define these in more generality.

Let

$$A = [a_1\ a_2\ \ldots\ a_n]$$

where $a_i \in [0, 1]$, i = 1, 2, …, n; A is a 1 × n fuzzy row matrix.

Let

$$B = \begin{bmatrix} b_1 \\ b_2 \\ \vdots \\ b_m \end{bmatrix},$$

$b_j \in [0, 1]$; j = 1, 2, …, m; B is a m × 1 fuzzy column matrix.

Let

$$C = \begin{bmatrix} a_{11} & a_{12} & \cdots & a_{1n} \\ a_{21} & a_{22} & \cdots & a_{2n} \\ \vdots & \vdots & & \vdots \\ a_{n1} & a_{n2} & \cdots & a_{nn} \end{bmatrix}$$

$a_{ij} \in [0, 1]$; $1 \le i, j \le n$. C is a $n \times n$ fuzzy square matrix.

Let

$$D = \begin{bmatrix} a_{11} & a_{12} & \cdots & a_{1n} \\ a_{21} & a_{22} & \cdots & a_{2n} \\ \vdots & \vdots & & \vdots \\ a_{m1} & a_{m2} & \cdots & a_{mn} \end{bmatrix}$$

$a_{ij} \in [0, 1]$, $1 \le i \le m$; $1 \le j \le n$. D is a $m \times n$ rectangular fuzzy matrix $m \neq n$.

*Note:* We may work with any fuzzy matrix.

It is still important to note while doing fuzzy mathematical models the fuzzy matrix may take its entries from the interval [–1, 1] then also they are known as fuzzy matrices. We mention here that Fuzzy Cognitive Maps (FCMs) model take their entries only from the interval [–1, 1].

The main operations performed on fuzzy matrices are usual max-min operations. Many cases the operation "+" i.e., usual addition may be replaced by 'max' i.e., (0.5 + 0.7) = max (0.5, 0.7) = 0.7.

"×" may be replaced by min i.e., min (0.5, 0.7) = 0.5.

But since matrix multiplication involves both addition and multiplication we make use of the min max or max min principle. Before we proceed on to define the very new notion of fuzzy supermatrix we just illustrate a few operations on the fuzzy matrices.

***Example 1.2.7:*** Let

$$A = \begin{bmatrix} 0.1 & 0 & 0.1 \\ 1 & 0.7 & 0 \\ 0.2 & 0.6 & 0.5 \end{bmatrix}$$

and

$$B = \begin{bmatrix} 0.9 & 1 & 0 & 0.8 \\ 0 & 0.2 & 1 & 0.5 \\ 0.1 & 0 & 0 & 0 \end{bmatrix}$$

i.e., if A = ($a_{ij}$) and b = ($b_{ij}$), then

max min {A, B} = max min {$a_{ij}$, $b_{jk}$} = max {min{($a_{11}$, $b_{11}$), min ($a_{12}$, $b_{21}$), min ($a_{13}$, $b_{31}$)} (for i = 1, k = 1) i.e.,

$$\begin{aligned}
r_{11} &= \max \{\min (0.1, 0.9), \min (0, 0), \min (0.1, 0.1)\} \\
&= \max \{0.1, 0, 0.1\} \\
&= 0.1 \\
r_{12} &= \max \{\min (0.1, 1), \min (0, 0.2), \min (0.1, 0)\} \\
&= \max \{0.1, 0, 0\} \\
&= 0.1 \\
r_{13} &= \max \{\min (0.1, 0), \min (0, 1), \min (0.1, 0)\} \\
&= \max \{0, 0, 0\} \\
&= 0 \\
r_{14} &= \max \{\min (0.1, 8), \min (0, 0.5), \min (0.1, 0)\} \\
&= \max \{0.1, 0, 0\} \\
&= 0.1
\end{aligned}$$

$$\begin{aligned}
r_{21} &= \max \{\min (1, 0.9), \min (0.7, 0), \min (0, 0.1)\} \\
&= \max \{0.9, 0, 0\} \\
&= 0.9 \\
r_{22} &= \max \{\min (1, 1), \min (0.7, 0.2), \min (0, 0)\} \\
&= \max \{1, 0.2, 0\} \\
&= 1 \\
r_{23} &= \max \{\min (1, 0), \min (0.7, 1), \min (0, 0)\} \\
&= \max \{0, 0.7, 0\} \\
&= 0.7
\end{aligned}$$

$$
\begin{aligned}
r_{24} &= \max\{\min(1, 0.8), \min(0.7, 0.5), \min(0, 0)\} \\
&= \max\{0.8, 0.5, 0\} \\
&= 0.8
\end{aligned}
$$

$$
\begin{aligned}
r_{31} &= \max\{\min(0.2, 0.9), \min(0.6, 0), \min(0.5, 0.1)\} \\
&= \max\{0.2, 0, 0.1\} \\
&= 0.2 \\
r_{32} &= \max\{\min(0.2, 1), \min(0.6, 0.2), \min(0.5, 0)\} \\
&= \max\{0.2, 0.2, 0\} \\
&= 0.2 \\
r_{34} &= \max\{\min(0.2, 0.8), \min(0.6, 0.5), \min(0.5, 0)\} \\
&= \max\{0.2, 0.5, 0\} \\
&= 0.5 \\
r_{33} &= \max\{\min(0.2, 0), \min(0.6, 1), \min(0.5, 0)\} \\
&= \max\{0, 0.6, 0\} \\
&= 0.6.
\end{aligned}
$$

Now

$$
R = \begin{bmatrix} 0.1 & 0.1 & 0 & 0.1 \\ 0.9 & 1 & 0.7 & 0.8 \\ 0.2 & 0.2 & 0.6 & 0.5 \end{bmatrix}.
$$

Thus when ever multiplication is compatible we can define on the fuzzy matrices the max min function or operation.

***Example 1.2.8:*** Let

$$
A = \begin{bmatrix} 0.7 & 0.5 & 0 \\ 1 & 0.2 & 0.1 \\ 0 & 0.3 & 1 \end{bmatrix} = (a_{ij})
$$

$A.A = \max\{\min(a_{ik}, a_{kj}) / 1 \le i, j, k \le 3\}$

i.e., $A.A = \{(a'_{ij}) / 1 \le i \le 3, 1 \le j \le 3\}$ using the min max operator.

$$A.A = \begin{bmatrix} 0.7 & 0.5 & 0 \\ 1 & 0.2 & 0.1 \\ 0 & 0.3 & 1 \end{bmatrix} . \begin{bmatrix} 0.7 & 0.5 & 0 \\ 1 & 0.2 & 0.1 \\ 0 & 0.3 & 1 \end{bmatrix}$$

$$= \begin{bmatrix} 0.7 & 0.5 & 0.1 \\ 0.7 & 0.5 & 0.1 \\ 0.3 & 0.3 & 1 \end{bmatrix}.$$

***Example 1.2.9:*** Let

$$A = \begin{bmatrix} 0.2 \\ 1 \\ 0.5 \\ 0 \\ 0.7 \\ 0 \\ 1 \\ 0.4 \end{bmatrix}_{8\times1}$$

be a 8 × 1 fuzzy column matrix and

$$B = [0\ 0.1\ 1\ 1\ 0.8\ 0.6\ 0.7\ 0.2]_{1\times8}$$

is a fuzzy row matrix.

Now

B.A = max {min (0, 0.2), min (0.1, 1), min (1, 0.5), min (1, 0), min (0.8, 0.7), min (0.6, 0), min (0.7, 1), min (0.2, 0.4)}
= max {0, 0.1, 0.5, 0, 0.7, 0, 0.7, 0.2}
= 0.7.

***Example 1.2.10:*** Let

$$A = \begin{bmatrix} 0.2 \\ 1 \\ 0.5 \\ 0 \\ 0.7 \\ 0 \\ 1 \\ 0.4 \end{bmatrix}$$

and

B = [0 0.1 1 1 0.8 0.6 0.7 0.2].

Now

$$\max\min\{A, B\} = \left\{ \begin{bmatrix} 0.2 \\ 1 \\ 0.5 \\ 0 \\ 0.7 \\ 0 \\ 1 \\ 0.4 \end{bmatrix}, \begin{bmatrix} 0 & 0.1 & 1 & 1 & 0.8 & 0.6 & 0.7 & 0.2 \end{bmatrix} \right\}$$

$$= \begin{bmatrix} 0 & 0.1 & 0.2 & 0.2 & 0.2 & 0.2 & 0.2 & 0.2 \\ 0 & 0.1 & 1 & 1 & 0.8 & 0.6 & 0.7 & 0.2 \\ 0 & 0.1 & 0.5 & 0.5 & 0.5 & 0.5 & 0.5 & 0.2 \\ 0 & 0 & 0 & 0 & 0 & 0 & 0 & 0 \\ 0 & 0.1 & 0.7 & 0.7 & 0.7 & 0.6 & 0.7 & 0.2 \\ 0 & 0 & 0 & 0 & 0 & 0 & 0 & 0 \\ 0 & 0.1 & 1 & 1 & 0.8 & 0.6 & 0.7 & 0.2 \\ 0 & 0.1 & 0.4 & 0.4 & 0.4 & 0.4 & 0.4 & 0.2 \end{bmatrix}.$$

***Example 1.2.11:*** Let

$$A = \begin{bmatrix} 0.2 \\ 0 \\ 1 \\ 0.7 \\ 0.3 \\ 1 \end{bmatrix}_{6\times 1}$$

and

$$B = \begin{bmatrix} 0.1 & 1 & 0 & 1 & 0 & 0.7 \\ 1 & 0.2 & 0.4 & 0 & 0.3 & 1 \\ 0.4 & 0 & 0 & 0.1 & 0 & 0.8 \\ 0 & 0 & 1 & 0.3 & 1 & 0.1 \\ 1 & 0 & 0.2 & 1 & 0.1 & 0 \\ 1 & 1 & 0 & 1 & 0 & 1 \end{bmatrix}_{6\times 6}$$

be two fuzzy matrices. The max min {B, A} = $[0.7\ 1\ 0.8\ 1\ 0.7\ 1]_{1\times 6}$ is a fuzzy row matrix.

***Example 1.2.12:*** Let

$$B = [0\ 0.1\ 0.2\ 1\ 0\ 0.5\ 1\ 0.7]_{1\times 8}$$

be a column fuzzy matrix and

$$A = \begin{bmatrix} 0.1 & 1 & 0 \\ 1 & 0 & 0.3 \\ 0.4 & 0.2 & 1 \\ 1 & 1 & 1 \\ 0 & 0 & 0 \\ 0.2 & 0.9 & 0.3 \\ 0.5 & 0.1 & 0.8 \\ 0.9 & 0.6 & 0.5 \end{bmatrix}_{8\times 3}$$

be a fuzzy rectangular matrix. max min {B, A} = $[0.7\ 1\ 1]_{1\times3}$ is a fuzzy row matrix.

***Example 1.2.13:*** Let

$$A = \begin{bmatrix} 0.3 & 1 & 0 \\ 1 & 0.2 & 1 \\ 0.2 & 1 & 0 \\ 0 & 0 & 0.3 \\ 1 & 1 & 0.2 \\ 0.1 & 0.3 & 0 \end{bmatrix}_{6\times3}$$

and

$$B = \begin{bmatrix} 0 & 0.2 & 1 & 0.3 & 0 \\ 0 & 0.1 & 0 & 1 & 0 \\ 0 & 1 & 0.1 & 0 & 1 \end{bmatrix}_{3\times5}$$

be two fuzzy rectangular matrices.

$$\text{max min } \{A, B\} = \begin{bmatrix} 0 & 0.2 & 0.3 & 1 & 0 \\ 0 & 1 & 1 & 0.3 & 1 \\ 0 & 0.2 & 0.2 & 1 & 0 \\ 0 & 0.3 & 0.1 & 0 & 0.3 \\ 0 & 0.2 & 1 & 1 & 0.2 \\ 0 & 0.1 & 0.1 & 0.3 & 0 \end{bmatrix}_{6\times5}$$

is a fuzzy rectangular matrix. Thus we have given some illustrative examples of operations on fuzzy matrices.

**Chapter Two**

# FUZZY SUPERMATRICES AND THEIR PROPERTIES

This chapter has three sections. In the first section we for the first time introduce the notion of fuzzy supermatrices and some operations on them which are essential for the fuzzy supermodels. In the second section we introduce the new notion of pseudo symmetric partition and pseudo symmetric supermatrices. In the third section special operations of fuzzy super special matrices are given.

## 2.1 Fuzzy Supermatrices and their Properties

Now we proceed on to introduce the notion of fuzzy supermatrices and operations on them. Throughout this chapter we consider matrices with entries only from the fuzzy interval [0, 1]. Thus all matrices in this chapter unless we make a specific mention of them will be fuzzy matrices.

**DEFINITION 2.1.1:** *Let $A_s = [A_1 \mid A_2 \mid \ldots \mid A_t]$ $(t > 1)$ where each $A_i$ is a fuzzy row vector (fuzzy row matrix), $i = 1, 2, \ldots, t$. We call A as the fuzzy super row vector or fuzzy super row matrix.*

***Example 2.1.1:*** Let $A_s$ = [0 0.1 0.3 | 0.6 0.7 0.5 1 0 1 0.9], $A_s$ is a fuzzy super row matrix, here t = 2 i.e., A = [$A_1$ | $A_2$] where $A_1$ = [0 0.1 0.3] and $A_2$ = [0.6 0.7 0.5 1 0 1 0.9].

*Note:* Given any fuzzy row matrix A = [$a_1$, $a_2$, …, $a_n$] where $a_i \in$ [0, 1]; i = 1, 2, …, n, we can partition the fuzzy row matrix A to form a fuzzy super row matrix. For instance if we have only two partitions i.e., A = [$a_1$ $a_2$ … $a_t$ | $a_{t+1}$ … $a_n$];

A is a fuzzy super row matrix with

$$
\begin{aligned}
A &= [A_1 \mid A_2] = [A_1\ A_2] \text{ where} \\
A_1 &= [a_1\ a_2 \dots a_t] \text{ and} \\
A_2 &= [a_{t+1}\ a_{t+2} \dots a_n]\,. \\
A &= [a_1\ a_2 \dots a_r \mid a_{r+1} \dots a_{r+s} \mid a_{r+s+1} \dots a_n] \\
&= [A_1 \mid A_2 \mid A_3] = [A_1\ A_2\ A_3]
\end{aligned}
$$

is also a fuzzy super row matrix got by partitioning the fuzzy row vector A, here $A_1$ = [$a_1$ … $a_r$] is a fuzzy row matrix, $A_2$ = [$a_{r+1}$ … $a_{r+s}$] and $A_3$ = [$a_{r+s+1}$ … $a_n$] are all fuzzy row matrices.

Thus a super fuzzy row matrix can be obtained by partitioning a fuzzy row matrix.

***Example 2.1.2:*** Let A = [0 1 0.3 1 1 0 0.5 0.3 0.4 1 1 0 0.3 0.2] be a fuzzy row matrix.

A = [0 1 0.3 | 1 1 0 0.5 0.3 | 0.4 1 1 0 | 0.3 0.2]

is a super fuzzy row matrix.

Now we proceed on to define super fuzzy column matrix.

**DEFINITION 2.1.2:** *Let*

$$A_s = \begin{bmatrix} A_1 \\ A_2 \\ \vdots \\ A_N \end{bmatrix}$$

*where N > 1 and each $A_i$ be a fuzzy column matrix for i = 1, 2, …, N; then $A_s$ is defined to be the super fuzzy column matrix.*

***Example 2.1.3:*** Let

$$A_s = \begin{bmatrix} A_1 \\ A_2 \\ A_3 \end{bmatrix} = \begin{bmatrix} 0.2 \\ 1 \\ 0.3 \\ \hline 0 \\ 0.14 \\ 1 \\ 0.7 \\ \hline 0.9 \\ 0.8 \end{bmatrix}$$

where

$$A_1 = \begin{bmatrix} 0.2 \\ 1 \\ 0.3 \end{bmatrix}, A_2 = \begin{bmatrix} 0 \\ 0.14 \\ 1 \\ 0.7 \end{bmatrix} \text{ and } A_3 = \begin{bmatrix} 0.9 \\ 0.8 \end{bmatrix}$$

be the three fuzzy column matrices.

$A_s$ is a super fuzzy column matrix. Here N = 3. Now as in case of super fuzzy row matrix, we can obtain super fuzzy column matrices by partitioning a fuzzy column matrix.

Suppose

$$B = \begin{bmatrix} b_1 \\ b_2 \\ b_3 \\ \vdots \\ b_m \end{bmatrix}$$

be a fuzzy column matrix.

Let

$$B_s = \left[\begin{array}{c} b_1 \\ b_2 \\ b_3 \\ b_4 \\ \hline b_5 \\ b_6 \\ \vdots \\ b_t \\ \hline b_{t+1} \\ \vdots \\ b_n \end{array}\right] = \left[\begin{array}{c} B_1 \\ \hline B_2 \\ \hline B_3 \end{array}\right] = \left[\begin{array}{c} B_1 \\ B_2 \\ B_3 \end{array}\right]$$

is a super fuzzy column matrix obtained by partitioning the fuzzy column matrix B.

***Example 2.1.4:*** Let

$$B = \begin{bmatrix} 0.2 \\ 1 \\ 0 \\ 0.3 \\ 0.7 \\ 0.9 \\ 0 \\ 0.1 \\ 1 \\ 0.3 \end{bmatrix}$$

be a fuzzy column matrix.

$$B_1 = \left[\begin{array}{c} 0.2 \\ 1 \\ 0 \\ \hline 0.3 \\ 0.7 \\ 0.9 \\ 0 \\ 0.1 \\ \hline 1 \\ 0.3 \end{array}\right]$$

is a super fuzzy column matrix. We can find super fuzzy column matrices according to our wish.

Now for instance

$$B_2 = \left[\begin{array}{c} 0.2 \\ 1 \\ \hline 0 \\ 0.3 \\ 0.7 \\ \hline 0.9 \\ 0 \\ \hline 0.1 \\ 0.3 \end{array}\right]$$

is also a super fuzzy column matrix.

$$B_3 = \begin{bmatrix} 0.2 \\ 1 \\ 0 \\ 0.3 \\ \hline 0.7 \\ \hline 0.9 \\ 0 \\ \hline 0.1 \\ 1 \\ 0.3 \end{bmatrix}$$

is also a super fuzzy column matrix. Clearly $B_1$, $B_2$ and $B_3$ are three distinct super fuzzy column matrices obtained by partitioning the fuzzy column matrix B.

Now we define super fuzzy matrices.

**DEFINITION 2.1.3:** *Let us consider a fuzzy matrix*

$$A = \begin{bmatrix} A_{11} & A_{12} & A_{13} \\ A_{21} & A_{22} & A_{23} \end{bmatrix}$$

*where $A_{11}$, $A_{12}$, $A_{13}$, $A_{21}$, $A_{22}$ and $A_{23}$ be fuzzy submatrices where number of columns in the fuzzy submatrices $A_{11}$ and $A_{21}$ are equal.*

*Similarly the columns in fuzzy submatrices of $A_{12}$ and $A_{22}$ are equal and columns of fuzzy matrices $A_{13}$ and $A_{23}$ are equal. This is evident from the second index of the fuzzy submatrices. One can also see, the number of row in fuzzy submatrices $A_{11}$, $A_{12}$ and $A_{13}$ are equal.*

*Similarly for fuzzy submatrices of $A_{21}$, $A_{22}$ and $A_{23}$ the number of rows are equal.*

*Thus a general super fuzzy matrix,*

$$A = \begin{bmatrix} A_{11} & A_{12} & \cdots & A_{1n} \\ A_{21} & A_{22} & \cdots & A_{2n} \\ \vdots & \vdots & & \vdots \\ A_{m1} & A_{m2} & \cdots & A_{mn} \end{bmatrix}$$

*where $A_{ij}$'s are fuzzy submatrices; i = 1, 2, …, m and j = 1, 2, …, n.*

We illustrate this by the following examples.

***Example 2.1.5:*** Let

$$A = \begin{bmatrix} A_{11} & A_{12} & A_{13} \\ A_{21} & A_{22} & A_{23} \\ A_{31} & A_{32} & A_{33} \\ A_{41} & A_{42} & A_{43} \end{bmatrix}$$

where $A_{ij}$ are fuzzy submatrices, $A_{k1}$, $A_{k2}$ and $A_{k3}$ have same number of rows for k = 1, 2, 3, 4. $A_{1m}$, $A_{2m}$, $A_{3m}$ and $A_{4m}$ have same number of columns for m = 1, 2, 3.

***Example 2.1.6:*** A is a super fuzzy matrix given below:

$$A = \left[\begin{array}{cc|cccc|ccc} 0.2 & 0.3 & 1 & 0 & 0.5 & 0.6 & 1 & 0.7 & 0.2 \\ 1 & 0 & 0.2 & 1 & 0.3 & 0.9 & 0 & 0.2 & 0 \\ \hline 0 & 1 & 0.3 & 0 & 0.4 & 1 & 0 & 0.2 & 0 \\ 0.9 & 0.3 & 0 & 1 & 1 & 0 & 0 & 1 & 0 \\ 0.6 & 1 & 0 & 0.4 & 1 & 0.3 & 1 & 0 & 0 \\ \hline 0.2 & 0 & 1 & 0 & 1 & 0.4 & 0 & 1 & 0.3 \\ 0.3 & 1 & 0.3 & 0.1 & 0 & 0.5 & 1 & 0.1 & 0 \\ 1 & 0.4 & 0.2 & 1 & 0 & 1 & 0.3 & 0 & 0.3 \\ \hline 0 & 1 & 0.1 & 0 & 0 & 1 & 0.2 & 0.5 & 0.3 \\ 0.4 & 0.3 & 1 & 0.2 & 0.5 & 0.7 & 1 & 1 & 0 \end{array}\right]$$

where

$$A_{11} = \begin{bmatrix} 0.2 & 0.3 \\ 1 & 0 \end{bmatrix},$$

$$A_{12} = \begin{bmatrix} 1 & 0 & 0.5 & 0.6 \\ 0.2 & 1 & 0.3 & 0.9 \end{bmatrix},$$

$$A_{13} = \begin{bmatrix} 1 & 0.7 & 0.2 \\ 0 & 0.2 & 0 \end{bmatrix},$$

$$A_{21} = \begin{bmatrix} 0 & 1 \\ 0.9 & 0.3 \\ 0.6 & 1 \end{bmatrix},$$

$$A_{22} = \begin{bmatrix} 0.3 & 0 & 0.4 & 1 \\ 0 & 1 & 1 & 0 \\ 0 & 0.4 & 1 & 0.3 \end{bmatrix},$$

$$A_{23} = \begin{bmatrix} 0 & 0.2 & 0 \\ 0 & 1 & 0 \\ 1 & 0 & 0 \end{bmatrix},$$

$$A_{31} = \begin{bmatrix} 0.2 & 0 \\ 0.3 & 1 \\ 1 & 0.4 \end{bmatrix},$$

$$A_{32} = \begin{bmatrix} 1 & 0 & 1 & 0.4 \\ 0.3 & 0.1 & 0 & 0.5 \\ 0.2 & 1 & 0 & 1 \end{bmatrix},$$

$$A_{33} = \begin{bmatrix} 0 & 1 & 0.3 \\ 1 & 0.1 & 0 \\ 0.3 & 0 & 0.3 \end{bmatrix},$$

$$A_{41} = \begin{bmatrix} 0 & 1 \\ 0.4 & 0.3 \end{bmatrix},$$

$$A_{42} = \begin{bmatrix} 0.1 & 0 & 0 & 1 \\ 1 & 0.2 & 0.5 & 0.7 \end{bmatrix}$$

and

$$A_{43} = \begin{bmatrix} 0.2 & 0.5 & 0.3 \\ 1 & 1 & 0 \end{bmatrix}.$$

Thus we see $A_{11}$, $A_{12}$ and $A_{13}$ has each two rows. $A_{21}$, $A_{22}$ and $A_{23}$ has each 3 rows. $A_{31}$, $A_{32}$ and $A_{33}$ has each three rows. $A_{41}$, $A_{42}$ and $A_{43}$ each has two rows. $A_{11}$, $A_{21}$, $A_{31}$ and $A_{41}$ each has two columns. $A_{12}$, $A_{22}$, $A_{32}$ and $A_{42}$ each has 4 columns. $A_{13}$, $A_{23}$, $A_{33}$ and $A_{43}$ each has 3 columns.

We do not call them as type 3 or type 4 or any type of super fuzzy matrices. We give yet another example of a super fuzzy matrix.

***Example 2.1.7:*** Let

$$A = \left[\begin{array}{cccc|cc} 0.1 & 0.4 & 1 & 0 & 0.5 & 0.3 \\ 1 & 0 & 0.7 & 0.3 & 1 & 0 \\ 0 & 1 & 0.2 & 0.4 & 0.6 & 0.7 \\ 0.4 & 0 & 0.3 & 0.2 & 0.2 & 0.3 \\ 0.6 & 1 & 0.5 & 0.4 & 1 & 0 \\ 0.7 & 1 & 0.2 & 1 & 0 & 1 \end{array}\right]$$

$= [A_1 \mid A_2]$ where

$$A_1 = \begin{bmatrix} 0.1 & 0.4 & 1 & 0 \\ 1 & 0 & 0.7 & 0.3 \\ 0 & 1 & 0.2 & 0.4 \\ 0.4 & 0 & 0.3 & 0.2 \\ 0.6 & 1 & 0.5 & 0.4 \\ 0.7 & 1 & 0.2 & 1 \end{bmatrix}$$

and

$$A_2 = \begin{bmatrix} 0.5 & 0.3 \\ 1 & 0 \\ 0.6 & 0.7 \\ 0.2 & 0.3 \\ 1 & 0 \\ 0 & 1 \end{bmatrix}$$

is a super fuzzy matrix.

***Example 2.1.8:*** Let

$$A = \left[\begin{array}{cccccc} 0.3 & 1 & 0 & 0.2 & 1 & 0.5 \\ 0 & 0.2 & 1 & 0 & 1 & 0.6 \\ 1 & 1 & 0 & 1 & 0.3 & 1 \\ \hline 0.3 & 1 & 0 & 0.1 & 0.2 & 1 \\ 0 & 0.3 & 0.2 & 0.1 & 1 & 0 \end{array}\right]$$

$$= \left[\begin{array}{c} A_1 \\ \hline A_2 \end{array}\right]$$

where

$$A_1 = \begin{bmatrix} 0.3 & 1 & 0 & 0.2 & 1 & 0.5 \\ 0 & 0.2 & 1 & 0 & 1 & 0.6 \\ 1 & 1 & 0 & 1 & 0.3 & 1 \end{bmatrix}$$

and

$$A_2 = \begin{bmatrix} 0.3 & 1 & 0 & 0.1 & 0.2 & 1 \\ 0 & 0.3 & 0.2 & 0.1 & 1 & 0 \end{bmatrix}$$

is a fuzzy supermatrix.

***Example 2.1.9:*** Let

$$A = \left[\begin{array}{ccc|cc} 0.3 & 0.5 & 1 & 0 & 0.2 \\ 0.6 & 1 & 0 & 0.7 & 0.1 \\ 0.4 & 0.3 & 1 & 0.3 & 0.5 \\ \hline 1 & 0.2 & 0 & 0.2 & 1 \\ 0 & 1 & 0.7 & 0.6 & 0 \\ 0.7 & 0.5 & 0.6 & 0.3 & 0.5 \\ 0.2 & 0.1 & 0.3 & 1 & 0.2 \end{array}\right]$$

$$= \begin{bmatrix} A_{11} & A_{12} \\ A_{21} & A_{22} \end{bmatrix}$$

where

$$A_{11} = \begin{bmatrix} 0.3 & 0.5 & 1 \\ 0.6 & 1 & 0 \\ 0.4 & 0.3 & 1 \end{bmatrix},$$

$$A_{12} = \begin{bmatrix} 0 & 0.2 \\ 0.7 & 0.1 \\ 0.3 & 0.5 \end{bmatrix},$$

$$A_{21} = \begin{bmatrix} 1 & 0.2 & 0 \\ 0 & 1 & 0.7 \\ 0.7 & 0.5 & 0.6 \\ 0.2 & 0.1 & 0.3 \end{bmatrix}$$

and

$$A_{22} = \begin{bmatrix} 0.2 & 1 \\ 0.6 & 0 \\ 0.3 & 0.5 \\ 1 & 0.2 \end{bmatrix}.$$

Now we proceed on to define the notion of transpose of a fuzzy supermatrix.

**DEFINITION 2.2.4:** *Let*

$$A = \begin{bmatrix} A_{11} & A_{12} & \cdots & A_{1n} \\ A_{21} & A_{22} & \cdots & A_{2n} \\ \vdots & \vdots & & \vdots \\ A_{m1} & A_{m2} & \cdots & A_{mn} \end{bmatrix}$$

*where $A_{ij}$ are fuzzy submatrices of A; $1 \leq i \leq m$ and $1 \leq j \leq n$, with entries from [0, 1]. A is a fuzzy supermatrix for elements of A belong to [0, 1].*

*The transpose of the fuzzy supermatrix A is defined to be $A^t$ denoted by A'.*

$$A^t = A' = \begin{bmatrix} A'_{11} & A'_{21} & \cdots & A'_{m1} \\ A'_{12} & A'_{22} & \cdots & A'_{m2} \\ \vdots & \vdots & & \vdots \\ A'_{m} & A'_{2n} & \cdots & A'_{mn} \end{bmatrix}$$

*where $A^t_{ij} = A'_{ij}$ ; $1 \leq i \leq m$ and $1 \leq j \leq n$.*

We illustrate them by the following examples.

***Example 2.1.10:*** Let A = [$A_1$ | $A_2$ | $A_3$ |$A_4$] be a super fuzzy row vector, where $A_1$ = [0 1 0.3], $A_2$ = [0.5 1 0 0.2], $A_3$ = [1 0] and $A_4$ = [0.3 1 0.6 0 1 0.7].

Now the transpose of the super fuzzy row vector A is given by

$$A' = \begin{bmatrix} A'_1 \\ A'_2 \\ A'_3 \\ A'_4 \end{bmatrix}$$

where

$$A'_1 = \begin{bmatrix} 0 \\ 1 \\ 0.3 \end{bmatrix}, A'_2 = \begin{bmatrix} 0.5 \\ 1 \\ 0 \\ 0.2 \end{bmatrix}, A'_3 = \begin{bmatrix} 1 \\ 0 \end{bmatrix} \text{ and } A'_4 = \begin{bmatrix} 0.3 \\ 1 \\ 0.6 \\ 0 \\ 1 \\ 0.7 \end{bmatrix}.$$

Now we proceed on to define some operations on them,

***Example 2.1.11:*** Let

$$A = \begin{bmatrix} A_1 \\ A_2 \\ A_3 \end{bmatrix}$$

be a fuzzy super column vector where

$$A_1 = \begin{bmatrix} 0.3 \\ 0.2 \\ 1 \\ 0 \end{bmatrix}, A_2 = \begin{bmatrix} 1 \\ 0 \\ 0.1 \\ 0.4 \end{bmatrix} \text{ and } A_3 = \begin{bmatrix} 0.5 \\ 0.1 \\ 0.5 \\ 0.2 \\ 1 \end{bmatrix}.$$

Now $A' = [A'_1\ A'_2\ A'_3] = [0.3\ 0.2\ 1\ 0 \mid 1\ 0\ 0.1\ 0.4 \mid 0.5\ 0.1\ 0.5\ 0.2\ 1]$.

***Example 2.1.12:*** Let

$$A = \begin{bmatrix} A_1 \\ A_2 \\ A_3 \end{bmatrix}$$

be a fuzzy super column vector where

$$A_1 = \begin{bmatrix} 0 \\ 1 \\ 0 \end{bmatrix}, A_2 = [\ 0.3\ ] \text{ and } A_3 = \begin{bmatrix} 1 \\ 0.2 \\ 0.5 \\ 0 \end{bmatrix}.$$

$$\begin{aligned} A' &= [A'_1\ A'_2\ A'_3] \\ &= [0\ 1\ 0 \mid 0.3 \mid 1\ 0.2\ 0.5\ 0] \end{aligned}$$

is a fuzzy super row vector.

$$\{\max \min\{A', A\}\} = A'.\ A.$$

(This is only a notational convenience)

$$\max\left[\max\min\left\{\left[A'_1\ A'_2\ A'_3\right], \begin{bmatrix} A_1 \\ A_2 \\ A_3 \end{bmatrix}\right\}\right]$$

$$= \max\left\{\max\min\left\{\left[0\ \ 1\ \ 0 \mid 0.3 \mid 1\ \ 0.2\ \ 0.5\ \ 0\right].\begin{bmatrix} 0 \\ 1 \\ 0 \\ \hline 0.3 \\ \hline 1 \\ 0.2 \\ 0.5 \\ 0 \end{bmatrix}\right\}\right\}$$

$$= \max\{A'_1\ .\ A_1, A'_2\ .\ A_2, A'_3\ .\ A_3\}$$

which will be known as the minor product of two fuzzy supermatrices.

max{A' . A}

$$= \max \left\{ \begin{bmatrix} 0 & 1 & 0 \end{bmatrix} . \begin{bmatrix} 0 \\ 1 \\ 0 \end{bmatrix}, [0.3] . [0.3], [1 \ 0.2 \ 0.5 \ 0] . \begin{bmatrix} 1 \\ 0.2 \\ 0.5 \\ 0 \end{bmatrix} \right\}$$

= max [{max {min (0, 0), min (1, 1), min (0, 0)}, max {min (0.3, 0.3)}, max{min (1, 1), min (0.2, 0.2), min (0.5, 0.5), min (0, 0)}]
= max {max (0, 1, 0), max (0.3), max {1, 0.2, 0.5, 0}}
= max {1, 0.3, 1}
= 1;

the way in which max{A' . A} = max {max{min(A', A)}} is defined is peculiar which may be defined as the super pseudo product of the transpose of the fuzzy supermatrix A' and the fuzzy supermatrix A.

Before we define product or min operation we give the illustration of it with

$$\min (A', A) = \min \left\{ \begin{bmatrix} A_1' & A_2' & A_3' \end{bmatrix}, \begin{bmatrix} A_1 \\ A_2 \\ A_3 \end{bmatrix} \right\}$$

we first make a note as submatrix multiplication the product is not compatible. So we multiply in this special way.

min (A', A)

$$= \min\left\{\left[0 \quad 1 \quad 0 \mid 0.3 \mid 1 \quad 0.2 \quad 0.5 \quad 0\right], \left[\begin{array}{c} 0 \\ 1 \\ 0 \\ \hline 0.3 \\ \hline 1 \\ 0.2 \\ 0.5 \\ 0 \end{array}\right]\right\}$$

$$= \min\left\{\left[a_1 \quad a_2 \quad a_3 \mid a_4 \mid a_5 \quad a_6 \quad a_7 \quad a_8\right], \left[\begin{array}{c} a_1 \\ a_2 \\ a_3 \\ \hline a_4 \\ \hline a_5 \\ a_6 \\ a_7 \\ a_8 \end{array}\right]\right\}$$

$$= \begin{bmatrix} \min(a_1a_1) & \min(a_1a_2) & \cdots & \min(a_1a_8) \\ \min(a_2a_1) & \min(a_2a_2) & \cdots & \min(a_2a_8) \\ \vdots & \vdots & & \vdots \\ \min(a_8a_1) & \min(a_8a_2) & \cdots & \min(a_8a_8) \end{bmatrix}$$

$$= \left[\begin{array}{ccc|c|cccc} 0 & 0 & 0 & 0 & 0 & 0 & 0 & 0 \\ 0 & 1 & 0 & 0.3 & 1 & 0.2 & 0.5 & 0 \\ 0 & 0 & 0 & 0 & 0 & 0 & 0 & 0 \\ \hline 0 & 0.3 & 0 & 0.3 & 0.3 & 0.2 & 0.3 & 0 \\ \hline 0 & 1 & 0 & 0.3 & 1 & 0.2 & 0.5 & 0 \\ 0 & 0.2 & 0 & 0.2 & 0.2 & 0.2 & 0.2 & 0 \\ 0 & 0.5 & 0 & 0.3 & 0.5 & 0.2 & 0.5 & 0 \\ 0 & 0 & 0 & 0 & 0 & 0 & 0 & 0 \end{array}\right].$$

Now we make three observations from the example.

1) The product is not compatible as usual matrices yielding a single element.
2) The partition is carried by the super pseudo product of type I in a very special way.
3) The resultant fuzzy matrix is symmetric. In fact it can be defined as a fuzzy symmetric supermatrix.

**DEFINITION 2.1.5:** *Let $A = [A_1 \mid A_2 \mid \ldots \mid A_n]$ be a fuzzy super row matrix i.e., each $A_i$ is a $1 \times t_i$ fuzzy row submatrix of A, i = 1, 2, …, n. Now we define using the transpose of A;*

$$A' = \begin{bmatrix} A'_1 \\ \hline A'_2 \\ \hline \vdots \\ \hline A'_n \end{bmatrix}$$

*the two types of products of these fuzzy super column and row matrices. i.e., A . A' and A' . A.*

$$\max\{A \,.\, A'\} = \max\left\{ \begin{bmatrix} A'_1 & A'_2 & \ldots & A'_n \end{bmatrix} . \begin{bmatrix} A'_1 \\ A'_2 \\ \vdots \\ A'_n \end{bmatrix} \right\}$$

*$\max \{A_1 \,.\, A'_1, A_2 \,.\, A'_2, \ldots, A_n \,.\, A'_n\} = \max \{\max \min (a_{i1}, a'_{11}), \max \min (a_{i2}, a'_{i2}), \ldots, \max \min (a_{in}, a'_{in})\}; 1 < i_1, i'_1 < t_1, 1 < i_2, i'_2 < t_2, \ldots, 1 < i_n, i'_n < t_n$ where*

$$A_1 = (a_1, \ldots, a_{t_1}) \text{ and } A'_1 = \begin{bmatrix} a'_1 \\ \vdots \\ a'_{t_1} \end{bmatrix}$$

*i.e., $A'_1 = A_1^t$, $A_2 = (a_1, a_2, \ldots, a_{t_2})$ and $A'_2 = A_2^t$ and so on. Thus $\max\{A \,.\, A'\} = \max\{A \,.\, A^t\} = a \in [0, 1]$.*

*This fuzzy supermatrix operation is the usual operation with product of $a_i$. $a'_j$ replaced by minimum ($a_i$, $a'_j$) and the sum of the elements $a_{i_1} + \ldots + a_{i_n}$ replaced by the maximum of ($a_{i_1} \ldots a_{i_n}$).*

Thus given any fuzzy super row matrix A we can always define the product where A' is the transpose of A the fuzzy super column matrix.

Now A . A' is defined as the product of two fuzzy super row matrix (vector) and fuzzy super column vector and the product is always an element from the fuzzy interval [0, 1].

Now how is A'A defined A' is a n × 1 fuzzy super column vector and A is a 1 × n fuzzy super row vector. How to define this product at the same time not fully destroying the partition or the submatrix structure.

Thus we define as a super pseudo product which is a super fuzzy n × n matrix.

**DEFINITION 2.1.6:** *Let A = ($a_1$ $a_2$ … $a_n$) = ($A_1$ | $A_2$ | …| $A_n$) be a super fuzzy row matrix and A' be the transpose of A.*

$$min\{A', A\} = \begin{bmatrix} min(a'_1a_1) & min(a'_1a_2) & \cdots & min(a'_1a_n) \\ min(a'_2a_1) & min(a'_2a_2) & \cdots & min(a'_2a_n) \\ \vdots & \vdots & & \vdots \\ min(a'_na_1) & min(a'_na_1) & \cdots & min(a'_na_n) \end{bmatrix}.$$

*Now min{A' . A} is a n × n square fuzzy matrix. It is partitioned as per the division of rows and columns of A and A' respectively. It is important to note min{A' . A} is a symmetric matrix about the diagonal*

$$min\{A', A\} = \begin{bmatrix} min(a_1,a_1) & min(a_1,a_2) & \cdots & min(a_1,a_n) \\ min(a_2,a_1) & min(a_2,a_2) & \cdots & min(a_2,a_n) \\ \vdots & \vdots & & \vdots \\ min(a_n,a_1) & min(a_n,a_2) & \cdots & min(a_n,a_n) \end{bmatrix},$$

*where $a'_i = a_i$ for the elements remain as it is while transposing the elements clearly since min $(a_2, a_1)$ = min $(a_1, a_2)$ we get the min{A', A} matrix to be a symmetric matrix. Further if $A = (A_1 A_2 \ldots A_n)$ with number of elements in $A_i$ is $t_i$ then we see $1 < i < n$ and $1 < t_i < n$. min{A', A} is a super fuzzy matrix with $i \times i$ fuzzy submatrices $i = 1, 2, \ldots, n$. and min{A', A} is a $n \times n$ fuzzy matrix.*

We first illustrate this by the following example.

***Example 2.1.13:*** Let

$$\begin{aligned} A &= [0.1\ 0\ 0.5 \mid 0.2\ 1 \mid 0.3\ 0.1\ 1\ 1] \\ &= [A_1 \mid A_2 \mid A_3] \end{aligned}$$

be a fuzzy super row matrix with

$$A_1 = [0.1\ 0\ 0.5],\ A_2 = [0.2\ 1]$$

and

$$A_3 = [0.3\ 0.1\ 1\ 1].$$

$$A' = \begin{bmatrix} A_1 \\ A_2 \\ A_3 \end{bmatrix}$$

we find

$$\min(A', A) = \min\left\{\begin{bmatrix} 0.1 \\ 0 \\ 0.5 \\ \hline 0.2 \\ 1 \\ \hline 0.3 \\ 0.1 \\ 1 \\ 1 \end{bmatrix}, \begin{bmatrix} 0.1 & 0 & 0.5 & | & 0.2 & 1 & | & 0.3 & 0.1 & 1 & 1 \end{bmatrix}\right\}$$

$$= \begin{bmatrix} \min(0.1,0.1) & \min(0.1,0) & \cdots & \min(0.1,1) \\ \min(0,0.1) & \min(0,0) & \cdots & \min(0,1) \\ \vdots & \vdots & & \vdots \\ \min(1,0.1) & \min(1,0) & \cdots & \min(1,1) \\ \min(1,0.1) & \min(1,0) & \cdots & \min(1,1) \end{bmatrix}.$$

We see min{A', A} is a 9 × 9 fuzzy supermatrix partitioned between 3rd and 4th row. Also between 5th and 6th row. Similarly min{A', A} is a fuzzy supermatrix partitioned between 3rd and 4th column and 5th and 6th column.

$$\min\{A', A\} = \left[\begin{array}{ccc|cc|cccc} 0.1 & 0 & 0.1 & 0.1 & 0.1 & 0.1 & 0.1 & 0.1 & 0.1 \\ 0 & 0 & 0 & 0 & 0 & 0 & 0 & 0 & 0 \\ 0.1 & 0 & 0.5 & 0.2 & 0.5 & 0.3 & 0.1 & 0.5 & 0.5 \\ \hline 0.1 & 0 & 0.2 & 0.2 & 0.2 & 0.2 & 0.1 & 0.2 & 0.2 \\ 0.1 & 0 & 0.5 & 0.2 & 1 & 0.3 & 0.1 & 1 & 1 \\ \hline 0.1 & 0 & 0.3 & 0.2 & 0.3 & 0.3 & 0.1 & 0.3 & 0.3 \\ 0.1 & 0 & 0.1 & 0.1 & 0.1 & 0.1 & 0.1 & 0.1 & 0.1 \\ 0.1 & 0 & 0.5 & 0.2 & 1 & 0.3 & 0.1 & 1 & 1 \\ 0.1 & 0 & 0.5 & 0.2 & 1 & 0.3 & 0.1 & 1 & 1 \end{array}\right]$$

$$= \begin{bmatrix} B_{11} & B_{12} & B_{13} \\ B'_{12} & B_{22} & B_{23} \\ B'_{13} & B'_{23} & B_{33} \end{bmatrix}$$

where

$$B_{11} = \begin{bmatrix} 0.1 & 0 & 0.1 \\ 0 & 0 & 0 \\ 0.1 & 0 & 0.5 \end{bmatrix} = \min\{A'_1, A_1\};$$

$$B_{12} = \begin{bmatrix} 0.1 & 0.1 \\ 0 & 0 \\ 0.2 & 0.5 \end{bmatrix} = \min\{A'_1, A_2\}$$

$$= \min\left\{\begin{bmatrix} 0.1 \\ 0 \\ 0.5 \end{bmatrix}, \begin{bmatrix} 0.2 & 1 \end{bmatrix}\right\}$$

$$= \begin{bmatrix} \min(0.1, 0.2) & \min(0.1, 1) \\ \min(0, 0.2) & \min(0, 1) \\ \min(0.5, 0.2) & \min(0.5, 1) \end{bmatrix}.$$

$$B_{13} = \min\{A'_1, A_3\} = \min\left\{\begin{bmatrix} 0.1 \\ 0 \\ 0.5 \end{bmatrix}, \begin{bmatrix} 0.3 & 0.1 & 1 & 1 \end{bmatrix}\right\}$$

$$= \begin{bmatrix} \min(0.1, 0.3) & \min(0.1, 0.1) & \min(0.1, 1) & \min(0.1, 1) \\ \min(0, 0.3) & \min(0, 0.1) & \min(0, 1) & \min(0, 1) \\ \min(0.5, 0.3) & \min(0.5, 0.1) & \min(0.5, 1) & \min(0.5, 1) \end{bmatrix}$$

$$= \begin{bmatrix} 0.1 & 0.1 & 0.1 & 0.1 \\ 0 & 0 & 0 & 0 \\ 0.3 & 0.1 & 0.5 & 0.5 \end{bmatrix} = B_{13}.$$

$$B'_{12} = \min\{A'_2, A_1\} = \min\left\{\begin{bmatrix} 0.2 \\ 1 \end{bmatrix}, \begin{bmatrix} 0.1 & 0 & 0.5 \end{bmatrix}\right\}$$

$$= \begin{bmatrix} \min(0.2, 0.1) & \min(0.2, 0) & \min(0.2, 0.5) \\ \min(1, 0.1) & \min(1, 0) & \min(1, 0.5) \end{bmatrix}$$

$$= \begin{bmatrix} 0.1 & 0 & 0.2 \\ 0.1 & 0 & 0.5 \end{bmatrix} = B'_{12} = B^t_{12}.$$

$$B_{22} = \min\{A'_2, A_2\} = \min \left\{ \begin{bmatrix} 0.2 \\ 1 \end{bmatrix}, \begin{bmatrix} 0.2 & 1 \end{bmatrix} \right\}$$

$$= \begin{bmatrix} \min(0.2, 0.2) & \min(0.2, 1) \\ \min(1, 0.2) & \min(1, 1) \end{bmatrix}$$

$$= \begin{bmatrix} 0.2 & 0.2 \\ 0.2 & 1 \end{bmatrix}.$$

$$B_{23} = \min\{A'_2, A_3\} = \min \left\{ \begin{bmatrix} 0.2 \\ 1 \end{bmatrix}, \begin{bmatrix} 0.3 & 0.1 & 1 & 1 \end{bmatrix} \right\}$$

$$= \begin{bmatrix} \min(0.2, 0.3) & \min(0.2, 0.1) & \min(0.2, 1) & \min(0.2, 1) \\ \min(1, 0.3) & \min(1, 0.1) & \min(1, 1) & \min(1, 1) \end{bmatrix}$$

$$= \begin{bmatrix} 0.2 & 0.1 & 0.2 & 0.2 \\ 0.3 & 0.1 & 1 & 1 \end{bmatrix}.$$

$$B'_{13} = \min\{A'_3, A_1\} = \min\{(A'_1. A_3)'\}$$

$$= \min \left\{ \begin{bmatrix} 0.3 \\ 0.1 \\ 1 \\ 1 \end{bmatrix}, \begin{bmatrix} 0.1 & 0 & 0.5 \end{bmatrix} \right\}$$

$$= \begin{bmatrix} \min(0.3, 0.1) & \min(0.3, 0) & \min(0.3, 0.5) \\ \min(0.1, 0.1) & \min(0.1, 0) & \min(0.1, 0.5) \\ \min(1, 0.1) & \min(1, 0) & \min(1, 0.5) \\ \min(1, 0.1) & \min(1, 0) & \min(1, 0.5) \end{bmatrix}$$

$$= \begin{bmatrix} 0.1 & 0 & 0.3 \\ 0.1 & 0 & 0.1 \\ 0.1 & 0 & 0.5 \\ 0.1 & 0 & 0.5 \end{bmatrix}.$$

$$B'_{23} = \min\{A'_3 , A_2\} = \min\{(A'_2, A_3)'\}$$

$$= \min \left\{ \begin{bmatrix} 0.3 \\ 0.1 \\ 1 \\ 1 \end{bmatrix}, \begin{bmatrix} 0.2 & 1 \end{bmatrix} \right\}$$

$$= \begin{bmatrix} \min(0.3, 0.2) & \min(0.3, 1) \\ \min(0.1, 0.2) & \min(0.1, 1) \\ \min(1, 0.2) & \min(1, 1) \\ \min(1, 0.2) & \min(1, 1) \end{bmatrix}$$

$$= \begin{bmatrix} 0.2 & 0.3 \\ 0.1 & 0.1 \\ 0.2 & 1 \\ 0.2 & 1 \end{bmatrix}.$$

$$B_{33} = \min\{A'_3 , A_3\}$$

$$= \min \left\{ \begin{bmatrix} 0.3 \\ 0.1 \\ 1 \\ 1 \end{bmatrix}, \begin{bmatrix} 0.3 & 0.1 & 1 & 1 \end{bmatrix} \right\}$$

$$= \begin{bmatrix} \min(0.3, 0.3) & \min(0.3, 0.1) & \min(0.3, 1) & \min(0.3, 1) \\ \min(0.1, 0.3) & \min(0.1, 0.1) & \min(0.1, 1) & \min(0.1, 1) \\ \min(1, 0.3) & \min(1, 0.1) & \min(1, 1) & \min(1, 1) \\ \min(1, 0.3) & \min(1, 0.1) & \min(1, 1) & \min(1, 1) \end{bmatrix}$$

$$= \begin{bmatrix} 0.3 & 0.1 & 0.3 & 0.3 \\ 0.1 & 0.1 & 0.1 & 0.1 \\ 0.3 & 0.1 & 1 & 1 \\ 0.3 & 0.1 & 1 & 1 \end{bmatrix}.$$

Thus the fuzzy supermatrix of A'A is a fuzzy supermatrix got by the super pseudo product of the transpose of A with A. i.e., of A' with A where A is a 1 × n fuzzy row vector. A'A is a n × n super fuzzy matrix here n = 9.

Next we define the notion of symmetric fuzzy supermatrix.

**DEFINITION 2.1.7:** *Let*

$$A = \begin{bmatrix} A_{11} & A_{12} & \cdots & A_{1n} \\ A_{21} & A_{22} & \cdots & A_{2n} \\ \vdots & \vdots & & \vdots \\ A_{n1} & A_{n2} & \cdots & A_{nn} \end{bmatrix}$$

*be a fuzzy supermatrix, we say A is a fuzzy super square matrix or supermatrix or fuzzy square supermatrix if the number of columns and number of rows of A are the same and the number of fuzzy submatrices along the rows is equal to the number of fuzzy column submatrices.*

*If the number of rows and columns of the fuzzy supermatrix is different then we call the fuzzy supermatrix to be a rectangular fuzzy supermatrix or fuzzy rectangular supermatrix or fuzzy super rectangular matrix.*

All these terms mean one and the same. Before we go for further investigations we give an example of each.

***Example 2.1.14:*** Let

$$A = \begin{bmatrix} A_{11} & A_{12} \\ A_{21} & A_{22} \\ A_{31} & A_{32} \end{bmatrix}$$

$$= \left[\begin{array}{cccc|cc} 0.1 & 1 & 0 & 0.2 & 1 & 0 \\ 1 & 1 & 1 & 0.5 & 0 & 0.6 \\ \hline 0.6 & 0.8 & 0.5 & 1 & 1 & 1 \\ \hline 0 & 1 & 0 & 0.6 & 1 & 0.4 \\ 0.2 & 0 & 1 & 0.2 & 0 & 0.1 \\ 1 & 0.2 & 0.3 & 0.3 & 0.2 & 0.2 \end{array}\right]$$

$$A_{11} = \begin{bmatrix} 0.1 & 1 & 0 & 0.2 \\ 1 & 1 & 1 & 0.5 \end{bmatrix},$$

$$A_{12} = \begin{bmatrix} 1 & 0 \\ 0 & 0.6 \end{bmatrix}$$

$$A_{21} = [0.6\ 0.8\ 0.5\ 1],$$

$$A_{22} = [1\ 1]$$

$$A_{31} = \begin{bmatrix} 0 & 1 & 0 & 0.6 \\ 0.2 & 0 & 1 & 0.2 \\ 1 & 0.2 & 0.3 & 0.3 \end{bmatrix}$$

and

$$A_{32} = \begin{bmatrix} 1 & 0.4 \\ 0 & 0.1 \\ 0.2 & 0.2 \end{bmatrix}.$$

This is a ordinary fuzzy supermatrix, but is not a fuzzy square supermatrix.

We see the transpose of A is given by

$$A^t = \begin{bmatrix} A_{11}^t & A_{21}^t & A_{31}^t \\ A_{12}^t & A_{22}^t & A_{32}^t \end{bmatrix}$$

$$= \left[\begin{array}{cc|c|ccc} 0.1 & 1 & 0.6 & 0 & 0.2 & 1 \\ 1 & 1 & 0.8 & 1 & 0 & 0.2 \\ 0 & 1 & 0.5 & 0 & 1 & 0.3 \\ 0.2 & 0.5 & 1 & 0.6 & 0.2 & 0.3 \\ \hline 1 & 0 & 1 & 1 & 0 & 0.2 \\ 0 & 0.6 & 1 & 0.4 & 0.1 & 0.2 \end{array}\right].$$

We can easily observe that the fuzzy matrix $A^t$ is partitioned between the rows 2 and 3 and 3 and 4. It is partitioned between the columns 4 and 5.

***Example 2.1.15:*** Let

$$A = \begin{bmatrix} a_{11} & a_{12} & \cdots & a_{1n} \\ a_{21} & a_{22} & \cdots & a_{2n} \\ \vdots & \vdots & & \vdots \\ a_{n1} & a_{n2} & \cdots & a_{nn} \end{bmatrix}$$

be a fuzzy square matrix. The fuzzy supermatrix got from A by partitioning A between the columns $a_{1i}$ and $a_{1i+1}$

$$A_s = \left[\begin{array}{cccc|ccc} a_{11} & a_{12} & \cdots & a_{1i} & a_{1i+1} & \cdots & a_{1n} \\ a_{21} & a_{22} & \cdots & a_{2i} & a_{2i+1} & \cdots & a_{2n} \\ \vdots & \vdots & & \vdots & \vdots & & \vdots \\ a_{n1} & a_{n2} & \cdots & a_{ni} & a_{ni+1} & \cdots & a_{nn} \end{array}\right]$$

is a fuzzy supermatrix but is not a square fuzzy supermatrix. We see the transpose of $A_S$ is

$$A_s^t = \left[\begin{array}{cccc} a_{11} & a_{21} & \cdots & a_{n1} \\ a_{12} & a_{22} & \cdots & a_{n2} \\ \vdots & \vdots & & \vdots \\ a_{1i} & a_{2i} & \cdots & a_{ni+1} \\ \hline a_{1i+1} & a_{2i+1} & \cdots & a_{ni+1} \\ \vdots & \vdots & & \vdots \\ a_{1n} & a_{2n} & \cdots & a_{nn} \end{array}\right].$$

We see yet another example.

***Example 2.1.16:*** Let A be any fuzzy square matrix.

$$A = \begin{bmatrix} a_{11} & a_{12} & \cdots & a_{1n} \\ a_{21} & a_{22} & \cdots & a_{n2} \\ \vdots & \vdots & & \vdots \\ a_{n1} & a_{n2} & \cdots & a_{nn} \end{bmatrix}.$$

We see A is partitioned between the rows $a_{i1}$ and $a_{i+11}$ and between the rows $a_{j1}$ and $a_{j+11}$.

The fuzzy supermatrix $A_s$ obtained from partitioning A is

$$A_s = \left[\begin{array}{cccc} a_{11} & a_{12} & \cdots & a_{1n} \\ a_{21} & a_{22} & \cdots & a_{2n} \\ a_{i1} & a_{i2} & \cdots & a_{in} \\ \hline a_{i+11} & a_{i+12} & \cdots & a_{i+1n} \\ \vdots & & & \\ a_{j1} & a_{j2} & \cdots & a_{jn} \\ \hline a_{j+11} & a_{j+12} & \cdots & a_{j+2n} \\ \vdots & & & \\ a_{n1} & a_{n2} & \cdots & a_{nn} \end{array}\right].$$

We see A is partitioned between the rows $a_{i1}$ and $a_{i+11}$ and between the rows $a_{ji}$ and $a_{j+11}$.

Clearly

$$A_s^t = \left[\begin{array}{cccc|ccc|ccc} a_{11} & a_{21} & \cdots & a_{i1} & a_{i+11} & \cdots & a_{j1} & a_{j+11} & \cdots & a_{n1} \\ a_{12} & a_{22} & \cdots & a_{i+2} & a_{i+12} & \cdots & a_{j2} & a_{j+12} & \cdots & a_{n2} \\ \vdots & \vdots & & \vdots & \vdots & & \vdots & \vdots & & \vdots \\ a_{1n} & a_{2n} & \cdots & a_{in} & a_{i+1n} & \cdots & a_{jn} & a_{j+1n} & \cdots & a_{nn} \end{array}\right].$$

Now $A_s^t$ is got by partitioning A between the rows $a_{i1}$ and $a_{i+11}$ and $a_{j1}$ and $a_{j+11}$. $A_s$ is not a fuzzy square supermatrix it is just a fuzzy supermatrix.

Now we still give another example.

***Example 2.1.17:*** Let A be a fuzzy square matrix.

$$A_s = \left[\begin{array}{cccc|cc} a_{11} & a_{12} & \cdots & a_{1i} & a_{1i+1} & a_{1n} \\ a_{21} & a_{22} & \cdots & a_{2i} & a_{2i+1} & a_{2n} \\ \vdots & \vdots & & \vdots & \vdots & \vdots \\ a_{j1} & a_{j2} & & a_{ji} & a_{ji+1} & a_{jn} \\ \hline a_{j+11} & a_{j+i2} & & a_{j+1i} & a_{j+1\,i+1} & a_{j+1n} \\ \vdots & \vdots & & \vdots & \vdots & \vdots \\ a_{n1} & a_{n2} & \cdots & a_{ni} & a_{ni+1} & a_{nn} \end{array}\right]$$

A partitioned between the rows j1 and j + 11 and between the columns 1i and 1i + 1 where $i \neq j$ is a fuzzy supermatrix but is not a square fuzzy supermatrix.

We give yet another example.

***Example 2.1.18:*** Let A be a fuzzy square matrix.

$$A_s = \left[\begin{array}{ccc|cc} 1 & 0 & 0.3 & 0.5 & 0.2 \\ 0.3 & 1 & 0.6 & 0.5 & 0.1 \\ 0.1 & 0 & 1 & 0 & 1 \\ \hline 0.2 & 0.8 & 0.5 & 1 & 0.8 \\ 0.5 & 1 & 1 & 0 & 0.2 \end{array}\right].$$

The 5 × 5 square fuzzy matrix has been partitioned between the 3rd column and the 4th column and also partitioned between the 3rd row and the 4th row. Thus $A_s$ is a square fuzzy supermatrix.

Now having seen an example of a square fuzzy supermatrix we proceed on to define the notion of it.

**DEFINITION 2.1.8:** *Let A be a n × n fuzzy square matrix. We call A a super fuzzy square matrix or a square fuzzy supermatrix if A is partitioned in the following way.*

*If the matrix A is partitioned between the columns say i and i + 1, j and j+1 and t and t+1 then A is also partitioned between the rows i and i+1, j and j+1 and t and t+1, now to be more general it can be partitioned arbitrarily between the columns* $i_1$ *and* $i_{1+1}$, $i_2$ *and* $i_{2+1}$, …, $i_r$ *and* $i_{r+1}$ *and between the rows* $i_1$ *and* $i_{1+1}$, $i_2$ *and* $i_{2+1}$, …, $i_r$ *and* $i_{r+1}$ *(* $r+1 < n$*).*

Now having defined a fuzzy square supermatrix we proceed on to define the notion of symmetric square fuzzy matrix.

**DEFINITION 2.1.9:** *Let* $A_s$ *be a fuzzy super square matrix or fuzzy square supermatrix. We say* $A_s$ *is a symmetric fuzzy super square matrix or a fuzzy symmetric super square matrix or a symmetric square super fuzzy matrix. If*

$$A_s = \begin{bmatrix} A_{11} & A_{12} & \cdots & A_{1N} \\ A_{21} & A_{22} & \cdots & A_{2N} \\ \vdots & \vdots & & \vdots \\ A_{N1} & A_{N2} & \cdots & A_{NN} \end{bmatrix},$$

*then* $A_{11}, A_{22}, \ldots, A_{NN}$ *are square fuzzy matrices and each of these fuzzy square matrices are symmetric square matrices and* $A^t_{ij} = A_{ji}$ *for* $1 \leq i \leq N$ *and* $1 \leq j \leq N$*.*

We illustrate this by the following example.

***Example 2.1.19:*** Let $A_s$ be a fuzzy supermatrix.

$$
A_s = \left[\begin{array}{ccc|cc|ccc|c}
0.8 & 1 & 0.8 & 0 & 1 & 0.1 & 0.9 & 1 & 0 \\
1 & 0 & 0.6 & 1 & 0 & 0.3 & 0.8 & 0 & 0.6 \\
0.6 & 1 & 0.2 & 0 & 1 & 0.2 & 1 & 0.2 & 1 \\
\hline
0.5 & 0.4 & 1 & 0.5 & 0.7 & 1 & 0 & 0 & 0.2 \\
0.3 & 1 & 0.4 & 1 & 1 & 0.6 & 1 & 1 & 0.3 \\
\hline
0.6 & 0.7 & 1 & 0 & 0 & 1 & 1 & 1 & 0 \\
1.0 & 1 & 0.3 & 0.3 & 0.1 & 0.4 & 0.8 & 0 & 0 \\
0.9 & 0 & 1 & 1 & 1 & 1 & 1 & 0.6 & 0.1 \\
\hline
1.0 & 0.9 & 1 & 0.4 & 0.4 & 0 & 0.6 & 1 & 1
\end{array}\right]
$$

$$
= \begin{bmatrix}
A_{11} & A_{12} & A_{13} & A_{14} \\
A_{21} & A_{22} & A_{23} & A_{24} \\
A_{31} & A_{32} & A_{33} & A_{34} \\
A_{41} & A_{42} & A_{43} & A_{44}
\end{bmatrix}.
$$

Clearly $A_s$ is a fuzzy super square matrix. For it is a 9 × 9 fuzzy matrix, partitioned between the 3rd and 4th column and row. It is partitioned between the 5th and 6th row and column and between the 8th and 9th row and column. Clearly $A_s$ is not a symmetric fuzzy supermatrix for we see $A^t_{ij} \neq A_{ji}$.
We illustrate this by the following example.

***Example 2.1.20:*** Let $A_s$ be a fuzzy supermatrix i.e.,

$$
A_s = \left[\begin{array}{ccc|cc|cccc}
0 & 1 & 0.3 & 0.7 & 1 & 0 & 1 & 0.9 & 0.3 \\
1 & 0.2 & 1 & 0 & 0 & 1 & 0 & 1 & 0 \\
0.3 & 1 & 0.5 & 0.3 & 0.5 & 0 & 1 & 0 & 1 \\
\hline
0.7 & 0 & 0.3 & 1 & 0.2 & 0.3 & 0.4 & 0.5 & 0.6 \\
1 & 0 & 0.5 & 0.2 & 0.6 & 0.6 & 0.5 & 0.4 & 0.3 \\
\hline
0 & 1 & 0 & 0.3 & 0.6 & 0.7 & 1 & 0 & 0.9 \\
1 & 0 & 1 & 0.4 & 0.5 & 1 & 0.6 & 1 & 0.2 \\
0.9 & 1 & 0 & 0.5 & 0.4 & 0 & 1 & 0 & 0.8 \\
0.3 & 0 & 1 & 0.6 & 0.3 & 0.9 & 0.2 & 0.8 & 1
\end{array}\right]
$$

$$= \begin{bmatrix} A_{11} & A_{12} & A_{13} \\ A_{21} & A_{22} & A_{23} \\ A_{31} & A_{32} & A_{33} \end{bmatrix}.$$

We see $A_{11}$, $A_{22}$ and $A_{33}$ are symmetric fuzzy matrices. Further $A_{12} = A^t_{21}$ and $A_{13} = A^t_{31}$ and $A_{23} = A^t_{32}$. Thus $A_s$ is a symmetric fuzzy supermatrix, as $A_s$ is also a square matrix.

We have just shown examples that if A is a fuzzy super row vector then A'A (A' the transpose of A) under the super pseudo product is a symmetric fuzzy square supermatrix i.e., A'A is a symmetric fuzzy supermatrix. Thus we can say if one wants to construct fuzzy symmetric supermatrices then one can take a fuzzy super row vector A and find the pseudo super product of A and A'. A'A will always be a symmetric fuzzy supermatrix.

Now we are interested in the following problem which we first illustrate and the propose it as a problem.

***Example 2.1.21:*** Let

$$A = \begin{bmatrix} a_{11} & a_{12} \\ a_{21} & a_{22} \end{bmatrix}$$

be a fuzzy 2 × 2 square matrix.

In how many ways can A be partitioned so that $A_s$ the fuzzy supermatrix got from A is distinct.

1. $A_s$ = A no partitioning.

2. $A_s = \left[\begin{array}{c|c} a_{11} & a_{12} \\ a_{21} & a_{22} \end{array}\right]$ partitioned between first and second column.

3. $A_s = \left[\begin{array}{cc} a_{11} & a_{12} \\ \hline a_{21} & a_{22} \end{array}\right]$ partitioned between the first and second row.

4. $A_s = \left[\begin{array}{c|c} a_{11} & a_{12} \\ \hline a_{21} & a_{22} \end{array}\right]$

partitioned between the first row and second row and between the first column and second column, which will also be known as the cell partition of A. A 2 × 2 matrix has 4 cell and a n × m matrix A will have mn cells by a cell partition of A. Thus we have four ways by which a 2 × 2 fuzzy matrix can be partitioned so that A is a super fuzzy matrix; the first of course no partition or trivial only a fuzzy matrix; so we say we have only 3 partitions. So if no partition is made on the fuzzy matrix we will not call it as a fuzzy supermatrix.

***Example 2.1.22:*** Let us consider the matrix

$$A = \begin{bmatrix} a_{11} & a_{12} & a_{13} \\ a_{21} & a_{22} & a_{23} \\ a_{31} & a_{32} & a_{33} \end{bmatrix}$$

where $a_{ij} \in [0, 1]$; $1 \leq i \leq 3$ and $1 \leq j \leq 3$. We now physically enumerate all the possible partitions of A so that A becomes a fuzzy supermatrix.

1. $A_s = \left[\begin{array}{c|cc} a_{11} & a_{12} & a_{13} \\ a_{21} & a_{22} & a_{23} \\ a_{31} & a_{32} & a_{33} \end{array}\right]$, partitioned between first and second column.

2. $A_s = \left[\begin{array}{cc|c} a_{11} & a_{12} & a_{13} \\ a_{21} & a_{22} & a_{23} \\ a_{31} & a_{32} & a_{33} \end{array}\right]$, partitioned between column two and three.

3. $A_s = \left[\begin{array}{c|c|c} a_{11} & a_{12} & a_{13} \\ a_{21} & a_{22} & a_{23} \\ a_{31} & a_{32} & a_{33} \end{array}\right]$, partitioned between the first and second column and between the second and third column.

4. $A_s = \left[\begin{array}{ccc} a_{11} & a_{12} & a_{13} \\ \hline a_{21} & a_{22} & a_{23} \\ a_{31} & a_{32} & a_{33} \end{array}\right]$, partitioned between the first and second row.

5. $A_s = \left[\begin{array}{ccc} a_{11} & a_{12} & a_{13} \\ a_{21} & a_{22} & a_{23} \\ \hline a_{31} & a_{32} & a_{33} \end{array}\right]$ partitioned between the second row and third row.

6. $A_s = \left[\begin{array}{ccc} a_{11} & a_{12} & a_{13} \\ \hline a_{21} & a_{22} & a_{23} \\ \hline a_{31} & a_{32} & a_{33} \end{array}\right]$ partitioned between the first and second row and between the second and third row.

7. $A_s = \left[\begin{array}{c|cc} a_{11} & a_{12} & a_{13} \\ \hline a_{21} & a_{22} & a_{23} \\ a_{31} & a_{32} & a_{33} \end{array}\right]$ partitioned between the first row and second row and between the first column and the second column.

8. $A_s = \left[\begin{array}{c|cc} a_{11} & a_{12} & a_{13} \\ a_{21} & a_{22} & a_{23} \\ \hline a_{31} & a_{32} & a_{33} \end{array}\right]$ partitioned between the first and second column and partitioned between the second and third row.

9. $A_s = \left[\begin{array}{c|cc} a_{11} & a_{12} & a_{13} \\ \hline a_{21} & a_{22} & a_{23} \\ \hline a_{31} & a_{32} & a_{33} \end{array}\right]$ partitioned, between the first and second column. Partitioned between the first and second row and between the second and third row.

10. $A_s = \left[\begin{array}{cc|c} a_{11} & a_{12} & a_{13} \\ \hline a_{21} & a_{22} & a_{23} \\ a_{31} & a_{32} & a_{33} \end{array}\right]$ is partitioned between the $2^{nd}$ and $3^{rd}$ column and partitioned between the first and second row.

11. $A_s = \left[\begin{array}{cc|c} a_{11} & a_{12} & a_{13} \\ a_{21} & a_{22} & a_{23} \\ \hline a_{31} & a_{32} & a_{33} \end{array}\right]$ partitioned between the $2^{nd}$ and $3^{rd}$ column and second and $3^{rd}$ row.

12. $A_s = \left[\begin{array}{cc|c} a_{11} & a_{12} & a_{13} \\ \hline a_{21} & a_{22} & a_{23} \\ \hline a_{31} & a_{32} & a_{33} \end{array}\right]$ is partitioned between the $2^{nd}$ and $3^{rd}$ column and between the first and second row and second and third row.

13. $A_s = \left[\begin{array}{c|c|c} a_{11} & a_{12} & a_{13} \\ \hline a_{21} & a_{22} & a_{23} \\ a_{31} & a_{32} & a_{33} \end{array}\right]$ is partitioned between the second and $3^{rd}$ column and first and second column. Also $A_s$ is partitioned between the first and the $2^{nd}$ row.

14. Now $A_s = \left[\begin{array}{c|c|c} a_{11} & a_{12} & a_{13} \\ \hline a_{21} & a_{22} & a_{23} \\ \hline a_{31} & a_{32} & a_{33} \end{array}\right]$ is partitioned between the first row and second row and between the second row and third

row. Partitioned between the first and second column and between second and third column, the cell partition of A.

We call this partition as the cell partition of the fuzzy matrix.

We see we have 14 types of 3 × 3 fuzzy supermatrix but only 3 in case of 2 × 2 fuzzy supermatrix. Thus we now pose the following problem.

1. Find the number of distinct fuzzy supermatrix obtained from partitioning a square fuzzy n × n matrix.

2. Find the number of symmetric fuzzy supermatrix using a n × n matrix.

*Note:* In case of 3 × 3 fuzzy supermatrices we have only two partitions leading to symmetric fuzzy supermatrix given by

$$A_s = \left[\begin{array}{cc|c} a_{11} & a_{12} & a_{13} \\ a_{21} & a_{22} & a_{23} \\ \hline a_{31} & a_{32} & a_{33} \end{array}\right] = \begin{bmatrix} A_{11} & A_{12} \\ A_{21} & A_{22} \end{bmatrix}$$

and

$$A_s = \left[\begin{array}{c|cc} a_{11} & a_{12} & a_{13} \\ \hline a_{21} & a_{22} & a_{23} \\ a_{31} & a_{32} & a_{33} \end{array}\right] = \begin{bmatrix} A_{11} & A_{12} \\ A_{21} & A_{22} \end{bmatrix}.$$

With $A_{11}$ and $A_{22}$ square symmetric fuzzy matrices and $A_{12} = (A_{21})^t$ or $A_{21} = (A_{12})^t$.

Now when we take the collection of all 4 × 4 fuzzy supermatrices.

$$A = \begin{bmatrix} a_{11} & a_{12} & a_{13} & a_{14} \\ a_{21} & a_{22} & a_{23} & a_{24} \\ a_{31} & a_{32} & a_{33} & a_{34} \\ a_{41} & a_{42} & a_{43} & a_{44} \end{bmatrix}$$

where $a_{ij} \in [0, 1]$. How many partition of A will lead to fuzzy symmetric supermatrices?

***Example 2.1.23:*** We have six partitions which can lead to fuzzy symmetric matrices which is as follows:

1. $$A_s = \left[\begin{array}{c|ccc} a_{11} & a_{12} & a_{13} & a_{14} \\ \hline a_{21} & a_{22} & a_{23} & a_{24} \\ a_{31} & a_{32} & a_{33} & a_{34} \\ a_{41} & a_{42} & a_{43} & a_{44} \end{array}\right]$$

where we have partitioned between first row and second row and partitioned between first column and second column.

2. $$A_s = \left[\begin{array}{ccc|c} a_{11} & a_{12} & a_{13} & a_{14} \\ a_{21} & a_{22} & a_{23} & a_{24} \\ a_{31} & a_{32} & a_{33} & a_{34} \\ \hline a_{41} & a_{42} & a_{43} & a_{44} \end{array}\right];$$

we have a partition between the $3^{rd}$ and $4^{th}$ row as well as between $3^{rd}$ and $4^{th}$ column.

3. $$A_s = \left[\begin{array}{cc|cc} a_{11} & a_{12} & a_{13} & a_{14} \\ a_{21} & a_{22} & a_{23} & a_{24} \\ \hline a_{31} & a_{32} & a_{33} & a_{34} \\ a_{41} & a_{42} & a_{43} & a_{44} \end{array}\right];$$

where we partition between the $2^{nd}$ and $3^{rd}$ row and also between the $2^{nd}$ and $3^{rd}$ column.

4. $$A_s = \left[\begin{array}{c|cc|c} a_{11} & a_{12} & a_{13} & a_{14} \\ \hline a_{21} & a_{22} & a_{23} & a_{24} \\ a_{31} & a_{32} & a_{33} & a_{34} \\ \hline a_{41} & a_{42} & a_{43} & a_{44} \end{array}\right],$$

we have partitioned between the first and second row and first and second column. We have also partitioned between the 3rd and 4th row and the 3rd and 4th column of A to have a symmetric super fuzzy matrix.

5. $$A_s = \left[\begin{array}{c|c|cc} a_{11} & a_{12} & a_{13} & a_{14} \\ \hline a_{21} & a_{22} & a_{23} & a_{24} \\ \hline a_{31} & a_{32} & a_{33} & a_{34} \\ a_{41} & a_{42} & a_{43} & a_{44} \end{array}\right],$$

we have partitioned the fuzzy matrix between the 1st and 2nd column and row respectively and between the 2nd and 3rd row and column respectively to obtain a supermatrix which may lead to the construction of symmetric fuzzy supermatrices.

6. $$A_s = \left[\begin{array}{cc|c|c} a_{11} & a_{12} & a_{13} & a_{14} \\ a_{21} & a_{22} & a_{23} & a_{24} \\ \hline a_{31} & a_{32} & a_{33} & a_{34} \\ \hline a_{41} & a_{42} & a_{43} & a_{44} \end{array}\right]$$

The 4 × 4 fuzzy matrix A is partitioned between the 2nd and 3rd row and column respectively and between the 3rd and 4th row and column respectively to obtain a fuzzy supermatrix which can lead to a fuzzy symmetric supermatrix.

We illustrate this by an example. We take for this basically a fuzzy symmetric 4 × 4 supermatrix.

***Example 2.1.24:*** Let

$$A = \begin{bmatrix} 1 & 0.6 & 1 & 0.4 \\ 0.6 & 0.2 & 0.4 & 0 \\ 1 & 0.4 & 0 & 0.5 \\ 0.4 & 0 & 0.5 & 0.9 \end{bmatrix}$$

be a fuzzy symmetric 4 × 4 matrix.

Now the six super symmetric matrices got from A are as follows.

1. $$A_s^1 = \left[\begin{array}{c|ccc} 1 & 0.6 & 1 & 0.4 \\ \hline 0.6 & 0.2 & 0.4 & 0 \\ 1 & 0.4 & 0 & 0.5 \\ 0.4 & 0 & 0.5 & 0.9 \end{array}\right] = \begin{bmatrix} A_{11} & A_{12} \\ A_{21} & A_{22} \end{bmatrix}.$$

$A_{11}$ and $A_{22}$ are square symmetric fuzzy matrices and

$$A_{12} = (A_{21})^t = \begin{bmatrix} 0.6 \\ 1 \\ 0.4 \end{bmatrix}^t = [0.6\ 1\ 0.4].$$

2. $$A_s^2 = \left[\begin{array}{cc|cc} 1 & 0.6 & 1 & 0.4 \\ 0.6 & 0.2 & 0.4 & 0 \\ \hline 1 & 0.4 & 0 & 0.5 \\ 0.4 & 0 & 0.5 & 0.9 \end{array}\right] = \begin{bmatrix} A_{11} & A_{12} \\ A_{21} & A_{22} \end{bmatrix}.$$

All the four fuzzy submatrices of the fuzzy supermatrix $A_s^2$ are symmetric and are square fuzzy submatrices of $A^2_s$.

3. $$A_s^3 = \left[\begin{array}{ccc|c} 1 & 0.6 & 1 & 0.4 \\ 0.6 & 0.2 & 0.4 & 0 \\ 1 & 0.4 & 0 & 0.5 \\ \hline 0.4 & 0 & 0.5 & 0.9 \end{array}\right] = \begin{bmatrix} A_{11} & A_{12} \\ A_{21} & A_{22} \end{bmatrix},$$

the four submatrices $A_{11}$, $A_{12}$, $A_{21}$ and $A_{22}$ are such that $A_{11}$ and $A_{22}$ are square fuzzy matrices which are symmetric and

$$A_{12} = (A_{21})^t = (0.4\ 0\ 0.5)^t = \begin{bmatrix} 0.4 \\ 0 \\ 0.5 \end{bmatrix}.$$

$A_s^3$ is yet another symmetric super fuzzy matrix.

4. $$A_s^4 = \left[\begin{array}{c|c|cc} 1 & 0.6 & 1 & 0.4 \\ \hline 0.6 & 0.2 & 0.4 & 0 \\ \hline 1 & 0.4 & 0 & 0.5 \\ 0.4 & 0 & 0.5 & 0.9 \end{array}\right] = \begin{bmatrix} A_{11} & A_{12} & A_{13} \\ A_{21} & A_{22} & A_{23} \\ A_{31} & A_{32} & A_{33} \end{bmatrix}$$

where $A_{11}$, $A_{22}$ and $A_{33}$ are fuzzy square matrices which are symmetric.

$$A_{12} = A^t_{21}.$$

$$A_{13} = (A_{31})^t = [1\ 0.4] = \begin{bmatrix} 1 \\ 0.4 \end{bmatrix}^t.$$

Similarly

$$A_{23} = [0.4\ 0] = [A_{32}]^t = \begin{bmatrix} 0.4 \\ 0 \end{bmatrix}^t$$

Thus $A_s^4$ is a symmetric fuzzy supermatrix with submatrices.

5. $$A_s^5 = \left[\begin{array}{cc|c|c} 1 & 0.6 & 1 & 0.4 \\ 0.6 & 0.2 & 0.4 & 0 \\ \hline 1 & 0.4 & 0 & 0.5 \\ \hline 0.4 & 0 & 0.5 & 0.9 \end{array}\right] = \begin{bmatrix} A_{11} & A_{12} & A_{13} \\ A_{21} & A_{22} & A_{23} \\ A_{31} & A_{32} & A_{33} \end{bmatrix}$$

where $A_{11}$, $A_{12}$, $A_{13}$ , …, $A_{32}$ and $A_{33}$ are fuzzy submatrices with $A_{11}$, $A_{12}$ and $A_{33}$ fuzzy symmetric matrices.

$$A_{12} = \begin{bmatrix} 1 \\ 0.4 \end{bmatrix} \text{ and } A_{13} = \begin{bmatrix} 0.4 \\ 0 \end{bmatrix}$$

$(A_{21})^t = A_{12} = (1\ 0.4)^t$

$(A_{31})^t = A_{13} = (0.4\ 0)^t$.

Next we consider the super symmetric fuzzy matrix.

6. $A_s^6 = \left[\begin{array}{c|cc|c} 1 & 0.6 & 1 & 0.4 \\ \hline 0.6 & 0.2 & 0.4 & 0 \\ 1 & 0.4 & 0 & 0.5 \\ \hline 0.4 & 0 & 0.5 & 0.9 \end{array}\right] = \begin{bmatrix} A_{11} & A_{12} & A_{13} \\ A_{21} & A_{22} & A_{23} \\ A_{31} & A_{32} & A_{33} \end{bmatrix}$

where $A_{11}$, $A_{22}$ and $A_{33}$ are symmetric fuzzy matrices and

$$A_{12} = [0.6\ 1] = \begin{bmatrix} 0.6 \\ 1 \end{bmatrix}^t$$

$$A_{13} = [0.4] = A_{31}$$

$$A_{23} = \begin{bmatrix} 0 \\ 0.5 \end{bmatrix} = [0\ 0.5]^t.$$

Thus we have seen a single symmetric 4 × 4 fuzzy matrix can lead to 6 symmetric super fuzzy matrices.

***Example 2.1.25:*** Given any n × n fuzzy symmetric matrix A how many partitions be made on A so that $A_s$ is a super fuzzy symmetric matrix. The advantage of symmetric super fuzzy matrices we can arrive at several symmetric super fuzzy matrices given any symmetric fuzzy matrix. We have seen just now given a 4 × 4 symmetric fuzzy matrix, we can construct six super symmetric fuzzy matrices. We just see how many symmetric super fuzzy matrices can be constructed using a 5 × 5 symmetric fuzzy matrix.

$$A = \begin{bmatrix} a_{11} & a_{12} & a_{13} & a_{14} & a_{15} \\ a_{21} & a_{22} & a_{23} & a_{24} & a_{25} \\ a_{31} & a_{32} & a_{33} & a_{34} & a_{35} \\ a_{41} & a_{42} & a_{43} & a_{44} & a_{45} \\ a_{51} & a_{52} & a_{53} & a_{54} & a_{55} \end{bmatrix}$$

$$A'_s = \left[\begin{array}{c|ccc} a_{11} & a_{12} & \cdots & a_{15} \\ \hline a_{21} & a_{22} & \cdots & a_{25} \\ \vdots & \vdots & & \vdots \\ a_{51} & a_{52} & \cdots & a_{55} \end{array}\right] = \begin{bmatrix} A_{11} & A_{12} \\ A_{21} & A_{22} \end{bmatrix}.$$

A is partitioned between the first row and second row and first column and second column.

$$A_s^2 = \left[\begin{array}{ccc|c} a_{11} & a_{12} & \cdots & a_{15} \\ a_{21} & a_{22} & \cdots & a_{25} \\ \vdots & \vdots & & \vdots \\ \hline a_{51} & a_{52} & \cdots & a_{55} \end{array}\right].$$

A is partitioned between the fourth and fifth row and fourth and fifth column respectively.

$$A_s^3 = \left[\begin{array}{cc|ccc} a_{11} & a_{12} & a_{13} & \cdots & a_{15} \\ a_{21} & a_{22} & a_{23} & \cdots & a_{25} \\ \hline a_{31} & a_{32} & a_{33} & \cdots & \\ \vdots & \vdots & & & \vdots \\ a_{51} & a_{52} & a_{53} & \cdots & a_{55} \end{array}\right].$$

A is partitioned between the second and third row and between the second and third column.

$$A_s^4 = \left[\begin{array}{ccc|cc} a_{11} & a_{12} & a_{13} & a_{14} & a_{15} \\ a_{21} & a_{22} & a_{23} & a_{24} & a_{25} \\ a_{31} & a_{32} & a_{33} & a_{34} & a_{35} \\ \hline a_{41} & a_{42} & a_{43} & a_{44} & a_{45} \\ a_{51} & a_{52} & a_{53} & a_{54} & a_{55} \end{array}\right].$$

A is partitioned between the $3^{rd}$ and $4^{th}$ row and between $3^{rd}$ and $4^{th}$ column.

Now we proceed on to describe the symmetric super fuzzy matrix $A_s^5$

$$A_s^5 = \left[\begin{array}{c|ccc|c} a_{11} & a_{12} & a_{13} & \cdots & a_{15} \\ \hline a_{21} & a_{22} & a_{23} & \cdots & a_{25} \\ a_{31} & a_{32} & a_{33} & \cdots & a_{35} \\ \vdots & \vdots & \vdots & & \vdots \\ \hline a_{51} & a_{52} & a_{53} & \cdots & a_{55} \end{array}\right] = \begin{bmatrix} A_{11} & A_{12} & A_{13} \\ A_{21} & A_{22} & A_{23} \\ A_{31} & A_{32} & A_{33} \end{bmatrix},$$

where $A_{11}$, $A_{22}$ and $A_{33}$ are symmetric fuzzy matrices.

$$A_s^6 = \left[\begin{array}{ccc|c|c} a_{11} & a_{12} & a_{13} & a_{14} & a_{15} \\ a_{21} & a_{22} & a_{23} & a_{24} & a_{25} \\ a_{31} & a_{32} & a_{33} & a_{34} & a_{35} \\ \hline a_{41} & a_{42} & a_{43} & a_{44} & a_{45} \\ \hline a_{51} & a_{52} & a_{53} & a_{54} & a_{55} \end{array}\right] = \begin{bmatrix} A_{11} & A_{12} & A_{13} \\ A_{21} & A_{22} & A_{23} \\ A_{31} & A_{32} & A_{33} \end{bmatrix}.$$

$A_{11}$, $A_{22}$ and $A_{33}$ are symmetric fuzzy matrices.

$$A_s^7 = \left[\begin{array}{c|c|ccc} a_{11} & a_{12} & a_{13} & \cdots & a_{15} \\ \hline a_{21} & a_{22} & a_{23} & \cdots & a_{25} \\ \hline a_{31} & a_{32} & a_{33} & \cdots & \\ \vdots & \vdots & \vdots & & \vdots \\ a_{51} & a_{52} & a_{53} & \cdots & a_{55} \end{array}\right] = \begin{bmatrix} A_{11} & A_{12} & A_{13} \\ A_{21} & A_{22} & A_{23} \\ A_{31} & A_{32} & A_{33} \end{bmatrix}.$$

Here also $A_s^6$ is different from $A_s^7$ and $A_s^8$.

$$A_s^8 = \left[\begin{array}{cc|c|c|c} a_{11} & a_{12} & a_{13} & a_{14} & a_{15} \\ a_{21} & a_{22} & a_{23} & a_{24} & a_{25} \\ \hline a_{31} & a_{32} & a_{33} & a_{34} & a_{35} \\ \hline a_{41} & a_{42} & a_{43} & a_{44} & a_{45} \\ \hline a_{51} & a_{52} & a_{53} & a_{54} & a_{55} \end{array}\right] = \begin{bmatrix} A_{11} & A_{12} & A_{13} & A_{14} \\ A_{21} & A_{22} & A_{23} & A_{24} \\ A_{31} & A_{32} & A_{33} & A_{34} \\ A_{41} & A_{42} & A_{43} & A_{44} \end{bmatrix}$$

where $A_{11}$ , $A_{12}$ $A_{13}$ and $A_{44}$ are symmetric fuzzy matrices.

Thus we have more symmetric super fuzzy matrices obtained from a fuzzy symmetric 5 × 5 matrix.

## 2.2 Pseudo Symmetric Supermatrices

In this section we for the first time define the notion of pseudo diagonal, pseudo symmetric matrix and pseudo super symmetric matrix or pseudo symmetric supermatrix.

Now we proceed onto define the notion of a pseudo symmetric supermatrix.

**DEFINITION 2.2.1:** *We just call a fuzzy matrix to be pseudo symmetric if it is symmetric about the opposite diagonal.*

That is if $A = \begin{bmatrix} a_{11} & a_{12} \\ a_{21} & a_{22} \end{bmatrix}$ we see $a_{11} = a_{22}$ so that pseudo symmetry is about the pseudo diagonal $a_{12}$ and $a_{12}$. The diagonal vertically opposite to the usual main diagonal of a square matrix A will be known as the pseudo diagonal of A.

***Example 2.2.1:*** Let

$$A = \begin{bmatrix} 0.4 & 1 \\ 0 & 0.4 \end{bmatrix}$$

i.e., in case of a 2 × 2 matrix we need the diagonal elements to be equal. The pseudo diagonal elements are 1 and 0.

***Example 2.2.2:*** Let

$$A = \begin{bmatrix} 0.3 & 0.1 & 0.5 \\ 1 & 0.4 & 0.1 \\ 0.2 & 1 & 0.3 \end{bmatrix}.$$

A is a fuzzy pseudo symmetric matrix for the pseudo diagonal is 0.5, 0.4 and 0.2 and the elements of this matrix are symmetrically distributed about the pseudo diagonal.

One natural question would be, can any symmetric matrix be pseudo symmetric. The answer is yes.

For A = $\begin{bmatrix} 0.2 & 0.5 \\ 0.5 & 0.2 \end{bmatrix}$ is a pseudo symmetric fuzzy matrix as well as symmetric fuzzy matrix.

***Example 2.2.3:*** Let

$$A = \begin{bmatrix} 0.5 & 0.2 & 0.7 \\ 0.2 & 0.1 & 0.2 \\ 0.7 & 0.2 & 0.5 \end{bmatrix},$$

A is both a pseudo symmetric fuzzy matrix as well as a symmetric fuzzy matrix.

*Note:* For a 3 × 3 fuzzy matrix

$$A = \begin{bmatrix} a_{11} & a_{12} & a_{13} \\ a_{21} & a_{22} & a_{23} \\ a_{31} & a_{32} & a_{33} \end{bmatrix},$$

to be both pseudo symmetric and symmetric we must have

1. $a_{11} = a_{33}$
2. $a_{13} = a_{31}$
3. $a_{12} = a_{21} = a_{23} = a_{32}$.

***Example 2.2.4:*** Let A be a fuzzy 4 × 4 matrix

$$A = \begin{bmatrix} 0.3 & 0.4 & 0.5 & 0.7 \\ 0.4 & 0.1 & 0.2 & 0.5 \\ 0.5 & 0.2 & 0.1 & 0.4 \\ 0.7 & 0.5 & 0.4 & 0.3 \end{bmatrix},$$

A is both symmetric and pseudo symmetric. Let

$$A = \begin{bmatrix} a_{11} & a_{12} & a_{13} & a_{14} \\ a_{21} & a_{22} & a_{23} & a_{24} \\ a_{31} & a_{32} & a_{33} & a_{34} \\ a_{41} & a_{42} & a_{43} & a_{44} \end{bmatrix}$$

is both pseudo symmetric and a symmetric fuzzy matrix, if

$$a_{11} = a_{44}$$
$$a_{14} = a_{41}$$
$$a_{33} = a_{22}$$
$$a_{12} = a_{21} = a_{34} = a_{43}$$
$$a_{13} = a_{31} = a_{24} = a_{42}$$
$$a_{23} = a_{32}.$$

It is to be noted that in general a fuzzy supermatrix which is pseudo symmetric need not be symmetric. We can have several examples to prove our claim.

***Example 2.2.5:***

$$A_s = \left[\begin{array}{cc|cccc} 0.3 & 1 & 1 & 0 & 0.2 & 0.5 \\ 1 & 0 & 0.1 & 0.9 & 0.7 & 0.2 \\ 0.5 & 1 & 0.8 & 0.3 & 0.9 & 0 \\ 0 & 0.7 & 0.4 & 0.8 & 0.1 & 1 \\ \hline 0.6 & 0.2 & 0.7 & 1 & 0 & 1 \\ 0.7 & 0.6 & 0 & 0.5 & 1 & 0.3 \end{array}\right]$$

is a pseudo symmetric fuzzy supermatrix. Clearly $A_s$ is not a symmetric fuzzy supermatrix. We see

$$A_s = \begin{bmatrix} A_{11} & A_{12} \\ A_{21} & A_{22} \end{bmatrix}$$

Clearly $A_{12}$ and $A_{21}$ are square fuzzy matrices which are pseudo symmetric about the pseudo diagonal. For

$$A_{12} = \begin{bmatrix} 1 & 0 & 0.2 & 0.5 \\ 0.1 & 0.9 & 0.7 & 0.2 \\ 0.8 & 0.3 & 0.9 & 0 \\ 0.4 & 0.8 & 0.1 & 1 \end{bmatrix}$$

Clearly $A_{12}$ is pseudo symmetric fuzzy matrix about the pseudo diagonal.

$$A_{21} = \begin{bmatrix} 0.6 & 0.2 \\ 0.7 & 0.6 \end{bmatrix}$$

is clearly not a symmetric fuzzy matrix but only a pseudo symmetric fuzzy matrix. Further $A_{11}$ and $A_{22}$ are not even square matrices. Further

$$A_{11} = \begin{bmatrix} 0.3 & 1 \\ 1 & 0 \\ 0.5 & 1 \\ 0 & 0.7 \end{bmatrix} \text{ and } A_{22} = \begin{bmatrix} 0.7 & 1 & 0 & 1 \\ 0 & 0.5 & 1 & 0.3 \end{bmatrix}$$

We clearly see $A_{11}^t \neq A_{22}$ or $A_{22}^t \neq A_{11}$. It is interesting to observe that $A_{11}$ and $A_{22}$ are not square fuzzy matrices so they can never be symmetric or pseudo symmetric we give yet another example of a fuzzy pseudo symmetric supermatrix but we see the entries in $A_{11}$ and $A_{22}$ some way related!

***Example 2.2.6:*** Let

$$A_s = \left[\begin{array}{ccc|ccc} 0.2 & 1 & 0 & 0 & 0.1 & 0.4 \\ 0.1 & 0.5 & 0.4 & 0.2 & 1 & 0.1 \\ 0.7 & 0.8 & 0.6 & 0.8 & 0.2 & 0 \\ \hline 1 & 0.5 & 0.2 & 0.6 & 0.4 & 0 \\ 0.7 & 0 & 0.5 & 0.8 & 0.5 & 1 \\ 0.4 & 0.7 & 1 & 0.7 & 0.1 & 0.2 \end{array}\right] = \begin{bmatrix} A_{11} & A_{12} \\ A_{21} & A_{22} \end{bmatrix}$$

where $A_{11}$, $A_{12}$, $A_{21}$ and $A_{22}$ are fuzzy submatrices of the super fuzzy matrix $A_s$. We see $A_{11}$ is neither symmetric nor a pseudo symmetric fuzzy matrix $A_{12}$ is a pseudo symmetric fuzzy matrix $A_{21}$ is also a pseudo symmetric fuzzy matrix.

We also see $A_{11} \neq A_{22}^t$ or $A_{22} \neq A_{11}^t$. To this end we define a new notion called pseudo transpose of a matrix.

**DEFINITION 2.2.2:** *Let*

$$A = \begin{bmatrix} a_{11} & a_{12} & \cdots & a_{1n} \\ a_{21} & a_{22} & \cdots & a_{2n} \\ \vdots & \vdots & & \vdots \\ a_{m1} & a_{m2} & \cdots & a_{mn} \end{bmatrix}$$

*be a rectangular m × n matrix. The pseudo transpose of A denoted by*

$$A^{t_p} = \begin{bmatrix} a_{mn} & \cdots & a_{m2} & a_{m1} \\ \vdots & & \vdots & \vdots \\ a_{2n} & \cdots & a_{22} & a_{21} \\ a_{1n} & \cdots & a_{12} & a_{11} \end{bmatrix}.$$

$A^{t_p}$ *is n × m matrix called the pseudo transpose of A.*

We illustrate this by the following examples.

***Example 2.2.7:*** Let

$$A = \begin{bmatrix} 0.3 & 1 & 0.5 \\ 0.2 & 0 & 0.6 \\ 0.7 & 0.4 & 0 \\ 0 & 0.9 & 1 \\ 0.8 & 1 & 0.4 \\ 0.9 & 0.5 & 0.2 \end{bmatrix}$$

be any 6 × 3 fuzzy rectangular matrix. The pseudo transpose of A denoted by $A^{t_p}$ is given by

$$A^{t_p} = \begin{bmatrix} 0.2 & 0.4 & 1 & 0 & 0.6 & 0.5 \\ 0.5 & 1 & 0.9 & 0.4 & 0 & 1 \\ 0.9 & 0.8 & 0 & 0.7 & 0.2 & 0.3 \end{bmatrix}.$$

$A^{t_p}$ is a 3 × 6 matrix, clearly $A^t \neq A^{t_p}$. For

$$A^t = \begin{bmatrix} 0.3 & 0.2 & 0.7 & 0 & 0.8 & 0.9 \\ 1 & 0 & 0.4 & 0.9 & 1 & 0.5 \\ 0.5 & 0.6 & 0 & 1 & 0.4 & 0.2 \end{bmatrix}.$$

Now when we define the notion of pseudo symmetric matrix we need the notion of pseudo transpose. We give yet another example of a pseudo symmetric super fuzzy matrix.

***Example 2.2.8:*** Let

$$A = \left[\begin{array}{c|ccc} 0.2 & 0.1 & 0.5 & 0.4 \\ 0.7 & 0.9 & 0.2 & 0.5 \\ 0.3 & 0.8 & 0.9 & 0.1 \\ \hline 0.6 & 0.3 & 0.7 & 0.2 \end{array}\right] = \begin{bmatrix} A_{11} & A_{12} \\ A_{21} & A_{22} \end{bmatrix}$$

be a fuzzy pseudo symmetric supermatrix with $A_{11}$, $A_{12}$, $A_{21}$ and $A_{22}$ as its submatrices. Clearly $A_{12}$ and $A_{21}$ are fuzzy square submatrices of A which are pseudo symmetric. Now

$$A_{11} = \begin{bmatrix} 0.2 \\ 0.7 \\ 0.3 \end{bmatrix}.$$

Then the pseudo transpose of $A_{11}$ is $A_{22}$ i.e., $A_{22} = A_{11}^{t_p} = [0.3\ 0.7\ 0.2]$. Clearly $A_{11}^{t} = [0.2\ 0.7\ 0.3] \neq A_{11}^{t_p}$.

Now we see how the pseudo transpose of a row matrix looks like.

**DEFINITION 2.2.3:** *Let $a = [a_1\ a_2 \ldots a_n]$ be a $1 \times n$ row matrix. The pseudo transpose of a denoted by*

$$a^{t_p} = \begin{bmatrix} a_n \\ a_{n-1} \\ \vdots \\ a_2 \\ a_1 \end{bmatrix}$$

*is a column matrix but is not the same as*

$$a^t = \begin{bmatrix} a_1 \\ a_2 \\ \vdots \\ a_n \end{bmatrix}.$$

*Let*

$$b = \begin{bmatrix} b_1 \\ b_2 \\ \vdots \\ b_m \end{bmatrix}$$

*be a $m \times 1$ column matrix. The pseudo transpose of b defined as $b^{t_p} = [b_m\ b_{m-1} \ldots b_2\ b_1]$. Clearly $b^{t_p} \neq b^t = [b_1\ b_2 \ldots b_m]$. It is important to note if*

$$A = \begin{bmatrix} a_{11} & a_{12} & \cdots & a_{1n} \\ a_{21} & a_{22} & \cdots & a_{2n} \\ \vdots & \vdots & & \vdots \\ a_{n1} & a_{n2} & \cdots & a_{nn} \end{bmatrix}$$

*be a $n \times n$ square matrix its pseudo transpose,*

$$A^{t_p} = \begin{bmatrix} a_{nn} & \cdots & a_{n2} & a_{n1} \\ \vdots & & \vdots & \vdots \\ a_{2n} & \cdots & a_{22} & a_{21} \\ a_{1n} & \cdots & a_{12} & a_{11} \end{bmatrix}$$

***Example 2.2.9:*** Let

$$A = \begin{bmatrix} 0.3 & 0.2 & 1 \\ 0.8 & 0.4 & 0.2 \\ 0.5 & 0.8 & 0.3 \end{bmatrix}$$

be a pseudo symmetric fuzzy matrix.
The pseudo transpose of A is given by

$$A^{t_p} = \begin{bmatrix} 0.3 & 0.2 & 1 \\ 0.8 & 0.4 & 0.2 \\ 0.5 & 0.8 & 0.3 \end{bmatrix}.$$

We see $A^{t_p} = A$. Thus we can define a pseudo symmetric matrix A to be a matrix in which its pseudo transpose is the same as the matrix i.e., $A = A^{t_p}$.

We now give in the following example the pseudo transpose of a pseudo symmetric matrix A.

***Example 2.2.10:*** Let

$$A = \begin{bmatrix} 0.3 & 0.1 & 1 & 0 & 0.2 \\ 1 & 0.2 & 0.6 & 0.9 & 0 \\ 0.7 & 1 & 0.3 & 0.6 & 1 \\ 0.5 & 1 & 1 & 0.2 & 0.1 \\ 0 & 0.5 & 0.7 & 1 & 0.3 \end{bmatrix}.$$

$$A^{t_p} = \begin{bmatrix} 0.3 & 0.1 & 1 & 0 & 0.2 \\ 1 & 0.2 & 0.6 & 0.9 & 0 \\ 0.7 & 1 & 0.3 & 0.6 & 1 \\ 0.5 & 1 & 1 & 0.2 & 0.1 \\ 0 & 0.5 & 0.7 & 1 & 0.3 \end{bmatrix}.$$

The pseudo transpose of A is $A^{t_p}$ so A is a pseudo symmetric fuzzy matrix.

Let us now find the pseudo transpose of a 2 × 5 fuzzy matrix.

***Example 2.2.11:*** Let

$$A = \begin{bmatrix} 0.1 & 0.8 & 1 & 0 & 0.7 \\ 0.4 & 0 & 1 & 0.7 & 0.6 \end{bmatrix}$$

be a 2 × 5 fuzzy rectangular matrix. The pseudo transpose of A is

$$A^{t_p} = \begin{bmatrix} 0.6 & 0.7 \\ 0.7 & 0 \\ 1 & 1 \\ 0 & 0.8 \\ 0.4 & 0.1 \end{bmatrix}.$$

$A^{t_p}$ is a 5 × 2 fuzzy rectangular matrix.

What is the pseudo transpose of a symmetric fuzzy matrix A? To this we give an example

***Example 2.2.12:*** Let

$$A = \begin{bmatrix} 0.2 & 1 & 0.6 & 0.3 \\ 1 & 0 & 0.7 & 0.5 \\ 0.6 & 0.7 & 1 & 0.8 \\ 0.3 & 0.5 & 0.8 & 0.4 \end{bmatrix}$$

be a 4 × 4 symmetric fuzzy matrix. The pseudo transpose of A;

$$A^{t_p} = \begin{bmatrix} 0.4 & 0.8 & 0.5 & 0.3 \\ 0.8 & 1 & 0.7 & 0.6 \\ 0.5 & 0.7 & 0 & 1 \\ 0.3 & 0.6 & 1 & 0.2 \end{bmatrix}.$$

Clearly we see $A^{t_p}$ is also a symmetric matrix but it is not equal to A. Thus one can by pseudo transpose obtain from a symmetric matrix another new symmetric matrix.

Now using the notion of pseudo transpose we can now define a symmetric matrix.

**DEFINITION 2.2.4:** *Let A be a square matrix if $A^{t_p}$ the pseudo transpose is a symmetric matrix then A is also a symmetric matrix, but different from A.*

***Example 2.2.13:*** Let

$$A = \begin{bmatrix} 0.5 & 0.3 & 0.7 \\ 0.4 & 1 & 0 \\ 0.8 & 0.9 & 0.1 \end{bmatrix}$$

be a 3 × 3 fuzzy matrix.

$$A^{t_p} = \begin{bmatrix} 0.1 & 0 & 0.7 \\ 0.9 & 1 & 0.3 \\ 0.8 & 0.4 & 0.5 \end{bmatrix}$$

is its pseudo transpose of A.
Now let

$$A^{t_p} = \begin{bmatrix} 0.1 & 0 & 0.7 \\ 0.9 & 1 & 0.3 \\ 0.8 & 0.4 & 0.5 \end{bmatrix} = B.$$

To find

$$B^{t_p} = \begin{bmatrix} 0.5 & 0.3 & 0.7 \\ 0.4 & 1 & 0 \\ 0.8 & 0.9 & 0.1 \end{bmatrix}.$$

Thus we see $\left(A^{t_p}\right)^{t_p}$ = A we can prove the following theorem for square matrices.

**THEOREM 2.2.1:** *Let*

$$A = \begin{bmatrix} a_{11} & a_{12} & \cdots & a_{1n} \\ a_{21} & a_{22} & \cdots & a_{2n} \\ \vdots & \vdots & & \vdots \\ a_{n1} & a_{n2} & \cdots & a_{nn} \end{bmatrix}$$

*be a n × n square matrix. Then $\left(A^{t_p}\right)^{t_p}$ = A.*

*Proof:* Given

$$A = \begin{bmatrix} a_{11} & a_{12} & \cdots & a_{1n-1} & a_{1n} \\ a_{21} & a_{22} & \cdots & a_{2n-1} & a_{2n} \\ \vdots & \vdots & & \vdots & \vdots \\ a_{n1} & a_{n2} & \cdots & a_{nn-1} & a_{nn} \end{bmatrix}$$

$$A^{t_p} = \begin{bmatrix} a_{nn} & \cdots & a_{2n} & a_{1n} \\ a_{nn-1} & \cdots & a_{2n-1} & a_{1n-1} \\ \vdots & & \vdots & \vdots \\ a_{n2} & \cdots & a_{22} & a_{12} \\ a_{n1} & \cdots & a_{21} & a_{11} \end{bmatrix}.$$

Now consider

$$\left(A^{t_p}\right)^{t_p} = \begin{bmatrix} a_{11} & a_{12} & \cdots & a_{1n-1} & a_{n1} \\ a_{21} & a_{22} & \cdots & a_{2n-1} & a_{n2} \\ \vdots & \vdots & & \vdots & \vdots \\ a_{n1} & a_{n2} & \cdots & a_{nn-1} & a_{nn} \end{bmatrix} = A.$$

Hence the claim.

**THEOREM 2.2.2:** *Let A be a n × n symmetric matrix then $A^{t_p}$ is another symmetric matrix.*

*Proof:* Given

$$A = \begin{bmatrix} a_{11} & a_{12} & \cdots & a_{1n} \\ a_{21} & a_{22} & \cdots & a_{2n} \\ \vdots & \vdots & & \vdots \\ a_{n1} & a_{n2} & \cdots & a_{nn} \end{bmatrix}$$

is a symmetric matrix as $a_{ij} = a_{ji}$ $(i \neq j)$ Now

$$A^{t_p} = \begin{bmatrix} a_{nn} & a_{n-1n} & \cdots & a_{2n} & a_{1n} \\ \vdots & \vdots & & \vdots & \vdots \\ a_{2n} & a_{2n-1} & \cdots & a_{22} & a_{12} \\ a_{1n} & a_{1n-1} & \cdots & a_{12} & a_{11} \end{bmatrix}.$$

Clearly as $a_{ij} = a_{ji}$ this is also a symmetric matrix.
Clearly it is easy to verify that if we see $A^{t_p}$ is also a symmetric matrix different from A.

Now having defined the notion of pseudo symmetry and obtained the notion of symmetry in terms of pseudo symmetry we proceed on to apply these concepts in case of fuzzy supermatrices. We have already defined the notion of symmetric fuzzy supermatrix; for which we said we have to partition only a fuzzy symmetric matrix and obtain several fuzzy super symmetric matrices from a given fuzzy symmetric matrix.

Now to define the fuzzy pseudo symmetric supermatrix we need to define the notion of pseudo partition of a square matrix for that alone can lead to the notion of pseudo symmetric supermatrix. We follow some sort of order in the pseudo partition and when the matrix is a rectangular one we cannot define the notion of pseudo partition for we do not have the notion of symmetric matrix or pseudo symmetric matrix in case of rectangular matrices.

**DEFINITION 2.2.5:** *Let A be a square n × n matrix i.e., A = ($a_{ij}$), $1 \leq i \leq n$ and $1 \leq j \leq n$.*

$$A = \begin{bmatrix} a_{11} & a_{12} & \cdots & a_{1n-1} & a_{1n} \\ a_{21} & a_{22} & \cdots & a_{2n-1} & a_{2n} \\ \vdots & \vdots & & \vdots & \vdots \\ a_{n-11} & a_{n-12} & & a_{n-1n-1} & a_{n-1n} \\ a_{n1} & a_{n2} & \cdots & a_{nn-1} & a_{nn} \end{bmatrix}$$

*be a n × n square matrix. The pseudo partition is carried out in this manner. Suppose we partition between r and [r+1]th column then we partition the row between the n – (r – 1)th row and (n – r)th row. Conversely if we partition the row between the rows s and (s+1)th row then we partition the column between [n – (s – 1)]$^{th}$ column and (n – s)$^{th}$ column. Such form of partitioning a square matrix is known as the pseudo partition.*

First we illustrate this by few examples.

***Example 2.2.14:*** Let A be a 7 × 7 matrix.

$$A = \left[\begin{array}{cc|ccc|cc} 0 & 1 & 2 & 5 & 0 & 2 & 1 \\ 1 & 1 & 2 & 0 & 7 & 4 & 3 \\ \hline 0 & 1 & 1 & 0 & 1 & 0 & 1 \\ 2 & 1 & 2 & 1 & 2 & 1 & 0 \\ 3 & 5 & 1 & 0 & 0 & 1 & 1 \\ \hline 1 & 0 & 1 & 1 & 2 & 2 & 0 \\ 0 & 1 & 0 & 1 & 2 & 0 & 2 \end{array}\right]$$

We partition between the 2$^{nd}$ and 3$^{rd}$ column; then we have to carry out the partition between the 5$^{th}$ and 6$^{th}$ row. Similarly for partitioning between 2$^{nd}$ and 3$^{rd}$ row we have to carry out a partition of the 5$^{th}$ and 6$^{th}$ column. Only this partition we call as a pseudo partition of the square matrix. Now we show yet another example to show that in general all partitions are not pseudo partitions but all pseudo partitions are partitions seen by the following example.

***Example 2.2.15:*** Let A be a 5 × 5 matrix

$$A = \left[\begin{array}{cc|ccc} 0 & 1 & 2 & 3 & 4 \\ 4 & 2 & 1 & 4 & 7 \\ \hline 0 & 1 & 0 & 1 & 0 \\ 2 & -1 & 2 & 8 & 6 \\ 1 & 2 & 0 & 1 & 2 \end{array}\right]$$

This clearly a partition which is not a pseudo partition of A.
Now we proceed on to give yet another example to show that certain partition of a square matrix can always be both a partition as well as a pseudo partition.

***Example 2.2.16:*** Let A be a 6 × 6 matrix.

$$A = \left[\begin{array}{ccc|ccc} 3 & 0 & 1 & 2 & 5 & 0 \\ 2 & 2 & 0 & 1 & 0 & 1 \\ 1 & 0 & 1 & 0 & 1 & 0 \\ \hline 0 & 1 & -1 & 1 & -1 & 1 \\ -1 & 0 & 1 & 2 & 3 & 4 \\ 5 & 6 & 7 & 8 & 9 & 1 \end{array}\right].$$

This is both a partition and a pseudo partition i.e., it can be performed as a partition which will result in a pseudo partition. A condition for it is stated below as a theorem.

**THEOREM 2.2.3:** *Let A be a 2n × 2n matrix. A partition of A by dividing the n and (n+1)th row or (n and (n+1)th column) and dividing (2n – (n – 1))th and (2n – n)th column or (2n –(n – 1)th and (2n – n)th row) is always a pseudo partition.*

*Proof:* Obvious by the very construction of the partition.

The interested reader can find nice characterization theorem to show conditions under which a partition is a pseudo partition.

***Example 2.2.17:*** Consider the matrix A;

$$A = \left[\begin{array}{ccc|cc} 1 & 2 & 3 & 0 & 4 \\ 5 & 1 & 0 & 2 & -1 \\ \hline 2 & 4 & 1 & 0 & 2 \\ -1 & 0 & 5 & 1 & 0 \\ 0 & 2 & 1 & 1 & 5 \end{array}\right]$$

This is a partition for we have partitioned between the 3rd and 4th column and between the 2nd and 3rd row. This is also a pseudo partition. We see if we pseudo partition a symmetric matrix it need not in general be a pseudo symmetric matrix.

**THEOREM 2.2.4:** *Let A be a symmetric matrix. A pseudo partition of A in general does not make A; a pseudo symmetric supermatrix.*

*Proof:* We prove this by an example first, though we can also prove the result for any n × n symmetric matrix.

Take the symmetric matrix.

$$A = \left[\begin{array}{cc|ccc} 0 & 1 & 2 & 3 & 4 \\ 1 & -1 & 0 & 1 & 0 \\ 2 & 0 & 2 & 1 & 2 \\ \hline 3 & 1 & 1 & 3 & 1 \\ 4 & 0 & 2 & 1 & 5 \end{array}\right]$$

Now make a pseudo partition on A by dividing (or partitioning) at the 2nd and 3rd column and at the [5 – (2 – 1)] = 4th row and the [5 – (3 – 1)] = 3rd row. Clearly A is not a pseudo symmetric matrix only a supermatrix which is neither super symmetric nor pseudo super symmetric. We make yet another theorem.

**THEOREM 2.2.5:** *Let A be any pseudo symmetric matrix. A symmetric partition of A does not in general make A to be super symmetric matrix.*

*Proof:* We show this by an example.

$$A = \begin{bmatrix} 0 & 1 & 2 & 3 & 4 \\ -1 & 2 & 4 & 5 & 3 \\ 7 & 6 & -1 & 4 & 2 \\ 3 & 9 & 6 & 2 & 1 \\ 1 & 3 & 7 & -1 & 0 \end{bmatrix}$$

is a pseudo symmetric matrix. We give a symmetric partition of A i.e.,

$$A_s = \left[\begin{array}{cc|ccc} 0 & 1 & 2 & 3 & 4 \\ -1 & 2 & 4 & 5 & 3 \\ \hline 7 & 6 & -1 & 4 & 2 \\ 3 & 9 & 6 & 2 & 1 \\ 1 & 3 & 7 & -1 & 0 \end{array}\right]$$

$A_s$ is not a symmetric supermatrix with the given partition hence the claim.

Thus we have seen in general a symmetric matrix cannot be made by a pseudo symmetric partition into a symmetric matrix or equally a pseudo symmetric matrix cannot be made into a symmetric matrix by a symmetric partition. Thus we see a pseudo symmetric matrix can be easily made into a pseudo symmetric supermatrix by a pseudo symmetric partition and a symmetric matrix can be made into a symmetric supermatrix by a symmetric partition.

Having made such observations we proceed onto illustrate this by some simple examples.

***Example 2.2.18:*** Let A be a pseudo symmetric matrix.

$$A = \begin{bmatrix} 0 & 1 & 2 & 3 & 4 & 5 & 6 \\ 1 & 0 & -1 & 2 & 7 & -1 & 5 \\ 2 & 5 & 6 & 8 & 9 & 7 & 4 \\ -1 & 4 & 2 & 3 & 8 & 2 & 3 \\ 6 & 2 & 7 & 2 & 6 & -1 & 2 \\ 5 & 3 & 2 & 4 & 5 & 0 & 1 \\ 0 & 5 & 6 & -1 & 2 & 1 & 0 \end{bmatrix}.$$

We have some of the pseudo partitions on A so that A is a pseudo symmetric supermatrix.

$$A_s^1 = \left[\begin{array}{c|cccccc} 0 & 1 & 2 & 3 & 4 & 5 & 6 \\ 1 & 0 & -1 & 2 & 7 & -1 & 5 \\ 2 & 5 & 6 & 8 & 9 & 7 & 4 \\ -1 & 4 & 2 & 3 & 8 & 2 & 3 \\ 6 & 2 & 7 & 2 & 6 & -1 & 2 \\ 5 & 3 & 2 & 4 & 5 & 0 & 1 \\ \hline 0 & 5 & 6 & -1 & 2 & 1 & 0 \end{array}\right]$$

$$A_s^2 = \left[\begin{array}{cc|ccccc} 0 & 1 & 2 & 3 & 4 & 5 & 6 \\ 1 & 0 & -1 & 2 & 7 & -1 & 5 \\ 2 & 5 & 6 & 8 & 9 & 7 & 4 \\ -1 & 4 & 2 & 3 & 8 & 2 & 3 \\ 6 & 2 & 7 & 2 & 6 & -1 & 2 \\ \hline 5 & 3 & 2 & 4 & 5 & 0 & 1 \\ 0 & 5 & 6 & -1 & 2 & 1 & 0 \end{array}\right]$$

$$A_s^3 = \left[\begin{array}{ccc|cccc} 0 & 1 & 2 & 3 & 4 & 5 & 6 \\ 1 & 0 & -1 & 2 & 7 & -1 & 5 \\ 2 & 5 & 6 & 8 & 9 & 7 & 4 \\ -1 & 4 & 2 & 3 & 8 & 2 & 3 \\ \hline 6 & 2 & 7 & 2 & 6 & -1 & 2 \\ 5 & 3 & 2 & 4 & 5 & 0 & 1 \\ 0 & 5 & 6 & -1 & 2 & 1 & 0 \end{array}\right]$$

$$A_s^4 = \left[\begin{array}{cccc|ccc} 0 & 1 & 2 & 3 & 4 & 5 & 6 \\ 1 & 0 & -1 & 2 & 7 & -1 & 5 \\ 2 & 5 & 6 & 8 & 9 & 7 & 4 \\ \hline -1 & 4 & 2 & 3 & 8 & 2 & 3 \\ 6 & 2 & 7 & 2 & 6 & -1 & 2 \\ 5 & 3 & 2 & 4 & 5 & 0 & 1 \\ 0 & 5 & 6 & -1 & 2 & 1 & 0 \end{array}\right]$$

$$A_s^5 = \left[\begin{array}{ccccc|cc} 0 & 1 & 2 & 3 & 4 & 5 & 6 \\ 1 & 0 & -1 & 2 & 7 & -1 & 5 \\ \hline 2 & 5 & 6 & 8 & 9 & 7 & 4 \\ -1 & 4 & 2 & 3 & 8 & 2 & 3 \\ 6 & 2 & 7 & 2 & 6 & -1 & 2 \\ 5 & 3 & 2 & 4 & 5 & 0 & 1 \\ 0 & 5 & 6 & -1 & 2 & 1 & 0 \end{array}\right]$$

$$A_s^6 = \left[\begin{array}{cccccc|c} 0 & 1 & 2 & 3 & 4 & 5 & 6 \\ \hline 1 & 0 & -1 & 2 & 7 & -1 & 5 \\ 2 & 5 & 6 & 8 & 9 & 7 & 4 \\ -1 & 4 & 2 & 3 & 8 & 2 & 3 \\ 6 & 2 & 7 & 2 & 6 & -1 & 2 \\ 5 & 3 & 2 & 4 & 5 & 0 & 1 \\ 0 & 5 & 6 & -1 & 2 & 1 & 0 \end{array}\right]$$

$$A_s^7 = \left[\begin{array}{cc|cc|ccc} 0 & 1 & 2 & 3 & 4 & 5 & 6 \\ 1 & 0 & -1 & 2 & 7 & -1 & 5 \\ 2 & 5 & 6 & 8 & 9 & 7 & 4 \\ \hline -1 & 4 & 2 & 3 & 8 & 2 & 3 \\ 6 & 2 & 7 & 2 & 6 & -1 & 2 \\ \hline 5 & 3 & 2 & 4 & 5 & 0 & 1 \\ 0 & 5 & 6 & -1 & 2 & 1 & 0 \end{array}\right].$$

Thus we have given in this example 7 pseudo symmetric supermatrices.

Next we proceed onto give an example of a symmetric supermatrix.

***Example 2.2.19:*** Let A be a 6 × 6 symmetric matrix.

$$A = \left[\begin{array}{cccccc} 2 & 1 & 0 & 1 & -2 & 5 \\ 1 & 6 & 8 & 2 & 1 & 0 \\ 0 & 8 & 4 & 3 & 9 & 7 \\ 1 & 2 & 3 & 0 & 6 & 4 \\ -2 & 1 & 9 & 6 & 10 & 5 \\ 5 & 0 & 7 & 4 & 5 & 1 \end{array}\right]$$

$$A_s^1 = \left[\begin{array}{c|ccccc} 2 & 1 & 0 & 1 & -2 & 5 \\ \hline 1 & 6 & 8 & 2 & 1 & 0 \\ 0 & 8 & 4 & 3 & 9 & 7 \\ 1 & 2 & 3 & 0 & 6 & 4 \\ -2 & 1 & 9 & 6 & 10 & 5 \\ 5 & 0 & 7 & 4 & 5 & 1 \end{array}\right].$$

This is a symmetric supermatrix of order 6.

$$A_s^2 = \left[\begin{array}{cc|cccc} 2 & 1 & 0 & 1 & -2 & 5 \\ 1 & 6 & 8 & 2 & 1 & 0 \\ \hline 0 & 8 & 4 & 3 & 9 & 7 \\ 1 & 2 & 3 & 0 & 6 & 4 \\ -2 & 1 & 9 & 6 & 10 & 5 \\ 5 & 0 & 7 & 4 & 5 & 1 \end{array}\right]$$

$$A_s^3 = \left[\begin{array}{cccc|cc} 2 & 1 & 0 & 1 & -2 & 5 \\ 1 & 6 & 8 & 2 & 1 & 0 \\ 0 & 8 & 4 & 3 & 9 & 7 \\ 1 & 2 & 3 & 0 & 6 & 4 \\ \hline -2 & 1 & 9 & 6 & 10 & 5 \\ 5 & 0 & 7 & 4 & 5 & 1 \end{array}\right]$$

$$A_s^4 = \left[\begin{array}{ccc|ccc} 2 & 1 & 0 & 1 & -2 & 5 \\ 1 & 6 & 8 & 2 & 1 & 0 \\ 0 & 8 & 4 & 3 & 9 & 7 \\ \hline 1 & 2 & 3 & 0 & 6 & 4 \\ -2 & 1 & 9 & 6 & 10 & 5 \\ 5 & 0 & 7 & 4 & 5 & 1 \end{array}\right].$$

We can make several such symmetric supermatrices.

The following interesting problems are proposed for the reader.

1. Given any n × n symmetric matrix A find the number of partition on it so that A is a symmetric supermatrix.

2. Given any n × n pseudo symmetric matrix $A_p$ find the number of pseudo partition on $A_p$ so that $A_p$ is a pseudo symmetric supermatrix.

***Example 2.2.20:*** Let us consider the pseudo symmetric matrix

$$A = \begin{bmatrix} 1 & 0 & 2 & 5 & 7 \\ 3 & 4 & 9 & 3 & 5 \\ -1 & 1 & 8 & 9 & 2 \\ 5 & 3 & 1 & 4 & 0 \\ 0 & 5 & -1 & 3 & 1 \end{bmatrix}.$$

Any pseudo symmetric partition of A will make A always a pseudo symmetric supermatrix.

But a partition will not always make a symmetric matrix into a symmetric supermatrix; this is illustrated by the following example

***Example 2.2.21:*** Let A be a symmetric 7 × 7 matrix.

$$A = \begin{bmatrix} 6 & 2 & 1 & 0 & 5 & 4 & 9 \\ 2 & 3 & 0 & 1 & 0 & 1 & 0 \\ 1 & 0 & 4 & 1 & 1 & 1 & 2 \\ 0 & 1 & 1 & 5 & 2 & 2 & 3 \\ 5 & 0 & 1 & 2 & 8 & 3 & 1 \\ 4 & 1 & 1 & 2 & 3 & 9 & 7 \\ 9 & 0 & 2 & 3 & 1 & 7 & 0 \end{bmatrix}.$$

Now consider a partition of A between the rows 2 and 3 and between the columns 4 and 5. Let $A_s$ be the resulting supermatrix.

$$A_s = \left[\begin{array}{cccc|ccc} 6 & 2 & 1 & 0 & 5 & 4 & 9 \\ 2 & 3 & 0 & 1 & 0 & 1 & 0 \\ \hline 1 & 0 & 4 & 1 & 1 & 1 & 2 \\ 0 & 1 & 1 & 5 & 2 & 2 & 3 \\ 5 & 0 & 1 & 2 & 8 & 3 & 1 \\ 4 & 1 & 1 & 2 & 3 & 9 & 7 \\ 9 & 0 & 2 & 3 & 1 & 7 & 0 \end{array}\right]$$

Clearly $A_s$ is not symmetric through A is a symmetric matrix. Thus every partition of a symmetric matrix need not in general lead to a symmetric supermatrix. Thus we call a partition to be a symmetric partition if A is itself only a square matrix. For a symmetric partition in general cannot be defined on a rectangular matrix. So if A is any square matrix, we say a partition of A is said to be a symmetric partition of A if the r and (r+1)th row is partitioned for A then it is a necessity that r and (r + 1)th column is partitioned $2 \leq r \leq n - 1$. Only such partitions will be known as a symmetric partitions.

It is important to note that all symmetric partitions are partitions but every partition in general need not be a symmetric partition.

We just illustrate this by an example

***Example 2.2.22:*** Let A be a 6 × 6 matrix

$$A = \begin{bmatrix} 1 & 2 & 3 & 4 & 5 & 6 \\ 7 & 8 & 9 & 0 & 1 & 2 \\ 1 & 1 & 0 & 1 & 0 & 2 \\ 3 & 3 & 1 & 0 & 3 & 2 \\ 5 & 7 & 9 & 1 & 4 & 7 \\ 6 & 8 & 0 & 2 & 5 & 8 \end{bmatrix}.$$

Let us partition A between the 3rd and 4th row and between the 2nd and 3rd column. Let the resulting supermatrix be denoted by $A_s$.

$$A_s = \left[\begin{array}{cc|cccc} 1 & 2 & 3 & 4 & 5 & 6 \\ 7 & 8 & 9 & 0 & 1 & 2 \\ 1 & 1 & 0 & 1 & 0 & 2 \\ 3 & 3 & 1 & 0 & 3 & 2 \\ \hline 5 & 7 & 9 & 1 & 4 & 7 \\ 6 & 8 & 0 & 2 & 5 & 8 \end{array}\right]$$

Clearly the partition is not symmetric by the very definition of the partition to be a symmetric partition.

Thus we have the following nice theorem.

**THEOREM 2.2.6:** *Let A be a n × n symmetric matrix. If P is any symmetric partition of A then the supermatrix $A_s$ is a symmetric supermatrix.*

*Proof:* Given A is a symmetric matrix. So if A = ($a_{ij}$) then $a_{ij}$ = $a_{ji}$, $1 \leq i \leq n$ and $1 \leq j \leq n$; $i \neq j$.

$$A = \begin{bmatrix} a_{11} & a_{12} & a_{13} & \cdots & a_{1n} \\ a_{21} & a_{22} & a_{23} & \cdots & a_{2n} \\ \vdots & \vdots & \vdots & & \vdots \\ a_{n1} & a_{n2} & a_{n3} & \cdots & a_{nn} \end{bmatrix}$$

Any symmetric partition of A given by $A_s$ where

$$A_s = \begin{bmatrix} A_{11} & A_{12} & \cdots & A_{1t} \\ A_{21} & A_{22} & \cdots & A_{2t} \\ \vdots & \vdots & & \vdots \\ A_{t1} & A_{t2} & \cdots & A_{tt} \end{bmatrix}.$$

Here $A_{11}$, $A_{12}$, …, $A_{tt}$ are submatrices of $A_s$. Clearly $A_{ii}$ are symmetric square matrices on the diagonal of the supermatrix $A_s$; i = 1, 2, …, t.

Further $A_{ij}$ = ($A_{ji}$) $i \neq j$, $1 \leq i, j \leq t$. Thus $A_s$ is a super symmetric matrix as A is given to be a symmetric matrix. Hence the claim.

All these concepts defined in the case of matrices can be carried out in case of fuzzy matrices as we have not used any form of multiplication or addition or any form of operations using two matrices. Thus the reader can easily define the notion of symmetric fuzzy matrix pseudo symmetric fuzzy matrix as well as the notion of super symmetric fuzzy matrix and the pseudo symmetric super fuzzy matrix.

## 2.3 Special Operations on Fuzzy Super Special Row and Column Matrix

Now we proceed onto define the minor product of a special type of fuzzy super special row matrix and its transpose. To this end we have to define the notion of fuzzy super special row matrix and fuzzy super special column matrix.

**DEFINITION 2.3.1:** *A fuzzy matrix of the form*

$$X = [X_1 \mid X_2 \ldots \mid X_n]$$

*will be called a fuzzy super special row matrix if $X_1, X_2, \ldots, X_n$ are submatrices having only the same number of rows but may have different number of columns.*

We illustrate this by the following example.

***Example 2.3.1:*** Let X = $[X_1 \mid X_2 \mid X_3 \mid X_4] = [X_1\ X_2\ X_3\ X_4]$ be a fuzzy super row special matrix where

$$X_1 = \begin{bmatrix} 0.2 & 1 & 0 & 1 & 0.7 \\ 0.5 & 0 & 1 & 0 & 0.5 \end{bmatrix},$$

$$X_2 = \begin{bmatrix} 1 & 0 \\ 0 & 0.5 \end{bmatrix},\ X_3 = \begin{bmatrix} 0.1 & 0.7 & 0.3 \\ 0 & 1 & 0.8 \end{bmatrix}$$

and

$$X_4 = \begin{bmatrix} 1 & 0 & 0.2 & 0.6 & 1 & 0.2 \\ 1 & 0 & 0.5 & 0.7 & 1 & 0.3 \end{bmatrix}.$$

***Example 2.3.2:*** Let X = $[X_1 \mid X_2 \mid X_3]$ where

$$X_1 = \begin{bmatrix} 0 & 1 \\ 0.5 & 0.2 \\ 0.3 & 0.7 \end{bmatrix},\ X_2 = \begin{bmatrix} 0 & 1 \\ 0.5 & 0.2 \\ 0.3 & 0.7 \end{bmatrix}$$

and

$$X_3 = \begin{bmatrix} 0 & 1 & 0 & 0.8 & 1 \\ 0.3 & 0 & 0.9 & 0.2 & 0.4 \\ 0.6 & 1 & 0.5 & 0.6 & 0.8 \end{bmatrix}.$$

X is a fuzzy super special row matrix.

***Example 2.3.3:*** $X = [X_1 \mid X_2]$ where

$$X_1 = \begin{bmatrix} 0.3 & 0.2 \\ 0.4 & 1 \\ 0.2 & 0 \\ 0.7 & 0.8 \end{bmatrix}$$

and

$$X_2 = \begin{bmatrix} 1 \\ 0.8 \\ 0.1 \\ 0.9 \end{bmatrix}.$$

X is a fuzzy super special row matrix.

Now we proceed onto define the notion of fuzzy super special column matrix.

**DEFINITION 2.3.2:** *Let*

$$Y = \begin{bmatrix} Y_1 \\ \hline Y_2 \\ \hline \vdots \\ \hline Y_m \end{bmatrix}$$

*be a fuzzy supermatrix where $Y_1$, $Y_2$, …, $Y_m$ are fuzzy submatrices of Y and each of the $Y_i$'s have the same number of columns but with many different rows. We call Y a fuzzy special super column matrix.*

We illustrate this by the following example.

***Example 2.3.4:*** Let

$$Y = \left[\begin{array}{cc} 0.2 & 0.7 \\ 0.1 & 0.5 \\ \hline 0.6 & 0.2 \\ 1 & 1 \\ 0.1 & 0.5 \\ 0.9 & 0.3 \\ 0.1 & 1 \\ \hline 0 & 0 \\ 0.7 & 1 \end{array}\right]$$

$$= \left[\begin{array}{c} Y_1 \\ \hline Y_2 \\ \hline Y_3 \end{array}\right]$$

is a super fuzzy special column matrix. Y has only two columns but $Y_1$ is a fuzzy matrix with 2 columns and 2 rows. $Y_2$ is a fuzzy matrix with 2 columns but with five rows and so on.

***Example 2.3.5:*** Let

$$Y = \left[\begin{array}{cccc} 0.3 & 0.2 & 0.1 & 0.5 \\ 1 & 1 & 1 & 1 \\ \hline 0.2 & 0.2 & 0.3 & 0.7 \\ 0 & 0 & 0 & 0 \\ 1 & 0 & 0.7 & 0.2 \\ \hline 0.7 & 0.6 & 0.8 & 1 \\ 0.1 & 1 & 0.9 & 1 \\ 1 & 0 & 1 & 1 \\ 0 & 0.5 & 0.2 & 1 \end{array}\right]$$

where Y is a fuzzy super special column matrix. This has only four columns and 9 rows contributed by the three fuzzy submatrices $Y_1$, $Y_2$ and $Y_3$ where

$$Y = \begin{bmatrix} Y_1 \\ Y_2 \\ Y_3 \end{bmatrix}.$$

***Example 2.3.6:*** Let

$$Y = \left[\begin{array}{ccccc} 0.3 & 0.2 & 0.1 & 1 & 0 \\ 0.6 & 0.7 & 0.8 & 0.9 & 1 \\ \hline 1 & 0.2 & 1 & 0 & 1 \end{array}\right]$$

$$= \begin{bmatrix} Y_1 \\ Y_2 \end{bmatrix}$$

where

$$Y_1 = \begin{bmatrix} 0.3 & 0.2 & 0.1 & 1 & 0 \\ 0.6 & 0.7 & 0.8 & 0.9 & 1 \end{bmatrix}$$

and

$$Y_2 = [1\ 0.2\ 1\ 0\ 1].$$

Y has five columns but only 3 rows. Y is a fuzzy super special column matrix.

Now we see if X is a fuzzy super special row matrix with say n rows and

$$X_s = [X_1\ X_2\ \ldots\ X_t]$$

where $X_1, X_2, \ldots, X_t$ are fuzzy submatrices each of them having only n rows but varying number of columns.

Now the transpose of X denoted by

$$X^t = [X_1\ X_2\ \ldots\ X_t]^t = \begin{bmatrix} X_1^t \\ X_2^t \\ X_3^t \\ \vdots \\ X_t^t \end{bmatrix}.$$

Now $X^t$ is a fuzzy super special column matrix with n columns but different number of rows.

Thus we see the transpose of a fuzzy super special row matrix is a fuzzy super special column matrix and vice versa.

Now we illustrate this by the following example

***Example 2.3.7:*** Let

$$X = \left[\begin{array}{cc|cccc|ccccc} 0.3 & 1 & 0.4 & 1 & 0.3 & 0.2 & 0.1 & 0 & 0 & 1 & 0.1 \\ 0.1 & 1 & 0.6 & 0.7 & 0.8 & 1 & 0 & 1 & 0 & 0.3 & 1 \\ 0.2 & 0 & 1 & 0.5 & 1 & 0.9 & 1 & 0.2 & 0 & 1 & 0 \end{array}\right]$$

be a fuzzy super special row matrix.

Now the transpose of X is

$$X^t = \left[\begin{array}{ccc} 0.3 & 0.1 & 0.2 \\ 1 & 1 & 0 \\ \hline 0.4 & 0.6 & 1 \\ 1 & 0.7 & 0.5 \\ 0.3 & 0.8 & 1 \\ 0.2 & 1 & 0.9 \\ \hline 0.1 & 0 & 1 \\ 0 & 1 & 0.2 \\ 0 & 0 & 0 \\ 1 & 0.3 & 1 \\ 0.1 & 1 & 0 \end{array}\right]$$

is a fuzzy super special column matrix.

***Example 2.3.8:*** Let

$$Y = \left[\begin{array}{ccc} 0.3 & 0.1 & 0.2 \\ 0.3 & 0.5 & 1 \\ \hline 1 & 0 & 1 \\ 0 & 1 & 0 \\ 0.3 & 0.8 & 0.9 \\ 0.1 & 0 & 1 \\ \hline 1 & 1 & 1 \\ 0 & 0 & 0 \\ 0.3 & 0.2 & 0.5 \end{array}\right]$$

be a fuzzy super special column matrix. Let $Y^t$ be the transpose of Y

$$Y^t = \left[\begin{array}{cc|cccc|ccc} 0.3 & 0.3 & 1 & 0 & 0.3 & 0.1 & 1 & 0 & 0.3 \\ 0.1 & 0.5 & 0 & 1 & 0.8 & 0 & 1 & 0 & 0.2 \\ 0.2 & 1 & 1 & 0 & 0.9 & 1 & 1 & 0 & 0.5 \end{array}\right].$$

Now we define a special product called the minor product moment of a fuzzy super special column matrix and fuzzy super special row matrix.

***Example 2.3.9:*** Let

$$X = \left[\begin{array}{cc|c|ccc} 0.2 & 0.3 & 0.4 & 0.3 & 0.2 & 0.3 \\ 0.1 & 0.4 & 1 & 1 & 0.1 & 0.2 \\ 0.2 & 0.1 & 0.2 & 0.4 & 0.3 & 0.2 \end{array}\right]$$

be a fuzzy super special row matrix. Now

$$X^t = \left[\begin{array}{ccc} 0.2 & 0.1 & 0.2 \\ 0.3 & 0.4 & 0.1 \\ \hline 0.4 & 1 & 0.2 \\ \hline 0.3 & 1 & 0.4 \\ 0.2 & 0.1 & 0.3 \\ 0.3 & 0.2 & 0.2 \end{array}\right]$$

is its transpose a fuzzy super special column matrix.

Now

$$X.X^t = \left[\begin{array}{cc|c|ccc} 0.2 & 0.3 & 0.4 & 0.3 & 0.2 & 0.3 \\ 0.1 & 0.4 & 1 & 1 & 0.1 & 0.2 \\ 0.2 & 0.1 & 0.2 & 0.4 & 0.3 & 0.2 \end{array}\right] . \left[\begin{array}{ccc} 0.2 & 0.1 & 0.2 \\ 0.3 & 0.4 & 0.1 \\ \hline 0.4 & 1 & 0.2 \\ \hline 0.3 & 1 & 0.4 \\ 0.2 & 0.1 & 0.3 \\ 0.3 & 0.2 & 0.2 \end{array}\right]$$

$$= \max \left[ \max\min \left\{ \begin{bmatrix} 0.2 & 0.3 \\ 0.1 & 0.4 \\ 0.2 & 0.1 \end{bmatrix} \begin{bmatrix} 0.2 & 0.1 & 0.2 \\ 0.3 & 0.4 & 0.1 \end{bmatrix} \right. , \right.$$

$$\max\min \begin{bmatrix} 0.4 \\ 1 \\ 0.2 \end{bmatrix} \begin{bmatrix} 0.4 & 1 & 0.2 \end{bmatrix} ,$$

$$\left. \left. \max\min \begin{bmatrix} 0.3 & 0.2 & 0.3 \\ 1 & 0.1 & 0.2 \\ 0.4 & 0.3 & 0.2 \end{bmatrix} \begin{bmatrix} 0.3 & 1 & 0.4 \\ 0.2 & 0.1 & 0.3 \\ 0.3 & 0.2 & 0.2 \end{bmatrix} \right\} \right]$$

= max {max min ($a_{ij}$, $b_{jk}$), max min ($c_{it}$ $d_{tk}$), max min ($p_{ij}$ $q_{jk}$)} (defined earlier) (pages 88-9 of this book).

$$= \max\left\{\begin{bmatrix} 0.3 & 0.3 & 0.2 \\ 0.3 & 0.4 & 0.1 \\ 0.2 & 0.1 & 0.2 \end{bmatrix}, \begin{bmatrix} 0.4 & 0.4 & 0.2 \\ 0.4 & 1 & 0.2 \\ 0.2 & 0.2 & 0.2 \end{bmatrix}, \begin{bmatrix} 0.3 & 0.3 & 0.3 \\ 0.3 & 1 & 0.4 \\ 0.3 & 0.4 & 0.4 \end{bmatrix}\right\}$$

$$= \begin{bmatrix} 0.4 & 0.4 & 0.3 \\ 0.4 & 1 & 0.4 \\ 0.3 & 0.4 & 0.4 \end{bmatrix}.$$

$XX^t$ is defined as the minor product moment of X.

***Example 2.3.10:*** Let

$$Y = \left[\begin{array}{cccc} 0.3 & 0.1 & 0.4 & 0.5 \\ \hline 0.1 & 1 & 1 & 0 \\ 1 & 0 & 0.7 & 0.2 \\ \hline 0.3 & 1 & 0.1 & 0.5 \\ 0.2 & 0 & 0.2 & 0.6 \\ 0.1 & 0.8 & 1 & 0.7 \\ 0 & 0.2 & 1 & 0.2 \end{array}\right]$$

be a fuzzy super special column matrix. $Y^t$ is the fuzzy super special row matrix.

$$\max\{\max\min\{Y^t, Y\}\}$$

$$= \max\{\max\min\left\{\left[\begin{array}{c|cc|cccc} 0.3 & 0.1 & 1 & 0.3 & 0.2 & 0.1 & 0 \\ 0.1 & 1 & 0 & 1 & 0 & 0.8 & 0.2 \\ 0.4 & 1 & 0.7 & 0.1 & 0.2 & 1 & 1 \\ 0.5 & 0 & 0.2 & 0.5 & 0.6 & 0.7 & 0.2 \end{array}\right] \times\right.$$

$$\left.\begin{bmatrix} 0.3 & 0.1 & 0.4 & 0.5 \\ \hline 0.1 & 1 & 1 & 0 \\ 1 & 0 & 0.7 & 0.2 \\ \hline 0.3 & 1 & 0.1 & 0.5 \\ 0.2 & 0 & 0.2 & 0.6 \\ 0.1 & 0.8 & 1 & 0.7 \\ 0 & 0.2 & 1 & 0.2 \end{bmatrix}\right\}$$

$$= \max\{\max\min\left\{\begin{bmatrix} 0.3 \\ 0.1 \\ 0.4 \\ 0.5 \end{bmatrix}, \begin{bmatrix} 0.3 & 0.1 & 0.4 & 0.5 \end{bmatrix}\right\},$$

$$\max\min\left\{\begin{bmatrix} 0.1 & 1 \\ 1 & 0 \\ 1 & 0.7 \\ 0 & 0.2 \end{bmatrix}, \begin{bmatrix} 0.1 & 1 & 1 & 0 \\ 1 & 0 & 0.7 & 0.2 \end{bmatrix}\right\},$$

$$\max\min\left\{\begin{bmatrix} 0.3 & 0.2 & 0.1 & 0 \\ 1 & 0 & 0.8 & 0.2 \\ 0.1 & 0.2 & 1 & 1 \\ 0.5 & 0.6 & 0.7 & 0.2 \end{bmatrix}, \begin{bmatrix} 0.3 & 1 & 0.1 & 0.5 \\ 0.2 & 0 & 0.2 & 0.6 \\ 0.1 & 0.8 & 1 & 0.7 \\ 0 & 0.2 & 1 & 0.2 \end{bmatrix}\right\},$$

$$= \max$$

$$\left\{\begin{bmatrix} 0.3 & 0.1 & 0.3 & 0.3 \\ 0.1 & 0.1 & 0.1 & 0.1 \\ 0.3 & 0.1 & 0.4 & 0.4 \\ 0.3 & 0.1 & 0.4 & 0.5 \end{bmatrix}, \begin{bmatrix} 1 & 0.1 & 0.7 & 0.2 \\ 0.1 & 1 & 1 & 0 \\ 0.7 & 1 & 1 & 0.2 \\ 0.2 & 0 & 0.2 & 0.2 \end{bmatrix}, \begin{bmatrix} 0.3 & 0.3 & 0.2 & 0.3 \\ 0.3 & 1 & 0.8 & 0.7 \\ 0.2 & 0.8 & 1 & 0.7 \\ 0.3 & 0.7 & 0.7 & 0.7 \end{bmatrix}\right\}$$

$$= \begin{bmatrix} 1 & 0.3 & 0.7 & 0.3 \\ 0.3 & 1 & 1 & 0.7 \\ 0.7 & 1 & 1 & 0.7 \\ 0.3 & 0.7 & 0.7 & 0.7 \end{bmatrix},$$

is the minor product moment of the transpose super fuzzy special column matrix Y, with the super fuzzy special column matrix. It is very important to note that neither Y nor $Y^t$ is a symmetric super fuzzy special column matrix or row matrix but $Y^t$ Y is a symmetric fuzzy matrix. It is nice to note the resultant of the product is not super symmetric fuzzy matrix only a symmetric fuzzy matrix. Likewise we see the minor product moment of a super fuzzy special row vector X with its transpose is a symmetric fuzzy matrix which is not a supermatrix. Thus we see if X is a super fuzzy special row matrix with n rows then $XX^t$ is just a $n \times n$ symmetric fuzzy matrix. Further if Y is a super fuzzy special column matrix with m rows then $Y^tY$ is a $m \times m$ symmetric fuzzy matrix. This special product yields more and more symmetric fuzzy matrices of desired order.

Now suppose we have super fuzzy special column matrix Y with m columns and a super fuzzy special row matrix X with m-rows then we can define their product YX. Let

$$Y = \begin{bmatrix} Y_1 \\ Y_2 \\ Y_3 \\ Y_t \end{bmatrix}$$

where Y is a super fuzzy special column matrix with m columns and the total number of rows in Y will be number of rows in Y i.e., it will be the number of rows in $Y_1$ + number of rows in $Y_2$ + … + number of rows in $Y_t$. Suppose the total number of rows in Y before partition be m. Suppose $X = [X_1\ X_2\ \ldots\ X_r]$ where X is a super fuzzy special row matrix with m rows then the number of columns in X before partitioning will be the number of columns in $X_1$ + number of columns in $X_2$ + … + number of columns in $X_r$. Suppose the number of columns in X be n, then the product YX is a super fuzzy matrix with m rows and n columns with tr number of fuzzy submatrices. First we illustrate this by an example.

***Example 2.3.11:*** Let

$$Y = \left[\begin{array}{ccc} 0.3 & 1 & 0 \\ 0.2 & 0.4 & 1 \\ \hline 1 & 0 & 0 \\ 0 & 1 & 0 \\ 0.2 & 0.5 & 1 \\ \hline 0.7 & 0 & 0.8 \\ 0.9 & 1 & 0.4 \end{array}\right]$$

be a fuzzy super special column matrix which has 3 fuzzy sub matrices and 7 rows before partitioning. Further Y has only 3 columns. Let

$$X = \left[\begin{array}{c|ccc} 0.2 & 0.6 & 1 & 0.5 \\ 1 & 0.7 & 0 & 0.6 \\ 0 & 1 & 0.2 & 0.8 \end{array}\right]$$

be a super fuzzy special row matrix with two fuzzy submatrices with 3 rows and X having 4 columns before partitioning. Major special product of Y with X given by

$$YX = \left[\begin{array}{ccc} 0.3 & 1 & 0 \\ 0.2 & 0.4 & 1 \\ \hline 1 & 0 & 0 \\ 0 & 1 & 0 \\ 0.2 & 0.5 & 1 \\ \hline 0.7 & 0 & 0.8 \\ 0.9 & 1 & 0.4 \end{array}\right] \left[\begin{array}{c|ccc} 0.2 & 0.6 & 1 & 0.5 \\ 1 & 0.7 & 0 & 0.6 \\ 0 & 1 & 0.2 & 0.8 \end{array}\right] =$$

$$
= \left[\begin{array}{c|c}
\begin{bmatrix} 0.3 & 1 & 0 \\ 0.2 & 0.4 & 1 \end{bmatrix}\begin{bmatrix} 0.2 \\ 1 \\ 0 \end{bmatrix} & \begin{bmatrix} 0.3 & 1 & 0 \\ 0.2 & 0.4 & 1 \end{bmatrix}\begin{bmatrix} 0.6 & 1 & 0.5 \\ 0.7 & 0 & 0.6 \\ 1 & 0.2 & 0.8 \end{bmatrix} \\
\hline
\begin{bmatrix} 1 & 0 & 0 \\ 0 & 1 & 0 \\ 0.2 & 0.5 & 1 \end{bmatrix}\begin{bmatrix} 0.2 \\ 1 \\ 0 \end{bmatrix} & \begin{bmatrix} 1 & 0 & 0 \\ 0 & 1 & 0 \\ 0.2 & 0.5 & 1 \end{bmatrix}\begin{bmatrix} 0.6 & 1 & 0.5 \\ 0.7 & 0 & 0.6 \\ 1 & 0.2 & 0.8 \end{bmatrix} \\
\hline
\begin{bmatrix} 0.7 & 0 & 0.8 \\ 0.9 & 1 & 0.4 \end{bmatrix}\begin{bmatrix} 0.2 \\ 1 \\ 0 \end{bmatrix} & \begin{bmatrix} 0.7 & 0 & 0.8 \\ 0.9 & 1 & 0.4 \end{bmatrix}\begin{bmatrix} 0.6 & 1 & 0.5 \\ 0.7 & 0 & 0.6 \\ 1 & 0.2 & 0.8 \end{bmatrix}
\end{array}\right]
$$

$$
= \left[\begin{array}{c|ccc}
1 & 0.7 & 0.3 & 0.6 \\
0.4 & 1 & 0.2 & 0.8 \\
\hline
0.2 & 0.6 & 1 & 0.5 \\
1 & 0.7 & 0 & 0.6 \\
0.5 & 1 & 0.2 & 0.8 \\
\hline
0.2 & 0.8 & 0.7 & 0.8 \\
1 & 0.7 & 0.9 & 0.6
\end{array}\right].
$$

We see YX is a fuzzy super matrix having six fuzzy sub matrices and its order is 7 × 4 before partitioning. Thus using major special product we can multiply fuzzy super special column matrix X with the fuzzy super special row matrix provided the number columns in X is equal to the number of rows in Y. Now we can also multiply $X^t$ with $Y^t$. We take the same example.

***Example 2.3.12:*** Taking X and Y from above example we get

$$
X^t = \left[\begin{array}{ccc}
0.2 & 1 & 0 \\
\hline
0.6 & 0.7 & 1 \\
1 & 0 & 0.2 \\
0.5 & 0.6 & 0.8
\end{array}\right]
$$

and

$$Y^t = \left[\begin{array}{cc|ccc|cc} 0.3 & 0.2 & 1 & 0 & 0.2 & 0.7 & 0.9 \\ 1 & 0.4 & 0 & 1 & 0.5 & 0 & 1 \\ 0 & 1 & 0 & 0 & 1 & 0.8 & 0.4 \end{array}\right].$$

$$X^tY^t = \left[\begin{array}{ccc} 0.2 & 1 & 0 \\ \hline 0.6 & 0.7 & 1 \\ 1 & 0 & 0.2 \\ 0.5 & 0.6 & 0.8 \end{array}\right] \left[\begin{array}{cc|ccc|cc} 0.3 & 0.2 & 1 & 0 & 0.2 & 0.7 & 0.9 \\ 1 & 0.4 & 0 & 1 & 0.5 & 0 & 1 \\ 0 & 1 & 0 & 0 & 1 & 0.8 & 0.4 \end{array}\right]$$

$$\left[\begin{array}{c|c|c} [0.2\ 1\ 0]\begin{bmatrix} 0.3 & 0.2 \\ 1 & 0.4 \\ 0 & 1 \end{bmatrix} & [0.2\ 1\ 0]\begin{bmatrix} 1 & 0 & 0.2 \\ 0 & 1 & 0.5 \\ 0 & 0 & 1 \end{bmatrix} & [0.2\ 1\ 0]\begin{bmatrix} 0.7 & 0.9 \\ 0 & 1 \\ 0.8 & 0.4 \end{bmatrix} \\ \hline \begin{bmatrix} 0.6 & 0.7 & 1 \\ 1 & 0 & 0.2 \\ 0.5 & 0.6 & 0.8 \end{bmatrix}\begin{bmatrix} 0.3 & 0.2 \\ 1 & 0.4 \\ 0 & 1 \end{bmatrix} & \begin{bmatrix} 0.6 & 0.7 & 1 \\ 1 & 0 & 0.2 \\ 0.5 & 0.6 & 0.8 \end{bmatrix}\begin{bmatrix} 1 & 0 & 0.2 \\ 0 & 1 & 0.5 \\ 0 & 0 & 1 \end{bmatrix} & \begin{bmatrix} 0.6 & 0.7 & 1 \\ 1 & 0 & 0.2 \\ 0.5 & 0.6 & 0.8 \end{bmatrix}\begin{bmatrix} 0.7 & 0.9 \\ 0 & 1 \\ 0.8 & 0.4 \end{bmatrix} \end{array}\right]$$

$$= \left[\begin{array}{cc|ccc|cc} 1 & 0.4 & 0.2 & 1 & 0.5 & 0.2 & 1 \\ \hline 0.7 & 1 & 0.6 & 0.7 & 1 & 0.8 & 0.7 \\ 0.3 & 0.2 & 1 & 0 & 0.2 & 0.7 & 0.9 \\ 0.6 & 0.8 & 0.5 & 0.6 & 0.8 & 0.8 & 0.6 \end{array}\right]$$

$$= (YX)^t.$$

[We have used max min operation in the above product].

Next we see when the product of fuzzy supermatrices are compatible.

Before we make a formal definition we give an illustrative example.

***Example 2.3.13:*** Let A and B be any two fuzzy super matrices where

$$A = \left[\begin{array}{c|cc|ccc} 1 & 0.2 & 1 & 0.3 & 0.2 & 1 \\ 0.2 & 0.3 & 1 & 0.2 & 1 & 0.2 \\ \hline 1 & 0.4 & 0.2 & 0.3 & 0.2 & 0.2 \\ 0.4 & 1 & 0.3 & 0.2 & 1 & 1 \\ 0.2 & 0.3 & 0.2 & 0.3 & 0.2 & 0.3 \\ 0.3 & 0.4 & 1 & 1 & 0.4 & 0.2 \\ \hline 0.2 & 1 & 0.2 & 0.2 & 1 & 0.3 \end{array}\right]$$

and

$$B = \left[\begin{array}{cccc|cc|cccc} 1 & 0.3 & 0.2 & 1 & 0.3 & 0.8 & 0.1 & 0.2 & 0.5 & 1 \\ \hline 0.3 & 0.4 & 1 & 0.2 & 1 & 0 & 1 & 0 & 0.2 & 0 \\ 1 & 1 & 0 & 0 & 0.4 & 0.5 & 0 & 1 & 0.4 & 0.7 \\ \hline 1 & 0.2 & 1 & 0.7 & 0.7 & 0.2 & 0.5 & 1 & 0.2 & 0.3 \\ 0 & 0.1 & 0.5 & 1 & 0.3 & 0.2 & 0.3 & 0.3 & 0.5 & 0.4 \\ 1 & 0.4 & 0.2 & 0.5 & 0.5 & 0.4 & 1 & 0 & 1 & 0 \end{array}\right].$$

Now we define the minor product of these fuzzy super matrices.

$$\text{Max } \{\text{max min A.B}\} =$$

$$\text{max } \{\text{max min}$$

$$\left\{ \left[\begin{array}{c} 1 \\ 0.2 \\ \hline 1 \\ 0.4 \\ 0.2 \\ 0.3 \\ \hline 0.2 \end{array}\right] \left[\begin{array}{cccc|cc|cccc} 1 & 0.3 & 0.2 & 1 & 0.3 & 0.8 & 0.1 & 0.2 & 0.5 & 1 \end{array}\right] \right\},$$

$$\left\{\left[\begin{array}{cc} 0.2 & 1 \\ 0.3 & 1 \\ \hline 0.4 & 0.2 \\ 1 & 0.3 \\ 0.3 & 0.2 \\ 0.4 & 1 \\ \hline 1 & 0.2 \end{array}\right]\left[\begin{array}{cccc|cc|cccc} 0.3 & 0.4 & 1 & 0.2 & 1 & 0 & 1 & 0 & 0.2 & 0 \\ 1 & 1 & 0 & 0 & 0.4 & 0.5 & 0 & 1 & 0.4 & 0.7 \end{array}\right]\right\},$$

$$\left\{\left[\begin{array}{ccc} 0.3 & 0.2 & 1 \\ 0.2 & 1 & 0.2 \\ \hline 0.3 & 0.2 & 0.2 \\ 0.2 & 1 & 1 \\ 0.3 & 0.2 & 0.3 \\ 1 & 0.4 & 0.2 \\ \hline 0.2 & 1 & 0.3 \end{array}\right]\left[\begin{array}{cccc|cc|cccc} 1 & 0.2 & 1 & 0.7 & 0.7 & 0.2 & 0.5 & 1 & 0.2 & 0.3 \\ 0 & 0.1 & 0.5 & 1 & 1 & 0.3 & 0.2 & 0.3 & 0.5 & 0.4 \\ 1 & 0.4 & 0.2 & 0.5 & 0.5 & 0.4 & 1 & 0 & 1 & 0 \end{array}\right]\right\} =$$

$$\left\{\left[\begin{array}{c|c|c} \begin{bmatrix} 1 \\ 0.2 \end{bmatrix}\begin{bmatrix} 1 & 0.3 & 0.2 & 1 \end{bmatrix} & \begin{bmatrix} 1 \\ 0.2 \end{bmatrix}\begin{bmatrix} 0.3 & 0.8 \end{bmatrix} & \begin{bmatrix} 1 \\ 0.2 \end{bmatrix}\begin{bmatrix} 0.1 & 0.2 & 0.5 & 1 \end{bmatrix} \\ \hline \begin{bmatrix} 1 \\ 0.4 \\ 0.2 \\ 0.3 \end{bmatrix}\begin{bmatrix} 1 & 0.3 & 0.2 & 1 \end{bmatrix} & \begin{bmatrix} 1 \\ 0.4 \\ 0.2 \\ 0.3 \end{bmatrix}\begin{bmatrix} 0.3 & 0.8 \end{bmatrix} & \begin{bmatrix} 1 \\ 0.4 \\ 0.2 \\ 0.3 \end{bmatrix}\begin{bmatrix} 0.1 & 0.2 & 0.5 & 1 \end{bmatrix} \\ \hline \begin{bmatrix} 0.2 \end{bmatrix}\begin{bmatrix} 1 & 0.3 & 0.2 & 1 \end{bmatrix} & \begin{bmatrix} 0.2 \end{bmatrix}\begin{bmatrix} 0.3 & 0.8 \end{bmatrix} & \begin{bmatrix} 0.2 \end{bmatrix}\begin{bmatrix} 0.1 & 0.2 & 0.5 & 1 \end{bmatrix} \end{array}\right],\right.$$

$$
\left[\begin{array}{c|c|c}
\begin{bmatrix}.2 & 1\\ .3 & 1\end{bmatrix}\begin{bmatrix}.3 & .4 & 1 & .2\\ 1 & 1 & 0 & 0\end{bmatrix} & \begin{bmatrix}.2 & 1\\ .3 & 1\end{bmatrix}\begin{bmatrix}1 & 0\\ .4 & .5\end{bmatrix} & \begin{bmatrix}.2 & 1\\ .3 & 1\end{bmatrix}\begin{bmatrix}1 & 0 & .2 & 0\\ 0 & 1 & .4 & .7\end{bmatrix}\\ \hline
\begin{bmatrix}.4 & .2\\ 1 & .3\\ .3 & .2\\ .4 & 1\end{bmatrix}\begin{bmatrix}.3 & .4 & 1 & .2\\ 1 & 1 & 0 & 0\end{bmatrix} & \begin{bmatrix}.4 & .2\\ 1 & .3\\ .3 & .2\\ .4 & 1\end{bmatrix}\begin{bmatrix}1 & 0\\ .4 & .5\end{bmatrix} & \begin{bmatrix}.4 & .2\\ 1 & .3\\ .3 & .2\\ .4 & 1\end{bmatrix}\begin{bmatrix}1 & 0 & .2 & 0\\ 0 & 1 & .4 & .7\end{bmatrix}\\ \hline
\begin{bmatrix}1 & .2\end{bmatrix}\begin{bmatrix}.3 & .4 & 1 & .2\\ 1 & 1 & 0 & 0\end{bmatrix} & \begin{bmatrix}1 & .2\end{bmatrix}\begin{bmatrix}1 & 0\\ .4 & .5\end{bmatrix} & \begin{bmatrix}1 & .2\end{bmatrix}\begin{bmatrix}1 & 0 & .2 & 0\\ 0 & 1 & .4 & .7\end{bmatrix}
\end{array}\right],
$$

$$
\left[\begin{array}{c|c|c}
\begin{bmatrix}.3 & .2 & 1\\ .2 & 1 & .2\end{bmatrix}\begin{bmatrix}1 & .2 & 1 & .7\\ 0 & .1 & .5 & 1\\ 1 & .4 & .2 & .5\end{bmatrix} & \begin{bmatrix}.3 & .2 & 1\\ .2 & 1 & .2\end{bmatrix}\begin{bmatrix}.7 & .2\\ 1 & .3\\ .5 & .4\end{bmatrix} & \begin{bmatrix}.3 & .2 & 1\\ .2 & 1 & .2\end{bmatrix}\begin{bmatrix}.5 & 1 & .2 & .3\\ .2 & .3 & .5 & .4\\ 1 & 0 & 1 & 0\end{bmatrix}\\ \hline
\begin{bmatrix}.3 & .2 & .2\\ .2 & 1 & 1\\ .3 & .2 & .3\\ 1 & .4 & .2\end{bmatrix}\begin{bmatrix}1 & .2 & 1 & .7\\ 0 & .1 & .5 & 1\\ 1 & .4 & .2 & .5\end{bmatrix} & \begin{bmatrix}.3 & .2 & .2\\ .2 & 1 & 1\\ .3 & .2 & .3\\ 1 & .4 & .2\end{bmatrix}\begin{bmatrix}.7 & .2\\ 1 & .3\\ .5 & .4\end{bmatrix} & \begin{bmatrix}.3 & .2 & .2\\ .2 & 1 & 1\\ .3 & .2 & .3\\ 1 & .4 & .2\end{bmatrix}\begin{bmatrix}.5 & 1 & .2 & .3\\ .2 & .3 & .5 & .4\\ 1 & 0 & 1 & 0\end{bmatrix}\\ \hline
\begin{bmatrix}.2 & 1 & .3\end{bmatrix}\begin{bmatrix}1 & .2 & 1 & .7\\ 0 & .1 & .5 & 1\\ 1 & .4 & .2 & .5\end{bmatrix} & \begin{bmatrix}.2 & 1 & .3\end{bmatrix}\begin{bmatrix}.7 & .2\\ 1 & .3\\ .5 & .4\end{bmatrix} & \begin{bmatrix}.2 & 1 & .3\end{bmatrix}\begin{bmatrix}.5 & 1 & .2 & .3\\ .2 & .3 & .5 & .4\\ 1 & 0 & 1 & 0\end{bmatrix}
\end{array}\right]\Bigg\}
$$

$$
= \max\left\{\left[\begin{array}{cccc|cc|cccc}
1 & 0.3 & 0.2 & 1 & 0.3 & 0.8 & 0.1 & 0.2 & 0.5 & 1\\
0.2 & 0.2 & 0.2 & 0.2 & 0.2 & 0.2 & 0.1 & 0.2 & 0.2 & 0.2\\ \hline
1 & 0.3 & 0.2 & 1 & 0.3 & 0.8 & 0.1 & 0.2 & 0.5 & 1\\
0.4 & 0.3 & 0.2 & 0.4 & 0.3 & 0.4 & 0.1 & 0.2 & 0.4 & 0.4\\
0.2 & 0.2 & 0.2 & 0.2 & 0.2 & 0.2 & 0.1 & 0.2 & 0.2 & 0.2\\
0.3 & 0.3 & 0.2 & 0.3 & 0.3 & 0.3 & 0.1 & 0.2 & 0.3 & 0.3\\ \hline
0.2 & 0.2 & 0.2 & 0.2 & 0.2 & 0.2 & 0.1 & 0.2 & 0.2 & 0.2
\end{array}\right]\right.,
$$

$$\left[\begin{array}{cccc|cc|cccc}
1 & 1 & 0.2 & 0.2 & 0.4 & 0.5 & 0.2 & 1 & 0.4 & 0.7 \\
1 & 1 & 0.3 & 0.2 & 0.4 & 0.5 & 0.3 & 1 & 0.4 & 0.7 \\ \hline
0.3 & 0.4 & 0.4 & 0.2 & 0.4 & 0.2 & 0.4 & 0.2 & 0.4 & 0.2 \\
0.3 & 0.4 & 1 & 0.2 & 1 & 0.3 & 1 & 0.3 & 0.3 & 0.3 \\
0.3 & 0.3 & 0.3 & 0.2 & 0.3 & 0.2 & 0.3 & 0.2 & 0.2 & 0.2 \\
1 & 1 & 0.4 & 0.2 & 0.4 & 0.5 & 0.4 & 1 & 0.4 & 0.7 \\ \hline
0.3 & 0.4 & 1 & 0.2 & 1 & 0.2 & 1 & 0.2 & 0.2 & 0.2
\end{array}\right],$$

$$\left.\left[\begin{array}{cccc|cc|cccc}
1 & 0.4 & 0.3 & 0.5 & 0.5 & 0.4 & 1 & 0.3 & 1 & 0.3 \\
0.2 & 0.2 & 0.5 & 1 & 1 & 0.3 & 0.2 & 0.3 & 0.5 & 0.4 \\ \hline
0.3 & 0.2 & 0.3 & 0.3 & 0.3 & 0.2 & 0.3 & 0.3 & 0.2 & 0.3 \\
1 & 0.4 & 0.5 & 1 & 1 & 0.4 & 1 & 0.3 & 1 & 0.4 \\
0.3 & 0.3 & 0.3 & 0.3 & 0.3 & 0.3 & 0.3 & 0.3 & 0.3 & 0.3 \\
1 & 0.2 & 1 & 0.7 & 0.4 & 0.3 & 0.5 & 1 & 0.3 & 0.4 \\ \hline
0.3 & 0.3 & 0.5 & 1 & 1 & 0.3 & 0.3 & 0.3 & 0.5 & 0.4
\end{array}\right]\right\}$$

$$= \left[\begin{array}{cccc|cc|cccc}
1 & 1 & 0.3 & 1 & 0.5 & 0.8 & 1 & 1 & 1 & 1 \\
1 & 1 & 0.5 & 1 & 1 & 0.5 & 0.3 & 1 & 0.5 & 0.7 \\ \hline
1 & 0.4 & 0.4 & 1 & 0.4 & 0.8 & 0.4 & 0.3 & 0.5 & 1 \\
1 & 0.4 & 1 & 1 & 1 & 0.4 & 1 & 0.3 & 1 & 0.4 \\
0.3 & 0.3 & 0.3 & 0.3 & 0.3 & 0.3 & 0.3 & 0.3 & 0.3 & 0.3 \\
1 & 1 & 1 & 0.7 & 0.4 & 0.5 & 0.5 & 1 & 0.5 & 0.7 \\ \hline
0.3 & 0.4 & 1 & 1 & 1 & 0.3 & 1 & 0.3 & 0.5 & 0.4
\end{array}\right].$$

Thus using these max min rules we operate with the super fuzzy matrices. All super matrix properties true for fuzzy super matrices can be carried out for the only simple difference between a super matrix and the fuzzy super matrix is that fuzzy super matrices take their entries from the unit interval [0,1] or at times from the set {–1,0,1}. We should keep in mind all fuzzy super matrices are super matrices only all super matrices in general are not fuzzy super matrices.

***Example 2.3.14:*** For instance if

$$S = \left[\begin{array}{cccc|ccc} 3 & 7 & 8 & 9 & 1 & 2 & 3 \\ 0 & 1 & 5 & 8 & 4 & 2 & 9 \\ \hline 9 & 0 & -4 & 2 & 1 & 0 & 9 \\ 7 & 4 & 0 & 1 & 2 & 7 & 4 \\ 3 & 7 & 8 & -4 & 19 & 1 & 14 \\ \hline 10 & 4 & 9 & 11 & 25 & 8 & 19 \end{array}\right]$$

is a super matrix but is not a fuzzy super matrix. We give yet another working with super fuzzy matrices.

***Example 2.3.15:*** Let A be a special row fuzzy matrix given by

$$A = \left[\begin{array}{cccc|cc|ccc|ccccc} 0 & 1 & 0 & 1 & 1 & 0 & 0 & 1 & 1 & 1 & 1 & 0 & 0 & 0 \\ 1 & 1 & 0 & 1 & 0 & 1 & 1 & 0 & 0 & 1 & 0 & 1 & 0 & 0 \\ 0 & 0 & 0 & 1 & 1 & 1 & 0 & 1 & 0 & 0 & 0 & 0 & 0 & 1 \\ 1 & 1 & 0 & 0 & 0 & 1 & 0 & 0 & 1 & 1 & 0 & 0 & 0 & 0 \end{array}\right].$$

Let

$$X = [1\ 0\ 1\ 1].$$
$$X \text{ o } A = [1\ 2\ 0\ 2 \mid 2\ 2 \mid 0\ 2\ 2 \mid 2\ 1\ 0\ 0\ 1] = B'.$$

(This ‘o’ operator is nothing but usual matrix multiplication). Clearly we see B' is not a fuzzy super row vector so we threshold B' to make B' a fuzzy super row vector by replacing every entry by 1 if that entry is a positive integer and by 0 if it is a negative integer or zero. Thus B' after thresholding becomes

$$[1\ 1\ 0\ 1 \mid 1\ 1 \mid 0\ 1\ 1 \mid 1\ 1\ 0\ 0\ 1] \quad = \quad B(\text{say}),$$
$$B \text{ o } A^T \hookrightarrow [1\ 1\ 1\ 1] \quad = \quad X_1\ .$$

'$\hookrightarrow$ denotes the resultant vector that has been thresholded'.
Now

$$X_1 \text{ o } A \hookrightarrow [1\ 1\ 1\ 1 \mid 1\ 1 \mid 1\ 1\ 1 \mid 1\ 1\ 1\ 1\ 1] = B_1.$$

.

Once again $B_1$ o $A^T$ will give only $X_1$. This sort of multiplication of fuzzy matrices is being carried out in the working of the super FRM models dealt in chapter three of this book.

Next we show how the product of a special fuzzy super column matrix is done which is also used in super FRM models.

***Example 2.3.16:*** Let P be a special super fuzzy column matrix given as follows :

$$P = \left[\begin{array}{ccccc} 1 & 0 & 1 & 1 & 1 \\ 0 & 1 & 0 & 0 & 1 \\ 1 & 1 & 0 & 0 & 0 \\ \hline 0 & 1 & 1 & 1 & 0 \\ 1 & 0 & 0 & 0 & 1 \\ 0 & 1 & 1 & 1 & 0 \\ 0 & 0 & 0 & 0 & 1 \\ \hline 1 & 0 & 0 & 0 & 0 \\ 0 & 1 & 0 & 0 & 0 \\ \hline 0 & 0 & 1 & 0 & 0 \\ 1 & 0 & 0 & 0 & 1 \end{array}\right].$$

Suppose

$$A = [1\ 0\ 1 \mid 1\ 1\ 1\ 0 \mid 1\ 0 \mid 1\ 0],\ AP = [4\ 3\ 4\ 3\ 2]$$

after thresholding we get

$$B = [1\ 1\ 1\ 1].\ BP^T = [4\ 2\ 2 \mid 3\ 2\ 3\ 1 \mid 1\ 1 \mid 1\ 2]$$

after thresholding we get

$$A_1 = [1\ 1\ 1 \mid 1\ 1\ 1\ 1 \mid 1\ 1 \mid 1\ 1].$$

Now

$$A_1P = [5\ 5\ 4\ 3\ 5] \hookrightarrow [1\ 1\ 1\ 1\ 1]$$

and so on.

We proceed on to work now with a diagonal super fuzzy matrix.

***Example 2.3.17:*** Let us consider $F_D$ a super diagonal fuzzy matrix

$$\left[\begin{array}{ccc|cccc|ccccc}
1 & 0 & 1 & 0 & 0 & 0 & 0 & 0 & 0 & 0 & 0 & 0 \\
1 & 0 & 0 & 0 & 0 & 0 & 0 & 0 & 0 & 0 & 0 & 0 \\
0 & 0 & 1 & 0 & 0 & 0 & 0 & 0 & 0 & 0 & 0 & 0 \\
1 & 1 & 0 & 0 & 0 & 0 & 0 & 0 & 0 & 0 & 0 & 0 \\
\hline
0 & 0 & 0 & 1 & 1 & 1 & 0 & 0 & 0 & 0 & 0 & 0 \\
0 & 0 & 0 & 0 & 1 & 0 & 1 & 0 & 0 & 0 & 0 & 0 \\
0 & 0 & 0 & 0 & 0 & 1 & 1 & 0 & 0 & 0 & 0 & 0 \\
0 & 0 & 0 & 0 & 0 & 0 & 1 & 0 & 0 & 0 & 0 & 0 \\
\hline
0 & 0 & 0 & 0 & 0 & 0 & 0 & 1 & 0 & 0 & 0 & 1 \\
0 & 0 & 0 & 0 & 0 & 0 & 0 & 0 & 1 & 0 & 1 & 0 \\
0 & 0 & 0 & 0 & 0 & 0 & 0 & 0 & 0 & 1 & 0 & 0 \\
0 & 0 & 0 & 0 & 0 & 0 & 0 & 0 & 0 & 0 & 1 & 1 \\
0 & 0 & 0 & 0 & 0 & 0 & 0 & 1 & 1 & 0 & 0 & 0 \\
0 & 0 & 0 & 0 & 0 & 0 & 0 & 0 & 0 & 1 & 0 & 1
\end{array}\right].$$

Now we see we have only diagonal matrix entries and the other entries are zero also we observe that the diagonal matrix are also of arbitrary order. How does this fuzzy super diagonal matrix $F_D$ function?

Suppose

$$A = [1\ 0\ 0\ 1 \mid 0\ 1\ 0\ 1 \mid 0\ 0\ 1\ 0\ 0\ 1]$$

then

$$A \text{ o } F_D = [2\ 1\ 1 \mid 0\ 1\ 0\ 2 \mid 0\ 0\ 2\ 0\ 1] = B'(\text{say})$$

after thresholding we get

$$B = [1\ 1\ 1 \mid 0\ 1\ 0\ 1 \mid 0\ 0\ 1\ 0\ 1]\ .$$

Now

$$B \text{ o } F_D^T = [2\ 1\ 1\ 2 \mid 1\ 2\ 1\ 1 \mid 1\ 0\ 1\ 1\ 0\ 2]$$

after thresholding we get

$$A_1 = [1\ 1\ 1\ 1 \mid 1\ 1\ 1\ 1 \mid 1\ 0\ 1\ 1\ 0\ 1].$$

$$A_1 \text{ o } F_D = [3\ 1\ 2 \mid 1\ 2\ 2\ 3 \mid 1\ 0\ 2\ 1\ 3]$$

after thresholding we get

$$B_1 = [1\ 1\ 1 \mid 1\ 1\ 1\ 1 \mid 1\ 0\ 1\ 1\ 1].$$

$$B_1 \text{ o } F_D^T = [2\ 1\ 1\ 2 \mid 3\ 2\ 2\ 1 \mid 2\ 1\ 1\ 2\ 1\ 2],$$

after thresholding we get

$$A_2 = [1\ 1\ 1\ 1 \mid 1\ 1\ 1\ 1 \mid 1\ 1\ 1\ 1\ 1\ 1].$$

$$A_2 \text{ o } F_D = [3\ 1\ 2 \mid 1\ 2\ 2\ 3 \mid 2\ 2\ 2\ 2\ 3]$$

after thresholding we get

$$B_2 = [1\ 1\ 1 \mid 1\ 1\ 1\ 1 \mid 1\ 1\ 1\ 1\ 1].$$

Thus our operation ends at this stage.

Now we proceed on to illustrate how a super fuzzy matrix functions.

***Example 2.3.18:*** Let X be a super fuzzy matrix given in the following;

$$X = \left[\begin{array}{ccc|ccccc|cccccc}
1 & 0 & 1 & 0 & 0 & 0 & 1 & 0 & 1 & 1 & 1 & 1 & 0 & 0 \\
0 & 1 & 0 & 1 & 0 & 0 & 0 & 0 & 0 & 0 & 0 & 0 & 1 & 1 \\
\hline
1 & 1 & 0 & 1 & 0 & 0 & 0 & 1 & 1 & 1 & 0 & 0 & 1 & 1 \\
1 & 0 & 1 & 0 & 0 & 0 & 1 & 0 & 0 & 0 & 1 & 1 & 0 & 0 \\
0 & 1 & 0 & 0 & 0 & 1 & 0 & 0 & 1 & 1 & 0 & 0 & 0 & 0 \\
1 & 0 & 0 & 0 & 1 & 0 & 0 & 0 & 0 & 0 & 1 & 1 & 1 & 1 \\
0 & 0 & 1 & 0 & 1 & 0 & 1 & 0 & 1 & 0 & 1 & 0 & 1 & 0 \\
1 & 0 & 0 & 1 & 0 & 0 & 0 & 0 & 1 & 0 & 0 & 0 & 0 & 0 \\
\hline
0 & 1 & 0 & 0 & 1 & 0 & 0 & 0 & 0 & 0 & 0 & 0 & 1 & 0 \\
0 & 0 & 1 & 0 & 0 & 1 & 0 & 0 & 0 & 0 & 0 & 0 & 0 & 1 \\
1 & 1 & 1 & 0 & 0 & 1 & 1 & 1 & 0 & 1 & 0 & 0 & 0 & 0 \\
1 & 0 & 1 & 1 & 1 & 0 & 0 & 0 & 0 & 0 & 1 & 0 & 0 & 0 \\
0 & 1 & 1 & 0 & 1 & 1 & 0 & 0 & 1 & 1 & 1 & 0 & 0 & 0
\end{array}\right].$$

We show how this super fuzzy matrix operation is carried out. Let

$$A = [0\ 1 \mid 0\ 0\ 1\ 1\ 0\ 0 \mid 1\ 0\ 1\ 0\ 0]$$

$$A \text{ o } X = [2\ 4\ 1 \mid 1\ 2\ 2\ 1\ 1 \mid 1\ 1\ 1\ 1\ 2\ 2];$$

after thresholding we get

$$B = [1\ 1\ 1 \mid 1\ 1\ 1\ 1\ 1 \mid 1\ 1\ 1\ 1\ 1\ 1]\ .$$

$$B \text{ o } X^T = [7\ 4 \mid 8\ 5\ 4\ 6\ 6 \mid 3\ 3\ 3\ 7\ 5\ 7]$$

after thresholding we get

$$A = [1\ 1 \mid 1\ 1\ 1\ 1\ 1\ 1 \mid 1\ 1\ 1\ 1\ 1].$$
$$A_1 \text{ o } X = [7\ 6\ 7 \mid 4\ 5\ 4\ 4\ 2 \mid 6\ 5\ 6\ 3\ 5\ 4]$$

after thresholding we get

$$B_1 = [1\ 1\ 1 \mid 1\ 1\ 1\ 1\ 1 \mid 1\ 1\ 1\ 1\ 1\ 1].$$

Thus we complete our calculations as we would obtain the same fuzzy row vector. Now we see this type of usual multiplication is carried out in case of all super FRM models. The only difference from the usual multiplication is that after each multiplication the resultant super row vector must be thresholded so that it is again a fuzzy super vector.

Next we proceed on to give yet another type of multiplication of fuzzy super matrices using the max min operator.

***Example 2.3.19:*** Let us consider a special fuzzy row matrix;

$$M = \left[\begin{array}{ccc|cc|ccccc} 0.1 & 1 & 0.3 & 1 & 0.7 & 0.7 & 0.2 & 0.3 & 1 & 0.4 \\ 1 & 0.2 & 0 & 0.3 & 0.2 & 0 & 0.5 & 0 & 0.3 & 1 \\ 0 & 0.7 & 1 & 0 & 0.1 & 1 & 0 & 0.4 & 1 & 0.8 \\ 0.4 & 0 & 0.6 & 0.5 & 0 & 0 & 0.4 & 0 & 0.5 & 0 \end{array}\right].$$

Let us consider the fuzzy row vector A = [0.3 0 0.1 0.7]. Using max min operation we find

$$\begin{aligned} \max\min\{A, M\} &= \max\min(a_i, m_{ij}) \\ &= [0.4\ 0.3\ 0.6 \mid 0.5\ 0.3 \mid 0.3\ 0.4\ 0.3\ 0.5\ 0.3] \\ &= B. \end{aligned}$$

Now max min $\{B, A^T\}$ or max min $\{A, B^T\}$ gives

$$[0.5\ 0.4\ 0.6\ 0.6] \quad = \quad A_1$$

using the max-min operation. Now using max min $\{A_1, M\}$ operation gives

$$[0.4\ 0.6\ 0.6 \mid 0.5\ 0.5 \mid 0.6\ 0.4\ 0.4\ 0.6\ 0.6] = B_1 \text{ (say)}.$$

This procedure is repeated.

Now we proceed on to show how a special super fuzzy column matrix works under min max operation.

***Example 2.3.20:*** Let T be a special column super fuzzy matrix which is given below;

$$T = \left[\begin{array}{ccccc} 0.1 & 1 & 0.3 & 0.2 & 0 \\ 1 & 0.4 & 1 & 1 & 0.5 \\ 0.6 & 0 & 0.7 & 0 & 0.1 \\ \hline 1 & 0.3 & 0 & 0.5 & 1 \\ 0 & 0.4 & 0.3 & 0.2 & 0.5 \\ \hline 0.7 & 1 & 0.5 & 1 & 0 \\ 0.1 & 0 & 0.3 & 0.1 & 0.2 \\ 0.7 & 0.2 & 0.1 & 0.6 & 0.5 \\ 0.6 & 0.3 & 1 & 0 & 1 \\ \hline 0.7 & 0.3 & 0.4 & 0.8 & 1 \\ 1 & 0 & 0.7 & 1 & 0.3 \\ 0.2 & 1 & 0.2 & 0.3 & 1 \\ 0.5 & 0.6 & 1 & 0.3 & 0.2 \\ 1 & 0 & 0.1 & 0 & 0.2 \end{array}\right].$$

Given A is fuzzy super row vector

$$\text{i.e., } A = [1\ 0.3\ 0 \mid 0.2\ 1 \mid 0\ 0.1\ 0\ 1 \mid 1\ 0.3\ 0.2\ 0.5\ 1].$$

We calculate max min $\{A, T\}$ =

$$[1\ 1\ 1\ 0.8\ 1] = B.$$

Now we calculate same max min $\{T, B^T\}$ to be

$$[1\ 1\ 0.7 \mid 1\ 0.5 \mid 1\ 0.3\ 0.7\ 1 \mid 1\ 1\ 1\ 1\ 1] = A_1.$$

We can calculate using the fuzzy super row vector $A_1$ the value max min $\{A_1, T\}$ and so on.

Now we proceed on to show how the special fuzzy super diagonal matrix operation is carried out using the max min principle.

***Example 2.3.21:*** Let S be a special fuzzy super diagonal matrix.

$$S = \left[\begin{array}{cccc|ccccc|ccc}
0.3 & 1 & 0.2 & 0.4 & 0 & 0 & 0 & 0 & 0 & 0 & 0 & 0 \\
1 & 0.3 & 0 & 0.7 & 0 & 0 & 0 & 0 & 0 & 0 & 0 & 0 \\
0.1 & 0.1 & 0.4 & 0 & 0 & 0 & 0 & 0 & 0 & 0 & 0 & 0 \\
1 & 1 & 1 & 0.8 & 0 & 0 & 0 & 0 & 0 & 0 & 0 & 0 \\
\hline
0 & 0 & 0 & 0 & 1 & 0 & 0.3 & 0.2 & 0.5 & 0 & 0 & 0 \\
0 & 0 & 0 & 0 & 0.7 & 1 & 0.1 & 0.3 & 1 & 0 & 0 & 0 \\
0 & 0 & 0 & 0 & 0.2 & 0.3 & 1 & 0.4 & 0.6 & 0 & 0 & 0 \\
\hline
0 & 0 & 0 & 0 & 0 & 0 & 0 & 0 & 0 & 1 & 0.7 & 0.3 \\
0 & 0 & 0 & 0 & 0 & 0 & 0 & 0 & 0 & 0.1 & 1 & 0.2 \\
0 & 0 & 0 & 0 & 0 & 0 & 0 & 0 & 0 & 0.7 & 0.3 & 0.9 \\
0 & 0 & 0 & 0 & 0 & 0 & 0 & 0 & 0 & 1 & 0 & 1 \\
0 & 0 & 0 & 0 & 0 & 0 & 0 & 0 & 0 & 0.2 & 1 & 0.3 \\
0 & 0 & 0 & 0 & 0 & 0 & 0 & 0 & 0 & 0 & 0.7 & 1
\end{array}\right].$$

Suppose we are given the fuzzy super row vector

$$A = [0.1\ 0\ 0.5\ 1 \mid 0.7\ 0.2\ 1 \mid 0.6\ 1\ 0.3\ 0.2\ 0.5\ 0.7].$$

Now using the max.min operation max min $\{A, S\}$ gives

$$[1\ 1\ 1\ 0.8 \mid 0.7\ 0.3\ 1\ 0.4\ 0.6 \mid 0.6\ 1\ 0.7] = B(\text{say}).$$

Now max min{S, $B^T$} gives

$$[1\ 1\ 0.4\ 1\ 0.7\ 0.7\ 1\ 0.7\ 1\ 0.7\ 0.7\ 1\ 0.7] = A_1.$$

$A_1$ is also a fuzzy super row vector.

We can proceed on in this way, to arrive at a resultant.

Now we show how to work with a fuzzy super matrix using a min max operator.

***Example 2.3.22:*** Let V be a super fuzzy matrix given by

$$\left[\begin{array}{cccc|ccc|ccccc}
0.3 & 0.2 & 1 & 0.6 & 0.7 & 0.1 & 0.8 & 0.3 & 0.2 & 0.5 & 0.7 & 1 \\
0.8 & 0.5 & 0.3 & 1 & 1 & 0 & 0.7 & 0.2 & 0.4 & 0.3 & 1 & 0.3 \\
1 & 0.2 & 0.4 & 0.6 & 0 & 1 & 0.8 & 0.8 & 0.7 & 0.5 & 0.7 & 1 \\
0.5 & 0.3 & 1 & 0.5 & 0.9 & 0.3 & 0.1 & 0.9 & 0.4 & 0.1 & 0 & 0.8 \\
0.2 & 1 & 0 & 0.7 & 1 & 0 & 1 & 0.1 & 0.3 & 1 & 0.7 & 1 \\
\hline
1 & 0.3 & 0 & 0.8 & 0.7 & 0.3 & 0.6 & 1 & 0 & 0.3 & 0.2 & 0.9 \\
0.3 & 1 & 0.6 & 1 & 0.3 & 0.8 & 0.2 & 0 & 1 & 0.8 & 1 & 0.1 \\
0.2 & 0 & 0.8 & 0.7 & 1 & 1 & 0 & 0.4 & 1 & 0.5 & 0.7 & 1 \\
0.7 & 0.2 & 0.1 & 0.6 & 0.8 & 0.9 & 1 & 0.7 & 0.6 & 1 & 0.8 & 0.4 \\
0 & 1 & 0.7 & 1 & 0.2 & 0 & 0.7 & 0 & 0 & 0 & 0.9 & 0.6 \\
0.9 & 0.7 & 0 & 0.7 & 0.8 & 0.8 & 1 & 0.5 & 0.4 & 0.8 & 0.9 & 0.2 \\
0.1 & 1 & 0.9 & 0.5 & 1 & 0.6 & 0.5 & 0.1 & 0.2 & 0 & 0.6 & 0.1 \\
\hline
0.3 & 1 & 0 & 0.6 & 0.3 & 0.9 & 1 & 0.7 & 1 & 0.8 & 0.6 & 1 \\
1 & 0.7 & 0.6 & 1 & 0 & 0.6 & 0.3 & 0.8 & 0 & 0.5 & 1 & 0 \\
0 & 0.6 & 0.5 & 0.8 & 1 & 0.2 & 0.1 & 0 & 1 & 0.6 & 0 & 1 \\
0.8 & 1 & 0.2 & 0.1 & 0.3 & 0.6 & 1 & 0.7 & 0 & 0.2 & 0.3 & 1
\end{array}\right].$$

Let

X = [0.2 0.3 0.2 0.1 0 | 0.3 0.2 1 0.7 0.8 0.4 0.4 | 0.2 0.5 0.6 0.3].

One can find min max{X, V}.

If

$$\min \max \{X, V\} = W$$

then we find min max{W, $V^T$} which will be a fuzzy super row vector.

Now having described the working of super fuzzy matrices using min max operation we just state that these methods help in working with all types of super FAM described in chapter three of this book.

Now we just show how we work with a special type of fuzzy super matrices which will be used when we use the super FCM models.

We just show how a special type of fuzzy super diagonal matrix is defined, the authors say this matrix is a special type for the diagonal matrix are square matrix with a further condition that the diagonal elements of these matrices are zero.

***Example 2.3.23:*** Let W be the special fuzzy diagonal matrix

$$W = \left[\begin{array}{ccccc|cccc|cccccc}
0 & 1 & -1 & 0 & 1 & 0 & 0 & 0 & 0 & 0 & 0 & 0 & 0 & 0 & 0 \\
1 & 0 & 0 & 1 & 0 & 0 & 0 & 0 & 0 & 0 & 0 & 0 & 0 & 0 & 0 \\
0 & 0 & 0 & 0 & 1 & 0 & 0 & 0 & 0 & 0 & 0 & 0 & 0 & 0 & 0 \\
0 & 1 & 0 & 0 & 0 & 0 & 0 & 0 & 0 & 0 & 0 & 0 & 0 & 0 & 0 \\
0 & 0 & 1 & 0 & 0 & 0 & 0 & 0 & 0 & 0 & 0 & 0 & 0 & 0 & 0 \\
\hline
0 & 0 & 0 & 0 & 0 & 0 & 1 & 0 & 0 & 0 & 0 & 0 & 0 & 0 & 0 \\
0 & 0 & 0 & 0 & 0 & 0 & 0 & -1 & 1 & 0 & 0 & 0 & 0 & 0 & 0 \\
0 & 0 & 0 & 0 & 0 & 1 & 0 & 0 & -1 & 0 & 0 & 0 & 0 & 0 & 0 \\
0 & 0 & 0 & 0 & 0 & 0 & 1 & 1 & 0 & 0 & 0 & 0 & 0 & 0 & 0 \\
\hline
0 & 0 & 0 & 0 & 0 & 0 & 0 & 0 & 0 & 0 & 1 & 1 & 0 & 0 & 1 \\
0 & 0 & 0 & 0 & 0 & 0 & 0 & 0 & 0 & 0 & 0 & 0 & 1 & 0 & 0 \\
0 & 0 & 0 & 0 & 0 & 0 & 0 & 0 & 0 & 1 & 0 & 0 & 0 & 1 & 0 \\
0 & 0 & 0 & 0 & 0 & 0 & 0 & 0 & 0 & 0 & 0 & 0 & 0 & 0 & 1 \\
0 & 0 & 0 & 0 & 0 & 0 & 0 & 0 & 0 & 1 & 0 & 1 & 1 & 0 & 0 \\
0 & 0 & 0 & 0 & 0 & 0 & 0 & 0 & 0 & 0 & 1 & 0 & 0 & 1 & 0
\end{array}\right].$$

Suppose we wish to work with the state fuzzy row vector

$$X = [1\ 0\ 0\ 0\ 0 \mid 0\ 1\ 1\ 0 \mid 0\ 0\ 0\ 1\ 1\ 0].$$

Now

$$X \text{ o } W = [0\ 1\ {-1}\ 0\ 1 \mid 1\ 0\ {-1}\ 0 \mid 1\ 0\ 1\ 1\ 0\ 1];$$

After thresholding and updating X o W we get

$$[1\ 1\ 0\ 0\ 1 \mid 1\ 1\ 1\ 0 \mid 1\ 0\ 1\ 1\ 1\ 1] = Y.$$

Now we find Y o W and so on until we arrive at a fixed point.

**Chapter Three**

# INTRODUCTION TO NEW FUZZY SUPER MODELS

In this chapter we introduce five new fuzzy super models and illustrate how they exploit the concept of supermatrices or special supervectors. These models will be very much useful when we have multi experts working with varying attributes. This chapter has five sections. In section one we introduce the notion of Super Fuzzy Relational Maps model, in section two the notion of Super Bidirectional Associative Memories (SBAM) models are introduced and, section 3 introduces the new notion of Super Fuzzy Associative Memories (SFAM). Section four gives the applications of these new super fuzzy models by means of illustrations. The final section introduces the new notion of super FCMs model.

We just say a supermatrix is a fuzzy supermatrix if its entries are from the interval [0, 1]. The operations on these models are basically max min operations whenever the compatibility exists.

## 3.1 New Super Fuzzy Relational Maps (SFRM) Model

In this section for the first time we introduce the new notion of Super Fuzzy Relational Maps (SFRMs) models and they are applied to real world problems, which is suited for multi expert problems. When we in the place of fuzzy relational matrix of the FRM (using single expert) use fuzzy supermatrix with multi experts we call the model as Super Fuzzy Relational Maps (SFRMs) model.

This section has two subsections. In the first subsection we just recall the basic properties of fuzzy relational maps in subsection two we define new fuzzy super relational maps models.

### 3.1.1 Introduction to Fuzzy Relational Maps (FRMs)

In this section we just recall the definition and properties of FRM given in [231].

In this section, we introduce the notion of Fuzzy Relational Maps (FRMs); they are constructed analogous to FCMs described and discussed in the earlier sections of [231]. In FCMs we promote the correlations between causal associations among concurrently active units. But in FRMs we divide the very causal associations into two disjoint units, for example, the relation between a teacher and a student or relation between an employee and employer or a relation between doctor and patient and so on. Thus for us to define a FRM we need a domain space and a range space which are disjoint in the sense of concepts. We further assume no intermediate relation exists within the domain elements or node and the range spaces elements. The number of elements in the range space need not in general be equal to the number of elements in the domain space.

Thus throughout this section we assume the elements of the domain space are taken from the real vector space of dimension n and that of the range space are real vectors from the vector space of dimension m (m in general need not be equal to n). We denote by R the set of nodes $R_1,\ldots, R_m$ of the range space, where $R = \{(x_1,\ldots, x_m) \mid x_j = 0 \text{ or } 1\}$ for j = 1, 2, …, m. If $x_i = 1$ it means that the node $R_i$ is in the on state and if $x_i = 0$ it means that the node $R_i$ is in the off state. Similarly D denotes the nodes $D_1, D_2, \ldots, D_n$ of the domain space where $D = \{(x_1,\ldots, x_n) \mid x_j = 0 \text{ or } 1\}$ for i = 1, 2,…, n. If $x_i = 1$ it means that the node $D_i$ is in the on state and if $x_i = 0$ it means that the node $D_i$ is in the off state.

Now we proceed on to define a FRM.

**DEFINITION 3.1.1.1:** *A FRM is a directed graph or a map from D to R with concepts like policies or events etc, as nodes and causalities as edges. It represents causal relations between spaces D and R.*

*Let $D_i$ and $R_j$ denote that the two nodes of an FRM. The directed edge from $D_i$ to $R_j$ denotes the causality of $D_i$ on $R_j$ called relations. Every edge in the FRM is weighted with a number in the set {0, ±1}. Let $e_{ij}$ be the weight of the edge $D_iR_j$, $e_{ij} \in \{0, \pm 1\}$. The weight of the edge $D_i R_j$ is positive if increase in $D_i$ implies increase in $R_j$ or decrease in $D_i$ implies decrease in $R_j$, i.e., causality of $D_i$ on $R_j$ is 1. If $e_{ij} = 0$, then $D_i$ does not have any effect on $R_j$ . We do not discuss the cases when increase in $D_i$ implies decrease in $R_j$ or decrease in $D_i$ implies increase in $R_j$ .*

**DEFINITION 3.1.1.2:** *When the nodes of the FRM are fuzzy sets then they are called fuzzy nodes. FRMs with edge weights {0, ±1} are called simple FRMs.*

**DEFINITION 3.1.1.3:** *Let $D_1, \ldots, D_n$ be the nodes of the domain space D of an FRM and $R_1, \ldots, R_m$ be the nodes of the range space R of an FRM. Let the matrix E be defined as $E = (e_{ij})$ where $e_{ij}$ is the weight of the directed edge $D_iR_j$ (or $R_jD_i$), E is called the relational matrix of the FRM.*

***Note:*** It is pertinent to mention here that unlike the FCMs the FRMs can be a rectangular matrix with rows corresponding to the domain space and columns corresponding to the range space. This is one of the marked difference between FRMs and FCMs.

**DEFINITION 3.1.1.4:** *Let $D_1, \ldots, D_n$ and $R_1,\ldots, R_m$ denote the nodes of the FRM. Let $A = (a_1,\ldots,a_n)$, $a_i \in \{0, 1\}$. A is called the instantaneous state vector of the domain space and it denotes the on-off position of the nodes at any instant. Similarly let $B = (b_1,\ldots, b_m)$, $b_i \in \{0, 1\}$. B is called instantaneous state vector of* *the range space and it denotes the on-off position of the nodes at any instant $a_i = 0$ if $a_i$ is off and $a_i = 1$ if $a_i$ is on for i= 1, 2,…, n. Similarly, $b_i = 0$ if $b_i$ is off and $b_i = 1$ if $b_i$ is on, for i= 1, 2,…, m.*

**DEFINITION 3.1.1.5:** *Let $D_1, \ldots, D_n$ and $R_1, \ldots, R_m$ be the nodes of an FRM. Let $D_iR_j$ (or $R_j D_i$) be the edges of an FRM, j = 1, 2,…, m and i= 1, 2,…, n. Let the edges form a directed cycle. An FRM is said to be a cycle if it posses a directed cycle. An FRM is said to be acyclic if it does not posses any directed cycle.*

**DEFINITION 3.1.1.6:** *An FRM with cycles is said to be an FRM with feedback.*

**DEFINITION 3.1.1.7:** *When there is a feedback in the FRM, i.e. when the causal relations flow through a cycle in a revolutionary manner, the FRM is called a dynamical system.*

**DEFINITION 3.1.1.8:** *Let $D_i R_j$ (or $R_j D_i$), $1 \leq j \leq m$, $1 \leq i \leq n$. When $R_i$ (or $D_j$) is switched on and if causality flows through edges of the cycle and if it again causes $R_i$ (or $D_j$), we say that the dynamical system goes round and round. This is true for any node $R_j$ (or $D_i$) for $1 \leq i \leq n$, (or $1 \leq j \leq m$). The equilibrium state of this dynamical system is called the hidden pattern.*

**DEFINITION 3.1.1.9:** *If the equilibrium state of a dynamical system is a unique state vector, then it is called a fixed point. Consider an FRM with $R_1, R_2, \ldots, R_m$ and $D_1, D_2, \ldots, D_n$ as nodes. For example, let us start the dynamical system by switching on $R_1$ (or $D_1$). Let us assume that the FRM settles down with $R_1$ and $R_m$ (or $D_1$ and $D_n$) on, i.e. the state vector remains as (1, 0, …, 0, 1) in R (or 1, 0, 0, … , 0, 1) in D), This state vector is called the fixed point.*

**DEFINITION 3.1.1.10:** *If the FRM settles down with a state vector repeating in the form*

$$A_1 \to A_2 \to A_3 \to \ldots \to A_i \to A_1 \text{ (or } B_1 \to B_2 \to \ldots \to B_i \to B_1)$$

*then this equilibrium is called a limit cycle.*

**METHODS OF DETERMINING THE HIDDEN PATTERN**

Let $R_1$, $R_2$, …, $R_m$ and $D_1$, $D_2$, …, $D_n$ be the nodes of a FRM with feedback. Let E be the relational matrix. Let us find a hidden pattern when $D_1$ is switched on i.e. when an input is given as vector $A_1 = (1, 0, …, 0)$ in $D_1$, the data should pass through the relational matrix E. This is done by multiplying $A_1$ with the relational matrix E. Let $A_1E = (r_1, r_2, …, r_m)$, after thresholding and updating the resultant vector we get $A_1 E \in R$. Now let $B = A_1E$ we pass on B into $E^T$ and obtain $BE^T$. We update and threshold the vector $BE^T$ so that $BE^T \in D$. This procedure is repeated till we get a limit cycle or a fixed point.

**DEFINITION 3.1.1.11:** *Finite number of FRMs can be combined together to produce the joint effect of all the FRMs. Let $E_1$,…, $E_p$ be the relational matrices of the FRMs with nodes $R_1$, $R_2$,…, $R_m$ and $D_1$, $D_2$,…, $D_n$, then the combined FRM is represented by the relational matrix $E = E_1+…+ E_p$.*

## 3.1.2 New Super Fuzzy Relational Maps models

In this section we introduce four types of new Super Fuzzy Relational Maps models.

**DEFINITION 3.1.2.1:** *Suppose we have some n experts working on a real world model and give their opinion. They all agree upon to work with the same domain space elements / attributes / concepts; using FRM model but do not concur on the attributes from the range space then we can use the special super fuzzy row vector to model the problem using Domain Super Fuzzy Relational Maps (DSFRMs) Model.*

The DSFRM matrix associated with this model will be given by $S_M$

$$S_M = \begin{array}{c} \\ D_1 \\ D_2 \\ \vdots \\ D_m \end{array} \begin{array}{c} \begin{array}{cccc} t_1^1 \ldots t_{r_1}^1 & t_1^2 \ldots t_{r_2}^2 & & t_1^n t_2^n \ldots t_{r_n}^n \end{array} \\ \left[ \begin{array}{c|c|c|c} & & & \\ & & & \\ & & & \\ & & & \end{array} \right] \end{array}$$

$$= \left[ S_M^1 \mid S_M^2 \mid \ldots \mid S_M^n \right]$$

where each $S_M^i$ is a m × $t_{r_i}^i$ matrix associated with a FRM given by the $i^{th}$ expert having $D_1, \ldots, D_m$ to be the domain attributes and ( $t_1^i\ t_2^i\ \ldots\ t_{r_i}^i$ ) to be the range attributes of the $i^{th}$ expert, i = 1, 2, …, n and $S_M$ the DSFRM matrix will be a special super row vector / matrix ($1 \leq i \leq n$).

However if n is even a very large value using the mode of programming one can easily obtain the resultant vector or the super hidden pattern for any input supervector which is under investigation.

These DSFRMs will be known as Domain constant DSFRMs for all the experts choose to work with the same domain space attributes only the range space attributes are varying denoted by DSFRM models.

Next we proceed on to define the notion of super FRMs with constant range space attributes and varying domain space attributes.

**DEFINITION 3.1.2.2:** *Let some m experts give opinion on a real world problem who agree upon to make use of the same space of attributes / concepts from the range space using FRMs but want to use different concepts for the domain space then we make use of the newly constructed special super column vector as the matrix to construct this new model. Thus the associated special super column matrix $S^M$ is*

$$S^M = \begin{array}{c} \\ D_1^1 \\ D_2^1 \\ \vdots \\ D_{t_1}^1 \\ D_1^2 \\ \vdots \\ D_{t_2}^2 \\ \vdots \\ D_1^m \\ \vdots \\ D_{t_m}^m \end{array} \begin{array}{c} \begin{array}{cccc} R_1 & R_2 & \ldots & R_s \end{array} \\ \left[ \begin{array}{cccc} & & & \\ & & & \\ & & & \\ & & & \\ \hline & & & \\ & & & \\ & & & \\ \hline & & & \\ \hline & & & \\ & & & \\ & & & \end{array} \right] \end{array}$$

*The column attributes i.e. the range space remain constant as $R_1, \ldots, R_s$ for all m experts only the row attributes for any $i^{th}$ expert is $D_1^i, D_2^i, \ldots, D_{t_i}^i$; i = 1, 2, …, m; vary from expert to expert. This system will be known as the Range constant fuzzy super FRM or shortly denoted as RSFRM model.*

We will illustrate these two definitions by examples i.e. by live problems in the $4^{th}$ section of this chapter.

One may be interested in finding a model which has both the range and the column vector varying for some experts how to construct a new model in that case.

**DEFINITION 3.1.2.3:** *Suppose we have m experts who wish to work with different sets of both row and column attributes i.e. domain and range space using FRMs, then to accommodate or form a integrated matrix model to cater to this need. We make use of the super diagonal fuzzy matrix, to model such a problem. Suppose the first expert works with the domain attributes $D_1^1, \ldots, D_{t_1}^1$ and range attributes $R_1^1, \ldots, R_{n_1}^1$, The* *second expert works with domain attributes* $D_1^2, \ldots, D_{t_2}^2$ *and with range attributes* $R_1^2, \ldots, R_{n_2}^2$ *and so on. Thus the* $m^{th}$ *expert works with* $D_1^m, \ldots, D_{t_m}^m$ *domain attributes and* $R_1^m, \ldots, R_{n_m}^m$ *range attributes. We have the following diagonal fuzzy supermatrix to model the situation. We are under the assumption that all the attributes both from the domain space as well as the range space of the m experts are different. The super fuzzy matrix S associated with this new model is given by*

$$
\begin{array}{c} \\ D_1^1 \\ \vdots \\ D_{t_1}^1 \\ D_1^2 \\ \vdots \\ D_{t_2}^2 \\ \vdots \\ D_1^m \\ \vdots \\ D_{t_m}^m \end{array}
\begin{array}{c}
\begin{array}{cccc} R_1^1\, R_2^1 \ldots R_{n_1}^1 & R_1^2\, R_2^2 \ldots R_{n_2}^2 & \ldots & R_1^m\, R_2^m \ldots R_{t_m}^m \end{array} \\
\left[\begin{array}{c|c|c|c}
M_1 & (0) & & (0) \\ \hline
(0) & M_2 & & (0) \\ \hline
(0) & (0) & & (0) \\ \hline
(0) & (0) & & M_n
\end{array}\right]
\end{array}
$$

*where each* $M_i$ *is a* $t_i \times n_i$ *matrix associated with the FRM, we see except, the diagonal strip all other entries are zero. We call this matrix as a special diagonal super fuzzy matrix and this model will be known as the Special Diagonal Super FRM Model which will be denoted briefly as (SDSFRM).*

Now we define the multi expert super FRM model.

**DEFINITION 3.1.2.4:** *Suppose one is interested in finding a model where some mn number of experts work on the problem and some have both domain and range attributes to be not* *coinciding with any other expert and a set of experts have only the domain attributes to be in common and all the range attributes are different. Another set of experts are such that only the range attributes to be in common and all the domain attributes are different, and all of them wish to work with the FRM model only; then we model this problem using a super fuzzy matrix. We have mn experts working with the problem.*

*Let the $t_1$ expert wish to work with domain attributes $P_1^1$, $P_2^1$, …, $P_{m(t_1)}^1$ and range attributes $q_1^1$, $q_2^1$, …, $q_{n(t_1)}^1$.*

*The $t_2$ expert works with $P_1^1$, $P_2^1$, …, $P_{m(t_1)}^1$ as domain attributes and the range attributes $q_1^2$, $q_2^2$, …, $q_{n(t_2)}^2$ and so on. Thus for the $t_i$ expert works with $P_1^i$, $P_2^i$, …, $P_{m(t_i)}^i$ as domain space attributes and $q_1^i$, $q_2^i$, …, $q_{n(t_i)}^i$, as range attributes ($1 \leq i \leq m(t_i)$ and $i \leq n(t_i)$).*

*So with these mn experts we have an associated super FRM matrix. Thus the supermatrix associated with the Super FRM (SFRM) model is a supermatrix of the form S(m) =*

$$\left[\begin{array}{c|c|c|c}
A^{11}_{m(t1)n(t1)} & A^{12}_{m(t1)n(t2)} & & A^{1n}_{m(t1)n(tn)} \\ \hline
A^{21}_{m(t2)n(t1)} & A^{22}_{m(t2)n(t2)} & & A^{2n}_{m(t2)n(tn)} \\ \hline
 & & & \\ \hline
A^{m1}_{m(tm)n(t1)} & A^{m2}_{m(tm)n(t2)} & & A^{mn}_{m(tm)n(tn)}
\end{array}\right]$$

where

$$A^{ij}_{m(t_i)n(t_j)} = \begin{array}{c} \\ P_1^i \\ P_2^i \\ \vdots \\ P_{m(t_i)}^i \end{array} \begin{array}{c} q_1^j \quad q_2^j \quad \cdots \quad q_{n(t_j)}^j \\ \left[ \begin{array}{c} \\ \\ (a^{ij}_{m(t_i)n(t_j)}) \\ \\ \\ \end{array} \right] \end{array}$$

*$1 \le i \le m$ and $1 \le j \le n$. S(m) is called the super dynamical FRM or a super dynamical system.*

*This matrix $A^{ij}_{m(t_i)n(t_j)}$ corresponds to the FRM matrix of the (ij)$^{th}$ expert with domain space attributes $P^i_1$, $P^i_2$, …, $P^i_{m(t_i)}$ and range space attributes $q^j_1$, $q^j_2$, …, $q^j_{n(t_j)}$, $1 \le i \le m$ and $1 \le j \le n$.*

Thus we have four types of super FRM models viz. range constant super FRM model RSFRM model or the row (vector) super FRM model with a row fuzzy supervector (supermatrix) associated as a dynamical system, domain constant super FRM model or DSFRM model with the column fuzzy supervector associated matrix as the dynamical system, diagonal super FRM model or SDSFRM with only diagonal having the fuzzy FRM matrices and its related dynamical system to be a fuzzy super diagonal matrix and finally the fuzzy super FRM (SFRM) model which is depicted by a fuzzy supermatrix.

Now having defined four types of super FRM models we proceed on to define super BAM models.

## 3.2 New Fuzzy Super Bidirectional Associative Memories (BAM) model

This section has two subsections. In the first subsection we recall the definition of Bidirectional Associative Memories models from [112]. In the subsection two the new types of super Bidirectional associative memories model are introduced.

### 3.2.1 Introduction to BAM model

In this section we just recall the definition of BAM of model from [112].

Now we go forth to describe the mathematical structure of the Bidirectional Associative Memories (BAM) model. Neural networks recognize ill defined problems without an explicit set of rules. Neurons behave like functions, neurons transduce an unbounded input activation x(t) at time t into a bounded output signal S(x(t)) i.e. Neuronal activations change with time.

Artificial neural networks consists of numerous simple processing units or neurons which can be trained to estimate sampled functions when we do not know the form of the functions. A group of neurons form a field. Neural networks contain many field of neurons. In our text $F_x$ will denote a neuron field, which contains n neurons, and $F_y$ denotes a neuron field, which contains p neurons. The neuronal dynamical system is described by a system of first order differential equations that govern the time-evolution of the neuronal activations or which can be called also as membrane potentials.

$$\dot{x}_i = g_i(X, Y, ...)$$
$$\dot{y}_j = h_j(X, Y, ...)$$

where $\dot{x}_i$ and $\dot{y}_j$ denote respectively the activation time function of the $i^{th}$ neuron in $F_X$ and the $j^{th}$ neuron in $F_Y$. The over dot denotes time differentiation, $g_i$ and $h_j$ are some functions of X, Y, ... where $X(t) = (x_1(t), \ldots, x_n(t))$ and $Y(t) = (y_1(t), \ldots, y_p(t))$ define the state of the neuronal dynamical system at time t. The passive decay model is the simplest activation model, where in the absence of the external stimuli, the activation decays in its resting value

$$\dot{x}_i = x_i$$

and $$\dot{y}_j = y_j.$$

The passive decay rate $A_i > 0$ scales the rate of passive decay to the membranes resting potentials $\dot{x}_i = -A_i x_i$. The default rate is $A_i = 1$, i.e. $\dot{x}_i = -A_i x_i$. The membrane time constant $C_i > 0$ scales the time variables of the activation dynamical system. The default time constant is $C_i = 1$. Thus $C_i \dot{x}_i = -A_i x_i$.

The membrane resting potential $P_i$ is defined as the activation value to which the membrane potential equilibrates in the absence of external inputs. The resting potential is an additive constant and its default value is zero. It need not be positive.

$$P_i = C_i \dot{x}_i + A_i x_i$$
$$I_i = \dot{x}_i + x_i$$

is called the external input of the system. Neurons do not compute alone. Neurons modify their state activations with external input and with feed back from one another. Now, how do we transfer all these actions of neurons activated by inputs their resting potential etc. mathematically. We do this using what are called synaptic connection matrices.

Let us suppose that the field $F_X$ with n neurons is synaptically connected to the field $F_Y$ of p neurons. Let $m_{ij}$ be a synapse where the axon from the $i^{th}$ neuron in $F_X$ terminates. $M_{ij}$ can be positive, negative or zero. The synaptic matrix M is a n by p matrix of real numbers whose entries are the synaptic efficacies $m_{ij}$.

The matrix M describes the forward projections from the neuronal field $F_X$ to the neuronal field $F_Y$. Similarly a p by n synaptic matrix N describes the backward projections from $F_Y$ to $F_X$. Unidirectional networks occur when a neuron field synaptically intra connects to itself. The matrix M be a n by n square matrix. A Bidirectional network occur if $M = N^T$ and $N = M^T$. To describe this synaptic connection matrix more simply, suppose the n neurons in the field $F_X$ synaptically connect to the p-neurons in field $F_Y$. Imagine an axon from the $i^{th}$ neuron in $F_X$ that terminates in a synapse $m_{ij}$, that about the $j^{th}$ neuron in $F_Y$. We assume that the real number $m_{ij}$ summarizes the synapse and that $m_{ij}$ changes so slowly relative to activation fluctuations that is constant.

Thus we assume no learning if $m_{ij} = 0$ for all t. The synaptic value $m_{ij}$ might represent the average rate of release of a neuro-transmitter such as norepinephrine. So, as a rate, $m_{ij}$ can be positive, negative or zero.

When the activation dynamics of the neuronal fields $F_X$ and $F_Y$ lead to the overall stable behaviour the bidirectional networks are called as Bidirectional Associative Memories (BAM).

Further not only a Bidirectional network leads to BAM also a unidirectional network defines a BAM if M is symmetric i.e. $M = M^T$. We in our analysis mainly use BAM which are bidirectional networks. However we may also use unidirectional BAM chiefly depending on the problems under investigations. We briefly describe the BAM model more technically and mathematically.

An additive activation model is defined by a system of n + p coupled first order differential equations that inter connects the fields $F_X$ and $F_Y$ through the constant synaptic matrices M and N.

$$x_i = -A_i x_i + \sum_{j=1}^{p} S_j(y_j) n_{ji} + I_i \qquad (3.2.1.1)$$

$$y_i = -A_j y_j + \sum_{i=1}^{n} S_i(x_i) m_{ij} + J_j \qquad (3.2.1.2)$$

$S_i(x_i)$ and $S_j(y_j)$ denote respectively the signal function of the $i^{th}$ neuron in the field $F_X$ and the signal function of the $j^{th}$ neuron in the field $F_Y$.

Discrete additive activation models correspond to neurons with threshold signal functions.

The neurons can assume only two values ON and OFF. ON represents the signal +1, OFF represents 0 or – 1 (– 1 when the representation is bipolar). Additive bivalent models describe asynchronous and stochastic behaviour.

At each moment each neuron can randomly decide whether to change state or whether to emit a new signal given its current activation. The Bidirectional Associative Memory or BAM is a non adaptive additive bivalent neural network. In neural literature the discrete version of the equation (3.2.1.1) and (3.2.1.2) are often referred to as BAMs.

A discrete additive BAM with threshold signal functions arbitrary thresholds inputs an arbitrary but a constant synaptic connection matrix M and discrete time steps K are defined by the equations

$$x_i^{k+1} = \sum_{j=1}^{p} S_j(y_j^k) m_{ij} + I_i \qquad (3.2.1.3)$$

$$y_j^{k+1} = \sum_{i=1}^{n} S_i\left(x_i^k\right) m_{ij} + J_j \qquad (3.2.1.4)$$

where $m_{ij} \in M$ and $S_i$ and $S_j$ are signal functions. They represent binary or bipolar threshold functions. For arbitrary real valued thresholds $U = (U_1, ..., U_n)$ for $F_X$ neurons and $V = (V_1, ..., V_P)$ for $F_Y$ neurons the threshold binary signal functions corresponds to

$$S_i(x_i^k) = \begin{cases} 1 \quad \text{if} \quad x_i^k > U_i \\ S_i(x_i^{k-1}) \text{ if } x_i^k = U_i \\ 0 \quad \text{if} \quad x_i^k < U_i \end{cases} \qquad (3.2.1.5)$$

and

$$S_j(x_j^k) = \begin{cases} 1 \quad \text{if} \quad y_j^k > V_j \\ S_j(y_j^{k-1}) \text{ if } y_j^k = V_j \\ 0 \quad \text{if} \quad y_j^k < V_j \end{cases} \qquad (3.2.1.6)$$

The bipolar version of these equations yield the signal value –1 when $x_i < U_i$ or when $y_j < V_j$. The bivalent signal functions allow us to model complex asynchronous state change patterns. At any moment different neurons can decide whether to compare their activation to their threshold. At each moment any of the 2n subsets of $F_X$ neurons or 2p subsets of the $F_Y$ neurons can decide to change state. Each neuron may randomly decide whether to check the threshold conditions in the equations (3.2.1.5) and (3.2.1.6). At each moment each neuron defines a random variable that can assume the value ON(+1) or OFF(0 or -1). The network is often assumed to be deterministic and state changes are synchronous i.e. an entire field of neurons is updated at a time. In case of simple asynchrony only one neuron makes a state change decision at a time. When the subsets represent the entire fields $F_X$ and $F_Y$ synchronous state change results.

In a real life problem the entries of the constant synaptic matrix M depends upon the investigator's feelings. The synaptic matrix is given a weightage according to their feelings. If $x \in F_X$ and $y \in F_Y$ the forward projections from $F_X$ to $F_Y$ is defined by the matrix M. $\{F(x_i, y_j)\} = (m_{ij}) = M$, $1 \leq i \leq n$, $1 \leq j \leq p$.

The backward projections is defined by the matrix $M^T$. $\{F(y_i, x_i)\} = (m_{ji}) = M^T$, $1 \leq i \leq n$, $1 \leq j \leq p$. It is not always true that the backward projections from $F_Y$ to $F_X$ is defined by the matrix $M^T$.

Now we just recollect the notion of bidirectional stability. All BAM state changes lead to fixed point stability. The property holds for synchronous as well as asynchronous state changes. A BAM system ($F_X$, $F_Y$, M) is bidirectionally stable if all inputs converge to fixed point equilibria. Bidirectional stability is a dynamic equilibrium. The same signal information flows back and forth in a bidirectional fixed point. Let us suppose that A denotes a binary n-vector and B denotes a binary p-vector. Let A be the initial input to the BAM system. Then the BAM equilibrates to a bidirectional fixed point ($A_f$, $B_f$) as

$$\begin{array}{lllll}
A & \rightarrow & M & \rightarrow & B \\
A' & \leftarrow & M^T & \leftarrow & B \\
A' & \rightarrow & M & \rightarrow & B' \\
A'' & \leftarrow & M^T & \leftarrow & B' \text{ etc.} \\
A_f & \rightarrow & M & \rightarrow & B_f \\
A_f & \leftarrow & M^T & \leftarrow & B_f \text{ etc.}
\end{array}$$

where A', A'', ... and B', B'', ... represents intermediate or transient signal state vectors between respectively A and $A_f$ and B and $B_f$. The fixed point of a Bidirectional system is time dependent.

The fixed point for the initial input vectors can be attained at different times. Based on the synaptic matrix M which is developed by the investigators feelings the time at which bidirectional stability is attained also varies accordingly.

### 3.2.2 Description and definition of Super Bidirectional Associative Memories model

In this section we for the first time define four types of fuzzy super BAM models. The row vector fuzzy super BAM (SDBAM) model, column vector fuzzy super BAM (SRBAM) model, diagonal fuzzy super BAM (SDSBAM) model and fuzzy super BAM model. This model will be known as the super BAM (SBAM) model. Clearly this is the most generalization of SDBAM, SRBAM and SDSBAM.

**DEFINITION 3.2.2.1:** *Suppose a set of n experts choose to work with a problem using a BAM model in which they all agree upon the same number of attributes from the space $F_x$ which will form the rows of the dynamical system formed by this multi expert BAM. Now n distinct sets of attributes are given from the space $F_y$ which forms a super row vector and they form the columns of the BAM model.*

*Suppose all the n experts agree to work with the same set of t-attributes say ($x_1\ x_2 \ldots x_t$) which forms the rows of the synaptic connection matrix M. Suppose the first expert works with the $p_1$ set of attributes given by ($y_1^1\ y_2^1 \ldots y_{p_1}^1$), the second expert with $p_2$ set of attributes given by ($y_1^1\ y_2^1 \ldots y_{p_1}^1$) and so on. Let the $i^{th}$ expert with $p_i$ set of attributes given by ($y_1^i\ y_2^i \ldots y_{p_i}^i$) for i = 1, 2, …, n. Thus the new BAM model will have its elements from $F_y$ where any element in $F_y$ will be a super row vector, T = ($y_1^1\ y_2^1 \ldots y_{p_1}^1 \mid y_1^2\ y_2^2 \ldots y_{p_2}^2 \mid \ldots \mid y_1^n\ y_2^n \ldots y_{p_n}^n$). Now the synaptic projection matrix associated with this new BAM model is a special row supervector $M_r$ given by*

$$M_r = \begin{array}{c} \\ x_1 \\ x_2 \\ \vdots \\ x_t \end{array} \begin{array}{c} \begin{array}{ccc|ccc|c|ccc} y_1^1 & y_2^1 \ldots & y_{p_1}^1 & y_1^2 & y_2^2 \ldots & y_{p_2}^2 & \ldots & y_1^n & y_2^n \ldots & y_{p_n}^n \end{array} \\ \left[ \begin{array}{c|c|c|c} \quad\quad & \quad\quad & \quad & \quad\quad \\ & & & \\ & & & \\ & & & \end{array} \right] \end{array}$$

*Here the elements /attributes from $F_x$ is a simple row rector where as the elements from $F_y$ is a super row vector.*

*We call this model to be a multi expert Special Domain Supervector BAM (SDBAM) model and the associated matrix is a special row vector matrix denoted by $M_r$. Let $X = (x_1\ x_2 \ldots x_t) \in F_x$. $Y = [y_1^1\ y_2^1 \ldots y_{p_1}^1 \mid y_1^2\ y_2^2 \ldots y_{p_2}^2 \mid \ldots \mid y_1^n\ y_2^n \ldots y_{p_n}^n] \in F_y$. If $X = (x_1\ x_2 \ldots x_t) \in F_x$ is the state vector given by the expert we find*

$$XM_r = Y_1 \in F_y$$
$$YM_r^T = X_1 \in F_x \ldots$$

*and so on. This procedure is continued until a equilibrium is arrived. Similarly if the expert chooses to work with $Y = [y_1^1\ y_2^1 \ldots y_{p_1}^1 \mid y_1^2\ y_2^2 \ldots y_{p_2}^2 \mid \ldots \mid y_1^n\ y_2^n \ldots y_{p_n}^n] \in F_y$ then we find the resultant by finding.*

*$YM_r^T \hookrightarrow X$, then find $XM_r$ and proceed on till the system arrives at an equilibrium state. This model will serve the purpose when row vectors from $F_x$ are a simple row vectors and row vectors from $F_y$ are super row vectors.*

Now we proceed on to define the second model which has super row vectors from $F_x$ and just simple row vectors from $F_y$.

**DEFINITION 3.2.2.2:** *Suppose we have a problem in which all m experts want to work using a BAM model. If they agree to work having the simple vectors from $F_y$ i.e., for the columns of the synaptic connection matrix i.e. there is no perpendicular partition of their related models matrix.*

*The rows are partitioned horizontally in this synaptic connection matrix i.e., the m experts have distinct sets of attributes taken from the space $F_x$ i.e. elements of $F_x$ are super row vectors. The resulting synaptic connection matrix $M_c$ is a special super column matrix. Let the $1^{st}$ expert have the set of row attributes to be $(x_1^1\ x_2^1 \ldots x_{q_1}^1)$, the $2^{nd}$ expert have the set of row attributes given by $(x_1^2\ x_2^2 \ldots x_{q_2}^2)$ and so on. Let the $i^{th}$ expert have the related row attributes as $(x_1^i\ x_2^i \ldots x_{q_i}^i)$; $i = 1, 2, \ldots, m$. Let the column vector given by all them is $[y_1 \ldots y_n]$. The related super synaptic connection matrix $M_c$ =*

$$
\begin{array}{c}
\\
x_1^1 \\ x_2^1 \\ \vdots \\ x_{q_1}^1 \\ x_1^2 \\ x_2^2 \\ \vdots \\ x_{q_2}^2 \\ \\ \vdots \\ \\ x_1^m \\ x_2^m \\ \vdots \\ x_{q_m}^m
\end{array}
\begin{array}{c}
\begin{array}{cccc} y_1 & y_2 & \cdots & y_n \end{array} \\
\left[\begin{array}{cccc}
 & & & \\
 & & & \\
 & & & \\
 & & & \\ \hline
 & & & \\
 & & & \\
 & & & \\
 & & & \\ \hline
 & & & \\
 & & & \\
 & & & \\ \hline
 & & & \\
 & & & \\
 & & & \\
 & & &
\end{array}\right]
\end{array}.
$$

*$M_c$ is a special super column vector / matrix.*

*Suppose an expert wishes to work with a super row vector X from $F_x$ then $X = [\, x_1^1\ x_2^1 \ldots x_{q_1}^1 \mid x_1^2\ x_2^2 \ldots x_{q_2}^2 \mid \ldots \mid x_1^m\ x_2^m \ldots x_{q_m}^m \,]$ we find $X \text{ o } M_c \hookrightarrow Y \in F_y$ , $YM_c^T = X_1 \in F_x$, we repeat the same procedure till the system attains its equilibrium i.e., a fixed point or a limit cycle.*

*This model which performs using the dynamical system $M_c$ is defined as the Special Super Range BAM (SRBAM) model.*

Next we describe the special diagonal super BAM model.

**DEFINITION 3.2.2.3:** *Suppose we have n experts to work on a specific problem and each expert wishes to work with a set of*

*row and column attributes distinct from others using a BAM model.*

*Then how to obtain a suitable integrated dynamical system using them. Let the first expert work with* $(x_1^1\ x_2^1 \ldots x_{n_1}^1)$ *attributes along the row of the related synaptic connection matrix of the related BAM and* $(y_1^1\ y_2^1 \ldots y_{p_1}^1)$ *the attributes related to the column, let the second expert give the row attributes of the synaptic connection matrix of the BAM to be* $(x_1^2\ x_2^2 \ldots x_{n_2}^2)$ *and that of the column be* $(y_1^2\ y_2^2 \ldots y_{p_2}^2)$ *and so on.*

*Let the* $i^{th}$ *expert give the row attributes of the synaptic connection matrix of the BAM to be* $(x_1^i\ x_2^i \ldots x_{n_i}^i)$ *and that of the column to be* $(y_1^i\ y_2^i \ldots y_{p_i}^i)$ *for i = 1, 2, …, n, the supermatrix described by*

$$M_D = \begin{array}{c} \\ x_1^1 \\ \vdots \\ x_{n_1}^1 \\ x_1^2 \\ \vdots \\ x_{n_2}^2 \\ \vdots \\ x_1^n \\ \vdots \\ x_{n_n}^n \end{array} \begin{array}{c} \begin{array}{cccc} y_1^1\, y_2^1 \ldots y_{p_1}^1 & y_1^2\, y_2^2 \ldots y_{p_2}^2 & \ldots & y_1^n\, y_2^n \ldots y_{p_n}^n \end{array} \\ \left[\begin{array}{c|c|c|c} A_1^1 & (0) & & (0) \\ \hline (0) & A_2^2 & & (0) \\ \hline (0) & (0) & & (0) \\ \hline (0) & (0) & & A_n^n \end{array}\right] \end{array}$$

*where* $A_i^i$ *is the synaptic connection matrix using the BAM model of the* $i^{th}$ *expert, i = 1, 2, …, n where*

$$A^i{}_1 = \begin{array}{c} \\ x_1^i \\ x_2^i \\ \vdots \\ x_{n_i}^i \end{array} \begin{array}{c} \begin{array}{cccc} y_1^i & y_2^i & \ldots & y_{p_i}^i \end{array} \\ \left[ \begin{array}{cccc} & & & \\ & & & \\ & & & \\ & & & \end{array} \right] \end{array}$$

*(0) denotes the zero matrix.*

*Thus this model has only non zero BAM synaptic connection matrices along the main diagonal described by* $M_D$*. The rest are zero.*

*The dynamical system associated with this matrix* $M_D$ *is defined to be the Special Diagonal Super BAM (SDSBAM) model.*

Next we describe the super BAM model.

**DEFINITION 3.2.2.4:** *Suppose we have mn number of experts who are interested in working with a specific problem using a BAM model; a multi expert model which will work as a single dynamical system is given by the Super BAM (SBAM) model.*

*Here a few experts have both the row and columns of the synaptic connection matrix of the BAM to be distinct. Some accept for same row attributes or vectors of the synaptic connection matrix but with different column attributes.*

*Some accept for same column attributes of the synaptic connection matrix of the BAM model but with different row attributes to find the related supermatrix associated with the super BAM model.*

*The supermatrix related with this new model will be denoted by* $M_s$ *which is described in the following.*

$$
\begin{array}{c}
\qquad\quad y_1^1\, y_2^1 \ldots y_{q_1}^1 \quad y_1^2\, y_2^2 \ldots y_{q_2}^2 \quad \ldots \quad y_1^n\, y_2^n \ldots y_{q_n}^n \\
\begin{array}{c} x_1^1 \\ \vdots \\ x_{p_1}^1 \\ x_1^2 \\ \vdots \\ x_{p_2}^2 \\ \vdots \\ x_1^m \\ \vdots \\ x_{p_m}^m \end{array}
\left[\begin{array}{c|c|c|c}
A_1^1 & A_2^1 & & A_n^1 \\ \hline
A_1^2 & A_2^2 & & A_n^2 \\ \hline
 & & & \\ \hline
A_1^m & A_2^m & & A_n^m
\end{array}\right]
\end{array}
$$

*where* $A_j^i$ *is the synaptic connection matrix of an expert who chooses to work with (* $x_1^i\ x_2^i\ \ldots\ x_{p_i}^i$ *) along the row of the BAM model and with (* $y_1^j\ y_2^j\ \ldots\ y_{q_j}^j$ *) along the column of the BAM model i.e.*

$$
A_j^i = \begin{array}{c} \\ x_1^i \\ x_2^i \\ \vdots \\ x_{p_i}^i \end{array}
\begin{array}{c} y_1^j \quad y_2^j \quad \ldots \quad y_{q_j}^j \\ \left[\begin{array}{ccc} \phantom{y_1^j} & \phantom{y_2^j} & \phantom{y_{q_j}^j} \\ & & \\ & & \\ & & \end{array}\right] \end{array}
$$

*$1 \le i \le m$ and $1 \le j \le n$. Thus for this model both the attributes from the spaces $F_x$ and $F_y$ are super row vectors given by*

$$X = [\, x_1^1\ x_2^1\ \ldots\ x_{p_1}^1 \mid x_1^2\ x_2^2\ \ldots\ x_{p_2}^2 \mid \ldots \mid x_1^m\ x_2^m\ \ldots\ x_{p_m}^m \,]$$

*in $F_x$ and*

$$Y = [\, y_1^1\ y_2^1\ \ldots\ y_{q_1}^1 \mid y_1^2\ y_2^2\ \ldots\ y_{q_2}^2 \mid \ldots \mid y_1^n\ y_2^n\ \ldots\ y_{q_n}^n \,]$$

*from the space or the neuronal field $F_y$.*

*The supermatrix $M_s$ is called the synaptic connection supermatrix associated with the multi expert super BAM model (SBAM model). Now having defined the multi expert super BAM model we proceed on to describe the functioning of the super dynamical system.*

Let $X = [x_1^1\ x_2^1 \ldots x_{p_1}^1 \mid x_1^2\ x_2^2 \ldots x_{p_2}^2 \mid \ldots \mid x_1^m\ x_2^m \ldots x_{p_m}^m] \in F_x$ be the super row vector given by the expert, its effect on the multi super dynamical system $M_s$.

$$X \text{ o } M_s \hookrightarrow Y$$
$$= [y_1^1\ y_2^1 \ldots y_{q_1}^1 \mid y_1^2\ y_2^2 \ldots y_{q_2}^2 \mid \ldots \mid y_1^n\ y_2^n \ldots y_{q_n}^n] \in F_y$$
$$Y \text{ o } M_s^T \hookrightarrow X_1 \in F_x.$$
$$X_1 \text{ o } M_s \hookrightarrow Y_1 \in F_y;$$

and so on and this procedure is repeated until the system attains a equilibrium.

## 3.3 Description of Super Fuzzy Associative Memories

In this section we for the first time introduce the notion of super fuzzy associative memories model. This section has two subsections in the first subsection we recall the definition of FAM model from [112]. In second subsection four new types of super FAM models are defined.

### 3.3.1 Introduction to FAM

In this section the notion of Fuzzy Associative Memories (FAM) is recalled from [112]. For more refer [112].

A fuzzy set is a map $\mu : X \to [0, 1]$ where X is any set called the domain and [0, 1] the range i.e., $\mu$ is thought of as a membership function i.e., to every element $x \in X$, $\mu$ assigns a membership value in the interval [0, 1]. But very few try to visualize the geometry of fuzzy sets. It is not only of interest but is meaningful to see the geometry of fuzzy sets when we discuss fuzziness. Till date researchers over looked such visualization [Kosko, 108-112], instead they have interpreted fuzzy sets as generalized indicator or membership functions; i.e., mappings μ from domain X to range [0, 1]. But functions are hard to visualize. Fuzzy theorist often picture membership functions as two-dimensional graphs with the domain X represented as a one-dimensional axis.

The geometry of fuzzy sets involves both domain X = $(x_1,\ldots, x_n)$ and the range [0, 1] of mappings $\mu : X \to [0, 1]$. The geometry of fuzzy sets aids us when we describe fuzziness, define fuzzy concepts and prove fuzzy theorems. Visualizing this geometry may by itself provide the most powerful argument for fuzziness.

An odd question reveals the geometry of fuzzy sets. What does the fuzzy power set $F(2^X)$, the set of all fuzzy subsets of X, look like? It looks like a cube, What does a fuzzy set look like? A fuzzy subsets equals the unit hyper cube $I^n = [0, 1]^n$. The fuzzy set is a point in the cube $I^n$. Vertices of the cube $I^n$ define a non-fuzzy set. Now with in the unit hyper cube $I^n = [0, 1]^n$ we are interested in a distance between points, which led to measures of size and fuzziness of a fuzzy set and more fundamentally to a measure. Thus within cube theory directly extends to the continuous case when the space X is a subset of $R^n$.

The next step is to consider mappings between fuzzy cubes. This level of abstraction provides a surprising and fruitful alternative to the prepositional and predicate calculus reasoning techniques used in artificial intelligence (AI) expert systems. It allows us to reason with sets instead of propositions. The fuzzy set framework is numerical and multidimensional. The AI framework is symbolic and is one dimensional with usually only bivalent expert rules or propositions allowed. Both frameworks can encode structured knowledge in linguistic form. But the fuzzy approach translates the structured knowledge into a flexible numerical framework and processes it in a manner that resembles neural network processing. The numerical framework also allows us to adaptively infer and modify fuzzy systems perhaps with neural or statistical techniques directly from problem domain sample data.

Between cube theory is fuzzy-systems theory. A fuzzy set defines a point in a cube. A fuzzy system defines a mapping between cubes. A fuzzy system S maps fuzzy sets to fuzzy sets. Thus a fuzzy system S is a transformation S: $I^n \to I^P$. The n-dimensional unit hyper cube $I^n$ houses all the fuzzy subsets of the domain space or input universe of discourse X = {$x_1$, …, $x_n$}. $I^p$ houses all the fuzzy subsets of the range space or output universe of discourse, Y = {$y_1$, …, $y_p$}. X and Y can also denote subsets of $R^n$ and $R^p$. Then the fuzzy power sets F ($2^X$) and F ($2^Y$) replace $I^n$ and $I^p$.

In general a fuzzy system S maps families of fuzzy sets to families of fuzzy sets thus S: $I^{n_1} \times \ldots \times I^{n_r} \to I^{p_1} \times \ldots \times I^{p_s}$ Here too we can extend the definition of a fuzzy system to allow arbitrary products or arbitrary mathematical spaces to serve as the domain or range spaces of the fuzzy sets. We shall focus on fuzzy systems S: $I^n \to I^P$ that map balls of fuzzy sets in $I^n$ to balls of fuzzy set in $I^p$. These continuous fuzzy systems behave as associative memories. The map close inputs to close outputs. We shall refer to them as Fuzzy Associative Maps or FAMs.

The simplest FAM encodes the FAM rule or association ($A_i$, $B_i$), which associates the p-dimensional fuzzy set $B_i$ with the n-dimensional fuzzy set $A_i$. These minimal FAMs essentially map one ball in $I^n$ to one ball in $I^p$. They are comparable to simple neural networks. But we need not adaptively train the minimal FAMs. As discussed below, we can directly encode structured knowledge of the form, "If traffic is heavy in this direction then keep the stop light green longer" is a Hebbian-style FAM correlation matrix. In practice we sidestep this large numerical matrix with a virtual representation scheme. In the place of the matrix the user encodes the fuzzy set association (Heavy, longer) as a single linguistic entry in a FAM bank linguistic matrix. In general a FAM system F: $I^n \to I^b$ encodes the processes in parallel a FAM bank of m FAM rules ($A_1$, $B_1$), …, ($A_m$ $B_m$). Each input A to the FAM system activates each stored FAM rule to different degree. The minimal FAM that stores ($A_i$, $B_i$) maps input A to $B_i$' a partly activated version of $B_i$. The more A resembles $A_i$, the more $B_i$' resembles $B_i$. The corresponding output fuzzy set B combines these partially activated fuzzy sets $B_1^1, B_2^1, \ldots, B_m^1$. B equals a weighted average of the partially activated sets B = $w_1 B_1^1 + ... + w_m B_m^1$ where $w_i$ reflects the credibility frequency or strength of fuzzy association ($A_i$, $B_i$). In practice we usually defuzzify the output waveform B to a single numerical value $y_j$ in Y by computing the fuzzy centroid of B with respect to the output universe of discourse Y.

More generally a FAM system encodes a bank of compound FAM rules that associate multiple output or consequent fuzzy sets $B^1{}_i$, …, $B_i^s$ with multiple input or antecedent fuzzy sets $A_i^1$, …, $A_i^r$. We can treat compound FAM rules as compound linguistic conditionals. This allows us to naturally and in many cases easily to obtain structural knowledge. We combine antecedent and consequent sets with logical conjunction, disjunction or negation. For instance, we could interpret the compound association ($A^1$, $A^2$, B), linguistically as the compound conditional "IF $X^1$ is $A^1$ AND $X^2$ is $A^2$, THEN Y is B" if the comma is the fuzzy association ($A^1$, $A^2$, B) denotes conjunction instead of say disjunction.

We specify in advance the numerical universe of discourse for fuzzy variables $X^1$, $X^2$ and Y. For each universe of discourse or fuzzy variable X, we specify an appropriate library of fuzzy set values $A_1^r$, …, $A_k^2$ Contiguous fuzzy sets in a library overlap. In principle a neural network can estimate these libraries of fuzzy sets. In practice this is usually unnecessary. The library sets represent a weighted though overlapping quantization of the input space X. They represent the fuzzy set values assumed by a fuzzy variable. A different library of fuzzy sets similarly quantizes the output space Y. Once we define the library of fuzzy sets we construct the FAM by choosing appropriate combinations of input and output fuzzy sets Adaptive techniques can make, assist or modify these choices.

An Adaptive FAM (AFAM) is a time varying FAM system. System parameters gradually change as the FAM system samples and processes data. Here we discuss how natural network algorithms can adaptively infer FAM rules from training data. In principle, learning can modify other FAM system components, such as the libraries of fuzzy sets or the FAM-rule weights $w_i$.

In the following subsection we propose and illustrate an unsupervised adaptive clustering scheme based on competitive learning to blindly generate and refine the bank of FAM rules. In some cases we can use supervised learning techniques if we have additional information to accurately generate error estimates. Thus Fuzzy Associative Memories (FAMs) are transformation. FAMs map fuzzy sets to fuzzy sets. They map unit cubes to unit cubes. In simplest case the FAM system consists of a single association. In general the FAM system consists of a bank of different FAM association. Each association corresponds to a different numerical FAM matrix or a different entry in a linguistic FAM-bank matrix. We do not combine these matrices as we combine or superimpose neural-network associative memory matrices. We store the matrices and access them in parallel. We begin with single association FAMs. We proceed on to adopt this model to the problem.

## 3.3.2 Super Fuzzy Associative Memories models

In this section four new types of super fuzzy associative memories models are introduced. These models are multi expert models which can simultaneously work with many experts using FAM models.

In this section we for the first time construct four types of super fuzzy associative memories.

**DEFINITION 3.3.2.1:** *We have a problem P on which n experts wishes to work using a FAM model which can work as a single unit multi expert system. Suppose all the n-experts agree to work with the same set of attributes from the domain space and they want to work with different and distinct sets of attributes from the range space. Suppose all the n experts wish to work with the domain attributes ($x_1$ $x_2$ … $x_t$) from the cube* $I^t = \underbrace{[0, 1] \times \ldots \times [0, 1]}_{t-times}$. *Let the first expert work with the range attributes ($y_1^1$ $y_2^1$ … $y_{p_1}^1$) and the second expert works with the* *range attributes ( $y_1^2\ y_2^2 \ldots y_{p_2}^2$ ) and so on. Thus the $i^{th}$ expert works with the range attributes ( $y_1^i\ y_2^i \ldots y_{p_i}^i$ ), i = 1, 2, …, n. Thus the range attributes*

$$Y = (\, y_1^1\ y_2^1 \ldots y_{p_1}^1 \mid y_1^2\ y_2^2 \ldots y_{p_2}^2 \mid \ldots \mid y_1^n\ y_2^n \ldots y_{p_n}^n\,)$$

*are taken from the cube* $I^{p_1+p_2+\ldots+p_n} = \underbrace{[0, 1]\times[0, 1]\times\ldots\times[0, 1]}_{p_1+p_2+\ldots+p_n\ times}$.

*This we see the range attributes are super row fuzzy vectors.*

*Now the matrix which can serve as the dynamical systems for this FAM model is given by $F_R$.*

$$\begin{array}{c} \\ x_1 \\ \vdots \\ x_t \end{array} \begin{array}{c} \begin{array}{cccc} y_1^1\, y_2^1 \ldots y_{p_1}^1 & y_1^2\, y_2^2 \ldots y_{p_2}^2 & \ldots & y_1^n\, y_2^n \ldots y_{p_n}^n \end{array} \\ \left[\begin{array}{c|c|c|c} & & & \\ A_1^1 & A_2^1 & & A_n^1 \\ & & & \end{array}\right] \end{array}$$

*Clearly $F_R$ is a special super row fuzzy vector. Thus F: $I^t \hookrightarrow I^{p_1+p_2+\ldots+p_n}$. Suppose using an experts opinion we have a fit vector, A = ($a_1$, $a_2$, …, $a_t$) ; $a_i \in\{0,1\}$, then A o $F_R$ = max min ($a_i$, $f_{ij}$); $a_i \in$ A and $f_{ij} \in F_R$. Let A o $F_R$ = B = ($b_j$), then $F_R$ o B = max min ($f_{ij}$, $b^i_j$) and so on, till we arrive at a fixed point or a limit cycle. The resultant fit vectors give the solution. This $F_R$ gives the dynamical system of the new model which we call as the Fuzzy Special Super Row vector FAM model (SRFAM model).*

Next we proceed on to describe how the FAM functions when we work with varying domain space and fixed range space attributes.

**DEFINITION 3.3.2.2:** *Suppose we have n experts working on a problem and they agree upon to work with the same range attributes and wish to work with distinct domain attributes using a FAM model. We built a new FAM model called the* *special super column fuzzy vector FRM model (SCFAM) and its related matrix is denoted by $F_c$. The fit vectors of the domain space are simple fit vectors where as the fit vectors of the range space are super row fit vectors.*

*Now we describe the special column fuzzy vector FAM, $F_c$ by the following matrix. The column attributes of the super fuzzy dynamical system $F_c$ are given by*

$$(y_1\ y_2 \ldots y_s) \in I^s = \underbrace{[0, 1] \times \ldots \times [0, 1]}_{s-times}.$$

*The row attributes of the first expert is given by ($x_1^1, x_2^1, \ldots, x_{p_1}^1$), the row attributes of the second expert is given by ($x_1^2, x_2^2, \ldots, x_{p_2}^2$). Thus the row attributes of the $i^{th}$ expert is given by ($x_1^i\ x_2^i \ldots x_{p_i}^i$), $i = 1, 2, \ldots, n$.*

*We have*

$$F_c = \begin{array}{c} \\ x_1^1 \\ \vdots \\ x_{p_1}^1 \\ x_1^2 \\ \vdots \\ x_{p_2}^2 \\ \vdots \\ x_1^n \\ \vdots \\ x_{p_n}^n \end{array} \begin{array}{c} y_1\ y_2\ \cdots\ y_s \\ \left[\begin{array}{c} A_1 \\ \hline A_2 \\ \hline \vdots \\ \hline A_n \end{array}\right] \end{array}$$

*to be a special super column fuzzy vector / matrix, where*

$$A_i = \begin{array}{c} \\ x_1^i \\ x_2^i \\ \vdots \\ x_{p_i}^i \end{array} \begin{array}{c} y_1 \quad y_2 \quad \cdots \quad y_s \\ \left[ \begin{array}{ccc} & & \\ & & \\ & & \\ & & \end{array} \right] \end{array}$$

*i = 1, 2, …, n is a fuzzy $p_i \times s$ matrix. Suppose the expert wishes to work with a fit vector (say) [ $x_1^1\ x_2^1 \ldots x_{p_1}^1 \mid x_1^2\ x_2^2 \ldots x_{p_2}^2 \mid \ldots \mid x_1^n\ x_2^n \ldots x_{p_n}^n$ ]. Then X o $F_c$= B where B is a simple row vector we find $\mathrm{F_c}\ o\ B = \underbrace{[0,\ 1] \times \ldots \times [0,\ 1]}_{p_1 + p_2 + \ldots + p_n \text{ times}} = I_{p_1 + \ldots + p_n}$ ; we proceed on to work till we arrive at an equilibrium state of the system.*

Next we proceed on to define FAM when n expert give opinion having distinct set of domain attributes and distinct set of range attributes.

**DEFINITION 3.3.2.3:** *Let n experts give opinion on a problem P and wish to use a FAM model, to put this data as an integrated multi expert system. Let the first expert give his/her attributes along the column as ( $y_1^1\ y_2^1 \ldots y_{q_1}^1$ ) and those attributes along the row as ( $x_1^1\ x_2^1 \ldots x_{p_1}^1$ ).*

*Let ( $y_1^2\ y_2^2 \ldots y_{q_2}^2$ ) and ( $x_1^2\ x_2^2 \ldots x_{p_2}^2$ ) be the column and row attributes respectively given by the second expert and so on. Thus any $i^{th}$ expert gives the row and column attributes as ( $x_1^i\ x_2^i \ldots x_{p_i}^i$ ) and ( $y_1^i\ y_2^i \ldots y_{q_i}^i$ ) respectively, i =1, 2, 3, …, n.*

*So for any $i^{th}$ expert the associated matrix of the FAM would be denoted by $A^i$ where*

$$A^i = \begin{array}{c} \\ x_1^i \\ x_2^i \\ \vdots \\ x_{p_i}^i \end{array} \begin{array}{c} y_1^i \quad y_2^i \quad \cdots \quad y_{q_i}^i \\ \left[ \begin{array}{ccc} & & \\ & & \\ & & \\ & & \end{array} \right] \end{array}$$

*Now form the multi expert FAM model using these n FAM matrices $A^1, A^2, \ldots, A^n$ and get the multi expert system which is denoted by*

$$F_D = \begin{array}{c} \\ x_1^1 \\ \vdots \\ x_{p_1}^1 \\ x_1^2 \\ \vdots \\ x_{p_2}^2 \\ \vdots \\ x_1^n \\ \vdots \\ x_{p_n}^n \end{array} \begin{array}{c} \begin{array}{cccc} y_1^1\, y_2^1 \ldots y_{q_1}^1 & y_1^2\, y_2^2 \ldots y_{q_2}^2 & \ldots & y_1^n\, y_2^n \ldots y_{q_n}^n \end{array} \\ \left[\begin{array}{c|c|c|c} A^1 & (0) & & (0) \\ \hline (0) & A^2 & & (0) \\ \hline (0) & (0) & & (0) \\ \hline (0) & (0) & & A^n \end{array}\right] \end{array}.$$

*This fuzzy supermatrix $F_D$ will be known as the diagonal fuzzy supermatrix of the FAM and the multi expert system which makes use of this diagonal fuzzy supermatrix $F_D$ will be known as the Fuzzy Super Diagonal FAM (SDFAM) model. Now the related fit fuzzy supervectors of this model $F_x$ and $F_y$ are fuzzy super row vectors given by X= $(x_1^1\ x_2^1 \ldots x_{p_1}^1 \mid x_1^2\ x_2^2 \ldots x_{p_2}^2 \mid \ldots \mid x_1^n\ x_2^n \ldots x_{p_n}^n) \in F_x$ and Y = $(y_1^1\ y_2^1 \ldots y_{q_1}^1 \mid y_1^2\ y_2^2 \ldots y_{q_2}^2 \mid \ldots \mid y_1^n\ y_2^n \ldots y_{q_n}^n) \in F_y$.*

*Now this new FAM model functions in the following way.*

*Suppose the expert wishes to work with the fuzzy super state fit vector X = $(x_1^1\ x_2^1 \ldots x_{p_1}^1 \mid x_1^2\ x_2^2 \ldots x_{p_2}^2 \mid \ldots \mid x_1^n\ x_2^n \ldots x_{p_n}^n)$ then Y= $(y_1^1\ y_2^1 \ldots y_{q_1}^1 \mid y_1^2\ y_2^2 \ldots y_{q_2}^2 \mid \ldots \mid y_1^n\ y_2^n \ldots y_{q_n}^n) \in F_y$.*

*Now $F_D$ o $Y = X_1 \in F_x$ and $X_i$ o $F_D = Y_1 \in F_y$ and so on.*

*This procedure is repeated until the system equilibrium is reached.*

**DEFINITION 3.3.2.4:** *Let us suppose we have a problem for which mn experts want to give their opinion. Here some experts give distinct opinion both for the row attributes and column attributes. Some experts concur on the row attributes but give different column attributes and a few others have the same set of row attributes but have a different set of column attributes. All of them concur to work using the FAM model. To find a multi expert FAM model which can tackle and give solution to the problem simultaneously.*

*To this end we make use of a super fuzzy matrix $F_s$ which is described in the following. Let the mn experts give their domain and column attributes as follows. The first expert works with the domain attributes as $(x_1^1\ x_2^1 \ldots x_{p_1}^1)$ and column attributes as $(y_1^1\ y_2^1 \ldots y_{q_1}^1)$. The second expert works with the same domain attributes viz $(x_1^1\ x_2^1 \ldots x_{p_1}^1)$ and but column attributes as $(y_1^2\ y_2^2 \ldots y_{q_2}^2)$. The $i^{th}$ expert, $1 \leq i \leq n$ works with $(x_1^1\ x_2^1 \ldots x_{p_1}^1)$ as the domain attributes and $(y_1^i\ y_2^i \ldots y_{q_i}^i)$ as the column attributes.*

*The $(n + 1)^{th}$ experts works with the new set of domain attributes $(x_1^2\ x_2^2 \ldots x_{p_2}^2)$ but with the same set of column attributes viz. $(y_1^1\ y_2^1 \ldots y_{q_1}^1)$. Now the $n + j^{th}$ expert works with using $(x_1^2\ x_2^2 \ldots x_{p_2}^2)$ as the domain attribute and $(y_1^i\ y_2^i \ldots y_{q_i}^i)$ as the column attribute $1 \leq j \leq n$. The $(2n + 1)^{th}$ expert works with $(x_1^3\ x_2^3 \ldots x_{p_3}^3)$ as the row attribute and $(y_1^1\ y_2^1 \ldots y_{q_1}^1)$ as the column attribute.*

*Thus any $(2n + k)^{th}$ expert uses $(x_1^2\ x_2^2 \ldots x_{p_2}^2)$ to be the row attribute and $(y_1^k\ y_2^k \ldots y_{q_k}^k)$ to be the column attribute $1 \leq k \leq n$. Thus any $(t_n + r)^{th}$ expert works with $(x_1^t\ x_2^t \ldots x_{p_t}^t)$ as the row* *attribute (1 ≤ t ≤ m) and (* $y_1^r\ y_2^r \ldots y_{q_r}^r$ *) as the column attribute 1 ≤ r ≤ n.*

*Now as 1 ≤ t ≤ m and 1 ≤ r ≤ n we get the FAM matrices of all the mn experts which is given by the supermatrix* $F_s$*.*

$$F_s = \begin{array}{c} \\ x_1^1 \\ \vdots \\ x_{p_1}^1 \\ x_1^2 \\ \vdots \\ x_{p_2}^2 \\ \vdots \\ x_1^m \\ \vdots \\ x_{p_m}^m \end{array} \begin{array}{c|c|c|c} y_1^1\, y_2^1 \ldots y_{q_1}^1 & y_1^2\, y_2^2 \ldots y_{q_2}^2 & \ldots & y_1^n\, y_2^n \ldots y_{q_n}^n \\ \hline A_1^1 & A_2^1 & & A_n^1 \\ \hline A_1^2 & A_2^2 & & A_n^2 \\ \hline & & & \\ \hline A_1^m & A_2^m & & A_n^m \end{array}$$

*where* $\mathrm{A}_j^i$ *is a fuzzy matrix associated with* $ij^{th}$ *expert*

$$\mathrm{A}_j^i = \begin{array}{c} \\ x_1^i \\ x_2^i \\ \vdots \\ x_{p_i}^i \end{array} \begin{array}{cccc} y_1^j & y_2^j & \ldots & y_{q_j}^j \\ \left[ \right. & & & \left. \right] \end{array}$$

*$1 \le i \le m$ and $1 \le j \le n$. This model is known as the multi expert fuzzy Super FAM (SFAM) model. The fit vectors associated with them are super row vectors from $F_x$ and $F_y$.*

*The fit super row vector X from $F_x$ is*

$$X = [\, x_1^1 \; x_2^1 \ldots x_{p_1}^1 \mid x_1^2 \; x_2^2 \ldots x_{p_2}^2 \mid \ldots \mid x_1^m \; x_2^m \ldots x_{p_m}^m \,]$$

*and the*

$$X \in F_x = \mathrm{I}^{p_1+p_2+\ldots+p_m} = \underbrace{[0, 1]\times[0, 1]\times\ldots\times[0, 1]}_{p_1+p_2+\ldots+p_m \text{ times}}.$$

*The fit super row vector Y from $F_y$ is*

$$Y = (\, y_1^1 \; y_2^1 \ldots y_{q_1}^1 \mid y_1^2 \; y_2^2 \ldots y_{q_2}^2 \mid \ldots \mid y_1^n \; y_2^n \ldots y_{q_n}^n \,);$$

$$Y \in F_y = \mathrm{I}^{q_1+q_2+\ldots+q_n} = \underbrace{[0, 1]\times[0, 1]\times\ldots\times[0, 1]}_{q_1+q_2+\ldots+q_m \text{ times}}.$$

*Thus if*

$$X = [\, x_1^1 \; x_2^1 \ldots x_{p_1}^1 \mid x_1^2 \; x_2^2 \ldots x_{p_2}^2 \mid \ldots \mid x_1^m \; x_2^m \ldots x_{p_m}^m \,]$$

*is the fit vector given by an expert; its effect on $F_s$ is given by $X$ o $F_s = Y \in F_y$; now $F_s$ o $Y = X^1 \in F_x$ then find the effect of $X^1$ on $F_s$ ; $X^1$ o $F_s = Y^1 \in F_y$ and so on.*

*We repeat this procedure until we arrive at a equilibrium state of the system.*

## 3.4 Illustration of Super Fuzzy Models

In this section we give illustrations of super fuzzy models. This section has 13 subsections the first twelve subsections give illustrations of the super fuzzy models defined in section 3.3. However the final section gives the uses of these super fuzzy models.

### 3.4.1 Illustration of the super row FRM model using super row vectors.

In this section we give the working and illustration of the super row FRM model using super row vectors.

In this model we wish to analyze the problem, the cause of dropouts in schools. Here we have a multiset of attributes given by multi experts. The assumption made in this model is that the attributes related with the domain space is the same for all experts as well as all the multi set of attributes for some fixed number of concepts or attributes are taken along the domain space. But the range space has multi set of attributes which is the different category of schools taken in this problem for investigation.

The row super FRM model is described in the following. Suppose we have some six attributes say $D_1$, $D_2$, …, $D_6$ associated with the school dropouts from school education in primary and secondary level. Suppose we are also interested in the study of the schools and their surroundings which contribute or influence to the school dropouts. We consider say 5 types of schools say

$S_1$ – rural corporation school,
$S_2$ – rural missionary run school,
$S_3$ – semi urban government aided school,
$S_4$ – city schools run by government and
$S_5$ – private owned posh city schools.

Under these five schools $S_1$, $S_2$, $S_3$, $S_4$ and $S_5$ we will have various attributes say $R_1^1$, …, $R_{t_1}^1$, related with the school $S_1$, $R_1^2$, …, $R_{t_2}^2$ be the attributes related under the school $S_2$ and so on i.e., with the school $S_i$ we have the attributes $R_1^i$, …, $R_{t_i}^i$ ; $1 \leq i \leq 5$.

Now how to model this not only multi experts system but a multi model, we model this using super row matrix. We take the domain attributes along the rows and the column attributes are attributes related with $S_1$, $S_2$, …, $S_5$. Then the model is made to function as follows.

The related super row matrix R of the super domain constant i.e., DSFRM model is as follows.

$$
\begin{array}{c}
 \\ D_1 \\ D_2 \\ D_3 \\ D_4 \\ D_5 \\ D_6
\end{array}
\begin{array}{c}
\qquad S_1 \qquad\qquad\qquad S_2 \qquad\qquad \dots \qquad\qquad S_5 \\
\begin{array}{ccc|ccc|c|ccc}
R_1^1 & R_2^1 \dots & R_{t_1}^1 & R_1^2 & R_2^2 \dots & R_{t_2}^2 & \dots & R_1^5 & R_2^5 \dots & R_{t_5}^5 \\
d_{11}^1 & d_{12}^1 \dots & d_{1t_1}^1 & d_{11}^2 & d_{12}^2 \dots & d_{1t_2}^2 & \dots & d_{11}^5 & d_{12}^5 \dots & d_{1t_5}^5 \\
d_{21}^1 & d_{22}^1 \dots & d_{2t_1}^1 & d_{21}^2 & d_{22}^2 \dots & d_{2t_2}^2 & & d_{21}^5 & d_{22}^5 \dots & d_{2t_5}^5 \\
 & & & & & & & & & \\
\vdots & \vdots & \vdots & \vdots & \vdots & \vdots & & \vdots & \vdots & \vdots \\
 & & & & & & & & & \\
d_{61}^1 & d_{62}^1 \dots & d_{6t_1}^1 & d_{61}^2 & d_{62}^2 \dots & d_{6t_2}^2 & \dots & d_{61}^5 & d_{62}^5 \dots & d_{6t_5}^5
\end{array}
\end{array}
$$

It may so happen that $t_i = t_j$ for some i and j, $1 \leq i, j \leq 5$.

Now we describe how the model works. Suppose the expert has given the related super row FRM matrix model which in this case is a $6 \times (t_1 + t_2 + \dots + t_5)$ super row matrix.

Suppose an expert wants to find the effect of a state vector X = (1 0 0 0 0 0) on the DSFRM model R, which will result in an hidden pattern.

$$
\begin{aligned}
X \text{ o } R &= (y_1^1\ y_2^1 \dots y_{q_1}^1 \mid y_1^2\ y_2^2 \dots y_{q_2}^2 \mid \dots \mid y_1^5\ y_2^5 \dots y_{q_5}^5) \\
&= Y,
\end{aligned}
$$

after updating and thresholding we get $XR \hookrightarrow Y$ (Here updating is not required as vector from domain space is taken) now we find $YR^T = (x_1\ x_2\ x_3\ x_4\ x_5\ x_6)$ after updating and thresholding we get $YR^T \hookrightarrow X_1$.

Now find $X_1R \hookrightarrow Y_1$ then find $Y_1$ o $R^T$ repeat this procedure until we land in a super fixed point or a super limit cycle. The fixed point or a limit cycle is certain as we have taken only elements from the set {0, 1}. Thus we are guaranteed of a fixed point or a limit cycle. Thus when we have both multi set of attributes and multi experts the DSFRM is best suited to give the hidden pattern of any desired state vector.

Now we proceed on to study the super FRM model in which multi set of attributes is to be analyzed from the domain space and the range space attributes remain the same.

## 3.4.2 Illustration of RSFRM model

We illustrate this model by the following problem

Suppose we have to study reservation among the OBC/SC/ST in the higher education. We may have a fixed set of attributes associated with each class of students so the variable will only be the class of students viz, upper caste students, B.C. students, MBC students, SC- students, ST students and students belonging to minority religions like Christianity, Jainism or Islam. Thus now the upper caste students may spell out the attributes related with them like we are very less percentage or so. Likewise the OBC may have attributes like poverty and failure of agriculture and so on. The SC/ST's may have attributes like discrimination, untouchability, instability, ill-treatment and so on.

How to connect or interrelate all these concepts and pose it as an integrated problem?

This is precisely done by the RSFRM (Range constant Super FRM model) or to be more specific super column vector FRM model i.e., RSFRM model. In this situation we make use of the super column vector / matrix. Suppose $R_1, \ldots, R_n$ are the n fixed set of attributes related with the problem for each and every class of students. Let

$S_1$ – denote the class of SC students,
$S_2$ – denote the class of ST students,
$S_3$ – class of students belonging to minority religions,
$S_4$ – the MBC class of students,
$S_5$ – the class of OBC students and
$S_6$ – the class of upper caste students.

Now these six classes of students may have a collection of attributes say the class $S_1$ has $S_1^1, \ldots, S_{m_1}^1$, related to the claim for reservation, $S_2$ has $S_1^2, \ldots, S_{m_2}^2$ sets of attributes and so on. Thus the class $S_i$ has $S_1^i, \ldots, S_{m_i}^i$, sets of attributes $1 \le i \le 6$. Now the RSFRM model matrix $R_C$ takes the following special form. The supermatrix related with this model is a super column vector.

$$\begin{array}{c} \\ S_1^1 \\ S_2^1 \\ \vdots \\ S_{m_1}^1 \\ S_1^2 \\ S_2^2 \\ \vdots \\ S_{m_2}^2 \\ \\ \vdots \\ \\ S_1^6 \\ S_2^6 \\ \vdots \\ S_{m_6}^6 \end{array}
\begin{array}{c} \begin{array}{cccc} R_1 & R_2 & \ldots & R_m \end{array} \\
\left[\begin{array}{cccc}
S_{11}^1 & S_{12}^1 & \ldots & S_{1n}^1 \\
S_{21}^1 & S_{22}^1 & \ldots & S_{2n}^1 \\
\vdots & \vdots & & \vdots \\
S_{m_{11}}^1 & S_{m_{12}}^1 & \ldots & S_{m_{1n}}^1 \\ \hline
S_{11}^2 & S_{12}^2 & \ldots & S_{1n}^2 \\
S_{21}^2 & S_{22}^2 & \ldots & S_{2n}^2 \\
\vdots & \vdots & & \vdots \\
S_{m_{21}}^2 & S_{m_{22}}^2 & \ldots & S_{m_{2n}}^2 \\ \hline
 & & & \\
 & & & \\
 & & & \\ \hline
S_{11}^6 & S_{12}^6 & \ldots & S_{1n}^6 \\
S_{21}^6 & S_{22}^6 & & S_{2n}^6 \\
\vdots & \vdots & & \vdots \\
S_{m_{61}}^6 & S_{m_{62}}^6 & \ldots & S_{m_{6n}}^6
\end{array}\right] \end{array} = R_c .$$

We see the state vectors of the domain space are super row vectors where as the state vectors of the range space are just simple row vectors. Now we using the wishes of an expert can find the hidden pattern of any desired state vector from the domain space or the range space of the super model under investigation.

Suppose

$$X = \left[x_1^1 \ldots x_{m_1}^1 \mid x_1^2 \ldots x_{m_2}^2 \mid \ldots \mid x_1^6 \ldots x_{m_6}^6\right]$$

be the given state supervector for which the hidden pattern is to be determined.

X o $R_c$ = ($r_1$, …, $r_n$) we need to only threshold this state vector let

$X \text{ o } R_c = (r_1, \ldots, r_n) \hookrightarrow (r_1^1 \ldots r_n^1)$
$= Y$

('$\hookrightarrow$' denotes the vector has been updated and thresholded)

Now

$$Y \text{ o } R_c^t = \left[t_1^1 \ldots t_{m_1}^1 \mid t_1^2 \ldots t_{m_2}^2 \mid \ldots \mid t_1^6 \ldots t_{m_6}^6\right]$$
$$= X_1 \text{ (say)}.$$

$X_1$ o $R_C \hookrightarrow Y_1$, this process is repeated until one arrives at a super fixed point or a super limit cycle.

Thus from the super hidden pattern one can obtain the necessary conclusions.

### 3.4.3 Example of SDSFRM model

In this section we describe the third super FRM model which we call as the diagonal super FRM model i.e., SDSFRM model and the super fuzzy matrix associated with it will be known as the diagonal fuzzy supermatrix. First we illustrate a diagonal fuzzy supermatrix before we define it abstractly.

***Example 3.4.3.1:*** Let T be a super fuzzy matrix given by

| | $x_1$ | $x_2$ | $x_3$ | $x_4$ | $x_5$ | $x_6$ | $x_7$ | $x_8$ | $x_9$ | $x_{10}$ | $x_{11}$ | $x_{12}$ |
|---|---|---|---|---|---|---|---|---|---|---|---|---|
| $y_1$ | 0.1 | 0.2 | 0 | 0 | 0 | 0 | 0 | 0 | 0 | 0 | 0 | 0 |
| $y_2$ | 1 | 0.7 | 0 | 0 | 0 | 0 | 0 | 0 | 0 | 0 | 0 | 0 |
| $y_3$ | 0 | 0 | 0.1 | 0.3 | 0.4 | 0.7 | 0 | 0 | 0 | 0 | 0 | 0 |
| $y_4$ | 0 | 0 | 1 | 0 | 0.3 | 0.2 | 0 | 0 | 0 | 0 | 0 | 0 |
| $y_5$ | 0 | 0 | 0.1 | 0.5 | 0.8 | 0.9 | 0 | 0 | 0 | 0 | 0 | 0 |
| $y_6$ | 0 | 0 | 0 | 0 | 0 | 0 | 0.3 | 0.1 | 0.7 | 1 | 0.6 | 0 |
| $y_7$ | 0 | 0 | 0 | 0 | 0 | 0 | 0.5 | 0.7 | 1 | 0.8 | 0.1 | 0.9 |
| $y_8$ | 0 | 0 | 0 | 0 | 0 | 0 | 0.1 | 0.4 | 0.1 | 0.5 | 0.7 | 1 |
| $y_9$ | 0 | 0 | 0 | 0 | 0 | 0 | 0.5 | 0.1 | 0.4 | 0.6 | 1 | 0.8 |

A super fuzzy matrix of this form will be known as the diagonal super fuzzy matrix. We see all the diagonal matrices are not square or they are of same order.

***Example 3.4.3.2:*** Let S be a fuzzy supermatrix i.e.,

$$S = \begin{array}{c} \\ y_1 \\ y_2 \\ y_3 \\ y_4 \\ y_5 \\ y_6 \\ y_7 \\ y_8 \\ y_9 \\ y_{10} \\ y_{11} \\ y_{12} \\ y_{13} \\ y_{14} \end{array} \begin{array}{c} \begin{array}{ccccccccc} x_1 & x_2 & x_3 & x_4 & x_5 & x_6 & x_7 & x_8 & x_9 \end{array} \\ \left[\begin{array}{ccc|cccc|cc} 0.1 & 0.3 & 0.5 & 0 & 0 & 0 & 0 & 0 & 0 \\ 0.9 & 0.1 & 0 & 0 & 0 & 0 & 0 & 0 & 0 \\ \hline 0 & 0 & 0 & 0.3 & 0.7 & 0.1 & 0.5 & 0 & 0 \\ 0 & 0 & 0 & 0.7 & 1 & 0.3 & 0.8 & 0 & 0 \\ 0 & 0 & 0 & 0.6 & 0.5 & 0.8 & 0.4 & 0 & 0 \\ 0 & 0 & 0 & 0.4 & 0.6 & 0.1 & 0.5 & 0 & 0 \\ 0 & 0 & 0 & 0.9 & 0.5 & 0.4 & 0.6 & 0 & 0 \\ 0 & 0 & 0 & 0.3 & 0.4 & 0.8 & 0.6 & 0 & 0 \\ 0 & 0 & 0 & 0.2 & 0.1 & 0.5 & 0.1 & 0 & 0 \\ 0 & 0 & 0 & 0.1 & 0.9 & 0 & 0.3 & 0 & 0 \\ 0 & 0 & 0 & 1 & 0 & 0.2 & 0.8 & 0 & 0 \\ \hline 0 & 0 & 0 & 0 & 0 & 0 & 0 & 0.3 & 0.2 \\ 0 & 0 & 0 & 0 & 0 & 0 & 0 & 0.9 & 0.1 \\ 0 & 0 & 0 & 0 & 0 & 0 & 0 & 0.8 & 0.5 \end{array}\right] \end{array}.$$

S is a diagonal fuzzy supermatrix.

**DEFINITION 3.4.3.2:** *Let S be a some $m \times n$ fuzzy supermatrix where*

$$m = ( t_1^1 + \ldots + t_{p_1}^1 + t_1^2 + \ldots + t_{p_2}^2 + \ldots + t_1^r + \ldots + t_{p_r}^r )$$

*and*

$$n = ( r_1^1 + \ldots + r_{q_1}^1 + r_1^2 + \ldots + r_{q_2}^2 + \ldots + r_1^t + \ldots + r_{q_t}^t ).$$

*Now the diagonal elements which are rectangular or square fuzzy matrices and rest of the elements in the matrix S are zeros.*

*The diagonal elements i.e., the fuzzy matrices are given by $A_1$, $A_2$, …, $A_r$ where*

$$S = \begin{bmatrix} A_1 & (0) & \dots & (0) \\ (0) & A_2 & \dots & (0) \\ \vdots & \vdots & & \vdots \\ (0) & (0) & \dots & A_r \end{bmatrix}.$$

For instance we continue $A_1, \dots, A_r$ as follows:

$$A_1 = \begin{array}{c} \\ r_1^i \\ r_2^i \\ \vdots \\ r_{q_1}^i \end{array} \begin{array}{c} \begin{array}{cccc} t_1^1 & t_2^1 & \dots & t_{p_1}^1 \end{array} \\ \begin{bmatrix} 0.1 & 0.9 & \dots & 0.12 \\ 0.7 & 0 & \dots & 0.3 \\ \vdots & \vdots & & \vdots \\ 0.1 & 0 & \dots & 0.11 \end{bmatrix} \end{array},$$

$$A_2 = \begin{array}{c} \\ r_1^2 \\ r_2^2 \\ \vdots \\ r_{q_2}^2 \end{array} \begin{array}{c} \begin{array}{cccc} t_1^2 & t_2^2 & \dots & t_{p_2}^2 \end{array} \\ \begin{bmatrix} 0.3 & 0.7 & \dots & 0.5 \\ 0.1 & 0.3 & \dots & 0.4 \\ \vdots & \vdots & & \vdots \\ 0.9 & 0.2 & \dots & 0.41 \end{bmatrix} \end{array}, \dots,$$

$$A_r = \begin{array}{c} \\ r_1^t \\ r_2^t \\ \vdots \\ r_{q_t}^t \end{array} \begin{array}{c} \begin{array}{cccc} t_1^r & t_2^r & \dots & t_{p_r}^r \end{array} \\ \begin{bmatrix} 0.4 & 0.7 & \dots & 0.8 \\ 0.5 & 0.3 & \dots & 0.7 \\ \vdots & \vdots & & \vdots \\ 0.9 & 0.14 & \dots & 0.19 \end{bmatrix} \end{array}.$$

(0) denotes a 0 matrix of the relevant order.

The values of $A_1, \dots, A_r$ are given very arbitrarily. This matrix S is called as the diagonal fuzzy supermatrix.

We illustrate the use of this by the following model.

Suppose we have a set of n experts who want to give their opinion on a problem and if each one of them have a distinct FRM associated with it i.e., they have different set of domain attributes as well as different set of range attributes then it becomes impossible for them to use the super row FRM model or the super column FRM model.

In this case they make use of the diagonal super FRM model.

Suppose one wants to study the problem of unemployment and the social set up i.e., unemployment as a social problem with n experts. Let the first expert have $t_i$ concepts associated with domain space and $r_i$ concepts associated with the range space i = 1, 2, …, n. The related diagonal SDSFRM model $S_D$ given by 4 experts is as follows

$$\left[\begin{array}{cccc|ccc|ccccc|cccc}
.6 & .1 & .2 & .3 & 0 & 0 & 0 & 0 & 0 & 0 & 0 & 0 & 0 & 0 & 0 & 0 \\
1 & 0 & .4 & .6 & 0 & 0 & 0 & 0 & 0 & 0 & 0 & 0 & 0 & 0 & 0 & 0 \\
.8 & .9 & 1 & .7 & 0 & 0 & 0 & 0 & 0 & 0 & 0 & 0 & 0 & 0 & 0 & 0 \\
\hline
0 & 0 & 0 & 0 & .8 & .9 & 1 & 0 & 0 & 0 & 0 & 0 & 0 & 0 & 0 & 0 \\
0 & 0 & 0 & 0 & .5 & .6 & 0 & 0 & 0 & 0 & 0 & 0 & 0 & 0 & 0 & 0 \\
0 & 0 & 0 & 0 & .9 & .2 & .4 & 0 & 0 & 0 & 0 & 0 & 0 & 0 & 0 & 0 \\
0 & 0 & 0 & 0 & .1 & .7 & .8 & 0 & 0 & 0 & 0 & 0 & 0 & 0 & 0 & 0 \\
\hline
0 & 0 & 0 & 0 & 0 & 0 & 0 & 0 & .3 & .2 & .1 & .7 & 0 & 0 & 0 & 0 \\
0 & 0 & 0 & 0 & 0 & 0 & 0 & .5 & .8 & .8 & .4 & .8 & 0 & 0 & 0 & 0 \\
0 & 0 & 0 & 0 & 0 & 0 & 0 & 0 & .7 & .9 & .6 & .1 & 0 & 0 & 0 & 0 \\
\hline
0 & 0 & 0 & 0 & 0 & 0 & 0 & 0 & 0 & 0 & 0 & 0 & .3 & .4 & .1 & 1 \\
0 & 0 & 0 & 0 & 0 & 0 & 0 & 0 & 0 & 0 & 0 & 0 & 0 & .8 & .9 & 0 \\
0 & 0 & 0 & 0 & 0 & 0 & 0 & 0 & 0 & 0 & 0 & 0 & .1 & 0 & 1 & .8 \\
0 & 0 & 0 & 0 & 0 & 0 & 0 & 0 & 0 & 0 & 0 & 0 & .4 & .1 & .6 & 1 \\
0 & 0 & 0 & 0 & 0 & 0 & 0 & 0 & 0 & 0 & 0 & 0 & 0 & .2 & .3 & .3
\end{array}\right]$$

Let X be any super row vector of the domain space i.e.,

$$X = \left[t_1^1 \;\; t_2^1 \;\; t_3^1 \;\middle|\; t_1^2 \;\; t_2^2 \;\; t_3^2 \;\; t_4^2 \;\middle|\; t_1^3 \;\; t_2^3 \;\; t_3^3 \;\middle|\; t_1^4 \;\; t_2^4 \;\; t_3^4 \;\; t_4^4 \;\; t_5^4\right]$$

$t^i_{j_t} \in [0, 1]$, $1 \le i \le 4$ and $1 \le j_t \le 5$. Let Y be any super row vector of the range space i.e.,

$$Y = \left[ r^1_1 \;\; r^1_2 \;\; r^1_3 \;\; r^1_4 \,|\, r^2_1 \;\; r^2_2 \;\; r^2_3 \,|\, r^3_1 \;\; r^3_2 \;\; r^3_3 \;\; r^3_4 \;\; r^3_5 \,|\, r^4_1 \;\; r^4_2 \;\; r^4_3 \;\; r^4_4 \right],$$

$r^i_{j_t} \in [0, 1]$, $1 \le i \le 4$ and $1 \le j_t \le 5$.
Thus if

X = [1 0 1 | 0 1 1 1 | 0 0 1| 0 0 0 1 0].

Then we find X o $S_D \hookrightarrow$ Y, is obtained Y is a super row vector find Y o $S_D^T$ and so on until we get a fixed point or a limit cycle. Since the elements of the state vector is from the set {0, 1} we get the result certainly after a finite number of steps. If some experts have the domain attributes in common and some others have the range attributes in common what should be the model. We use the Super FRM (SFRM) model.

### 3.4.4 Super FRM (SFRM) model

A Super FRM (SFRM) model makes use of the fuzzy supermatrix. In this case we have all the three types of super FRM models viz. Super Domain FRM model (DSFRM), Super Range FRM (RSFRM) model and the Super Diagonal FRM (SDSFRM) models are integrated and used.

Suppose we have 6 = 3 × 2 experts of whom 2 of them choose to have both domain and range attribute to be distinct. The problem we choose to illustrate this model is the cause of migrant labourers becoming easy victims of HIV/AIDS.

The nodes / attributes taken by the first expert as the domain space.

$P^1_1$ - No education / no help by government

$P^1_2$ - Awareness program never reaches them

$P^1_3$ - No responsibility of parents to educate children.

$P_4^1$ - Girls at a very young age say even at 11 years are married

The range nodes given by the first expert

$M_1^1$ - Government help never reaches the rural poor illiterate

$M_2^1$ - Availability of cheap liquor

$M_3^1$ - Cheap availability of CSWs

The nodes given by the second expert. The domain nodes / attributes given by the 2nd expert

$P_1^2$ - Addiction to cheap liquor

$P_2^2$ - Addiction to smoke and visit of CSWs

$P_3^2$ - Very questionable living condition so to have food atleast once a day they migrate

The nodes / attributes given by the range space.

$M_1^2$ - No job opportunities in their native place

$M_2^2$ - No proper health center

$M_3^2$ - No school even for primary classes

$M_4^2$ - Acute poverty

The nodes given by the third expert who wishes to work with the domain nodes $P_1^1, P_2^1$, $P_3^1$ and $P_4^1$ and chooses the range nodes / attributes as $M_1^2$, $M_2^2$, $M_3^2$ and $M_4^2$.

The fourth expert chooses the domain nodes as $P_1^1, P_2^1$, $P_3^1$ and $P_4^1$ and range nodes as $M_1^3$, $M_2^3$ and $M_3^3$ which are given as

$M_1^3$ - No proper road or bus facilities

$M_2^3$ - Living conditions questionably poor

$M_3^3$ - Government unconcern over their development in any plan; so only introduction of machine for harvest etc. has crippled agricultural labourers, their labour opportunities

The fifth expert wishes to work with $P_1^2$, $P_2^2$ and $P_3^2$ as domain attributes and $M_1^1$, $M_2^1$ and $M_3^1$ as range attributes. The sixth expert works with $P_1^2$, $P_2^2$ and $P_3^2$ as domain nodes and $M_1^3$, $M_2^3$ and $M_3^3$ as the range nodes. The super FRM has the following super fuzzy matrix which is a 7 × 10 fuzzy supermatrix.

Let $M_s$ denote the super fuzzy FRM matrix.

$$
\begin{array}{c} \\ P_1^1 \\ P_2^1 \\ P_3^1 \\ P_4^1 \\ P_1^2 \\ P_2^2 \\ P_3^2 \end{array}
\begin{array}{c}
\begin{array}{cccccccccc} M_1^1 & M_2^1 & M_3^1 & M_1^2 & M_2^2 & M_3^2 & M_4^2 & M_1^3 & M_2^3 & M_3^3 \end{array} \\
\left[\begin{array}{ccc|cccc|ccc}
1 & 0 & 0 & 0 & 0 & 1 & 0 & 1 & 0 & 0 \\
0 & 0 & 1 & 0 & 1 & 0 & 0 & 0 & 0 & 1 \\
0 & 1 & 0 & 0 & 0 & 1 & 1 & 0 & 1 & 0 \\
1 & 0 & 0 & 1 & 0 & 0 & 1 & 0 & 1 & 0 \\ \hline
0 & 1 & 0 & 0 & 0 & 1 & 0 & 0 & 0 & 1 \\
0 & 1 & 1 & 1 & 0 & 0 & 0 & 0 & 1 & 0 \\
1 & 0 & 0 & 0 & 0 & 0 & 1 & 1 & 0 & 0
\end{array}\right]
\end{array}.
$$

The domain state vectors are super row vectors of the form; $A = [a_1^1\ a_2^1\ a_3^1\ a_4^1 \mid a_1^2\ a_2^2\ a_3^2]$, $a_i^1, a_j^2 \in \{0, 1\}$; $1 \leq i \leq 4$ and $1 \leq j \leq 3$.

The range state vectors are super row vectors of the form $B = [b_1^1\ b_2^1\ b_3^1 \mid b_1^2\ b_2^2\ b_3^2\ b_4^2 \mid b_1^3\ b_2^3\ b_3^3]$; $b_i^1, b_j^2, b_k^3 \in \{0, 1\}$; $1 \leq i \leq 3$, $1 \leq j \leq 4$ and $1 \leq k \leq 3$.

Suppose the experts want to work with $X = [0\ 1\ 0\ 0 \mid 0\ 0\ 1]$

$$
\begin{aligned}
XM_S &= Y \\
&= [1\ 0\ 1 \mid 0\ 1\ 0\ 1 \mid 1\ 0\ 1]
\end{aligned}
$$

$$Y M_S^T \hookrightarrow [1\ 1\ 1\ 1\ |\ 1\ 1\ 1]$$
$$= X'$$
$$X'M_S \hookrightarrow [1\ 1\ 1\ |\ 1\ 1\ 1\ 1\ |\ 1\ 1\ 1].$$

All the nodes come to on state. This is only an illustrative model to show how the system super FRM i.e., SFRM functions.

### 3.4.5 Special DSBAM model illustration

In this section we give illustration of the special super row bidirectional associative memories model. Suppose we have four experts who want to work using the super row BAM i.e. they choose to work with the same set of domain attributes but wish to work with varied column attributes. Suppose they have six domain attributes say $x_1, x_2, ..., x_6$ and the first expert has the range attributes as $y_1^1$, $y_2^1$, $y_3^1$, $y_4^1$ and $y_5^1$.

The second expert with $y_1^2$, $y_2^2$, $y_3^2$ and $y_4^2$, range attributes the third expert with the nodes $y_1^3$, $y_2^3$ and $y_3^3$ as range attributes and the fourth expert with the nodes $y_1^4$, $y_2^4$, $y_3^4$, $y_4^4$ and $y_5^4$ as range attributes.

The row constant super BAM model $M_R$ is given by

$$\begin{array}{c} \\ x_1 \\ x_2 \\ x_3 \\ x_4 \\ x_5 \\ x_6 \end{array} \begin{array}{c} \begin{array}{ccccc|cccc|ccc|ccccc} y_1^1 & y_2^1 & y_3^1 & y_4^1 & y_5^1 & y_1^2 & y_2^2 & y_3^2 & y_4^2 & y_1^3 & y_2^3 & y_3^3 & y_1^4 & y_2^4 & y_3^4 & y_4^4 & y_5^4 \end{array} \\ \left[ \begin{array}{c|c|c|c} \quad & \quad & \quad & \quad \\ & & & \\ & & & \\ & & & \\ & & & \\ & & & \end{array} \right] \end{array}$$

Suppose $A = [a_1\ a_2\ a_3\ a_4\ a_5\ a_6]$ is the state vector given or chosen by an expert the effect of A on $M_R$ is given by

$$A \circ M_R \hookrightarrow [b_1^1\ b_2^1\ b_3^1\ b_4^1\ b_5^1 \mid b_1^2\ b_2^2\ b_3^2\ b_4^2 \mid b_1^3\ b_2^3\ b_3^3 \mid b_1^4\ b_2^4\ b_3^4\ b_4^4\ b_5^4] = B.$$

$$B \circ M_R^T \hookrightarrow A_1 = [a'_1\ a'_2\ a'_3\ a'_4\ a'_5\ a'_6].$$

$$A_1 \circ M_R^T \hookrightarrow B_1$$

and so on. Until one arrives at a fixed point or a limit cycle. One can also start with the state vector

$$Y = [y_1^1\ y_2^1\ y_3^1\ y_4^1\ y_5^1 \mid y_1^2\ y_2^2\ y_3^2\ y_4^2 \mid y_1^3\ y_2^3\ y_3^3 \mid y_1^4\ y_2^4\ y_3^4\ y_4^4\ y_5^4]$$

and work in the similar manner using $M_R$.

### 3.4.6 RSBAM model illustration

In this section we give illustration of the Special Super Row BAM model (RSBAM model).

Let us suppose we have some 4 experts working with a problem with same column attributes but 4 distinct set of row attributes which has the synaptic connection matrix to be a column supermatrix.

Let the column attributes be given by $(y_1\ y_2 \ldots y_8)$. The first expert gives the attributes along the row to be $(x_1^1\ x_2^1\ \ldots\ x_5^1)$, the second experts row attributes for the same problem is given by $(x_1^2\ x_2^2\ \ldots\ x_6^2)$ and so on. Thus the $i^{th}$ expert gives the set of row attributes as $(x_1^i\ x_2^i\ \ldots\ x_{p_i}^i)$; i = 1, 2, …, 4. Now as all the four experts have agreed to work with the same set of column attributes viz., $(y_1\ y_2 \ldots\ y_8)$.

The related connection matrix of the RSBAM is given by $M_C$.

$$
\begin{array}{c}
\\
x_1^1 \\ x_2^1 \\ x_3^1 \\ x_4^1 \\ x_5^1 \\ x_1^2 \\ x_2^2 \\ x_3^2 \\ x_4^2 \\ x_5^2 \\ x_6^2 \\ x_1^3 \\ x_2^3 \\ x_3^3 \\ x_4^3 \\ x_1^4 \\ x_2^4 \\ x_3^4 \\ x_4^4 \\ x_5^4 \\ x_6^4
\end{array}
\begin{array}{c}
y_1 \quad y_2 \quad y_3 \quad y_4 \quad y_5 \quad y_6 \quad y_7 \quad y_8 \\
\left[\begin{array}{c}
M_1 \\ \hline
M_2 \\ \hline
M_3 \\ \hline
M_4
\end{array}\right]
\end{array}.
$$

We see $M_C$ is a special super column matrix. We see the fit vector of the domain space $F_X$ are super row vectors where as the fit vectors of the range space $F_Y$ are just ordinary or elementary row vectors.

Thus if

$$X = [x_1^1 \; x_2^1 \; \ldots \; x_5^1 \mid x_1^2 \; x_2^2 \; \ldots \; x_6^2 \mid x_1^3 \; x_2^3 \; x_3^3 \; x_4^3 \mid x_1^4 \; x_2^4 \; \ldots \; x_6^4]$$

$$XM_C \hookrightarrow Y$$
$$= [y_1\ y_2\ \ldots\ y_8] \in F_Y.$$
$$Y^T M_C \hookrightarrow X^1 \in F_X$$

now

$$X^1 M_C \hookrightarrow Y^1 \in F_Y;$$

and so on, till we arrive at a fixed point or a limit cycle.

Next we proceed on to give the illustration of a diagonal super BAM (SDBAM) model.

### 3.4.7 Special super diagonal BAM model illustration

In this section we give illustration of the special super diagonal BAM model with live illustration and show how it functions.

In this model we study the problems faced by the labourers working in garment industries and the flaws related with the garment industries. We only give 3 experts opinion for proper understanding. However it is possible to work with any number of experts once a proper programming is constructed. The attributes given by the first expert.

The problems related to the labourers working in garment industries.

$L_1^1$ - Minimum wage act not followed

$L_2^1$ - Payment of bonus not followed

$L_3^1$ - Cannot voice for safety hazards at work place

$L_4^1$ - Women discriminated

$L_5^1$ - Pay not proportional to work and hours of work

$L_6^1$ - Child labour as home can be their work place

The flaws related with the garment industries as given by the first expert.

$I_1^1$ - Working and living conditions of workers are bad

$I_2^1$ - Garment workers not organized in trade union for fear of losing their jobs

$I_3^1$ - Very poor or no health / life safety to workers

$I_4^1$ - Forced labour on migrant women and children

$I_5^1$ - Most of the industries acts are flouted

Next we proceed onto give the attributes of the second expert about the same problem using a BAM model.

Attributes given by the second expert on labourers working in the garment industries.

$L_1^2$ - No proper wages or PF benefit

$L_2^2$ - Child labour at its peak

$L_3^2$ - Workers loose their jobs if they associate themselves with the trade unions

$L_4^2$ - Migrant labourers as daily wagers because of easy availability of contractors

$L_5^2$ - Women discriminated

Attributes related with the garment industries as given by the second expert.

$I_1^2$ - Most workers are not aware of their right so are very easily exploited

$I_2^2$ - Outstanding in its performance very high turn out

$I_3^2$ - Only youngsters and children are employed in majority

$I_4^2$ - Forced labour on children and migrant women

$I_5^2$ - Garment workers not organized in any trade union

$I_6^2$ - Long-term workers turn into labour contractors or sub contractors

$I_7^2$ - Almost all industrial laws are flouted

Now we proceed on to give the views of the 3rd expert. The views of the 3rd expert related to the labourers working in garment industries.

$L_1^3$ - Pay not proportional to work

$L_2^3$ - No proper pay / bonus / PF or even gratuity

$L_3^3$ - Women and children discriminated badly

$L_4^3$ - No hygiene or safety at the work place

The concepts / attributes given by the third expert for the BAM related to the flaws related to these garment industries.

$I_1^3$ - Industry acts are not followed

$I_2^3$ - Very outstanding performance turnout in crores by these industries

$I_3^3$ - Living conditions of garment workers are very bad

$I_4^3$ - Atmosphere very hazardous to health

$I_5^3$ - No concern of the garment industry owners about the living conditions of their workers

Now using the 3 experts opinion we give the associated super diagonal fuzzy matrix model, which is a 15 × 17 super fuzzy matrix of the BAM which we choose to call as the special diagonal super fuzzy matrix. The nodes / attributes of $F_X$ are super row vectors given by

$$X = [\, L_1^1\ L_2^1\ L_3^1\ L_4^1\ L_5^1\ L_6^1 \mid L_1^2\ L_2^2\ L_3^2\ L_4^2\ L_5^2 \mid L_1^3\ L_2^3\ L_3^3\ L_4^3 \,].$$

The coordinates of super row vectors from $F_y$ which form the columns of the BAM are given by

$$Y = [\, I_1^1 \;\; I_2^1 \;\; I_3^1 \;\; I_4^1 \;\; I_5^1 \mid I_1^2 \; I_2^2 \; I_3^2 \; I_4^2 \; I_5^2 \; I_6^2 \; I_7^2 \mid I_1^3 \; I_2^3 \; I_3^3 \; I_4^3 \; I_5^3 \,].$$

Any fit vector given by the expert would be from the same scale which was used to construct the BAM model. In this problem the experts agreed to work on the scale [–4, 4]. When an expert gives the fit vector $X = [\, x_1^1 \; x_2^1 \ldots x_6^1 \mid x_1^2 \; x_2^2 \ldots x_5^2 \mid x_1^3 \ldots x_4^3 \,]$ where $x_{ij} \in [-4,\ 4\,]$; $1 \le i \le 3$, $1 \le j \le 6$ or 5 or 4.

Let $M_D$ denote the diagonal super fuzzy matrix; $M_D =$

$$\begin{array}{c} \\ L_1^1 \\ L_2^1 \\ L_3^1 \\ L_4^1 \\ L_5^1 \\ L_6^1 \\ L_1^2 \\ L_2^2 \\ L_3^2 \\ L_4^2 \\ L_5^2 \\ L_1^3 \\ L_2^3 \\ L_3^3 \\ L_4^3 \end{array}
\begin{array}{c} \begin{array}{ccccc|ccccccc|ccccc} I_1^1 & I_2^1 & I_3^1 & I_4^1 & I_5^1 & I_1^2 & I_2^2 & I_3^2 & I_4^2 & I_5^2 & I_6^2 & I_7^2 & I_1^3 & I_2^3 & I_3^3 & I_4^3 & I_5^3 \end{array} \\
\left[\begin{array}{ccccc|ccccccc|ccccc}
3 & 2 & 1 & 4 & 0 & & & & & & & & & & & & \\
1 & -2 & 0 & 3 & 2 & & & & & & & & & & & & \\
2 & 1 & -1 & 0 & 2 & & & & (0) & & & & & & (0) & & \\
3 & 1 & 0 & 3 & 1 & & & & & & & & & & & & \\
-1 & 2 & 1 & 2 & -1 & & & & & & & & & & & & \\
4 & 3 & -1 & 4 & 3 & & & & & & & & & & & & \\
\hline
 & & & & & 0 & 1 & 2 & 3 & -1 & 2 & 1 & & & & & \\
 & & & & & 3 & 1 & -1 & 0 & 3 & 1 & 0 & & & & & \\
 & & (0) & & & 4 & 2 & 0 & 1 & 4 & -4 & 2 & & & (0) & & \\
 & & & & & 1 & -1 & 3 & 0 & 1 & 0 & -1 & & & & & \\
 & & & & & -2 & 3 & -1 & 2 & 0 & -2 & 0 & & & & & \\
\hline
 & & & & & & & & & & & & 3 & 0 & -2 & 1 & 0 \\
 & & & & & & & & & & & & 2 & 2 & 1 & 0 & 1 \\
 & & (0) & & & & & & (0) & & & & -1 & 1 & 0 & 2 & 1 \\
 & & & & & & & & & & & & 1 & -1 & 0 & 3 & 0
\end{array}\right]
\end{array}$$

Suppose the expert gives the fit vector as

$$X_K \quad = \quad [2\; 3\; 0\; 1\; -1\; -2 \mid 0\; 2\; -1\; -3\; -5 \mid 3\; 0\; -2\; -4]$$

at the $k^{th}$ time period; using the activation function S we get

$S(X_K)$ = [1 1 0 1 0 0 | 0 1 0 0 0 | 1 0 0 0 ]
$S(X_K)\ M_D$ = [7 1 1 10 3 | 3 1 –1 0 3 1 0 | 3 0 –2 1 0]
= $Y_{K+1}$
$S(Y_{K+1})$ = [1 1 1 1 1 | 1 1 0 0 1 1 0 | 1 0 0 1 0]

$S(Y_{K+1})$ o $M_D^T$ and so on.

### 3.4.8 Super BAM model illustration

In this section we give illustration of super BAM model by a real world problem.

Now the super fuzzy matrix $M_s$ associated with four experts for the same problem of garment industries is given below.

$$
\begin{array}{c}
\\
I_1^1 \\ I_2^1 \\ I_3^1 \\ I_4^1 \\ I_5^1 \\ I_1^2 \\ I_2^2 \\ I_3^2 \\ I_4^2
\end{array}
\begin{array}{c}
\begin{array}{cccccccccc}
L_1^1 & L_2^1 & L_3^1 & L_4^1 & L_1^2 & L_1^2 & L_3^2 & L_4^2 & L_5^2 & L_6^2
\end{array} \\
\left[\begin{array}{cccc|cccccc}
2 & 1 & -1 & 3 & 0 & 1 & 0 & 2 & 3 & 4 \\
3 & 4 & 0 & 1 & 2 & 3 & 1 & 1 & -2 & -1 \\
-2 & 3 & 2 & -1 & 1 & 0 & 1 & -2 & 0 & 2 \\
1 & 1 & 1 & 2 & 2 & -1 & 4 & 0 & 0 & 3 \\
4 & 3 & -1 & 0 & 3 & 2 & -3 & 1 & 1 & -2 \\
\hline
3 & 0 & -1 & 2 & 0 & 1 & 0 & 2 & 1 & -3 \\
2 & -2 & 3 & 1 & 1 & 0 & 3 & -1 & 0 & 2 \\
-2 & 2 & 2 & -1 & 2 & 2 & 3 & 0 & 1 & 0 \\
-3 & 0 & 2 & 4 & 3 & 1 & -1 & 3 & 4 & -1
\end{array}\right]
\end{array}.
$$

The first expert works with $(L_1^1\ L_2^1\ L_3^1\ L_4^1)$ as the column attribute and row attribute $(I_1^1\ I_2^1\ I_3^1\ I_4^1\ I_5^1)$. The second expert works with $(L_1^2\ L_2^2\ L_3^2\ L_4^2\ L_5^2\ L_6^2)$ as the column attribute and $(I_1^1\ I_2^1\ I_3^1\ I_4^1\ I_5^1)$ as the row attribute.

The third expert works with $(L_1^1\ L_2^1\ L_3^1\ L_4^1)$ as the column attributes and $(I_1^2\ I_2^2\ I_3^2\ I_4^2)$ as the row attributes and the fourth expert use ($L_1^2$ $L_2^2$ $L_3^2$ $L_4^2$ $L_5^2$ $L_6^2$) as the column attributes and ($I_1^2$ $I_2^2$ $I_3^2$ $I_4^2$) as row attributes.

However all the four experts agree to work on the same interval viz. [–4, 4]. Now any given state / fit vector by the expert is synchronized and made to work on the super fuzzy matrix $M_S$. Clearly once again the fit vectors which are super row vectors take their values only from the interval [–4, 4]. This super BAM model also functions as the usual BAM model with a simple difference these resultant row vectors which is a fixed point or a limit cycle is just a super row vector.

### 3.4.9 Special Row FAM model illustration (SRFAM)

In this section we give illustration of special row FAM model and show how it functions.

Suppose we have some four experts who wish to work on a specific problem using a special row FAM. All of them agree upon to work with the same set of row attributes. We indicate how to study the multi expert problem to analyse the problems faced by the labourers in Garment industries. The attributes related to the labourers given by all the four experts is as follows:

$W_1$ - Minimum wages act not followed with no bonus or PF or gratuity
$W_2$ - Women discriminated
$W_3$ - More hours of work pay not proportional to work
$W_4$ - Cannot voice for safety hazards
$W_5$ - Child labour is at its peak

The attributes given by the first expert related to the flaws of garment industries.

$I_1^1$ - Living conditions of garment workers are bad
$I_2^1$ - Garment workers not organized in trade unions for fear of losing jobs

$I_3^1$ - Outstanding in its performance very high turn out

$I_4^1$ - Youngsters and children employed for better turn out

$I_5^1$ - No medical facilities or protection from health hazards

The attributes related to the garment industries given by the second expert.

$I_1^2$ - They run with very high profit yet deny proper pay to workers

$I_2^2$ - They do not allow the functioning of any workers union

$I_3^2$ - They employ only children and youngsters on contract basis so that when they become little middle aged or women get married they are sent home

$I_4^2$ - Forced labour on migrant women and children

$I_5^2$ - Living conditions of the workers is very poor

$I_6^2$ - Long term workers are turned into labour contractors or sub contractors (contracting inside the industry)

The attributes given by the third expert related to the flaws in running the garment industry.

$I_1^3$ - Garment workers not organized in trade unions for fear of losing job

$I_2^3$ - All there industries run with a very high profit

$I_3^3$ - Living conditions of the workers questionable

$I_4^3$ - Forced labour on children and women

$I_5^3$ - All industry acts flouted

Now we proceed on to give the attributes given by the fourth expert related to the flaws related with the industry.

$I_1^4$ - Industry acts are flouted

$I_2^4$ - All labour laws disobeyed

$I_3^4$ - Children employed for better turn out

$I_4^4$ - No medical / health protection for the workers

$I_5^4$ - All industries run with 100% profit

$I_6^4$ - The workers social conditions very poor and questionable

$I_7^4$ - Workers cannot be organized into trade unions for fear of losing job

The super row FAM is given by $F_R$

$$\begin{array}{c} \\ W_1 \\ W_2 \\ W_3 \\ W_4 \\ W_5 \end{array} \begin{array}{c} I_1^1\ I_2^1\ I_3^1\ I_4^1\ I_5^1\ I_1^2\ I_2^2\ I_3^2\ I_4^2\ I_5^2\ I_6^2\ I_1^3\ I_2^3\ I_3^3\ I_4^3\ I_5^3\ I_1^4\ I_2^4\ I_3^4\ I_4^4\ I_5^4\ I_6^4\ I_7^4 \\ \left[\begin{array}{c|c|c|c} \\ \\ A_1^1 & A_2^2 & A_3^3 & A_4^4 \\ \\ \\ \end{array}\right] \end{array}$$

where $A_i^i$ is a fuzzy matrix i.e. it takes the entries from the unit interval [0, 1], $1 \leq i \leq 4$.

If A = [0 1 0 1 0] is a fit vector given by an expert then

$$A \text{ o } F_R = \max \{\min (a_i\ a_{ij}^t)\};\ 1 \leq t \leq 4$$
$$= B \text{ (say)}.$$

B is a super fuzzy row vector.

$$F_R \text{ o } B = \max \min (a_{ij}^t\ b_j)$$
$$= a_i \in A^1.$$

This process is repeated until we arrive at a fixed point. The pair of resultant fit vectors $(A_k, B_t)$ would be a fuzzy row vector with

$$A_k = [a_k^1\ a_k^2\ a_k^3\ a_k^4\ a_k^5];$$

$a_k^t \in [0, 1]$ $1 \leq t \leq 5$ and

$$B_t = [b_1^1\ b_2^1\ b_3^1\ b_4^1\ b_5^1 \mid b_1^2\ b_2^2\ b_3^2\ b_4^2\ b_5^2\ b_6^2 \mid$$
$$b_1^3\ b_2^3\ b_3^3\ b_4^3\ b_5^3 \mid b_1^4\ b_2^4\ b_3^4\ b_4^4\ b_5^4\ b_6^4\ b_7^4],$$

$b_{jk}^s \in [0, 1]; 1 \leq s \leq 4$ is super fuzzy row vector. Now from the resultant we can see the entries in the resultant are values from [0, 1] so that we can get the gradation of importance of each of the attributes.

### 3.4.10 Special Column super FAM model

In this section we show how a special column super FAM model operates on a real world problem. Suppose we have a set of experts who wish to work with a problem having a same set of column attributes but want to use a distinct set of row attributes, the multi expert FAM model makes use of the special column supervector for this study. We illustrate this by an example, however the model has been described in 193-5. Suppose three experts wish to study using FAM model the socio economic problems of HIV/AIDS affected women patients. All of them agree upon the column nodes or attributes related with the women to be $W_1$, $W_2$, $W_3$, $W_4$, $W_5$ and $W_6$ where

- $W_1$ - Child marriage / widower marriage, child married to men thrice or four times their age.
- $W_2$ - Causes of women being infected with HIV/AIDS
- $W_3$ - Disease untreated till it is chronic or they are in last stages of life
- $W_4$ - Women not bread winners thus a traditional set up of our society
- $W_5$ - Deserted by family when they have HIV/AIDS
- $W_6$ - Poverty and ignorance a major reason for becoming HIV/AIDS victims

The concepts associated with society, men/husband leading women to become HIV/AIDS patients given by the first expert.

$R_1^1$ - Female child a burden so the sooner they give in marriage the better relief economically.

$R_2^1$ - Poverty / not owners of property.

$R_3^1$ - No moral responsibility on the part of husbands and they infect their wives willfully.

$R_4^1$ - STD/VD infected husbands leading wives to frequent natural abortion or death of born infants.

$R_5^1$ - Bad habits of men/husbands.

The concepts given by the second expert in relation to the causes of women becoming HIV/AIDS infected.

$R_1^2$ - Husbands hide their disease from their family members so the wife become HIV/AIDS infected

$R_2^2$ - STD/VD infected husbands

$R_3^2$ - Poverty, don't own poverty

$R_4^2$ - Bad habits of men/husbands

The 3rd experts views on how HIV/AIDS affect rural women.

$R_1^3$ - STD/VD infected husbands

$R_2^3$ - Female child a burden so they want to dispose off as soon as possible to achieve economic freedom

$R_3^3$ - Bad habits of men/husbands

$R_4^3$ - No moral responsibility on the part of husbands and they infect their wives willfully

$R_5^3$ - Poverty a major draw back

$R_6^3$ - Migration of men to earn, due to lack of employment in the rural areas

Let $F_c$ denote the related special super column fuzzy matrix of the super column FAM model.

$$
\begin{array}{c}
 \\ R_1^1 \\ R_2^1 \\ R_3^1 \\ R_4^1 \\ R_5^1 \\ R_1^2 \\ R_2^2 \\ R_3^2 \\ R_4^2 \\ R_1^3 \\ R_2^3 \\ R_3^3 \\ R_4^3 \\ R_5^3 \\ R_6^3
\end{array}
\begin{array}{c}
\begin{array}{cccccc} W_1 & W_2 & W_3 & W_4 & W_5 & W_6 \end{array} \\
\left[\begin{array}{c}
B_1^1 \\ \hline B_2^2 \\ \hline B_3^3
\end{array}\right]
\end{array}
$$

where each of the $B_i^i$ is a fuzzy matrix; i = 1, 2, 3. Suppose we have a fit vector

$$A = [a_1^1\ a_2^1\ a_3^1\ a_4^1\ a_5^1 \mid a_1^2\ a_2^2\ a_3^2\ a_4^2 \mid a_1^3\ a_2^3\ a_3^3\ a_4^3\ a_5^3\ a_6^3]$$

which is a super fuzzy row vector then

$$A \text{ o } F_c = B = [b_1\ b_2\ b_3\ b_4\ b_5\ b_6];$$

now B the fit vector is only a fuzzy row vector.

The effect of B on the system is given by $F_c$ o B = $A_1$, $A_1$ is a fit row vector which is a super fuzzy row vector, we proceed on until we arrive at a system stability.

### 3.4.11 Super fuzzy diagonal FAM model illustration

This section gives illustration of how a super diagonal FAM functions. Here we give a super diagonal FAM model using 3 experts opinion. Suppose some three experts are interested in the study of the problem of cause of migration and migrant labourers becoming victim to HIV/AIDS and all of them wish to work with distinct set of row and column attributes.

The row attributes given by the first expert related to attributes that lead to migrant labourers becoming victim of HIV/AIDS.

$M_1^1$ - Government help never reaches the rural poor illiterate

$M_2^1$ - Availability of cheap liquor

$M_3^1$ - Poverty

$M_4^1$ - No job opportunities in their home town / village

$M_5^1$ - Living conditions questionably poor

$M_6^1$ - No proper health center in their home town / village

Attributes given by the first expert related to the migrant labourers becoming victims of HIV/AIDS.

$P_1^1$ - No education / No help by government

$P_2^1$ - Smoke and visit of CSWs

$P_3^1$ - Addiction of cheap liquor

$P_4^1$ - Awareness program never reaches them

$P_5^1$ - No responsibility of parents to educate their children.

Now we enlist the row and column attributes given by the second expert in the following.

The attributes which promote migration.

$M_1^2$ - No job opportunities in the place where they live

$M_2^2$ - Availability of CSWs and cheap liquor

$M_3^2$ - Poverty

$M_4^2$ - No proper health center

The attributes related to the migrant labourers who become victims of HIV/AIDS.

$P_1^2$ - No education

$P_2^2$ - Very questionable living conditions so migrate

$P_3^2$ - Addiction to cheap liquor, CSWs and smoke

$P_4^2$ - Awareness program never reaches them

$P_5^2$ - No responsibility of parents to educate children.

Now we proceed on to work with the 3rd expert.
The attributes connected with the migrant labourers becoming HIV/AIDS.

$M_1^3$ - No proper health center

$M_2^3$ - No job opportunities at their home town

$M_3^3$ - Poverty

$M_4^3$ - Availability of CSWs and liquor at cheap rates

$M_5^3$ - Government programmes never reaches the poor illiterate.

$M_6^3$ - Government silence over their serious problems.

$M_7^3$ - No school even for primary classes.

The nodes related to the problems of migrant labourers who suffer and ultimately become victims of HIV/AIDS as given by the third expert.

| | | |
|---|---|---|
| $P_1^3$ | - | No education |
| $P_2^3$ | - | Very questionable living condition so migrate |
| $P_3^3$ | - | Awareness program never reaches them |
| $P_4^3$ | - | No responsibility of parents to educate children |
| $P_5^3$ | - | Addiction to cheap liquor and CSWs. |

Now we proceed on to give the diagonal super fuzzy matrix $F_D$ related to the super FAM using the three experts. This supermatrix will have 17 row elements and 15 column elements. The rest of the non diagonal elements are zero.

$$
\begin{array}{c}
 \\
M_1^1 \\ M_2^1 \\ M_3^1 \\ M_4^1 \\ M_5^1 \\ M_6^1 \\ M_1^2 \\ M_2^2 \\ M_3^2 \\ M_4^2 \\ M_1^3 \\ M_2^3 \\ M_3^3 \\ M_4^3 \\ M_5^3 \\ M_6^3 \\ M_7^3
\end{array}
\begin{array}{c}
\begin{array}{ccccc|ccccc|ccccc}
P_1^1 & P_2^1 & P_3^1 & P_4^1 & P_5^1 & P_1^2 & P_2^2 & P_3^2 & P_4^2 & P_5^2 & P_1^3 & P_2^3 & P_3^3 & P_4^3 & P_5^3
\end{array} \\
\left[\begin{array}{c|c|c}
(F_{ij}^1) & (0) & (0) \\ \hline
(0) & (F_{ij}^2) & (0) \\ \hline
(0) & (0) & (F_{ij}^3)
\end{array}\right]
\end{array}
$$

$$\left(F_{ij}^1\right) = \begin{array}{c} \\ M_1^1 \\ M_2^1 \\ M_3^1 \\ M_4^1 \\ M_5^1 \\ M_6^1 \end{array} \begin{array}{c} \begin{array}{ccccc} P_1^1 & P_2^1 & P_3^1 & P_4^1 & P_5^1 \end{array} \\ \left[\begin{array}{ccccc} & & & & \\ & & & & \\ & & & & \\ & & & & \\ & & & & \\ & & & & \end{array}\right] \end{array},$$

$$\left(F_{ij}^2\right) = \begin{array}{c} \\ M_1^2 \\ M_2^2 \\ M_3^2 \\ M_4^2 \end{array} \begin{array}{c} \begin{array}{ccccc} P_1^2 & P_2^2 & P_3^2 & P_4^2 & P_5^2 \end{array} \\ \left[\begin{array}{ccccc} & & & & \\ & & & & \\ & & & & \\ & & & & \end{array}\right] \end{array}$$

and

$$\left(F_{ij}^3\right) = \begin{array}{c} \\ M_1^3 \\ M_2^3 \\ M_3^3 \\ M_4^3 \\ M_5^3 \\ M_6^3 \\ M_7^3 \end{array} \begin{array}{c} \begin{array}{ccccc} P_1^3 & P_2^3 & P_3^3 & P_4^3 & P_5^3 \end{array} \\ \left[\begin{array}{ccccc} & & & & \\ & & & & \\ & & & & \\ & & & & \\ & & & & \\ & & & & \\ & & & & \end{array}\right] \end{array}.$$

The fit vectors from the domain space X would be a super row fuzzy vector where

$$X = [\, x_1^1 \; x_2^1 \; x_3^1 \; x_4^1 \; x_5^1 \; x_6^1 \mid x_1^2 \; x_2^2 \; x_3^2 \; x_4^2 \mid x_1^3 \; x_2^3 \; x_3^3 \; x_4^3 \; x_5^3 \; x_6^3 \; x_7^3 \,].$$

The elements $x_i^1$, $x_j^2$, $x_k^3$ $\in [0, 1]$; $1 \leq i \leq 6$, $1 \leq j \leq 4$ and $1 \leq k \leq 7$. The fit vector from the range space Y would be a super row fuzzy vector where

$$Y = [\, y_1^1 \; y_2^1 \; y_3^1 \; y_4^1 \; y_5^1 \mid y_1^2 \; y_2^2 \; y_3^2 \; y_4^2 \; y_5^2 \mid y_1^3 \; y_2^3 \; y_3^3 \; y_4^3 \; y_5^3 \,].$$

The elements $y_j^t \in [0, 1]$; $1 \leq t \leq 3$ and $1 \leq j \leq 5$. We find X o $F_D$ = Y now $F_D$ o Y = $X_1$ and so on until we arrive at the equilibrium of the system.

Now we proceed on to give the next model.

### 3.4.12 Fuzzy Super FAM model Illustration

In this section we illustrate a fuzzy super FAM model and its mode of functioning. Suppose we wish to study the Employee and Employer relationship model. We have 12 = 4 × 3 experts working on the problem. The concepts / nodes related with the employee given by the first expert.

$D_1^1$ - Pay with allowance to the employee
$D_2^1$ - Only pay to employee
$D_3^1$ - Average performance by the employee
$D_4^1$ - Poor performance by the employee
$D_5^1$ - Employee works more than 8 hours a day

The concepts related with the employer as given by the first expert is as follows:

$R_1^1$ - The industry runs with maximum profit
$R_2^1$ - No loss, no gain
$R_3^1$ - Trade unions / employee unions not encouraged
$R_4^1$ - Employee unions have a say in the industry
$R_5^1$ - Contract labourers form the majority of workers
$R_6^1$ - The profit/loss workers remain in the same state

The second expert wishes to work with the same domain attributes but however has given a set of range attributes which is as follows:

$R_1^2$ - The industry is unconcerned about the social and the economic conditions of their workers.
$R_2^2$ - Industry for the past one decade has always run with profit or no loss.
$R_3^2$ - Employee unions do not exist in their industries.
$R_4^2$ - Owners always have the trend to exploit the employees.
$R_5^2$ - Bonus is given every year.

The third expert wishes to work with the same set of domain attributes as the first expert but has a set of different range attributes.

$R_1^3$ - Industries least bothered about the employees health hazards.
$R_2^3$ - The industrialists have no mind set to look into the problems of the employees.
$R_3^3$ - The industry aims only for profit unconcerned about any other factors.
$R_4^3$ - Employees get proper medical care
$R_5^3$ - Employees enjoy all sorts of benefit including pension after retirement.
$R_6^3$ - Employees are made to work not following any of the industrial acts.

Now the fourth expert wishes to work with a new set of domain attributes but is happy to accept the first experts range attributes viz., $R^1_1$, $R^1_2$, …, $R^1_6$.

The set of domain attributes given by the 4th expert is as follows:

$D_1^2$ - No pension benefits or medical aid

$D_2^2$ - The employee do not perform well

$D_3^2$ - They work for more hours with less pay

$D_4^2$ - The employees are vexed, for profit or loss the workers are not given any benefit

$D_5^2$ - Workers are not informed of the health hazards for their labour conditions are questionable.

The fifth expert wishes to work with the domain attributes $D_1^2$, $D_2^2$, $D_3^2$, $D_4^2$ and $D_5^2$ the range attributes taken by him are $R_1^2$, $R_2^2$, $R_3^2$, $R_4^2$ and $R_5^2$.

Now the sixth expert wishes to work with the domain attributes as $D_1^2$, $D_2^2$, $D_3^2$, $D_4^2$ and $D_5^2$. The range attributes taken by him are $R_1^3$, $R_2^3$, $R_3^3$, $R_4^3$, $R_5^3$ and $R_6^3$.

The seventh expert works with the new set of domain attributes are taken as $D_1^3$, $D_2^3$, $D_3^3$, $D_4^3$, $D_5^3$ and $D_6^3$ which are as follows.

$D_1^3$ - Living conditions of the workers are questionable.

$D_2^3$ - The workers do not receive any good pay or medical aid or pension benefits.

$D_3^3$ - Employees cannot form any association or be members of any trade unions.

$D_4^3$ - Most of them do not have any job security.

$D_5^3$ - The employees are not informed about the health hazards when they do their job.

This expert however wishes to work with the range attributes as $R_1^1$, $R_2^1$, …, $R_6^1$.

The eighth expert works with the domain attributes $D_1^3$, $D_2^3$, $D_3^3$, $D_4^3$ and $D_5^3$ and range attributes as $R_1^2$, $R_2^2$, $R_3^2$, $R_4^2$ and $R_5^2$.

The ninth expert works with the domain attributes $D_1^3$, $D_2^3$, $D_3^3$, $D_4^3$ and $D_5^3$. The range attributes taken by him are $R_1^3$, $R_2^3$, $R_3^3$, $R_4^3$, $R_5^3$ and $R_6^3$.

The tenth expert works with the new set of domain attributes which is as follows.

| | | |
|---|---|---|
| $D_1^4$ | - | No pension or medical benefit |
| $D_2^4$ | - | The industrialist least bothered about the social or economic or health conditions of the workers |
| $D_3^4$ | - | The employee are not given pay proportion to their work |
| $D_4^4$ | - | The employee are not allowed to associate themselves with any trade union or allowed to form any union |
| $D_5^4$ | - | The employee suffer from several types of health hazards |
| $D_6^4$ | - | The employee feel, loss or gain to industry, it does not affect their economic status |

However the range attributes taken by him are $R_1^1$, $R_2^1$, …, $R_6^1$.

The eleventh expert works with the domain attributes as given by the tenth expert viz., $D_1^4$, $D_2^4$, $D_3^4$, $D_4^4$, $D_5^4$ and $D_6^4$. The range attributes are taken as $R_1^2$, $R_2^2$, …, $R_5^2$.

The twelfth expert works with the domain attributes as $D_1^4$, $D_2^4$, $D_3^4$, $D_4^4$, $D_5^4$ and $D_6^4$ and the range attributes are taken as $R_1^3$, $R_2^3$, …, $R_6^3$.

The associated dynamical system which is a super fuzzy matrix is as follows given by $F_s$.

$$
\begin{array}{c|c|c|c}
 & R_1^1\ R_2^1\ R_3^1\ R_4^1\ R_5^1\ R_6^1 & R_1^2\ R_2^2\ R_3^2\ R_4^2\ R_5^2 & R_1^3\ R_2^3\ R_3^3\ R_4^3\ R_5^3\ R_6^3 \\ \hline
\begin{matrix} D_1^1 \\ D_2^1 \\ D_3^1 \\ D_4^1 \\ D_5^1 \end{matrix} & F_1^1 & F_2^1 & F_3^1 \\ \hline
\begin{matrix} D_1^2 \\ D_2^2 \\ D_3^2 \\ D_4^2 \\ D_5^2 \end{matrix} & F_4^2 & F_5^2 & F_6^2 \\ \hline
\begin{matrix} D_1^3 \\ D_2^3 \\ D_3^3 \\ D_4^3 \\ D_5^3 \\ D_6^3 \end{matrix} & F_7^3 & F_8^3 & F_9^3 \\ \hline
\begin{matrix} D_1^4 \\ D_2^4 \\ D_3^4 \\ D_4^4 \\ D_5^4 \\ D_6^4 \end{matrix} & F_{10}^4 & F_{11}^4 & F_{12}^4
\end{array}
$$

Each of the $F_r^t$; $1 \leq t \leq 4$ and $1 \leq r \leq 12$ are fuzzy elementary supermatrix and $F_s$ forms the fuzzy supermatrix. Any fit vector in the domain space is a fuzzy super row vector given by

$$
\begin{aligned}
X = {} & [x_1^1\ x_2^1\ x_3^1\ x_4^1\ x_5^1 \mid x_1^2\ x_2^2\ x_3^2\ x_4^2\ x_5^2 \mid \\
& x_1^3\ x_2^3\ x_3^3\ x_4^3\ x_5^3\ x_6^3 \mid x_1^4\ x_2^4\ x_3^4\ x_4^4\ x_5^4\ x_6^4],
\end{aligned}
$$

$x_i^t$, $x_j^t$ $\in$ [0, 1] with $1 \le t \le 4$, $1 \le j \le 5$ and $1 \le i \le 6$. The fit vector of the range space is a fuzzy super row vector given by

$$Y = [\, y_1^1\ y_2^1\ y_3^1\ y_4^1\ y_5^1\ y_6^1 \mid y_1^2\ y_2^2\ y_3^2\ y_4^2\ y_5^2 \mid y_1^3\ y_2^3\ y_3^3\ y_4^3\ y_5^3\ y_6^3 \,],$$

$y_i^t$, $y_j^t$ $\in$ [0, 1]; $1 \le t \le 3$ and $1 \le i \le 6$, $1 \le j \le 5$.

Now X o $F_s$ is calculated. If the resultant is Y we find $F_s$ o Y and so on until we arrive at a equilibrium state of the system. Here the resultant fuzzy fit vectors help us in the results as they are graded and not like the usual fit vectors which take values as 0 and 1. These fit vectors take values from the interval [0, 1] so a gradation of preference is always possible and using the resultant we can predict the most important attribute and a less important node.

Now having described these super fuzzy structures which have been introduced in – we proceed in to mention some of its uses.

### 3.4.13 Uses of super fuzzy models

1. When we have more than one expert giving his/her opinion on the problem super fuzzy models come handy. As far as super fuzzy models are concerned the same program can be used only in the final step we need to make a partition of the resultant fuzzy row vector so that it is a super row vector depicting the views of each and every expert simultaneously. This helps the observers to compare them and study their relative effects.

2. When combined FRM is used we see that some times a negative value –1 cancels with the positive value +1 making the effect of a specific attribute over the other to zero or no impact but in actually we see the pair of nodes have an impact. In such cases we see the super FRMs or a diagonal super FRM or the special column FRM or the special row

FRM is not only handy but can give the opinion of each and every expert by a single super row vector distinctly.

3. The super BAM models likewise can act as a multi expert model we see the BAM models can at a time serve as a single expert system but super BAM models have the capacity to serve as an multi expert system. Also we can compare the n experts opinion as they would be given by a single super row vector. Thus these super fuzzy BAM models are new and are a powerful multi experts model using the same BAM model technique. Even if they choose to opt different scales or intervals still, they can serve as a single dynamical system. Only supermatrices help to model them.

4. Another advantage of these super models is they can function even when the experts agree upon to work with the same set of attributes along the rows (or columns) and distinct set of attributes along the columns (or rows). Even if each expert wishes to use a distinct set of attributes we have the special diagonal super FRM model or special diagonal super BAM model or the special diagonal super FAM model to serve the purpose.

These models are constructed more to the engineers, doctors and above all socio scientists for in our opinion the socio scientists are the ones who do lots of service to humanity.

The development of any tool or subject if reaches the last man or is useful to the last man the authors feel it is worthwhile. Thus this book mainly aims to built models which are with in the reach of a analyst with minimum mathematical knowledge. Using the computers the solutions can be easily programmed. Further the authors have made it more easy by such extensions in these super model. The same program can be used only a small command or a construction should be added to divide or partition the resultant row vectors to make it a supervector in keeping with the supermatrix used for the investigation. For more about programming see [228-231].

## 3.5 Super FCM Models

We have just seen super FRM model, super BAM model and super FAM models. All these three models makes use of special row fuzzy supermatrix or special column fuzzy supermatrix, fuzzy super special diagonal matrix and fuzzy supermatrix. But when we want to model Super Fuzzy Cognitive Maps i.e. super FCMs we cannot have the connection matrix to be a special fuzzy row supermatrix or a special fuzzy column supermatrix or the fuzzy supermatrix. Only the special diagonal fuzzy supermatrix alone can be used when we want to study a single supermatrix to depict a multi expert model using the FCMs. This section has two subsections.

We just briefly recall the properties of FCM models in the first section. In section two we define the new super FCM models.

### 3.5.1 Introduction to FCM models

This section has two subsections in the first subsection we just briefly recall the functioning of the FCM from [108, 112].

In this section we recall the notion of Fuzzy Cognitive Maps (FCMs), which was introduced by Bart Kosko [108] in the year 1986. We also give several of its interrelated definitions. FCMs have a major role to play mainly when the data concerned is an unsupervised one. Further this method is most simple and an effective one as it can analyse the data by directed graphs and connection matrices.

**DEFINITION 3.5.1.1:** *An FCM is a directed graph with concepts like policies, events etc. as nodes and causalities as edges. It represents causal relationship between concepts.*

***Example 3.5.1.1:*** In Tamil Nadu (a southern state in India) in the last decade several new engineering colleges have been approved and started. The resultant increase in the production of engineering graduates in these years is disproportionate with the need of engineering graduates. This has resulted in thousands of unemployed and underemployed graduate engineers. Using an expert's opinion we study the effect of such unemployed people on the society. An expert spells out the five major concepts relating to the unemployed graduated engineers as

| | | |
|---|---|---|
| $E_1$ | – | Frustration |
| $E_2$ | – | Unemployment |
| $E_3$ | – | Increase of educated criminals |
| $E_4$ | – | Under employment |
| $E_5$ | – | Taking up drugs etc. |

The directed graph where $E_1, \ldots, E_5$ are taken as the nodes and causalities as edges as given by an expert is given in the following Figure 3.5.1.1:

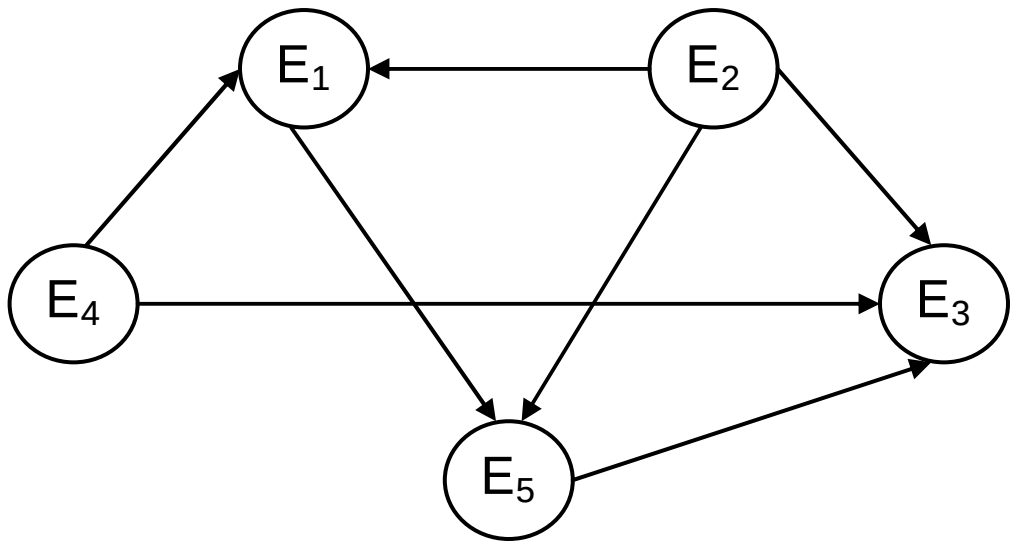

**FIGURE: 3.5.1.1**

According to this expert, increase in unemployment increases frustration. Increase in unemployment, increases the educated criminals. Frustration increases the graduates to take up to evils like drugs etc. Unemployment also leads to the increase in number of persons who take up to drugs, drinks etc. to forget their worries and unoccupied time. Under-employment forces then to do criminal acts like theft (leading to murder) for want of more money and so on. Thus one cannot actually get data for this but can use the expert's opinion for this unsupervised data to obtain some idea about the real plight of the situation. This is just an illustration to show how FCM is described by a directed graph.

{If increase (or decrease) in one concept leads to increase (or decrease) in another, then we give the value 1. If there exists no relation between two concepts the value 0 is given. If increase (or decrease) in one concept decreases (or increases) another, then we give the value –1. Thus FCMs are described in this way.}

**DEFINITION 3.5.1.2:** *When the nodes of the FCM are fuzzy sets then they are called as fuzzy nodes.*

**DEFINITION 3.5.1.3:** *FCMs with edge weights or causalities from the set {–1, 0, 1} are called simple FCMs.*

**DEFINITION 3.5.1.4:** *Consider the nodes / concepts $C_1, \ldots, C_n$ of the FCM. Suppose the directed graph is drawn using edge weight $e_{ij} \in \{0, 1, -1\}$. The matrix E be defined by $E = (e_{ij})$ where $e_{ij}$ is the weight of the directed edge $C_i C_j$. E is called the adjacency matrix of the FCM, also known as the connection matrix of the FCM.*

It is important to note that all matrices associated with an FCM are always square matrices with diagonal entries as zero.

**DEFINITION 3.5.1.5:** *Let $C_1, C_2, \ldots, C_n$ be the nodes of an FCM. $A = (a_1, a_2, \ldots, a_n)$ where $a_i \in \{0, 1\}$. A is called the instantaneous state vector and it denotes the on-off position of the node at an instant.*

$$a_i = 0 \text{ if } a_i \text{ is off and}$$
$$a_i = 1 \text{ if } a_i \text{ is on}$$

*for $i = 1, 2, \ldots, n$.*

**DEFINITION 3.5.1.6:** *Let $C_1, C_2, \ldots, C_n$ be the nodes of an FCM. Let $\overrightarrow{C_1C_2}, \overrightarrow{C_2C_3}, \overrightarrow{C_3C_4}, \ldots, \overrightarrow{C_iC_j}$ be the edges of the FCM ($i \neq j$). Then the edges form a directed cycle. An FCM is said to be cyclic if it possesses a directed cycle. An FCM is said to be acyclic if it does not possess any directed cycle.*

**DEFINITION 3.5.1.7:** *An FCM with cycles is said to have a feedback.*

**DEFINITION 3.5.1.8:** *When there is a feedback in an FCM, i.e., when the causal relations flow through a cycle in a revolutionary way, the FCM is called a dynamical system.*

**DEFINITION 3.5.1.9:** *Let* $\overline{C_1C_2}, \overline{C_2C_3}, \ldots, \overrightarrow{C_{n-1}C_n}$ *be a cycle. When* $C_i$ *is switched on and if the causality flows through the edges of a cycle and if it again causes* $C_i$ *, we say that the dynamical system goes round and round. This is true for any node* $C_i$ *, for i = 1, 2, … , n. The equilibrium state for this dynamical system is called the hidden pattern.*

**DEFINITION 3.5.1.10:** *If the equilibrium state of a dynamical system is a unique state vector, then it is called a fixed point.*

***Example 3.5.1.2:*** Consider a FCM with $C_1$, $C_2$, …, $C_n$ as nodes. For example let us start the dynamical system by switching on $C_1$. Let us assume that the FCM settles down with $C_1$ and $C_n$ on i.e. the state vector remains as (1, 0, 0, …, 0, 1) this state vector (1, 0, 0, …, 0, 1) is called the fixed point.

**DEFINITION 3.5.1.11:** *If the FCM settles down with a state vector repeating in the form*

$$A_1 \to A_2 \to \ldots \to A_i \to A_1$$

*then this equilibrium is called a limit cycle.*

Methods of finding the hidden pattern are discussed in the following Section 1.2.

**DEFINITION 3.5.1.12:** *Finite number of FCMs can be combined together to produce the joint effect of all the FCMs. Let* $E_1$, $E_2$, … , $E_p$ *be the adjacency matrices of the FCMs with nodes* $C_1$, $C_2$, …, $C_n$ *then the combined FCM is got by adding all the adjacency matrices* $E_1$, $E_2$, …, $E_p$.

*We denote the combined FCM adjacency matrix by $E = E_1 + E_2 + \ldots + E_p$.*

**NOTATION:** Suppose $A = (a_1, \ldots, a_n)$ is a vector which is passed into a dynamical system E. Then $AE = (a'_1, \ldots, a'_n)$ after thresholding and updating the vector suppose we get $(b_1, \ldots, b_n)$ we denote that by

$$(a'_1, a'_2, \ldots, a'_n) \hookrightarrow (b_1, b_2, \ldots, b_n).$$

Thus the symbol ' $\hookrightarrow$ ' means the resultant vector has been thresholded and updated.

FCMs have several advantages as well as some disadvantages. The main advantage of this method it is simple. It functions on expert's opinion. When the data happens to be an unsupervised one the FCM comes handy. This is the only known fuzzy technique that gives the hidden pattern of the situation. As we have a very well known theory, which states that the strength of the data depends on, the number of experts' opinion we can use combined FCMs with several experts' opinions.

At the same time the disadvantage of the combined FCM is when the weightages are 1 and –1 for the same $C_i$ $C_j$, we have the sum adding to zero thus at all times the connection matrices $E_1, \ldots, E_k$ may not be conformable for addition.

Combined conflicting opinions tend to cancel out and assisted by the strong law of large numbers, a consensus emerges as the sample opinion approximates the underlying population opinion. This problem will be easily overcome if the FCM entries are only 0 and 1. We have just briefly recalled the definitions. For more about FCMs please refer Kosko [108-112].

### 3.5.2 Description of super fuzzy cognitive maps models with illustration.

In this section we for the first time define a super FCM model and describe how it functions. Here we give the description of the multi expert super FCM model using the special super fuzzy diagonal matrix.

**DEFINITION 3.5.2.1***: Suppose n experts want to work with a problem P using a FCM model, then how to form an integrated dynamical system which can function simultaneously using the n experts opinion.*

*Suppose the first expert spells out the attributes of a problem as* $x_1^1, x_2^1, \ldots, x_{t_1}^1$, *the second expert gives the attributes as* $x_1^2, x_2^2, \ldots, x_{t_2}^2$ *and so on. Thus the* $i^{th}$ *expert gives the attributes with which he wishes to work as* $x_1^i, x_2^i, \ldots, x_{t_i}^i$; $i = 1, 2, 3, \ldots, n$. *Now we model the problem using the special diagonal super fuzzy matrix; this supermatrix will be called as the super connection matrix of the Super FCM (SFCM). We see the special feature of this special super diagonal fuzzy matrix would be all the diagonal matrices are square matrices and the main diagonal of each of these submatrices of the special fuzzy super diagonal matrix is zero. The special diagonal super fuzzy matrix for the problem P takes the following form and is denoted by* $M_D$.

$$
\begin{array}{c}
 \\
x_1^1 \\ x_2^1 \\ \vdots \\ x_{t_1}^1 \\ x_1^2 \\ x_2^2 \\ \vdots \\ x_{t_2}^2 \\ \vdots \\ x_1^n \\ x_2^n \\ \vdots \\ x_{t_n}^n
\end{array}
\begin{array}{c}
\begin{array}{cccc} x_1^1 \; x_2^1 \ldots x_{t_1}^1 & x_1^2 \; x_2^2 \ldots x_{t_2}^2 & \ldots & x_1^n \; x_2^n \ldots x_{t_n}^n \end{array} \\
\left[\begin{array}{c|c|c|c}
M_1^1 & (0) & \ldots & (0) \\
\hline
(0) & M_2^2 & & (0) \\
\hline
 & & & \\
\hline
(0) & (0) & & M_n^n
\end{array}\right]
\end{array}.
$$

*We see $M_i^i$ is a fuzzy matrix with main diagonal elements to be zero i.e.*

$$M_i^i = \begin{bmatrix} 0 & m_{12}^i & \dots & m_{1t_i}^i \\ m_{21}^i & 0 & \dots & m_{2t_i}^i \\ \vdots & \vdots & & \vdots \\ m_{t_i 2}^i & m_{t_i 2}^i & \dots & 0 \end{bmatrix}$$

*i = 1, 2, … , n.*

We illustrate the functioning of this model.

This model will be known as the multi expert super fuzzy cognitive maps model and the associated fuzzy supermatrix would be known as the special diagonal fuzzy supermatrix.

***Example 3.5.2.1:*** Suppose we have 3 experts who wish to work with a problem using FCM. The problem they wish to investigate is to analyze the Indian political situation to predict the possible electoral winner or how people tend to prefer a particular politician and so on or so forth. All of them choose to use the FCM model.

The first expert wishes to work with the following six nodes.

$x_1^1$ - Language
$x_2^1$ - Community
$x_3^1$ - Service to people
$x_4^1$ - Finance they have
$x_5^1$ - Media they can accesses to
$x_6^1$ - Party's strength and opponents strength.

The second expert wants to work with the following five nodes.

$x_1^2$ - Working members of the party

$x_2^2$ - Party's popularity in media

$x_3^2$ - The local communities strength and weakness in relative to the politicians community

$x_4^2$ - Media's accessibility

$x_5^2$ - Popularity of the politician in the context of public opinion.

The third expert wishes to work with the following attributes or nodes

$x_1^3$ - Language and caste of the public

$x_2^3$ - The finance the politician can spend in propaganda

$x_3^3$ - Opponents strength

$x_4^3$ - Parties weaknesses

$x_5^3$ - Service done by the party in that village.

$x_6^3$ - The party's popularity in media.

$x_7^3$ - Public figure configuration.

Now using these 3 experts we obtain the following super special diagonal fuzzy matrix $M_D$ which is as follows.

$$
\begin{array}{c}
\\
\begin{matrix} x_1^1 \\ x_2^1 \\ x_3^1 \\ x_4^1 \\ x_5^1 \\ x_6^1 \end{matrix} \\ \hline
\begin{matrix} x_1^2 \\ x_2^2 \\ x_3^2 \\ x_4^2 \\ x_5^2 \end{matrix} \\ \hline
\begin{matrix} x_1^3 \\ x_2^3 \\ x_3^3 \\ x_4^3 \\ x_5^3 \\ x_6^3 \\ x_7^3 \end{matrix}
\end{array}
\begin{array}{c}
\begin{array}{ccc}
x_1^1\ x_2^1\ x_3^1\ x_4^1\ x_5^1\ x_6^1 & x_1^2\ x_2^2\ x_3^2\ x_4^2\ x_5^2 & x_1^3\ x_2^3\ x_3^3\ x_4^3\ x_5^3\ x_6^3\ x_7^3
\end{array} \\
\left[\begin{array}{c|c|c}
(x_{ij}^1) & (0) & (0) \\ \hline
(0) & (x_{ij}^2) & (0) \\ \hline
(0) & (0) & (x_{ij}^3)
\end{array}\right]
\end{array}
$$

We see $(x_{ij}^1)$ is a $6 \times 6$ matrix in which $x_{ii}^1 = 0$ for i = 1, 2, 3, …, 6 i.e. the diagonal terms are zero and the $x_{ij}^1 \in \{0, 1, -1\}$; $1 \leq i, j \leq 6$. Similarly $(x_{ij}^2)$ is a $5 \times 5$ matrix in which $x_{kk}^2 = 0$, k = 1, 2, ..., 5 and $x_{ij}^2 \in \{0, 1, -1\}$; $1 \leq i, j \leq 5$. The fuzzy matrix $(x_{ij}^3)$ is such that $x_{ij}^3 = 0$, j = 1, 2, ..., 7 with $x_{ij}^3 \in \{0, 1, -1\}$. The diagonal elements of this special fuzzy supermatrix are all square matrices.

However we see this problem is for mere illustration.

| | $x_1^1$ | $x_2^1$ | $x_3^1$ | $x_4^1$ | $x_5^1$ | $x_6^1$ | $x_1^2$ | $x_2^2$ | $x_3^2$ | $x_4^2$ | $x_5^2$ | $x_1^3$ | $x_2^3$ | $x_3^3$ | $x_4^3$ | $x_5^3$ | $x_6^3$ | $x_7^3$ |
|---|---|---|---|---|---|---|---|---|---|---|---|---|---|---|---|---|---|---|
| $x_1^1$ | 0 | 1 | –1 | 0 | 0 | 1 | | | | | | | | | | | | |
| $x_2^1$ | 1 | 0 | 0 | 1 | 1 | 0 | | | | | | | | | | | | |
| $x_3^1$ | –1 | 0 | 0 | 1 | 0 | 0 | | | (0) | | | | | | (0) | | | |
| $x_4^1$ | 0 | 0 | 1 | 0 | –1 | 0 | | | | | | | | | | | | |
| $x_5^1$ | 1 | 0 | 0 | 0 | 0 | 1 | | | | | | | | | | | | |
| $x_6^1$ | 0 | 1 | –1 | 0 | 0 | 0 | | | | | | | | | | | | |
| $x_1^2$ | | | | | | | 0 | 1 | 0 | 1 | 0 | | | | | | | |
| $x_2^2$ | | | | | | | 1 | 0 | 1 | 0 | 1 | | | | | | | |
| $x_3^2$ | | | (0) | | | | 0 | 1 | 0 | –1 | 0 | | | | (0) | | | |
| $x_4^2$ | | | | | | | 1 | 0 | –1 | 0 | 1 | | | | | | | |
| $x_5^2$ | | | | | | | 0 | 0 | 1 | –1 | 0 | | | | | | | |
| $x_1^3$ | | | | | | | | | | | | 0 | 0 | 1 | 0 | 1 | 1 | 1 |
| $x_2^3$ | | | | | | | | | | | | 0 | 0 | 1 | 1 | 0 | 0 | 1 |
| $x_3^3$ | | | | | | | | | | | | 1 | 0 | 0 | 0 | 1 | 1 | 0 |
| $x_4^3$ | | | (0) | | | | | | (0) | | | 0 | 0 | 1 | 0 | 0 | 0 | 1 |
| $x_5^3$ | | | | | | | | | | | | 1 | 1 | 0 | 0 | 0 | 1 | 1 |
| $x_6^3$ | | | | | | | | | | | | 1 | –1 | 0 | 0 | 0 | 0 | –1 |
| $x_7^3$ | | | | | | | | | | | | –1 | 0 | 1 | 0 | 1 | 0 | 0 |

Now let us see how the dynamical multi expert system $M_D$ functions. Suppose the expert wishes to work with

$$X = [1\ 0\ 0\ 0\ 0\ 1 \mid 0\ 1\ 0\ 0\ 0 \mid 0\ 0\ 1\ 0\ 0\ 1\ 0\ ];$$

we want to obtain the hidden pattern of this state vector X.

$$X \text{ o } M_D = [0\ 2\ –2\ 0\ 0\ 1 \mid 1\ 0\ 1\ 0\ 1 \mid 2\ –1\ 0\ 0\ 1\ 1\ –1]$$

after updating and thresholding we get

$$X_1 = [1\ 1\ 0\ 0\ 0\ 1 \mid 1\ 1\ 1\ 0\ 1 \mid 1\ 0\ 1\ 0\ 1\ 1\ 0\ ].$$

$$X_1 \text{ o } M_D = [1\ 2\ –2\ 1\ 1\ 1 \mid 1\ 2\ 2\ –1\ 1 \mid 3\ 0\ 1\ 0\ 2\ 3\ 1]$$

after updating and thresholding we get the resultant $X_2$ to be

$$X_2 = [1\ 1\ 0\ 1\ 1\ 1 \mid 1\ 1\ 1\ 0\ 1 \mid 1\ 0\ 1\ 0\ 1\ 1\ 1]$$

and so on until we arrive at fixed point or a limit cycle. Thus we see we can use the same C-program used for FCMs only what we need is to partition the hidden pattern properly. We see because of the programs we can use any number of experts opinion. Further the demerit of combined FCM which at times tends to zero if we have 1 and –1 as a entry is over come by this method. Further as every ones feeling is given by a single super fuzzy row vector which helps in a very easy comparison. Further this is the easy or simple model which can be handled by any non mathematician.

# FURTHER READING

We have given a long list for further reading so that the interested reader can make use of them:

# INDEX

## A

## B

## C

## D

## E

## F

## H

## L

## M

## N

## O

## P

## R

## S

## T

# ABOUT THE AUTHORS

**Dr.W.B.Vasantha Kandasamy** is an Associate Professor in the Department of Mathematics, Indian Institute of Technology Madras, Chennai. In the past decade she has guided 11 Ph.D. scholars in the different fields of non-associative algebras, algebraic coding theory, transportation theory, fuzzy groups, and applications of fuzzy theory of the problems faced in chemical industries and cement industries.

She has to her credit 640 research papers. She has guided over 57 M.Sc. and M.Tech. projects. She has worked in collaboration projects with the Indian Space Research Organization and with the Tamil Nadu State AIDS Control Society. This is her 32nd book.

On India's 60th Independence Day, Dr.Vasantha was conferred the Kalpana Chawla Award for Courage and Daring Enterprise by the State Government of Tamil Nadu in recognition of her sustained fight for social justice in the Indian Institute of Technology (IIT) Madras and for her contribution to mathematics. (The award, instituted in the memory of Indian-American astronaut Kalpana Chawla who died aboard Space Shuttle Columbia). The award carried a cash prize of five lakh rupees (the highest prize-money for any Indian award) and a gold medal.

She can be contacted at vasanthakandasamy@gmail.com
You can visit her on the web at: http://mat.iitm.ac.in/~wbv or: http://www.vasantha.net

---

**Dr. Florentin Smarandache** is an Associate Professor of Mathematics at the University of New Mexico in USA. He published over 75 books and 100 articles and notes in mathematics, physics, philosophy, psychology, literature, rebus. In mathematics his research is in number theory, non-Euclidean geometry, synthetic geometry, algebraic structures, statistics, neutrosophic logic and set (generalizations of fuzzy logic and set respectively), neutrosophic probability (generalization of classical and imprecise probability). Also, small contributions to nuclear and particle physics, information fusion, neutrosophy (a generalization of dialectics), law of sensations and stimuli, etc.

He can be contacted at smarand@unm.edu

---

**Rev. Fr. K. Amal, SJ** is presently pursuing his Ph.D. in Mathematics in the Madurai Kamaraj University, Tamil Nadu, India. Having completed his bachelor's degree in mathematics at the St. Joseph’s

College, Trichy, he entered the Society of Jesus in July 1977. In 1985, he completed his master's degree in mathematics in Loyola College, Chennai and his master's degree in sociology in 1992. Ordained in 1988, he started his service as a parish priest in Mallihapuram, a Dalit parish, in the Chengelpet district. He later served the All India Catholic University Federation (AICUF ), a student movement, for 15 long years till 2003; in the latter seven-and-a-half years he occupied the position of National Director, AICUF. During his tenure, the movement witnessed newer areas of operation and scaled greater heights in its expansion and animation works, particularly with Dalit and Tribal students. Transferred to the Indian Social Institute, Bangalore, in 2004, as its Superior, he was invited to be the Director of Works of the Human Rights Desk and the Training Unit. Besides this, he has been an Advisor to the **Asia-Pacific** Chapter of the International Movement of Catholic Students (IMCS) since June 2003.

In this capacity he had traveled far and wide in the Asia Pacific region and elsewhere too. He has organized, participated, animated and served as a resource person in many workshops, seminars and meetings. He can be contacted at kamalsj@gmail.com